\documentclass[11pt, english]{smfart}
\usepackage{amsfonts, amssymb, amsmath, latexsym}
\usepackage{mathptmx}
\usepackage[english]{babel}
\usepackage{mysectionsen}
\usepackage[latin1]{inputenc}
\usepackage[T1]{fontenc}
\usepackage{textcomp}

\input xy
\xyoption{all}

\newcommand{\G}{\mathrm{G}}
\newcommand{\PGL}{\mathrm{PGL}}
\newcommand{\SL}{\mathrm{SL}}
\newcommand{\Par}{\mathrm{Par}}
\newcommand{\V}{\mathrm{V}}
\newcommand{\W}{\mathrm{W}}

\newcommand{\K}{\mathrm{K}}

\newcommand{\Qr}{\mathrm{Q}}
\newcommand{\an}{\mathrm{an}}
\newcommand{\rP}{\mathrm{P}}
\newcommand{\U}{\mathrm{U}}
\newcommand{\Sr}{\mathrm{S}}
\newcommand{\R}{\mathrm{R}}

\newcommand{\X}{\mathrm{X}}
\newcommand{\A}{\mathrm{A}}
\newcommand{\T}{\mathrm{T}}
\newcommand{\N}{\mathrm{N}}
\newcommand{\I}{\mathrm{I}}
\newcommand{\J}{\mathrm{J}}
\newcommand{\E}{\mathrm{E}}
\newcommand{\Y}{\mathrm{Y}}
\newcommand{\F}{\mathrm{F}}
\newcommand{\Z}{\mathrm{Z}}
\newcommand{\Hr}{\mathrm{H}}
\newcommand{\B}{\mathrm{B}}
\newcommand{\M}{\mathrm{M}}
\newcommand{\Lr}{\mathrm{L}}
\newcommand{\C}{\mathrm{C}}
\newcommand{\D}{\mathrm{D}}

\newcommand{\rad}{\mathrm{rad}}
\newcommand{\radu}{\mathrm{rad}^{\rm u}}

\newcommand{\inv}{^{-1}}

\begin{document}

\thispagestyle{empty}

\begin{center} {\Large \bfseries{\textsc{Bruhat-Tits theory from Berkovich's point of view.}}}

\vspace{0.1cm}
{\Large \bfseries{\textsc{I --- Realizations and compactifications of buildings}}}
\end{center}

\vspace{0.3cm}

\begin{center}
\textsc{Bertrand R\'emy, Amaury Thuillier and Annette Werner}
\end{center}

\vspace{0.3cm}

\begin{center} March 2009
\end{center}

\vspace{5cm}

\hrule

\vspace{0,5cm}

{\small 
\noindent
{\bf Abstract:}
We investigate Bruhat-Tits buildings and their compactifications by means of Berkovich analytic geometry over complete non-Archimedean fields. For every reductive group  $\G$ over a suitable non-Archimedean field $k$ we define a map
from the Bruhat-Tits building $\mathcal{B}(\G,k)$ to the Berkovich analytic space $\G^{\rm an}$ asscociated with $\G$.
Composing this map with the projection of $\G^{\rm an}$ to its flag varieties, we define a family of compactifications of $\mathcal{B}(\G,k)$. This generalizes results by Berkovich in the case of split groups. 

Moreover, we show that the boundary strata of the compactified buildings are precisely the  Bruhat-Tits buildings associated with a certain class of parabolics.  We also investigate the stabilizers of boundary points and prove a mixed Bruhat decomposition theorem for them. 

\vspace{0,2cm}


\vspace{0,1cm}

\noindent
{\bf AMS classification (2000):}
20E42, 
51E24, 
14L15,
14G22.
}

\vspace{0,5cm}

\hrule

\newpage

\tableofcontents

\newpage

\section*{Introduction}
\label{s - intro}
{\bf 1.}~In the mid 60ies, F. Bruhat and J. Tits initiated a theory which led to a deep understanding of reductive algebraic groups over valued fields~\cite{BT1a},~\cite{BT1b}.
The main tool (and a concise way to express the achievements) of this long-standing work is the notion of a \textit{building}.
Generally speaking, a building is a gluing of (poly)simplicial subcomplexes, all isomorphic to a given tiling naturally acted upon by a Coxeter group~\cite{AbramenkoBrown}.
The copies of this tiling in the building are called \textit{apartments} and must satisfy, by definition, strong incidence properties which make the whole space very symmetric.
The buildings considered by F. Bruhat and J. Tits are Euclidean ones, meaning that their apartments are Euclidean tilings (in fact, to cover the case of non-discretely valued fields, one has to replace Euclidean tilings by affine spaces acted upon by a Euclidean reflection group with a non-discrete, finite index, translation subgroup~\cite{TitsCome}).
A Euclidean building carries a natural non-positively curved metric, which allows one to classify in a geometric way maximal bounded subgroups in the rational points of a given non-Archimedean semisimple algebraic group.
This is only an instance of the strong analogy between the Riemannian symmetric spaces associated with semisimple real Lie groups and Bruhat-Tits buildings~\cite{TitsICM}.
This analogy is our guideline here.

Indeed, in this paper we investigate Bruhat-Tits buildings and their compactification by means of analytic geometry over non-Archimedean valued fields, as developed by V. Berkovich --- see~\cite{BerkoICM} for a survey.
Compactifications of symmetric spaces is now a very classical topic, with well-known applications to group theory (e.g., group cohomology~\cite{BorelSerre}) and to number theory (via the study of some relevant moduli spaces modeled on Hermitian symmetric spaces~\cite{DeligneShimura}).
For deeper motivation and a broader scope on compactifications of symmetric spaces, we refer to the recent book~\cite{BorelJi}, in which the case of locally symmetric varieties is also covered.
One of our main results is to construct for each semisimple group $\G$ over a suitable non-Archimedean valued field $k$, a family of compactifications of the Bruhat-Tits building $\mathcal{B}(\G,k)$ of $\G$ over $k$.
This family is finite, actually indexed by the conjugacy classes of proper parabolic $k$-subgroups in $\G$.
Such a family is of course the analogue of the family of Satake~\cite{Sa} or Furstenberg~\cite{Furst} compactifications of a given Riemannian non-compact symmetric space --- see~\cite{GJT} for a general exposition.

In fact, the third author had previously associated, with each Bruhat-Tits building, a family of compactifications also indexed by the conjugacy classes of proper parabolic $k$-subgroups~\cite{Wer2} and generalizing the "maximal" version constructed before by E. Landvogt~\cite{La}.
The Bruhat-Tits building $\mathcal{B}(\G,k)$ of $\G$ over $k$ is defined as the quotient for a suitable equivalence relation, say $\sim$, of the product of the rational points $\G(k)$ by a natural model, say $\Lambda$, of the apartment; we will refer to this kind of construction as a \textit{gluing procedure}.
The family of compactifications of~\cite{Wer2} was obtained by suitably compactifying $\Lambda$ to obtain a compact space $\overline{\Lambda}$ and extending $\sim$ to an equivalence relation on $\G(k) \times \overline{\Lambda}$.
As expected, for a given group $\G$ we eventually identify the latter family of compactifications with the one we construct here, see \cite{RTW2}.

Our compactification procedure makes use of embeddings of Bruhat-Tits buildings in the analytic versions of some well-known homogeneous varieties (in the context of algebraic transformation groups), namely flag manifolds.
The idea goes back to V. Berkovich in the case when $\G$ splits over its ground field $k$~\cite[\S 5]{Ber1}.
One aesthetical advantage of the embedding procedure is that it is similar to the historical ways to compactify symmetric spaces, e.g., by seeing them as topological subspaces of some projective spaces of Hermitian matrices or inside spaces of probability measures on a flag manifold.
More usefully (as we hope), the fact that we specifically embed buildings into compact spaces from Berkovich's theory may make these compactifications useful for a better understanding of non-Archimedean spaces relevant to number theory (in the case of Hermitian symmetric spaces).
For instance, the building of ${\rm GL}_n$ over a valued field $k$ is the "combinatorial skeleton" of the Drinfel'd half-space $\Omega^{n-1}$ over $k$~\cite{BoutotCarayol}, and it would be interesting to know whether the precise combinatorial description we obtain for our compactifications might be useful to describe other moduli spaces for suitable choices of groups and parabolic subgroups.
One other question about these compactifications was raised by V. Berkovich himself~\cite[5.5.2]{Ber1} and deals with the potential generalization of Drinfel'd half-spaces to non-Archimedean semisimple algebraic groups of arbitrary type.

\vspace{3pt}

{\bf 2.}~Let us now turn to the definition of the embedding maps that allow us to compactify Bruhat-Tits buildings.
Let $\G$ be a $k$-isotropic semisimple algebraic group defined over the non-Archimedean valued field $k$ and let $\mathcal{B}(\G,k)$ denote the Euclidean building provided by Bruhat-Tits theory
\cite{TitsCorvallis}.
We prove the following statement (see \ref{ss - buildings in flags} and Prop. \ref{prop.compact.theta-t}):
\textit{
assume that the valued field $k$ is a local field (i.e., is locally compact) and (for simplicity) that $\G$ is almost $k$-simple; then for any conjugacy class of proper parabolic $k$-subgroup, say $t$, there exists a continuous, $\G(k)$-equivariant map $\vartheta_t: \mathcal{B}(\G,k) \rightarrow \mathrm{Par}_t(\G)^{\rm an}$ which is a homeomorphism onto its image.
}
Here $\mathrm{Par}_t(\G)$ denotes the connected component of type $t$ in the proper variety $\mathrm{Par}(\G)$ of all parabolic subgroups in $\G$ (on which $\G$ acts by conjugation)
\cite[Expos\'e XXVI, Sect. 3]{SGA3}.
The superscript ${}^{\mathrm{an}}$ means that we pass from the $k$-variety $\mathrm{Par}_t(\G)$ to the Berkovich $k$-analytic space associated with it~\cite[3.4.1-2]{Ber1}; the space $\mathrm{Par}(\G)^{\mathrm{an}}$ is compact since $\mathrm{Par}(\G)$ is projective.
We denote by $\overline{\mathcal{B}}_t(\G,k)$ the closure of the image of $\vartheta_t$ and call it the \textit{Berkovich compactification} of type $t$ of the Bruhat-Tits building $\mathcal{B}(\G,k)$.

Roughly speaking, the definition of the maps $\vartheta_t$ takes up the first half of this paper, so let us provide some further information about it.
As a preliminary, we recall some basic but helpful analogies between (scheme-theoretic) algebraic geometry and $k$-analytic geometry (in the sense of Berkovich).
Firstly, the elementary blocks of $k$-analytic spaces in the latter theory are the so-called \textit{affinoid} spaces; they, by and large, correspond to affine schemes in algebraic geometry. Affinoid spaces can be glued together to define $k$-analytic spaces, examples of which are provided by analytifications of affine schemes: if $\X={\rm Spec}(\A)$ is given by a finitely generated $k$-algebra $\A$, then the set underlying the analytic space $\X^{\mathrm{an}}$ consists of multiplicative seminorms on $\A$ extending the given absolute value on $k$.
Let us simply add that it follows from the "spectral analytic side" of Berkovich theory that each affinoid space $\X$ admits a \textit{Shilov boundary}, namely a (finite) subset on which any element of the Banach $k$-algebra defining $\X$ achieves its minimum.
We have enough now to give a construction of the maps $\vartheta_t$ in three steps:

\begin{itemize}
\item[Step 1:] we attach to any point $x \in \mathcal{B}(\G,k)$ an affinoid subgroup $\G_x$ whose $k$-rational points coincide with the parahoric subgroup $\G_x(k)$ associated with $x$ by Bruhat-Tits theory (Th.
\ref{thm.affinoid.subgroups}).
\item[Step 2:] we attach to any so-obtained analytic subgroup $\G_x$ a point $\vartheta(x)$ in $\G^{\mathrm{an}}$ (in fact the unique point in the Shilov boundary of $\G_x$), which defines a map $\vartheta: \mathcal{B}(\G,k) \to \G^{\mathrm{an}}$ (Prop \ref{prop.theta}).
\item[Step 3:] we finally compose the map $\vartheta$ with an "orbit map" to the flag variety $\mathrm{Par}_t(\G)^{\mathrm{an}}$ of type $t$ (Def. \ref{def.theta-t}).
\end{itemize}

Forgetting provisionally that we wish to compactify the building $\mathcal{B}(\G,k)$ (in which case we have to assume that $\mathcal{B}(\G,k)$ is locally compact, or equivalently, that $k$ is local), this three-step construction of the map $\vartheta_t: \mathcal{B}(\G,k) \to \mathrm{Par}_t(\G)^{\mathrm{an}}$ works whenever the ground field $k$ allows the functorial existence of $\mathcal{B}(\G,k)$ (see \ref{ss - Bruhat-Tits} for a reminder of these conditions).
We note that in Step 2, the uniqueness of the point $\vartheta(x)$ in the Shilov boundary of $\G_x$ comes from the use of a field extension splitting $\G$ and allowing to see $x$ as a special point (see below) and from the fact that integral structures attached to special points in Bruhat-Tits theory are explicitly described by means of Chevalley bases.
At last, the point $\vartheta(x)$ determines $\G_x$ because the latter analytic subgroup is the holomorphic envelop of $\vartheta(x)$ in $\G^{\mathrm{an}}$.
Here is a precise statement for Step 1 (Th. \ref{thm.affinoid.subgroups}).

\vspace{3pt}
\noindent \emph{\textbf{Theorem 1}} ---
\textit{
For any point $x$ in $\mathcal{B}(\G,k)$, there is a unique affinoid subgroup $\G_x$ of $\G^{\rm an}$ satisfying the following condition: for any non-Archimedean extension $\K$ of $k$, we have $\G_x(\K) = \mathrm{Stab}_{\G(\K)}(x).$
}
\vspace{3pt}

This theorem (hence Step 1) improves an idea used for another compactification procedure, namely the one using the map attaching to each point $x \in \mathcal{B}(\G,k)$ the biggest parahoric subgroup of $\G(k)$ fixing it~\cite{GuiRem}.
The target space of the map $x \mapsto \G_x(k)$ in [loc. cit.] is the space of closed subgroups of $\G(k)$, which is compact for the Chabauty topology~\cite[VIII.5]{Integration78}.
This idea doesn't lead to a compactification of $\mathcal{B}(\G,k)$ but only of the set of vertices of it: if $k$ is discretely valued and if $\G$ is simply connected, any two points in a given facet of the Bruhat-Tits building $\mathcal{B}(\G,k)$ have the same stabilizer.
Roughly speaking, in the present paper we use Berkovich analytic geometry, among other things, to overcome these difficulties thanks to the fact that we can use arbitrarily large non-Archimedean extensions of the ground field.
More precisely, up to taking a suitable non-Archimedean extension $\K$ of $k$, any point $x \in \mathcal{B}(\G,k)$ can be seen as a special point in the bigger (split) building $\mathcal{B}(\G,\K)$, in which case we can attach to $x$ an affinoid subgroup of $(\G \otimes_k \K)^{\mathrm{an}}$.
As a counterpart, in order to obtain the affinoid subgroup $\G_x$ defined over $k$ as in the above theorem, we have to apply a Banach module avatar
of Grothendieck's faithfully flat descent formalism~\cite[VIII]{SGA1} (Appendix 1).

As an example, consider the case where $\G={\rm SL}(3)$ and the field $k$ is discretely valued.
The apartments of the building are then tilings of the Euclidean plane by regular triangles (\textit{alcoves} in the Bruhat-Tits terminology).
If the valuation $v$ of $k$ is normalized so that $v(k^\times)=\mathbb{Z}$, then in order to define the group $\G_x$ when $x$ is the barycenter of a triangle, we have to (provisionally) use a finite ramified extension $\K$ such that $v(\K^\times)={1\over 3}\mathbb{Z}$ (the apartments in $\mathcal{B}(\G,\K)$ have "three times more walls" and $x$ lies at the intersection of three of them).
The general case, when the barycentric coordinates of the point $x$ (in the closure of its facet) are not \textit{a priori} rational, requires an \textit{a priori} infinite extension.

As already mentioned, when $\G$ splits over the ground field $k$, our compactifications have already been defined by V. Berkovich~\cite[\S 5]{Ber1}.
His original procedure relies from the very beginning on the explicit construction of reductive group schemes over $\mathbb{Z}$ by means of Chevalley bases~\cite{ChevalleyBBK}.
If $\T$ denotes a maximal split torus (with character group $X^*({\T})$), then the model for an apartment in $\mathcal{B}(\G,k)$ is $\Lambda = {\rm Hom}(X^*({\T}),\mathbb{R}_+^\times)$ seen as a real affine space.
Choosing a suitable (special) maximal compact subgroup ${\bf P}$ in $\G^{\rm an}$, V. Berkovich identifies $\Lambda$ with the image of $\T^{\rm an}$ in the quotient variety $\G^{\rm an}/{\bf P}$.
The building $\mathcal{B}(\G,k)$ thus appears in $\G^{\rm an}/{\bf P}$ as the union of the transforms of $\Lambda$ by the proper action of the group of $k$-rational points $\G(k)$ in $\G^{\rm an}/{\bf P}$.
Then V. Berkovich uses the notion of a \textit{peaked point} (and other ideas related to holomorphic convexity) in order to construct a section map $\G^{\rm an}/{\bf P} \to \G^{\rm an}$.
This enables him to realize $\mathcal{B}(\G,k)$ as a subset of $\G^{\rm an}$, which is closed if $k$ is local.

The hypothesis that $\G$ is split is crucial for the choice of the compact subgroup ${\bf P}$. The construction in Step 1 and 2 is different from Berkovich's original approach and allows a generalization to the non-split case. We finally note that, in Step 3, the embedding map $\vartheta_t: \mathcal{B}(\G,k) \rightarrow \mathrm{Par}_t(\G)^{\mathrm{an}}$ only depends on the type $t$; in particular, it doesn't depend on the choice of a parabolic $k$-subgroup in the conjugacy class corresponding to $t$.

\vspace{3pt}

{\bf 3.}~Let us henceforth assume that the ground field $k$ is locally compact.
We fix a conjugacy class of parabolic $k$-subgroups in $\G$, which provides us with a $k$-rational type $t$.
The building $\mathcal{B}(\G,k)$ is the product of the buildings of all almost-simple factors of $\G$, and we let $\mathcal{B}_t(\G,k)$ denote the quotient of $\mathcal{B}(\G,k)$ obtained by removing each almost-simple factor of $\G$ on which $t$ is trivial.
The previous canonical, continuous and $\G(k)$-equivariant map $\vartheta_t: \mathcal{B}(\G,k) \rightarrow \mathrm{Par}_t(\G)^{\rm an}$ factors through an injection $\mathcal{B}_t(\G,k) \hookrightarrow \mathrm{Par}_t(\G)^{\rm an}$.
We then consider the question of describing as a $\G(k)$-space the so-obtained compactification $\overline{\mathcal{B}}_t(\G,k)$, that is the closure of ${\rm Im}(\vartheta_t) = \mathcal{B}_t(\G,k)$ in $\mathrm{Par}_t(\G)^{\mathrm{an}}$.

The type $t$ and the scheme-theoretic approach to flag varieties we adopt in Step 3 above (in order to see easily the uniqueness of $\vartheta_t$), lead us to distinguish some other types of conjugacy classes of parabolic $k$-subgroups (\ref{ss - relevant parabolics}).
These classes are called \textit{$t$-relevant} and are defined by means of flag varieties, but we note afterwards that $t$-relevancy amounts also to a combinatorial condition on roots (Prop.
\ref{prop.roots.relevant}) which we illustrate in Example \ref{ex - relevant SL} in the case of the groups ${\rm SL}(n)$.

Moreover each parabolic subgroup $\rP \in \mathrm{Par}(\G)$ defines a closed \textit{osculatory} subvariety $\mathrm{Osc}_t(\rP)$ of $\mathrm{Par}_t(\G)$, namely the one consisting of all parabolics of type $t$ whose intersection with $\rP$ is a parabolic subgroup (Prop. \ref{prop.osc}).
Then $\rP$ is $t$-relevant if it is maximal among all parabolic $k$-subgroups defining the same osculatory subvariety.
It is readily seen that each parabolic subgroup is contained in a unique $t$-relevant one.
For instance, if $\G = {\rm PGL}({\rm V})$ and if $\delta$ is the type of flags $(0 \subset {\rm H} \subset {\rm V})$ where ${\rm H}$ is a hyperplane of the $k$-vector space ${\rm V}$, then $\delta$-relevant parabolic $k$-subgroups are those corresponding to flags $(0 \subset {\rm W} \subset {\rm V})$, where ${\rm W}$ is a linear subspace of ${\rm V}$.
Moreover $\overline{\mathcal{B}}_{\delta}({\rm PGL}({\rm V}),k)$ is the seminorm compactification described in~\cite{Wer1}.
In general, we denote by $\R_{t}(\rP)$ the kernel of the algebraic action of $\rP$ on the variety $\mathrm{Osc}_t(\rP)$ and by $\pi_{t,\rP}$ the natural projection $\rP \twoheadrightarrow \rP/\R_t(\rP)$.
The following theorem sums up several of our statements describing $\overline{\mathcal{B}}_t(\G,k)$ as a $\G(k)$-space (see e.g., Th. \ref{thm.stratification}, Th. \ref{thm.point.stabilizer} and Prop. \ref{prop - mixed Bruhat dec}).

\vspace{3pt}
\noindent \emph{\textbf{Theorem 2}} ---
\textit{
Let $\G$ be a connected semisimple linear algebraic group defined over a non-Archimedean local field $k$ and let $t$ be the type of a proper parabolic $k$-subgroup in $\G$.
We denote by $\mathcal{B}(\G,k)$ its Bruhat-Tits building and by $\overline{\mathcal{B}}_t(\G,k)$ the Berkovich compactification of type $t$ of the latter space.
\begin{itemize}
\item[(i)] For any proper $t$-relevant parabolic $k$-subgroup $\rP$, there exists a natural continuous map
$\mathcal{B}_t(\rP/{\rm rad}(\rP),k) \to \overline{\mathcal{B}}_t(\G,k)$ whose image lies in the boundary.
These maps altogether provide the following stratification:
$$\overline{\mathcal{B}}_t(\G,k) = \bigsqcup_{\textrm{t-relevant  $\rP$'s}} \mathcal{B}_t(\rP/{\rm rad}(\rP),k),$$
where the union is indexed by the $t$-relevant parabolic $k$-subgroups in $G$.
\item[(ii)] Let $x$ be a point in a stratum $\mathcal{B}_t(\rP/{\rm rad}(\rP),k)$.
Then there is a $k$-analytic subgroup ${\rm Stab}_{\G}^t(x)$ of $\G^{\rm an}$ such that ${\rm Stab}_{\G}^t(x)(k)$ is the stabilizer of $x$ in $\G(k)$.
Moreover we have ${\rm Stab}_{\G}^t(x) = \pi_{t,P}^{-1}((\rP/\R_t(\rP))_x)$, where $(\rP/\R_t(\rP))_x$ is the $k$-affinoid subgroup of $(\rP/\R_t(\rP))^{\rm an}$ attached by theorem 1 to the point $x$ of
$\mathcal{B}_t(\rP/{\rm rad}(\rP),k)=\mathcal{B}(\rP/\R_t(\rP),k)$.
\item[(iii)] Any two points $x, y$ in $\overline{\mathcal{B}}_t(\G,k)$ lie in a common compactified apartment $\overline{A}_t(\Sr,k)$ and we have:
$$\G(k) = {\rm Stab}_{\G}^t(x)(k) \N(k) {\rm Stab}_{\G}^t(y)(k),$$ where $\N$ is the normalizer of the maximal split torus $\Sr$ defining the apartment $\A(\Sr,k)$.
\end{itemize}
}
\vspace{3pt}

Statement (i) in the above theorem says that the boundary added by taking a closure in the embedding procedure consists of Bruhat-Tits buildings, each of these being isomorphic to the Bruhat-Tits building of some suitable Levi factor (Prop. \ref{prop.strata.groups}).
This phenomenon is well-known in the context of symmetric spaces~\cite{Sa}.
Statement (ii) first says that a boundary point stabilizer is a subgroup of a suitable parabolic $k$-subgroup in which, roughly speaking, some almost simple factors of the Levi factor are replaced by parahoric subgroups geometrically determined by the point at infinity.
In the case $\G = {\rm PGL}({\rm V})$ with $\delta$ as above, the $\delta$-relevant parabolic $k$-subgroups (up to conjugacy) are those having exactly two diagonal blocks, and the boundary point stabilizers are simply obtained by replacing exactly one block by a parahoric subgroup of it.
At last, statement (iii) is often referred to as the \textit{mixed Bruhat decomposition}.


\vspace{3pt}

{\bf 4.}~At this stage,  we understand the finite family of Berkovich compactifications $\overline{\mathcal{B}}_t(\G,k)$, indexed by the $k$-rational types $t$. We describe in 4.2 the natural continuous and $\G(k)$-equivariant maps between these compactifications arising from fibrations between flag varieties and we show in Appendix C that no new compactification arises from non-rational types of parabolic subgroup. In a sequel to this article~\cite{RTW2}, we will (a) compare Berkovich compactifications with the ones defined by the third author in~\cite{Wer2}, relying on a gluing procedure and the combinatorics of weights of an absolutely irreducible linear representations of $\G$, and (b) as suggested in [loc.cit], show (from two different viewpoints) that these compactifications can also be described in a way reminiscent to Satake's original method for compactifying riemanniann symmetric spaces.

\vspace{3pt}

{\bf 5.}~Let us close this introduction by two remarks.
The first one simply consists in mentioning why it is interesting to have non-maximal compactifications of Bruhat-Tits buildings.
This is (at least) because in the case of Hermitian locally symmetric spaces, some interesting compactifications, namely the Baily-Borel ones~\cite{BailyBorel}, are obtained as quotients of \textit{minimal} compactifications (of a well-defined type) by arithmetic lattices.
The second remark deals with the Furstenberg embedding approach, consisting in sending a symmetric space into the space of probability measures on the various flag varieties of the isometry group~\cite{Furst}.
In the Bruhat-Tits case, this method seems to encounter new difficulties compared to the real case.
The main one is that not all maximal compact subgroups in a simple non-Archimedean Lie group act transitively on the maximal flag variety of the group.
This is well-known to specialists in harmonic analysis (e.g., one has to choose a special maximal compact subgroup to obtain a Gelfand pair).
The consequence for Furstenberg compactifications is that, given a non-special vertex $v$ with stabilizer $\G_v(k)$, it is not clear, in order to attach a $\G_v(k)$-invariant probability measure $\mu_v$ to $v$, how to distribute the mass of $\mu_v$ among the $\G_v(k)$-orbits in the flag variety.
We think that in the measure-theoretic approach, some subtle problems of this kind deserve to be investigated, though the expected compactifications are constructed in the present paper by the Berkovich approach.

\vspace{3pt}

{\bf Conventions.}~
Let us simply recall a standard convention (already used above): a \textit{local field} is a non-trivially and discretely valued field which is locally compact for the topology arising from the valuation; this amounts to saying that it is complete and that the residue field is finite.

Roughly speaking this paper applies some techniques from algebraic (and analytic) geometry in order to prove some group-theoretic statements.
Conventions in these two different fields are sometimes in conflict.
We tried to uniformly prefer the conventions from algebraic geometry since they are the ones which are technically used.
For instance, it is important for us to use varieties of parabolic subgroups~\cite{SGA3} rather than flag varieties, even though they don't have any rational point over the ground field
and the affine and projective spaces are those defined in~\cite{EGA}.

Accordingly, our notation for valued fields are that of V. Berkovich's book~\cite{Ber1}; in particular, the valuation ring of such a field $k$ is denoted by $k^\circ$ and its maximal ideal is denoted by $k^{\circ\circ}$ (1.2.1).

\vspace{0.1cm} {\bf Working hypothesis.}~ The basic idea underlying this work is to rely on functoriality of Bruhat-Tits buildings with respect to field extensions. The required assumptions on the group or on the base field are discussed in (1.3.4).  

\vspace{5pt}

{\bf Structure of the paper.}~
In the first section, we briefly introduce Berkovich's theory of analytic geometry over complete non-Archimedean fields and Bruhat-Tits theory of reductive algebraic groups over valued fields. 
The second section is devoted to realizing the Bruhat-Tits buildings of reductive groups over complete valued fields as subsets of several spaces relevant to analytic geometry, namely the analytic spaces attached to the groups themselves, as well as the analytic spaces associated with the various flag varieties of the groups.
The third section deals with the construction of the compactifications, which basically consists in taking the closures of the images of the previous maps; it has also a Lie theoretic part which provides in particular the tools useful to describe the compactifications in terms of root systems and convergence in Weyl chambers.
The fourth section is dedicated to describing the natural action of a non-Archimedean reductive group on Berkovich compactifications of its building.

At last, in one appendix we extend the faithfully flat descent formalism in the Berkovich context because it is needed in the second section, and in the other appendix we prove some useful technicalities on fans, in connection with compactifications.


\section{Berkovich geometry and Bruhat-Tits buildings}
\label{s - preliminaries}

The main goal of this section is to recall some basic facts about the main two topics which are "merged" in the core of the paper in order to compactify Euclidean buildings.
These topics are non-Archimedean analytic geometry according to Berkovich and the Bruhat-Tits theory of algebraic groups over valued fields, whose main tools are the geometry of buildings and integral structures for the latter groups.
This requires to fix first our basic conventions about algebraic groups.
Concerning Berkovich geometry, the only new result here is a criterion for being affinoid over the ground field; it uses Grothendieck's formalism of faithfully flat descent adapted (in the first appendix) to our framework.
Concerning Bruhat-Tits theory, we use (huge) extensions of the ground field to improve dramatically the transitivity properties of the actions of the rational points of a reductive group on the points of the associated building.

\subsection{Algebraic groups}
\label{ss - group schemes}
This subsection is merely dedicated to fix some conventions concerning the algebraic groups and related algebraic varieties we use in this paper.
The main point is that, in order to fully use Berkovich's approach to non-Archimedean analytic geometry, we need to adopt the framework of group schemes as developed in~\cite{SGA3} or~\cite{DemazureGabriel} --- an introductory reference is~\cite{Waterhouse}.
As an example, if $k$ is a field and $\G$ is a $k$-group scheme, a subgroup of $\G$ will always mean a $k$-subgroup scheme.

\vspace{0.2cm}
\noindent \textbf{(1.1.1)}
We use throughout the text the language of group schemes, which is arguably more precise and flexible than the somehow classical viewpoint adopted for example in~\cite{Borel}. Thus, a \emph{linear algebraic group} over a field $k$ is by definition a smooth connected affine group scheme over $k$ (see for example~\cite[21.12]{Involutions} for a translation). If $\G$ is a linear algebraic group over $k$ and $k'/k$ is any field extension, we denote by $\G(k')$ the group of $k'$-valued points of $\G$ and by $\G \otimes_k k'$ the linear algebraic group over $k'$ deduced from base change.

Let first $\overline k$ denote an algebraically closed field and let ${\rm H}$ denote a linear algebraic group over $\overline k$. The \textit{radical} (\textit{unipotent radical}, respectively) of ${\rm H}$ is the biggest connected normal solvable (unipotent, respectively) subgroup of ${\rm H}$; we denote these subgroups by $\rad({\rm H})$ and $\radu({\rm H})$, respectively.
We say that ${\rm H}$ is \textit{reductive} if the subgroup $\radu({\rm H})$ is trivial; this is so whenever $\rad({\rm H})$ is trivial, in which case we say that ${\rm H}$ is \textit{semisimple}. In general, we let $\Hr_{\rm ss}$ denote the semisimple group $\Hr/{\rm rad}(\G)$.

Now let $\Sr$ denote a scheme.
A group scheme over $\Sr$ is a group object in the category of schemes over $\Sr$~\cite[II, \S 1]{DemazureGabriel}.
Let $\G$ be such a group scheme.
For any point $s \in \Sr$ we denote by $\overline s$ the spectrum of an algebraic closure of the residue field $\kappa(s)$ of $s$.
Following~\cite[Expos\'e XIX, 2.7]{SGA3}, we say that $\G$ is a \textit{reductive} (\textit{semisimple}, respectively) $\Sr$-group scheme if it is smooth and affine over $\Sr$, and if all its geometric fibers $\G_{\overline s} = \G \times_{\Sr} \overline{s}$ are connected and reductive (connected and semisimple, respectively) in the above sense. A particular instance of reductive $\Sr$-group schemes are \emph{tori}: a smooth affine $\Sr$-group scheme $\T$ is a \emph{torus} if there exists an \'etale covering $\Sr' \rightarrow \Sr$ such that $\T \otimes_{\Sr} \Sr'$ is isomorphic to a diagonalizable group $\mathbb{G}_{\mathrm{m},\Sr'}^n$. A torus $\T$ is \emph{split} if it is isomorphic to some $\mathbb{G}_{\mathrm{m},\Sr}^n$. Finally, a \emph{maximal torus} of a $\Sr$-group scheme $\G$ is a closed subgroup $\T$ satisfying the following two conditions: (a) $\T$ is a torus, (b) for each $s \in \Sr$, $\T_{\overline s}$ is a maximal torus of $\G_{\overline{s}}$.

Reductive $\Sr$-group schemes are defined and thoroughly studied by M. Demazure in~\cite[Expos\'es XIX to XXVI]{SGA3}; an introductory article is~\cite{DemazureThese}.
We will use the results contained in these references throughout the present paper.

\vspace{0.2cm}
\noindent \textbf{(1.1.2)}
By the fundamental work of C. Chevalley, over a fixed algebraically closed field, a reductive algebraic group is completely characterized by its \textit{root datum}~\cite[9.6.2 and 10.1.1]{Springer}: this is basically a combinatorial datum encoding the Dynkin diagram, the maximal torus and the way roots appear as characters of the latter.

For any scheme $\Sr$, a \textit{Chevalley} (\textit{Demazure}, respectively) group scheme over $\Sr$ is a semisimple (reductive, respectively) $\Sr$-group scheme containing a split maximal torus; one speaks also of a split semisimple (reductive, respectively) $\Sr$-group scheme. To each such $\Sr$-group is attached a root datum and the main result of~\cite{DemazureThese} is the following generalization of Chevalley's theorem: for any non-empty scheme $\Sr$, there is a one-to-one correspondence between isomorphism classes of Demazure group schemes and isomorphism classes of root data. In particular, any Demazure group scheme over a non-empty scheme $\Sr$ comes from a Demazure group scheme over $\mathrm{Spec}(\mathbb{Z})$ by base change. For each root datum, the existence of the corresponding split semisimple group scheme $\G$ over $\mathrm{Spec}(\mathbb{Z})$ was essentially proved by Chevalley in~\cite{ChevalleyBBK}; starting with the semisimple complex group $\G(\mathbb{C})$, his strategy was to introduce a $\mathbb{Z}$-form of its Lie algebra $\mathfrak{g}_{\mathbb{C}}$ in terms of specific bases called \textit{Chevalley bases}~\cite{Steinberg}. We will use them.

\vspace{0.2cm}
\noindent \textbf{(1.1.3)}
One of the main tools we use for the compactifications is the variety of parabolic subgroups of a reductive group scheme. The reference for what follows is~\cite[Expos\'e XXVI]{SGA3}.

Let $k$ denote a field and let ${\rm G}$ denote a reductive group over $\overline k$.
A closed subgroup $\rP$ of ${\rm G}$ is called \textit{parabolic} if it is smooth and if the quotient variety ${\rm G}/\rP$ is proper over $\overline k$.
This last condition amounts to saying that $\rP$ contains a \textit{Borel subgroup} of ${\rm G}$, that is a closed connected and solvable smooth subgroup of ${\rm G}$, maximal for these properties.

More generally, if $\Sr$ is a scheme and $\G$ is a reductive group scheme over $\Sr$, then a subgroup $\rP$ is called \textit{parabolic} if it is smooth over $\Sr$ and if the quotient $\G_{\overline s}/\rP_{\overline s}$ is a proper $\overline s$-scheme for any $s \in \Sr$.
In this case, the (fppf) quotient-sheaf $\G/\rP$ is represented by an $\Sr$-scheme, which is smooth and proper over $\Sr$.

For any reductive group scheme $\G$ over $\Sr$, the functor from the category of $\Sr$-schemes to the category of sets, attaching to each $\Sr$-scheme $\Sr'$ the set of parabolic subgroups of $\G_{\Sr'}$, is representable by a smooth and projective $\Sr$-scheme $\mathrm{Par}(\G)$. If $\rP$ is a parabolic subgroup of $\G$, the natural morphism $\lambda_{\rP} : \G \rightarrow \mathrm{Par}(\G)$ defined by $$\G(\Sr') \rightarrow \mathrm{Par}(\G)(\Sr'), \ g \mapsto g (\rP \times_{\Sr} \Sr') g^{-1}$$ for any $\Sr$-scheme $\Sr'$, induces an isomorphism between $\G/\rP$ and an open and closed subscheme of $\mathrm{Par}(\G)$.

If $k$ denotes a field with algebraic closure $k^a$ and if $\Sr={\rm Spec}(k)$, then the connected components of ${\rm Par}(\G)$ are in natural one-to-one correspondence with the $\mathrm{Gal}(k^{a}|k)$-stable subsets of the set of vertices of the Dynkin diagram of $\G \otimes_k k^a$.

The \textit{type} of a parabolic subgroup, say $\rP$, of $\G$ is the connected component of ${\rm Par}(\G)$ containing it; it is denoted by $t(\rP)$.
The connected component corresponding to a given type $t$ is denoted by ${\rm Par}_t(\G)$; it is a projective smooth homogeneous space under the group $\G$, image of the canonical map $\lambda_{\rP} : \G/\rP \to {\rm Par}(\G)$.

The connected component ${\rm Par}_t(\G)$ of $\mathrm{Par}(\G)$ will occasionally be called the \textit{flag variety of type $t$} associated with $\G$; it need not contain a rational point over $k$.
When it does, such a point corresponds to a parabolic $k$-subgroup of $\G$ and we say that the type $t$ is \emph{k-rational} (or even \emph{rational}, if non confusion seems likely to arise).

Finally, if the type $t$ corresponds to the empty set of simple roots in the above description of connected components of $\mathrm{Par}(\G)$, the scheme ${\rm Par}_{\varnothing}(\G)$ is the \textit{variety of Borel subgroups} and we denote it by ${\rm Bor}(\G)$. Since ${\rm Bor}(\G)(k)$ is the set of Borel subgroups of $\G$, we have $\mathrm{Bor}(\G)(k) \neq \varnothing$ if and only if $\G$ has a Borel subgroup, i.e., if and only if $\G$ is \textit{quasi-split}. We will use ${\rm Bor}(\G)$ for an arbitrary reductive $k$-group $\G$.

\subsection{Non-Archimedean analytic geometry}
\label{ss - Berkovich}

We begin this subsection by a brief review of analytic geometry over a non-Archimedean field.
We then turn to the question of descending the base field, a technical device which lies at the core of Berkovich realizations of buildings (Sect. \ref{s - realizations}).

\vspace{0.2cm}
\noindent \textbf{(1.2.1)}
A \emph{non-Archimedean field} is a field $k$ equipped with a complete non-Archimedean absolute value $|.|$, which we assume to be non-trivial. Elements of $k$ with absolute value less or equal to $1$ define a subring $k^{\circ}$ of $k$; this ring is local and we let $\widetilde{k}$ denote its residue field. For the algebraic theory of valuations, we refer to~\cite[Chapitre 6]{BourbakiAlgCo56} and~\cite{BGR}.

\vspace{0.1cm}

A \emph{non-Archimedean extension} $\K/k$ is an isometric extension of non-Archimedean fields. For a non-Archimedean field $k$, we let $\mathbf{BMod}(k)$ denote the category of Banach $k$-modules with bounded $k$-homomorphisms. We call a bounded $k$-homomorphism $u: \M \rightarrow \N$ \emph{strict} if the quotient norm on $\M/\mathrm{ker}(u)$ is equivalent to the norm induced by $\N$ (in~\cite{Ber1} V. Berkovich calls such a homomorphism \emph{admissible}).
We let $\mathbf{BMod}^{\mathbf{st}}(k)$ denote the subcategory of $\mathbf{BMod}(k)$ in which morphisms are required to be strict.

For a non-Archimedean field $k$ and an $n$-tuple $\mathbf{r}=(r_1,\ldots, r_n)$ of positive real numbers, we let $k\{r_1^{-1}\xi_1, \ldots, r_n^{-1}\xi_n\}$ denote the $k$-algebra
$$\left\{f = \sum_{\nu=(\nu_1, \ldots, \nu_n) \in \mathbb{N}^{n}} a_{\nu} \xi_1^{\nu_1} \ldots \xi_n^{\nu_n} \ \Big| \ |a_{\nu}|r_1^{\nu_1} \ldots r_n^{\nu_n} \rightarrow 0 \textrm{ when } \nu_1 + \ldots + \nu_n \rightarrow \infty \right\}$$ equipped with the Banach norm $$||f|| = \max_{\nu \in \mathbb{N}^n} |a_{\nu}|r_1^{\nu_1} \ldots r_n^{\nu_n}.$$

When $\mathbf{r} = (r_1, \ldots, r_n)$ is a vector of positive real numbers which are \emph{linearly independent} in $(\mathbb{R}_{>0}/|k^{\times}|) \otimes_{\mathbb{Z}}\mathbb{Q}$, i.e., such that $r_1^{\nu_1} \ldots r_n^{\nu_n} \notin |k^{\times}|$ for any $\nu=(\nu_1, \ldots,\nu_n) \in \mathbb{Z}^n - \{0\}$, the Banach $k$-algebra $k\{r_1^{-1}\xi_1, r_1 \xi_1^{-1}, \ldots, r_n^{-1}\xi_n, r_n \xi_n^{-1}\}$ is a non-Archimedean field which we denote by $k_{\mathbf{r}}$.

Let ${\rm M}$ and ${\rm N}$ be Banach $k$-modules, all of whose norms are denoted by $\parallel . \parallel$.
Then we can consider on t0he classical (i.e., algebraic) tensor product ${\rm M} \otimes_k {\rm N}$ a norm, also denoted by $\parallel . \parallel$, and defined by
$\parallel f \parallel = \inf \max_i \parallel m_i \parallel \cdot \parallel n_i \parallel$, where the infimum is taken over all the representations $f = \sum_i m_i \otimes n_i$.
The completion of  ${\rm M} \otimes_k {\rm N}$ with respect to $\parallel . \parallel$ is called the \textit{completed tensor product} of ${\rm M}$ and ${\rm N}$, and is denoted by
${\rm M} \widehat{\otimes}_k {\rm N}$.
The notion of completed tensor product of homomorphisms is defined similarly.

\vspace{0.2cm}
\noindent \textbf{(1.2.2)}
Let ${\rm A}$ denote a commutative Banach ring with unit.
V. Berkovich calls \textit{spectrum} of ${\rm A}$, and denotes by $\mathcal{M}({\rm A})$, the set of all bounded multiplicative seminorms on ${\rm A}$; this is a non-empty set if $\A \neq 0$~\cite[1.2]{Ber1}. We adopt the following notation: for an element $x$ of $\mathcal{M}(\A)$ and an element $f$ of $\A$, we write $|f(x)|$ or $|f|(x)$ instead of $x(f)$. Equipped with the weakest topology for which all the functions ${\rm A} \to \mathbb{R}_+, f \mapsto \mid f \mid(x)$ are continuous, $\mathcal{M}(\A)$ is a Hausdorff and compact topological space [loc. cit.].

\vspace{0.1cm} Roughly speaking, this notion of spectrum for Banach commutative rings plays in Berkovich theory a role similar to the one played by the notion of spectrum (set of prime ideals) in algebraic geometry.
For instance, as in commutative algebra, any bounded homomorphism of commutative Banach rings $\varphi : {\rm A} \to {\rm B}$ gives rise to a continuous map $\mathcal{M}({\rm B}) \to \mathcal{M}({\rm A})$, which we denote by ${}^{a}\varphi$.
At last, if $x$ is a point of $\mathcal{M}({\rm A})$, then its kernel $\mathfrak{p}_x = \{ f \in {\rm A} \ ; \  \mid f \mid(x) = 0\}$ is a closed prime ideal in ${\rm A}$.
The completion of the fraction field of ${\rm A}/\mathfrak{p}_x$, with respect to the natural extension of $\mid.\mid (x)$, is called the \textit{complete residue field} of ${\rm A}$ at $x$; it is denoted by $\mathcal{H}(x)$.

Strictly speaking, the building blocks of algebraic geometry are spectra of commutative rings seen as locally ringed spaces, that is spectra endowed with a sheaf of local rings.
Analogously, one has to define a sheaf of local rings on each space $\X = \mathcal{M}({\rm A})$ where ${\rm A}$ is a commutative Banach ring with unit, in order to obtain the so-called $k$-affinoid spaces. However, since the building blocks are compact, one has first of all to single out a class of compact subspaces in $\mathcal{M}(\A)$. Here is a brief summary of their definition, given in several steps~\cite[2.1-2.3]{Ber1}.

\vspace{0.1cm}
First of all, we are interested in spectra of suitable Banach algebras over non-Archimedean complete valued fields: the \textit{affinoid algebras}.
Let us be more precise. We let $\mathbf{BAlg}(k)$ denote the category of Banach $k$-algebras and bounded $k$-homomorphisms.

\begin{Def}
\label{def.affinoid} \begin{itemize}
\item[(i)] A Banach $k$-algebra $\A$ is called \emph{$k$-affinoid} if there exists a strict epimorphism $k\{r_1^{-1}\xi_1, \ldots, r_n^{-1}\xi_n\} \twoheadrightarrow \A$.
It is called \emph{strictly $k$-affinoid} if we can take $r_i=1$ for each index $i$.

\item[(ii)] A Banach $k$-algebra $\A$ is called \emph{affinoid} if there exists a non-Archimedean extension $\K/k$ such that $\A$ is a $\K$-affinoid algebra.
\item[(iii)] We let $k-\mathbf{Aff}$ ($\mathbf{Aff}(k)$, respectivly) denote the full subcategory of $\mathbf{BAlg}(k)$ whose objects are $k$-affinoid algebras (of affinoid $k$-algebras, respectively).
\end{itemize}
\end{Def}

We henceforth fix a $k$-affinoid algebra ${\rm A}$ and consider its (Berkovich) spectrum $\X = \mathcal{M}({\rm A})$.
A \emph{$k$-affinoid domain} in $\X$  [loc. cit., 2.2 p. 27] is by definition a subset $\D$ of $\X$ such that the functor
$$\F_{\D}: \mathbf{Aff}(k) \rightarrow \mathbf{Sets}, \ \ \B \mapsto \{\varphi \in \mathrm{Hom}_{\mathbf{BAlg}(k)}(\A, \B) \ | \ \mathrm{Im}(^{a}\varphi) \subset \D\}$$
is representable by a pair $(\A_{\D}, \varphi_{\D})$ consisting of a $k$-affinoid algebra $\A_{\D}$ and a bounded $k$-homomorphism  $\varphi_{\D}: \A \rightarrow \A_{\D}$.
This pair is then unique up to a unique isomorphism and the morphism $^{a}\varphi_{\D}$ maps $\mathcal{M}(\A_{\D})$ homeomorphically onto $\D$.
The \textit{special subsets} of $\X$ are then defined to be the finite unions of $k$-affinoid domains contained in $\X$ [loc. cit., p. 30]; to such a space $\D$ is naturally attached a Banach $k$-algebra, say ${\rm A}_\D$, which can be computed as follows: if $\{ \D_i \}_{i \in I}$ is a finite covering of $\D$ by $k$-affinoid domains, then ${\rm A}_\D = {\rm Ker}( \prod_{i \in I} {\rm A}_{\D_i} \to \prod_{i,j \in I} {\rm A}_{\D_i \cap \V_j})$ --- Tate's acyclicity theorem implies that the kernel doesn't depend on the choice of the covering.
Thanks to the latter class of closed subsets in $\X$, a sheaf of local rings $\mathcal{O}_\X$ can finally be defined on $\X$ by setting for each open subset $\U$ of $\X$:
$$\Gamma(\mathcal{O}_\X,\U) = \underleftarrow{\lim}_{\D} \, {\rm A}_\V,$$
where the projective limit is taken over all special subsets $\D \subset \U$.
The so-obtained locally ringed spaces $(\X, \mathcal{O}_\X)$ are called \textit{$k$-affinoid spaces} [loc. cit., p. 32]; an \textit{affinoid space over $k$} is a $\K$-affinoid space for some non-Archimedean extension $\K/k$.

\vspace{0.1cm}
In Berkovich theory, the next step is then to define $k$-analytic spaces. Since we will not need this notion in full generality, let us simply mention that a $k$-analytic space is defined by gluing affinoid spaces along affinoid domains, and that the functorial definition of affinoid domains in an affinoid space given above extends to any $k$-analytic space; we refer to~\cite[\S 1]{Ber2} for a detailed exposition. A $k$-analytic space is simultaneously a locally ringed space and a locally ringed site with respect to the Grothendieck topology generated by its affinoid subspaces. One relies on the latter structure to define morphisms. The category $k$-$\mathbf{An}$ of $k$-analytic spaces has finite products and a $k$-analytic group is a group object in the category of $k$-analytic spaces. As for schemes, the underlying set of a $k$-analytic group is not a group in general.

We will need the notion of an analytic space $\X^{\rm an}$ associated with a scheme $\X$ locally of finite type over a non-Archimedean field $k$~\cite[3.4]{Ber1}. As in complex algebraic geometry, $\X^{\mathrm{an}}$ is equipped with a morphism of locally ringed spaces $\rho : \X^{\mathrm{an}} \rightarrow \X$ and $(\X^{\mathrm{an}}, \rho)$ represents the functor $$k-\mathbf{An} \rightarrow \mathbf{Set}, \ \ \Y \mapsto \mathrm{Hom}_{\mathbf{loc.rg.sp.}}(\Y,\X).$$ When $\X$ is affine, the analytic space $\X^{\rm an}$ is described set-theoretically as consisting of the multiplicative seminorms on the coordinate ring $\mathcal{O}(\X)$ of $\X$ extending the absolute value of $k$ and the map $\rho$ sends a seminorm on its kernel --- which is a (closed) prime ideal of $\mathcal{O}(\X)$.

In general, the underlying set of $\X^{\mathrm{an}}$ can be described as the quotient of the set $\mathcal{X} = \bigcup_{\K / k} \X(\K)$ --- where the union is over all non-Archimedean extensions $\K / k$ --- by the equivalence relation which identifies $x' \in \X(\K')$ and $x'' \in \X(\K'')$ whenever there exist embeddings of non-Archimedean fields $\K \to \K '$ and $\K \to \K''$ such that $x'$ and $x''$ come from the same point in $\X(\K)$ [loc. cit., 3.4.2].

\begin{Lemma} \label{lemma.compact} Let $\X$ be an affine algebraic $k$-scheme. Any compact subset of $\X^{{\rm an}}$ is contained in a $k$-affinoid domain.
\end{Lemma}

\vspace{0.1cm}
\noindent
\emph{\textbf{Proof}}.
Write $\X = \mathrm{Spec}(\A)$ and choose a $k$-epimorphism $\pi: k[\xi_1, \ldots, \xi_n] \rightarrow \A$. For a compact subset $\C$ of $\X^{\mathrm{an}}$, there exists a positive real number $r$ such that $\max_{1 \leqslant i \leqslant n} \sup_{\C} |\pi(\xi_i)| \leqslant r$ and, if we set $\B = k\{r^{-1}\xi_1, \ldots, r^{-1}\xi_n\}/\mathrm{ker}(\pi)$, then $\mathcal{M}(\B)$ is a $k$-affinoid domain in $\X^{\mathrm{an}}$ containing $\C$. \hfill $\Box$

\vspace{0.8cm}
\noindent \textbf{(1.2.3)}
Let ${\rm A}$ be a commutative Banach algebra.
Recall that any element $f \in {\rm A}$ has a \textit{spectral radius}~\cite[I, 2.3]{BBKspectral}:
$$\rho(f) = \lim_{n\to\infty} \parallel f^n \parallel^{1\over n} = \inf_{n \in \mathbb{N}} \parallel f^n \parallel^{1\over n}.$$
Then the subset ${\rm A}^\circ = \{ f \in {\rm A} \mid \rho(f) \leqslant 1\}$ is a subring of ${\rm A}$ and ${\rm A}^{\circ\circ} = \{ f \in {\rm A} \mid \rho(f) < 1\}$ is an ideal of ${\rm A}^\circ$; we denote by $\widetilde {\rm A} = {\rm A}^\circ / {\rm A}^{\circ\circ}$ the corresponding quotient ring.
Let $x \in \mathcal{M}({\rm A})$.
Then we have a sequence of bounded ring homomorphisms: ${\rm A} \to {\rm A}/\mathfrak{p}_x \hookrightarrow {\rm Quot}\left({\rm A}/\mathfrak{p}_x\right) \hookrightarrow \mathcal{H}(x)$, where Quot denotes the fraction field of an integral domain, and $\mathcal{H}(x)$ is the complete residue field of $x$ defined above.
This provides a ring homomorphism $\widetilde{\rm A} \to \widetilde{\mathcal{H}(x)}$ whose kernel is a prime ideal since $\widetilde{\mathcal{H}(x)}$ is a field.
We finally obtain a map $\pi : \mathcal{M}({\rm A}) \to {\rm Spec}(\widetilde{\rm A})$ by attaching to $x$ the prime ideal $\{f \in \A^{\circ} \ | \ |f|(x) < 1\}$.
This map is called the \textit{reduction map} of the Banach algebra ${\rm A}$~\cite[2.4]{Ber1}.

A useful notion from spectral theory is that of Shilov boundary of ${\rm A}$: a closed subset, say $\Gamma$, of $\mathcal{M}({\rm A})$ is called a boundary if any element of ${\rm A}$ achieves its maximum in $\Gamma$.
Minimal boundaries exist by Zorn's lemma, and when there is a unique minimal boundary, the latter is called the \textit{Shilov boundary} of {\rm A} and is denoted by $\Gamma({\rm A})$.
In the case when ${\rm A}$ is a $k$-affinoid algebra, the Shilov boundary exists and is a finite subset of $\X = \mathcal{M}({\rm A})$ [loc. cit., 2.4.5] such that $$\rho(f) = \max_{\Gamma({\rm A})}|f|$$ for any $f \in \A$ [loc. cit., 1.3.1]. If $\A$ is strictly $k$-affinoid, then there is a remarkable connection with the reduction map $\pi$ : the Shilov boundary $\Gamma({\rm A})$ is the preimage of the set of generic points (of the irreducible components) of ${\rm Spec}(\widetilde{\rm A})$ [loc. cit., 2.4.4].
We will make crucial use of arguments in this circle of ideas in section \ref{s - realizations}.

\vspace{0.8cm}

\noindent \textbf{(1.2.4)}
Let $\mathcal{A}$ be a finitely presented $k^{\circ}$-algebra whose spectrum we denote $\mathcal{X}$. The \emph{generic fibre} (the \emph{special fibre} respectively) of $\mathcal{X}$ is the $k$-scheme $\X = \mathrm{Spec}(\mathcal{A} \otimes_{k^{\circ}} k)$ (the $\widetilde{k}$-scheme $\mathcal{X}_s = \mathrm{Spec}(\mathcal{A} \otimes_{k^{\circ}} \widetilde{k})$ respectively). The map $$||.||_{\mathcal{A}} : \mathcal{A} \otimes_{k^{\circ}} k \rightarrow \mathbb{R}_{\geqslant 0}, \ \  a \mapsto \inf \{ |\lambda| \ ; \ \lambda \in k^{\times} \ \textrm{and} \ a \in \lambda(\mathcal{A} \otimes 1)\}$$ is a seminorm on $\A \otimes_{k^{\circ}} k$. The Banach algebra $\A$ obtained by completion is a strictly $k$-affinoid algebra whose spectrum we denote by $\mathcal{X}^{{\rm an}}$. This affinoid space is naturally an affinoid domain in $\X^{\mathrm{an}}$ whose points are multiplicative seminorms on $\mathcal{A} \otimes_{k^{\circ}} k$ which are bounded with respect to the seminorm $||.||_{\mathcal{A}}$. Moreover, there is a reduction map $\tau : \mathcal{X}^{\rm an} \rightarrow \mathcal{X}_s$ defined as follows: a point $x$ in $\mathcal{X}^{\rm an}$ gives a sequence of ring homomorphisms $$\mathcal{A} \rightarrow \mathcal{H}(x)^{\circ} \rightarrow \widetilde{\mathcal{H}(x)}$$ whose kernel $\tau(x)$ defines a prime ideal of $\mathcal{A} \otimes_{k^{\circ}} \widetilde{k}$, i.e., a point in $\mathcal{X}_s$.

\vspace{0.1cm} This construction and the one described in 1.2.3 are related as follows. The ring $\A^{\circ}$ of power-bounded elements in the affinoid algebra $\A$ is the integral closure of $\mathcal{A}$~\cite[6.1.2 and 6.3.4]{BGR} and we have a commutative diagram $$\xymatrix{& \mathrm{Spec}(\widetilde{\A}) \ar@{->}[dd] \\ \mathcal{X}^{\rm an} \ar@{->}[ru]^{\pi} \ar@{->}[rd]_{\tau} & \\ & \mathcal{X}_s}$$ in which the vertical arrow is a finite morphism. It follows that, if the scheme $\mathcal{X}$ is integrally closed in its generic fibre --- in particular if $\mathcal{X}$ is smooth --- then $\pi = \tau$.

\vspace{0.1cm} The construction above extends to any scheme $\mathcal{X}$ finitely presented over $k^{\circ}$. One defines a $k$-analytic space $\mathcal{X}^{\rm an}$ by gluing the affinoid spaces $(\mathcal{X}_{|\U})^{\rm an}$ associated with affine open subschemes $\mathcal{X}_{|U}$ of $\mathcal{X}$. This $k$-analytic space is equipped with a reduction map $\tau$ to the special fibre of $\mathcal{X}$. It comes also with a canonical morphism $\mathcal{X}^{\rm an} \rightarrow \X^{\rm an}$, where $\X=\mathcal{X} \otimes_{k^{\circ}} k$ denotes the generic fibre of $\mathcal{X}$ and $\X^{\rm an}$ its analytification (see 1.2.2); if $\mathcal{X}$ is proper, this map is an isomorphism~\cite[$\S$ 5]{Ber4}.

\vspace{0.8cm}
\noindent \textbf{(1.2.5)} Let $\X = \mathcal{M}(\A)$ be a $k$-affinoid space.  For any non-Archimedean extension $\K/k$, the preimage of a $k$-affinoid domain $\D \subset \X$ under the canonical projection ${\rm pr}_{\K/k}: \X_{\K} = \X \widehat{\otimes}_k \K \rightarrow \X$ is a $\K$-affinoid domain in $\X_\K$ since the functor $\F_{{\rm p}_{\K/k}^{-1}(\D)}$ is easily seen to be represented by the pair $(\A_{\D} \widehat{\otimes}_k \K, \varphi_{\D} \widehat{\otimes} \mathrm{id}_{\K})$. The converse assertion holds if the extension $\K/k$ is affinoid, i.e., if $\K$ is a $k$-affinoid algebra.

\begin{Prop} \label{prop.descent.analytic}
Let $\X$ be a $k$-affinoid space and let $\K/k$ be an affinoid extension. A subset $\D$ of $\X$ is a $k$-affinoid domain if and only if the subset ${\rm pr}_{\K/k}^{-1}(\D)$ of $\X_{\K}$ is a $\K$-affinoid domain.
\end{Prop}

\vspace{0.2cm}

This proposition, which gives the key to Berkovich realizations of Bruhat-Tits buildings, is a special case of faithfully flat descent in non-Archimedean analytic geometry.
Since we couldn't find a suitable reference, we provide in the first appendix a complete proof of this result (and of some related technical statements we will need).

\subsection{Bruhat-Tits theory}
\label{ss - Bruhat-Tits}
In this section, we sum up the main facts from Bruhat-Tits theory we need in this paper.
Concerning the hypotheses under which we will be using the theory, we need a weak version of the functoriality of Bruhat-Tits buildings with respect to extensions of the ground field (this is automatically satisfied when the ground field is locally compact).
Thanks to huge non-Archimedean extensions, we note that we can obtain interesting transitivity properties of the corresponding groups of rational points acting on their buildings.

\vspace{0.3cm}

\noindent
\textbf{(1.3.1)}~
We very quickly introduce the main terminology of Bruhat-Tits theory; we refer to ~\cite{RousseauGrenoble}, and particularly to Sect. 10 and 11 therein, for a reasonably detailed introduction to these notions.
The first two parts of this reference also contain a very useful geometric introduction to non-discrete Euclidean buildings.

\vspace{0.1cm}

Let $\G$ be a reductive group over a (by convention, complete) non-Archimedean field $k$.
We choose a maximal split torus $\Sr$ in $\G$ and denote by $\Phi(\G,\Sr)$ the corresponding set of roots. For every root $\alpha \in \Phi(\G,\Sr)$, the root group $\U_\alpha$ is the subgroup of $\G$ whose Lie algebra is the eigenspace associated with the character $\alpha$ or $2\alpha$ in the adjoint representation of $\Sr$ on ${\rm Lie}(\G)$.
Then $({\rm Cent}_\G(\Sr)(k), (\U_\alpha(k))_{\alpha \in \Phi(\Sr,\G)})$ is a \textit{generating root datum}~\cite[6.1.3, example c)]{BT1a}~of type $\Phi(\Sr,\G)$ in the sense of Bruhat-Tits.
This statement sums up a substantial part of Borel-Tits theory on the combinatorial structure of $\G(k)$~\cite{BoTi},~\cite{Borel} --- note that the fact that the field $k$ is valued hasn't been used so far.

\vspace{0.1cm}

We now take into account the ultrametric absolute value of $k$, which we write as $|\cdot| = e^{-\omega(\cdot)}$ for some valuation $\omega : k \to \mathbb{R}$.
Let us denote by ${\rm M}$ the centralizer ${\rm Cent}_\G(\Sr)$ of $\Sr$ in $\G$: it is a reductive group sometimes called the \textit{reductive anisotropic kernel}~of $\G$.
We also denote by $\V(\Sr)$, or simply $\V$ if no confusion is possible, the real vector space $\V(\Sr) = {\rm Hom}_{\mathbb{Ab}}(\X^*(\Sr), \mathbb{R})$ where $\X^*(\Sr)$ is the character group $\mathrm{Hom}_{k-\mathbf{Gr}}(\Sr, \mathbb{G}_{\mathrm{m},k})$.
By restriction from $\M$ to $\Sr$, each character $\chi \in \X^*(\M)$ is a linear form on $\V(\Sr)$ and there exists a homomorphism $\nu : {\rm M}(k) \to \V(\Sr)$ such that for any such $\chi$, we have the equality: $\chi(\nu(z)) = -\omega(\chi(z))$ for any $z \in \Sr(k)$~\cite[Prop. 11.1]{RousseauGrenoble}.

Moreover the homomorphism $\nu$ can be extended to a homomorphism from ${\rm Norm}_\G(\Sr)(k)$ to the group of affine automorphisms of $\V(\Sr)$, in such a way that the image of any $z \in \Sr(k)$ is the translation by $\nu(z) \in \V(\Sr)$ and the linear part of the image of any $n \in {\rm Norm}_\G(\Sr)(k)$ is given by the image of $n$ in the \textit{spherical Weyl group}~$\W^{v}={\rm Norm}_\G(\Sr)(k)/{\rm Cent}_\G(\Sr)(k)$ [loc. cit., Prop. 11.3].
It follows from the axioms of a root datum that for each $\alpha \in \Phi(\Sr,\G)$ and any non-trivial $u \in \U_\alpha(k)$ there exist non-trivial $u', u'' \in
\U_{-\alpha}(k)$ such that $m(u)=u'uu''$ normalizes $\Sr(k)$ and the image of $m(u)$ in $\W^v$ is the reflection associated with $\alpha$.
The group generated by all the so-obtained reflections is called the \textit{affine Weyl group}~of $\G$; we denote it by $\W$.
Finally we denote by $\A(\Sr,k)$ the \textit{apartment}~of $\Sr$, that is the affine space with underlying real vector space $\V(\Sr)$, endowed with the previously defined affine action by the group $\W$.

\vspace{0.1cm}

The main result of Bruhat-Tits theory concerning the combinatorial structure of $\G(k)$ is, under suitable assumptions on $\G$ and $k$, the existence of a \textit{valuation}~of the above root datum in $\G(k)$, in the sense of~\cite[6.2]{BT1a} --- we go back in (1.3.4) to the these assumptions, since we have to make our own (\textit{a priori}~stronger) assumptions for this paper.
Roughly speaking, a valuation is a collection $(\varphi_\alpha)_{\alpha \in \Phi(\Sr,\G)}$ of maps $\U_\alpha \to \mathbb{R} \cup \{ +\infty \}$ which corresponds, in the split case, to the valuation $\omega$ of $k$ when one chooses consistent additive parameterizations of the root groups.
In general, for each root $\alpha$ and each real number $m$, the preimage $\U_{\alpha, m}=\varphi_\alpha^{-1}([m, +\infty])$ is a subgroup of the root group $\U_\alpha(k)$; moreover the family $(\U_{\alpha, m})_{m \in \mathbb{R}}$ is a filtration of $\U_\alpha(k)$; the groups $\U_{\alpha, m}$ satisfy axioms requiring some consistency conditions, as well as a suitable behaviour with respect to commutators and to the above (well-defined) map $u \mapsto m(u)$ given by the root datum axioms.
In this framework, to each point $x \in \A(\Sr,k)$ is attached a well-defined subgroup $\U_{\alpha, x}=\U_{\alpha,-\alpha(x)}$ of $\U_\alpha(k)$ for each root $\alpha$.

\vspace{0.3cm}

\noindent
\textbf{(1.3.2)}~
Assuming the existence of a valuation for the root datum $(\M, (\U_\alpha(k))_{\alpha \in \Phi(\Sr,\G)})$, one attaches to each point $x \in \A(\Sr,k)$ two groups.
The first group is denoted by $\rP_x$, it is by definition generated by ${\rm Ker}(\nu\!\mid_{\M(k)})$ and the groups $\U_{\alpha, x}$ when $\alpha$ varies over $\Phi(\Sr,\G)$; the second group is denoted by
$\widehat{\rP}_x$, it is defined by $\widehat{\rP}_x = \rP_x\cdot \N(k)_x$ where $\N(k)_x$ denotes the stabilizer of $x$ in ${\rm Norm}_\G(\Sr)(k)$.
The \textit{Bruhat-Tits building} of $\G$ over $k$, denoted by $\mathcal{B}(\G,k)$, is defined~\cite[7.4]{BT1a} as the quotient of the product space $\G(k) \times \A(\Sr,k)$ by the equivalence (gluing) relation $\sim$ given by:

\vspace{0.1cm}

\centerline{$(g,x) \sim (h,y)$ if and only if there exists $n \in {\rm Norm}_\G(\Sr)(k)$ such that $y=\nu(n)x$ and $g^{-1}hn \in \widehat{\rP}_x$.}

\vspace{0.1cm}

\noindent
We obtain in this way a space $\mathcal{B}(\G,k)$ on which $\G(k)$ acts naturally; denoting by $[g,x]$ the class of $(g,x)$ for $\sim$, the action is described by $h[g,x] = [hg,x]$.
Each subgroup $\U_{\alpha, m}$ fixes a non-empty half-space of $\A(\Sr,k)$ whose boundary is the fixed-point set in $\A(\Sr,k)$ of the reflection $m(u)$ as above, for some suitable $u \in \U_{\alpha, m}$.

The Euclidean building $\mathcal{B}(\G,k)$ carries a natural non-positively curved metric~\cite[8.2]{BT1a}, which is complete since so is $k$ in the present paper; moreover, the action of any $g \in \G(k)$ is isometric.
The apartments, which are the $\G(k)$-transforms of the space $\A(\Sr,k)$ (embedded in $\mathcal{B}(\G,k)$ by the map $x \mapsto [1_{\G(k)},x]$), can be seen as the maximal subspaces of $\mathcal{B}(\G,k)$ isometric to some Euclidean space.
They are in one-to-one $\G(k)$-equivariant correspondence with the maximal $k$-split tori in $\G$; this follows from the fact that $\G(k)$ acts transitively by conjugation on the maximal $k$-split tori in $\G$ and the fact that ${\rm Norm}_\G(\Sr)(k)$ is exactly the stabilizer of $\A(\Sr,k)$ in $\G(k)$.

\vspace{0.1cm}

One point we would like to emphasize is that, though we are mainly interested in compactifying buildings, we must use Bruhat-Tits theory in full generality.
More precisely, a necessary condition for a metric space to admit a compactification is that the space be locally compact, which corresponds for Bruhat-Tits buildings like $\mathcal{B}(\G,k)$ to the case when $k$ is locally compact; still, the techniques we use lead us to deal with non-discretely valued fields: the geometric counterpart is the use of non-discrete buildings.
The non-discretely valued case is of course covered by~\cite{BT1a} and~\cite{BT1b}, but is less popular --- e.g., not covered by the survey~\cite{TitsCorvallis}.
In this case, Bruhat-Tits buildings are not cellular complexes and facets are filters of subsets in the building.
Still, a useful approach by incidence axioms close to those of buildings of affine Tits systems has been developped in~\cite{RousseauGrenoble}.

Roughly speaking, the model of an apartment is a Euclidean affine space acted upon by a group generated by reflections such that its vector quotient is a finite reflection group.
The \textit{walls}~are the fixed-point sets of the affine reflections and a \textit{half-apartment}~is a half-space bounded by some wall.
This defines also in the underlying vector space a partition into \textit{Weyl}~or \textit{vector chambers} and \textit{vector facets}.
Given an apartment $\A$, a point  $x \in \A$ and a vector facet ${\rm F}^v$, the \textit{facet}~${\rm F}_{x,{\rm F}^v}$ based at $x$ and of direction ${\rm F}^v$ is the filter of subsets of $\A$ which contain the intersection of finitely many half-apartments or walls containing a subset of the form $\Omega \cap (x+{\rm F}^v)$, where $\Omega$ is an open neighborhood of $x$ [loc. cit., \S 5]; an \textit{alcove} is a maximal facet.
With these definitions, a theory a buildings by means of incidence axioms of apartments close to that of discrete buildings~\cite{AbramenkoBrown} can be written~\cite[Part II]{RousseauGrenoble}; note that this is quite different from Tits' axiomatic introduced in~\cite{TitsCome} since for classifications purposes, J. Tits uses there the geometry at infinity (in particular Weyl chambers and spherical buildings at infinity) in order to define and investigate Euclidean (\textit{a priori}~not necessarily Bruhat-Tits) buildings.

With this terminology in mind, we can go back to the group action on $\mathcal{B}(\G,k)$ in order to formulate geometrically some well-known decompositions of $\G(k)$.
Recall that a point $x \in \A(\Sr,k)$ is called \textit{special}~if for any direction of wall there is a wall of $\A$ actually passing through $x$~\cite[1.3.7]{BT1a}.
The \textit{Cartan decomposition}~of $\G(k)$ says that if $x \in \A$ is special, then a fundamental domain for the ${\rm Stab}_W(x)$-action on $\A$ (i.e., a \textit{sector}~of tip $x$) is also a fundamental domain for the $\G(k)$-action on $\mathcal{B}(\G,k)$.
This decomposition implies that in order to describe a compactification, it is enough to describe sufficiently many converging sequences of points in the closure of a given sector.
It also implies that given any point $x \in \A$, the apartment $\A$ contains a complete set of representatives for the $\G(k)$-action on $\mathcal{B}(\G,k)$, that is: $\mathcal{B}(\G,k) = {\rm Stab}_{\G(k)}(x)\A$.
We also have to use \textit{Iwasawa decompositions}~\cite[7.3.1]{BT1a}.
Such a decompositon is associated with a point $x\in \A$ and with vector chambers ${\rm D}, {\rm D}'$ for $\A$.
It decomposes $\G(k)$ as: $\G(k) = \U_{\rm D}^+(k) \cdot {\rm Stab}_{\G(k)}(\A) \cdot  {\rm Stab}_{\G(k)}({\rm F}_{x,{\rm D'}})$, where $\U_{\rm D}^+(k)$ is the unipotent radical of the parabolic subgroup stabilizing the chamber at infinity $\partial_\infty {\rm D}$.

\vspace{0.1cm}

At last, let us say a few words about extended Bruhat-Tits buildings.
We denote by $\X^*(\G)$ the character group $\mathrm{Hom}_{k-\mathbf{Gr}}(\G,\mathbb{G}_{\mathrm{m},k})$.
If $\Z$ denotes the center of $\G$, the restriction homomorphism $\X^*(\G) \rightarrow \X^*(\Z^\mathrm{o})$ has finite kernel and its image has finite index in $ \X^*(\Z^\mathrm{o})$.
Therefore it induces an isomorphism between the dual $\mathbb{R}$-vector spaces.
When $\mathcal{B}(\G,k)$ exists, we let $\mathcal{B}^e(\G,k)$ denote the \emph{extended Bruhat-Tits}~building of $\G(k)$; it is simply the product of the building $\mathcal{B}(\G,k)$ by the real vector space $\V(\G) = \mathrm{Hom}_{\mathbf{Ab}}(\X^*(\G),\mathbb{R})$.
The space $\mathcal{B}^e(\G,k)$ is endowed with the $\G(k)$-action induced by the natural actions on each factor --- see~\cite[1.2]{TitsCorvallis} for the case of an apartment, and~\cite[2.1.15]{RousseauHDR} for the building case.

\vspace{0.4cm}

\noindent
\textbf{(1.3.3)}~
The first result we recall is an easy consequence of measure-theoretic arguments when $k$ is locally compact, i.e., when $\G(k)$ carries a Haar measure --- see e.g., the use of~\cite[2.5.3]{Margulis} in the proof of~\cite[Lemma 12]{GuiRem}.

\begin{Lemma}
\label{Iwahori.Zariski.dense}
For any $x \in \mathcal{B}(\G,k)$, the pointwise stabilizer $\G_x(k)$ is Zariski dense in $\G$.
\end{Lemma}

\vspace{0.1cm}
\noindent
\emph{\textbf{Proof}}.
We pick a maximal split torus $\Sr$ in $\G$ such that $x \in {\rm A}(\Sr,k)$.
We denote by ${\rm M} = {\rm Cent}_\G(\Sr)$ the corresponding Levi factor; its derived subgroup $[{\rm M},{\rm M}]$ is semisimple anisotropic and the bounded group $[{\rm M},{\rm M}](k)$ stabilizes pointwise the apartment ${\rm A}(\Sr,k)$.
Let $c$ be an alcove in ${\rm A}(\Sr,k)$ whose closure contains $x$.
It is enough to show that the Zariski closure $\overline{\G_c(k)}^Z$ of the pointwise stabilizer $\G_c(k)$ of $c$, is equal to $\G$.
First, since ${\rm Norm}_\G(\Sr)(k) = {\rm Stab}_{\G(k)}({\rm A}(\Sr,k))$, we deduce that ${\rm Norm}_\G(\Sr) \cap \G_c(k)$ is equal to the pointwise stabilizer of ${\rm A}(\Sr,k)$, hence contains $[{\rm M},{\rm M}](k)$.
Since ${\rm M}=[{\rm M},{\rm M}] \cdot \Sr$ is a reductive group, we already deduce that we have: ${\rm M} \subset \overline{\G_c(k)}^Z$\cite[18.3]{Borel}.

We pick now a special point in ${\rm A}(\Sr,k)$ and use the corresponding valuation $(\varphi_\alpha)_{ \alpha \in \Phi(\Sr,\G)}$ of the root datum $({\rm M},(\U_\alpha)_{ \alpha \in \Phi(\Sr,\G)})$ of $\G$ with respect to $\Sr$.
Let $\alpha$ be a root in $\Phi(\Sr,\G)$.
The group $\G_c(k)$ contains a suitable subgroup, say $\U_{\alpha,c}$, of the filtration given by $\varphi_\alpha$.
Using for instance the cocharacter of $\Sr$ corresponding to the coroot of $\alpha$, we see that we have $\U_\alpha(k) = \bigcup_{s \in \Sr(k)} s \U_{\alpha,c} s^{-1}$, which proves that $\U_\alpha(k) \subset \overline{\G_c(k)}^Z$ because $\overline{\G_c(k)}^Z$ contains $\Sr$.
Since $\U_\alpha(k)$ is Zariski dense in $\U_\alpha$~\cite[3.20]{BoTi}, we deduce that $\overline{\G_c(k)}^Z$ contains the root group $\U_\alpha$ for each $\alpha \in \Phi(\Sr,\G)$.
This proves our claim, since $\G$ is generated by $\M$ and the root groups $\U_\alpha$ for $\alpha$ varying in $\Phi(\Sr,\G)$.
\hfill$\Box$

\vspace{0.2cm}

Let $\G$ be a split connected reductive group over $k$.
With each special point $x \in \mathcal{B}^e(\G,k)$ is associated a well-defined Demazure group scheme $\mathcal{G}_x$ over $k^{\circ}$ with generic fibre $\G$ such that $\mathcal{G}_x(k^{\circ})$ is the stabilizer of $x$ in $\G(k)$~\cite[4.6.22]{BT1b}.
More generally, for any connected reductive group $\G$ over $k$, Bruhat-Tits theory attaches to each point $x$ of $\mathcal{B}(\G,k)$ some group scheme over $k^{\circ}$ with generic fibre $\G$.
In the present paper, we \emph{never} use these integral group schemes if $\G$ is not split or if $x$ is not special.
We replace these group schemes over $k^\circ$ by affinoid subgroups over $k$ which come from Berkovich theory.
The latter affinoid subgroups are defined thanks to the Demazure group schemes of the special points in the split case, and a faithfully flat descent argument.
To perform this, we have to use huge non-Archimedean extensions of $k$ and some weak functoriality property of $\mathcal{B}(\G,k)$.

\vspace{0.4cm}
\noindent
\textbf{(1.3.4)}~We can now be more precise about our working assumptions. First of all, the use of Berkovich theory implies that we systematically work with complete valued fields. Given a non-Archimedean field $k$ and a reductive group $\G$ over $k$, we need of course the \emph{existence} of a Bruhat-Tits building $\mathcal{B}(\G,k)$ for $\G$, as well as some functoriality with respect to non-Archimedean extensions. Ideally, we would assume full \emph{functoriality} of buildings with respect to non-Archimedean extensions, namely the existence of a functor $\mathcal{B}(\G, \cdot)$ from the category $\mathbf{E}(k)$ of non-Archimedean extensions of $k$ to the category of sets, mapping a field extension $\K'/\K$ to a $\G(\K)$-equivariant injection $\mathcal{B}(\G,\K) \rightarrow \mathcal{B}(\G,\K')$. This hypothesis is too strong to hold in full generality \cite[III.5]{RousseauHDR}. It is however fulfilled if $\G$ is split (easy) or only quasi-split; more generally, functoriality holds if the following two technical conditions are satisfied \cite[V.1.2 and errata]{RousseauHDR}
\begin{itemize}
\item[(i)] $\G$ quasi-splits over a tamely ramified (finite) Galois extension $k'/k$ ;
\item[(ii)] there is a maximal split torus $\T$ of $\G \otimes_k k'$ containing a maximal split torus of $\G$ and whose apartment $\A(\T,k')$ contains a Galois-fixed point.
\end{itemize}
Using functoriality in the quasi-split case, one remarks that both conditions are fulfilled over any non-Archimedean extension $\K$ of $k$ as soon as they are fulfilled over the base field $k$.

Condition (i) holds if the residue field of $k$ is perfect. Condition (ii) holds if the valuation of $k$ is discrete, in which case it follows from the so-called "descente \'etale" of Bruhat-Tits. As a consequence, our strong hypothesis is fulfilled if the non-Archimedean field is \emph{discretely valued with a perfect residue field}; this is in particular the case if $k$ is a \emph{local field}.

\vspace{0.2cm} \noindent \begin{Rk} According to \cite[II.4.14]{RousseauHDR}, condition (ii) holds if and only if the exists a maximal split torus $\Sr$ of $\G$ satisying condition (DE) in \cite[5.1.5]{BT1b}).
\end{Rk}

\vspace{0.2cm} However, a weaker form of functoriality suffices in order to perform our basic construction of affinoid groups in 2.1. It is enough to assume the existence of a functor $\mathcal{B}(\G, \cdot)$ on a full subcategory $\mathbf{E}_0(k)$ of $\mathbf{E}(k)$ which is \emph{cofinal}, i.e., each non-Archimedean extension of $k$ is contained in some extension $\K \in \mathbf{E}_0(k)$. Thanks to functoriality in the split case, this condition holds as soon as the building $\mathcal{B}(\G,k)$ sits inside the Galois-fixed point set of the building $\mathcal{B}(\G,k')$ of $\G$ over some finite Galois extension $k'/k$ splitting $\G$, in which case we can take for $\mathbf{E}_0(k)$ the full subcategory of $\E(k)$ consisting of non-Archimedean extensions of $\K$ containing $k'$.

It turns out that in the cases when the Bruhat-Tits building of $\G$ over $k$ is known to exist, existence follows from ``descending'' the valuation of the root datum of $\G$ over a splitting extension down to $k$. This is always possible when $k$ is a local field, and under much broader hypotheses of \cite[Introduction]{BT1b}, by the famous two-step descent argument of the whole latter article (which, by the way, justifies that we can use the machinery of \cite{BT1a} that we have just summed up). There is also a one-step descent available in G.Rousseau's habilitation \cite{RousseauHDR} and in the more recent papers \cite{Prasad} and \cite{PrasadYu}.

\vspace{0.4cm}

\noindent \textbf{(1.3.5)} We work under the functoriality hypothesis discussed in (1.3.4). For a point $x$ in $\mathcal{B}(\G,k)$ and a non-Archimedean extension $\K/k$, we let $x_{\K}$ denote the image of $x$ in $\mathcal{B}(\G,\K)$.

\begin{Prop} \label{prop.split.special}
Let $x$ be a point in $\mathcal{B}(\G,k)$. There exists an affinoid extension $\K/k$ satisfying the following two conditions:
\begin{itemize}
\item[(i)] the group $\G \otimes_k \K$ is split;
\item[(ii)] the canonical injection $\mathcal{B}(\G,k) \rightarrow \mathcal{B}(\G,\K)$ maps $x$ to a special point.
\end{itemize}
\end{Prop}

\vspace{0.1cm}
\noindent
\emph{\textbf{Proof}}. Let $k'$ be a finite extension of $k$ splitting the group $\G$ and set $x' = x_{k'}$. Pick a split maximal torus $\T$ in $\G' = \G \otimes_k k'$ such that $x'$ lies in the apartment $\A(\T,k')$ and set $\N = \mathrm{Norm}_{\G'}(\T)$; we recall the notation $\V(\T) = {\rm Hom}_{\mathbf{Ab}}(\X^*(\T),\mathbb{R})$. Finally, let $x_0$ be a special point of $\mathcal{B}(\G,k')$ contained in $\A(\T,k')$. Since $(x_0)_{\K}$ is a special point of $\mathcal{B}(\G,\K)$ for any non-Archimedean extension $\K$ of $k'$, the unique affine bijection with identical linear part $\V(\T \otimes_{k'} \K) \rightarrow \A(\T,\K)$ mapping $0$ to $(x_0)_{\K}$ is $\N(\K)$-equivariant. Indeed, the local Weyl group $\W_{x_0}$ coincides with the full spherical Weyl group $\W^v$ and $\N(\K) = \T(\K) \rtimes \W^v$.

The image of the map $$\T(\K) \rightarrow \V(\T \otimes_{k'} \K) = \V(\T), \ \ t \mapsto \left( \chi \mapsto \log |\chi(t)| \right)$$ consists of all linear forms $u$ on the vector space $\X^*(\T) \otimes_{\mathbb{Z}} \mathbb{R}$ such that $\langle u, \chi \rangle$ belongs to the subgroup $\log |\K^{\times}|$ of $\mathbb{R}$ for any character $\chi \in \X^*(\T)$. In the identification above, the point $x'$ of $\A(\T,k')$ corresponds to a linear form $u$ on $\X^*(\T)$ whose image is a finitely generated subgroup of $\mathbb{R}$ since $\X^*(\T)$ is a finitely generated abelian group. Now, if we consider any affinoid extension $\K$ of $k$ containing $k'$ and such that $\mathrm{im}(u) \subset \mathrm{log}|\K^{\times}|$, then $u$ and $0$ belong to the same orbit in $\V(\T)$ under the group $\T(\K)$, hence $x$ and $x_0$ belong to the same orbit under $\N(\K)$.
The point $x$ is therefore a special point of $\mathcal{B}(\G,\K)$.
\hfill$\Box$

\vspace{0.1cm}

\begin{Prop}
\label{prop.transitivity}
For any two points $x$ and $y$ in $\mathcal{B}(\G,k)$, there exist an affinoid extension $\K/k$ and an element $g$ of $\G(\K)$ such that
$g x_{K} = y_K$.
\end{Prop}

\vspace{0.1cm}
\noindent
\emph{\textbf{Proof}}.~
Let $\K/k$ be a non-Archimedean extension splitting $\G$ and such that the point $x_{\K}$ is a special vertex of $\mathcal{B}(\G,\K)$. Pick a split maximal torus $\T$ of $\G \otimes_k \K$ whose apartment in $\mathcal{B}(\G,\K)$ contains both $x_{\K}$ and $y_{\K}$. Taking $x_{\K}$ as a base point, the argument used in the proof above shows that there exists a non-Archimedean extension $\K'/\K$ such that $x_{\K'}$ and $y_{\K'}$ lie in the same orbit under $\G(\K')$.
\hfill$\Box$

\noindent
\textbf{(1.3.5)}~
As in~\cite{GuiRem}, we also need to see the buildings of Levi factors of $\G$ inside $\mathcal{B}(\G,k)$ --- see [loc. cit., 1.4] for further introductory details.
Let $\rP$ a parabolic subgroup of $\G$ containing $\Sr$.
The image $\overline{\Sr}$ of $\Sr$ under the canonical projection $p: \rP \rightarrow \rP_{\mathrm{ss}} = \rP/\rad (\rP)$ is a maximal split torus of the semisimple group $\rP_{\mathrm{ss}}$.
The map $\X^*(\overline{\Sr}) \rightarrow \X^*(\Sr), \ \chi \mapsto \chi \circ p$ is an injective homomorphism and we let $p^{\vee}: \V(\Sr) \rightarrow \V(\overline{\Sr})$ denote the dual projection.

Let $\Lr_{\rP}$ denote the Levi subgroup of $\rP$ containing $\mathrm{Cent}_{\G}(\Sr)$. The projection $p$ induces an isomorphism between the reductive groups $\Lr_{\rP}$ and $\rP/\radu (\rP)$, and between the semisimple groups $\Lr_{\rP}/\rad (\Lr_{\rP})$ and $\rP/\rad (\rP)$.
The set $\Phi(\rP,\Sr)$ of roots of $\rP$ with respect to $\Sr$ is the union of the disjoint closed subsets $\Phi(\Lr_{\rP},\Sr)$ and $\Phi(\radu (\rP),\Sr)$.
The subset $\Phi(\Lr_{\rP},\Sr)$ of $\X^*(\Sr)$ spans $\X^*(\overline{\Sr}) \otimes_{\mathbb{Z}} \mathbb{Q}$, hence the kernel of $p^{\vee}$ is the linear subspace
\begin{eqnarray*} \X^*(\overline{\Sr})^{\bot} & = & \left\{u \in \V(\Sr) \ | \ \forall \alpha \in \X^*(\overline{\Sr}), \ \alpha(u) = 0 \right\} \\ &
= & \left\{u \in \V(\Sr) \ | \ \forall \alpha \in \Phi(\Lr_{\rP}, \Sr), \ \alpha(u) = 0\right\}. \end{eqnarray*}
The following proposition is a particular instance of the results proved in~\cite[Sect. 7.6]{BT1a}.

\begin{Prop}
\label{prop.apartment.levi}
There is a natural affine map $\A(\Sr,k) \rightarrow \A(\overline{\Sr},k)$ with linear part $p^{\vee}$, mapping special points to special points and inducing an isomorphism
$$\A(\Sr,k)/\X^*(\overline{\Sr})^{\bot} \xymatrix{{} \ar@{->}[r]^{\sim} & {}} \A(\overline{\Sr},k).$$
\end{Prop}

Moreover let $k'/k$ be a finite Galois extension splitting $\G$ and $\T$ any maximal torus of $\G$. If we set $\Gamma = \mathrm{Gal}(k'/k)$ and let $\Sr$ denote the maximal split torus of $\T$, the map $\A(\T,k') \rightarrow \A(\overline{\T},k')$ defined as above is $\Gamma$-equivariant and the natural diagram $$\xymatrix{\A(\T,k') \ar@{->}[r] & \A(\overline{\T},k') \\ \A(\Sr,k) \ar@{^{(}->}[u] \ar@{->}[r] & \A(\overline{\Sr},k) \ar@{^{(}->}[u]},$$ in which the vertical maps are the canonical injections, is commutative.

\section{Realizations of buildings} \label{s - realizations}

In this section we define, for a given reductive group $\G$ over a complete non-Archimedean field $k$, various maps from the Bruhat-Tits building $\mathcal{B}(\G,k)$ (or its extended version) to analytic spaces over $k$.
The target spaces are first the Berkovich analytic space $\G^{\rm an}$ associated with $\G$ and then the ones associated with the connected components ${\rm Par}_t(\G)$ of the variety ${\rm Par}(\G)$ of parabolic subgroups of $\G$.

The construction of the fundamental (first) map $\vartheta : \mathcal{B}(\G,k) \to \G^{\rm an}$ relies on the idea to attach to each point $x \in \mathcal{B}(\G,k)$ an affinoid subgroup $\G_x$ such that $\G_x(k)$ is the stabilizer of $x$ in $\G(k)$ (Th. \ref{thm.affinoid.subgroups}). In the split case, this map was defined Berkovich \cite[5.4]{Ber1} in a different way. 
Our construction requires a faithfully flat descent result in the context of Berkovich geometry, which is proved in the first appendix.
The other maps are derived from $\vartheta$.
The analytic space ${\rm Par}_t(\G)^{\rm an}$ attached to the (projective) "flag variety of type $t$" ${\rm Par}_t(\G)$ is compact.
The so-obtained maps $\vartheta_t : \mathcal{B}(\G,k) \to {\rm Par}_t(\G)^{\rm an}$, which only depend on $t$, are used in the next section to define the compactifications.

\vspace{0.2cm} 
We consider a reductive group $\G$ over a non-Archimedean field $k$ and recall that our working hypothesis, detailed in (1.3.4), are fulfilled in particular if $k$ is a local field, or more generally if $\mathcal{B}^e (\G,k)$ is obtained by descent of the ground field from a splitting field down to $k$.

\subsection{Affinoid subgroups associated with points of a building}

\noindent \textbf{(2.1.1)} The fundamental fact underlying Berkovich's point of view on Bruhat-Tits theory is the following result.

\begin{Thm} \label{thm.affinoid.subgroups} Let $x$ be a point in $\mathcal{B}^e (\G,k)$. There exists a unique $k$-affinoid subgroup $\G_x$ of $\G^{{\rm an}}$ such that, for any non-Archimedean extension $\K/k$, we have: $$\G_x(\K) = \mathrm{Stab}_{\G(\K)}(x_{\K}).$$
\end{Thm}

\vspace{0.1cm}
\noindent
\emph{\textbf{Proof}}.
Given a non-Archimedean extension $\K/k$, we say that a $\K$-point $g \in \G(\K)$ of $\G$ is \emph{localized} in the point $z$ of $\G^{\mathrm{an}}$ if $\{z\}$ is the image of the morphism $g : \mathcal{M}(\K) \rightarrow \G^{\mathrm{an}}$.

\vspace{0.1cm} Define $\G_x$ as the subset of $\G^{\mathrm{an}}$ consisting of the points $z$ satisfying the following condition:
\begin{center}\emph{there exist a non-Archimedean extension $\K/k$ and a $\K$-point $g: \mathcal{M}(\K) \rightarrow G^{\rm an}$ of $\G$ localized in $z$ \\ such that $gx_{\K} = x_{\K}$.}
\end{center}

\vspace{0.1cm}
Let $k'/k$ be a non-Archimedean extension, denote by $p$ the canonical projection of $(\G \otimes_{k} k')^{\rm an} = \G^{\rm an} \widehat{\otimes}_k k'$ onto $\G^{\rm an}$ and set $x'=x_{k'}$. We claim that $p^{-1}(\G_x)=(\G \otimes_k k')_{x'}$.

By definition, a point $z$ of $\G^{\rm an}$ belongs to $\G_x$ if and only if there exist a non-Archimedean extension $\K/k$ and a point $g \in \G(\K)$ fixing $x_{\K}$ and sitting in a commutative diagram $$\xymatrix{k[\G] \ar@{->}[rr]^g \ar@{->}[rd]& & \K. \\ & \mathcal{H}(z) \ar@{->}[ru]&}$$
Given $z' \in (\G \otimes_k k')_x'$, such a diagram exists for the extension $\mathcal{H}(z')$ of $\mathcal{H}(z)$ and therefore $z$ belongs to $\G_x$. Conversely, if $z'$ is a point of $(\G \otimes_k k')^{\rm an}$ over $z \in \G_x$, there exists a non-Archimedean extension $\K'$ of $\mathcal{H}(z)$ covering both $\K$ and $\mathcal{H}(z')$; since $x'_{\K'} = x_{\K'}$, the element $g$ of $\G(\K)$ seen in $\G(\K')$ fixes $x'_{\K'}$ and therefore $z'$ belongs to $(\G \otimes_k k')_{x'}$.



\vspace{0.1cm} Let us temporarily assume that the group $\G$ is \emph{split} and that the point $x$ is a \emph{special} point of $\mathcal{B}^e(\G,k)$.
According to Bruhat-Tits theory, there exists a Demazure group scheme $\mathcal{G}_{x}$ over the ring $k^{\circ}$ with generic fibre $\G$ such that $\mathcal{G}_x(k^{\circ})$ is the subgroup of $\G(k)$ fixing the point $x$. More generally, for any non-Archimedean extension $\K/k$, the subgroup $\mathcal{G}_x(\K^{\circ})$ is the stabilizer of the point $x_{\K}$ in $\G(\K)$; indeed, $x_{\K}$ is still a special point of $\mathcal{B}^e(\G,\K)$ and $\mathcal{G}_{x_{\K}} = \mathcal{G}_x \otimes_{k^{\circ}} \K^{\circ}$.

Applying the construction described in 1.2.4, one gets an affinoid subgroup $\mathcal{G}_x^{\rm an}$ of $\G^{\mathrm{an}}$. We have the equality $\mathcal{G}_x^{\rm an}(\K) = \mathcal{G}_x(\K^{\circ})$ in $\G(\K)$. This amounts to saying that, for any non-Archimedean extension $\K/k$, a point $g: \mathcal{M}(\K) \rightarrow \G^{\mathrm{an}}$ is localized in $\mathcal{G}^{\rm an}_x$ if and only if $gx_{\K} = x_{\K}$, hence $\mathcal{G}_x^{\rm an} = \G_x$ as subsets of $\G^{\mathrm{an}}$. Hence $\G_x$ is in this case a $k$-affinoid domain of $\G^{\mathrm{an}}$ and, for any non-Archimedean extension $\K/k$, $\G_x(\K)$ is the subgroup of $\G(\K)$ fixing the point $x_{\K}$.

\vspace{0.1cm} We now remove the two assumptions above. Let $\K/k$ be an affinoid extension splitting $\G$ and such that $x_{\K}$ is a special point of $\mathcal{B}^e(\G, \K)$ (Proposition \ref{prop.split.special}). In view of what has been proved so far, ${\rm pr}_{\K/k}^{-1}(\G_x) = \G_{x_{\K}}$ is a $\K$-affinoid domain in $\G^{\mathrm{an}} \widehat{\otimes}_k \K = (\G \otimes_{k} \K)^{\mathrm{an}}$; in particular, $\G_x = {\rm pr}_{\K/k}(\G_{x_{\K}})$ is a compact subset of $\G^{\mathrm{an}}$. Since any compact subset of $\G^{\mathrm{an}}$ is contained in a $k$-affinoid domain (Lemma \ref{lemma.compact}), we conclude from Proposition \ref{prop.descent.analytic} that $\G_x$ is a $k$-affinoid domain in $\G^{\mathrm{an}}$.

Finally, let $\K/k$ be any non-Archimedean extension and pick an extension $\K'$ of $\K$ splitting $\G$ and such that $x_{\K'}$ is a special point. We have:
\begin{eqnarray*} \G_{x}(\K) & = & \G_x(\K') \cap \G(\K) \\ & = & \mathrm{Stab}_{\G(\K')}(x_{\K'}) \cap \G(\K) \\ & = & \mathrm{Stab}_{\G(\K)}(x_{\K}) \end{eqnarray*}
and $\G_x$ is a $k$-affinoid subgroup of $\G^{\mathrm{an}}$ in view of the next lemma. \hfill $\Box$

\begin{Lemma} \label{lemma.affinoid.subgroups} Let $\X$ be a $k$-analytic group. For any $k$-affinoid domain $\D$ of $\X$, the following conditions are equivalent:
\begin{itemize}
\item $\D$ is a $k$-affinoid subgroup of $\X$;
\item for any non-Archimedean extension $\K$ of $k$, the subset $\D(\K)$ of $\X(\K)$ is a subgroup. \end{itemize}
\end{Lemma}

\vspace{0.1cm}
\noindent
\emph{\textbf{Proof}}.
By definition, a non-empty $k$-affinoid domain $\D$ of $\X$ is a subgroup of $\X$ if the multiplication $m_{\X}: \X \times \X \rightarrow \X$ and the inversion $i_{\X}: \X \rightarrow \X$ factor through the canonical immersions $\iota_{\D}: \D \rightarrow \X$ and $\iota_{\D} \times \iota_{\D}: \D \times \D \rightarrow \X \times \X$: $$\xymatrix{\D \times \D  \ar@{-->}[r] \ar@{->}[d]_{\iota_{\D} \times \iota_{\D}} & \D \ar@{->}[d]^{\iota_{\D}} \\ \X \times \X \ar@{->}[r]_{m_\X} & \X} \ \hspace{2cm} \ \xymatrix{\D  \ar@{-->}[r] \ar@{->}[d]_{\iota_{\D}} & \D \ar@{->}[d]^{\iota_{\D}} \\ \X \ar@{->}[r]_{i_\X} & \X.}$$
Equivalently, both morphisms $m_{\X} \circ \iota_{\D} \times \iota_{\D}$ and $i_{\X} \circ \iota_{\D}$ are required to factor through $\iota_{\D}$. In view of the definition of $k$-affinoid domains in terms of representability of a functor (1.2.2), this is the case if and only if their images lie in the subset $\D$ of $\X$. Since each point of a $k$-analytic space $\Y$ is the image of a morphism $\mathcal{M}(\K) \rightarrow \Y$ for a suitable non-Archimedean extension $\K$ of $k$, the latter condition amounts exactly to saying that $\D(\K)$ is a subgroup of $\G(\K)$ for any non-Archimedean extension $\K/k$. \hfill $\Box$

\vspace{0.1cm} \noindent \begin{Rk} Let $x$ be a point in $\mathcal{B}(\G,k)$. The theorem above has the following consequence: given a non-Archimedean extension $\K$ of $k$ and a $\K$-point $g \in \G(\K)$ fixing $x_{\K}$ in $\mathcal{B}(\G,\K)$, any other $\K$-point $h \in \G(\K)$, inducing the same seminorm as $g$ on the coordinate algebra $k[\G]$ of $\G$, fixes $x_{\K}$.
\end{Rk}

\begin{Prop} \label{prop.theta}
Let $x$ be a point in $\mathcal{B}^e(\G,k)$.
\begin{itemize}
\item[(i)] The $k$-affinoid subgroup $\G_x$ of $\G^{{\rm an}}$ has a unique Shilov boundary point, which we denote by $\vartheta(x)$. It is a norm on the $k$-algebra $k[\G]$ extending the absolute value of $k$.
\item[(ii)] The $k$-affinoid group $\G_x$ is completely determined by the point $\vartheta(x)$: its $k$-affinoid algebra is the completion of $\left(k[\G], |.| (\vartheta(x))\right)$
and we have: $$\G_x = \{z \in \G^{\rm an} \ | \ |f|(z) \leqslant |f|(\vartheta(x)) \ \textrm{for all } f \in k[\G]\}.$$
\item[(iii)] If we let the group $\G(k)$ act on $\G^{{\rm an}}$ by conjugation, then the subgroup of $G(k)$ fixing the point $\vartheta(x)$ is $(\Z\G_x)(k)$, where $\Z = \mathrm{Cent}(\G)$.
\item[(iv)] If we let the group $\G(k)$ act on $\G^{{\rm an}}$ by translation (on the left or on the right), then $\G_x(k)$ is the subgroup of $\G(k)$ fixing the point $\vartheta(x)$.
\end{itemize}
\end{Prop}

\vspace{0.1cm}
\noindent
\emph{\textbf{Proof}}.
(i) Pick a non-Archimedean extension $\K/k$ splitting $\G$ and such that $x_{\K}$ is a special point of $\mathcal{B}^e(\G,\K)$; under these assumptions, the $\K$-affinoid group $\G_{x} \widehat{\otimes}_k \K = (\G \otimes_k \K)_{x_{\K}}$ is deduced from a Demazure group scheme $\mathcal{G}$ over $\K^{\circ}$. Since $\mathcal{G}$ is smooth, the reduction map (1.2.4) $\G_{x_{\K}} \rightarrow \mathcal{G} \otimes_{\K^{\circ}} \widetilde{\K}$ induces a bijection between the Shilov boundary of $\G_{x_{\K}}$ and the set of generic points in the special fibre of $\mathcal{G}$. Since Demazure group schemes are connected (by definition), $\mathcal{G} \otimes_{\K^{\circ}} \widetilde{\K}$ has only one generic point and therefore the Shilov boundary of $\G_{x} \widehat{\otimes}_k \K$ is reduced to a point.  This is \emph{a fortiori} true for $\G_x$ since the canonical projection $\G_x \widehat{\otimes}_k \K \rightarrow \G_x$ maps the Shilov boundary of the range \emph{onto} the Shilov boundary of the target~\cite[proof of Proposition 2.4.4]{Ber1}.

\vspace{0.1cm}
By definition, $\vartheta(x)$ is a multiplicative seminorm on the $k$-algebra $k[\G]$. That $\vartheta(x)$ is in fact a norm can be checked after any non-Archimedean extension $\K/k$ since $$|f|(\vartheta(x_{\K}) = \max_{\G_{x_\K}} |p^{*}(f)| = \max_{\G_x}|f| = |f|(\vartheta(x)),$$ where $p$ denotes the projection of $\G_{x_{\K}} = \G_x \widehat{\otimes}_k \K$ onto $\G_x$. We can therefore assume that $\G$ is split, in which case the conclusion follows most easily from the explicit formula (i) of Proposition \ref{prop.theta.explicit}, whose proof is independant of assertions (ii), (iii), (iv) below.

\vspace{0.1cm}
(ii) Denote by $\A(x)$ the completion of $\left(k[\G],|.|(\vartheta(x))\right)$ and let $\A_x$ be the $k$-affinoid algebra of $\G_x$. Since $\A_x$ is reduced, we may --- and shall --- assume that its norm coincides with its spectral norm \cite[Proposition 2.1.4]{Ber1}, hence with $|.|(\vartheta(x))$ as $\vartheta(x)$ is the only Shilov boundary point of $\G_x$. The immersion $\G_x \rightarrow \G^{\mathrm{an}}$ corresponds to an injective homomorphism of $k$-algebras $k[\G] \rightarrow \A_x$ and thus extends to an isometric embedding $i$ of $\A(x)$ into $\A_x$.

Consider a non-Archimedean extension splitting $\G$ and such that $x_{\K}$ is a special point of $\mathcal{B}(\G,\K)$. We let $\mathcal{G}_{x_{\K}}$ denote the Demazure group scheme over $\K^{\circ}$ attached to $x_{\K}$. By definition, we have an isomorphism of Banach algebras $\A_{x} \widehat{\otimes}_k \K \simeq \A_{x_{\K}}$; moreover, since $\A_{x_{\K}}$ is the completion of $\K[\G]$ with respect to the gauge norm coming from $\K^{\circ}[\mathcal{G}_{x_{\K}}]$ (see 1.2.4), $\K[\G]$ is dense in $\A_{x_{\K}}$. It follows that $i \widehat{\otimes} \K : \A(x) \widehat{\otimes}_k \K \rightarrow \A_x \widehat{\otimes}_k \K$ is an isomorphism of Banach algebras, hence $\A(x) = \A_x$ by descent (Lemma A.5).

\vspace{0.1cm}
(iii) We adapt the argument given by Berkovich in~\cite[Lemma 5.3.2]{Ber1}. Consider $g \in \G(k)$ such that $g \vartheta(x) g^{-1} = \vartheta(x)$. In view of (ii), we have $g\G_x g^{-1} = \G_x$ and thus $g\G_{x_{\K}}g^{-1} = \G_{x_{\K}}$ for any non-Archimedean extension $\K/k$ since $\G_{x_{\K}} = \G_x \widehat{\otimes}_k \K$. Choose such an extension splitting $\G$ and making $x$ a special point. Letting $\mathcal{G}_{x}$ denote the Demazure group scheme over $\K^{\circ}$ associated with $x_{\K}$, our element $g$ of $\G(k)$ induces a $\K$-automorphism $\gamma$ of $\mathcal{G}_{x} \widehat{\otimes}_{\K^{\circ}} \K$. By (ii), the affinoid algebra $\A_x$ of $\G_x$ is the completion of $K[\G]$ with respect to the gauge (semi)norm attached to $\K^{\circ}[\mathcal{G}_x]$. Since the $\K^{\circ}$-scheme $\mathcal{G}_x$ is smooth and connected, this norm is multiplicative hence coincides with the spectral norm $|.|(\vartheta(x))$; moreover, $\K^{\circ}[\mathcal{G}_x]$ is integrally closed in $\K[\G]$ and therefore $\K^{\circ}[\mathcal{G}_x] = \K[\G] \cap \A_x^{\circ}$ (see 1.2.4). By hypothesis, the automorphism ${\rm int}(g)^*$ of $\K[\G]$ is an isometry with respect to the norm $|.|(\vartheta(x))$, hence ${\rm int}(g)^*$ induces a $\K$-automorphism of $\K^{\circ}[\mathcal{G}_x]$, i.e., the automorphism ${\rm int}(g)$ of $\G \otimes_k \K$ induces an automorphism $\gamma$ of $\mathcal{G}_x$. Now we consider the following commutative diagram with exact rows $$\xymatrix{1 \ar@{->}[r] & \mathcal{G}_{x}(\K^{\circ})/\mathcal{Z}(\K^{\circ}) \ar@{->}[r] \ar@{^{(}->}[d] & \mathrm{Aut}(\mathcal{G}_x, \K^{\circ}) \ar@{->}[r] \ar@{^{(}->}[d] & \mathrm{Aut ext}(\mathcal{G}_x, \K^{\circ}) \ar@{->}[r] \ar@{=}[d] & 1 \\ 1 \ar@{->}[r] & \mathcal{G}_{x}(\K)/\mathcal{Z}(\K) \ar@{->}[r] & \mathrm{Aut}(\mathcal{G}_x, \K) \ar@{->}[r] & \mathrm{Aut ext}(\mathcal{G}_x, \K) \ar@{->}[r]  & 1}$$ where $\mathcal{Z}$ is the center of $\mathcal{G}_x$ and ${\rm Autext}$ denotes the group of outer automorphisms ~\cite[Expos\'e XXIV, Sect. 1]{SGA3}.

Since $\gamma$ has trivial image in $\mathrm{Aut ext}(\mathcal{G}_x, \K^{\circ})$, there exists $h \in \mathcal{G}_x (\K^{\circ})$ such that $\gamma = {\rm int}(h)$. It follows that $g = h z$ in $\G(\K)$ with $h \in \mathcal{G}_{x}(\K^{\circ}) = \G_x(\K)$ and $z \in \Z(\K)$, and therefore $g$ is a $k$-point of the group $\G_x\Z$.

\vspace{0.1cm} (iv) Consider $g \in \G(k)$ such that $g \vartheta(x) = \vartheta(x)$. In view of (ii), we have $g \G_x = \G_x$ and thus $g$ belongs to $\G(k) \cap \G_x = \G_x(k)$. \hfill $\Box$

\begin{Cor}
\label{cor.theta.equivariant}
For any $x \in \mathcal{B}^e(\G, k)$ and any $g \in \G(k)$, $$\G_{gx} = g \G_x g^{-1} \ \ \ \textrm{ and } \ \ \ \vartheta(gx) = g\vartheta(x)g^{-1}.$$
\end{Cor}

\vspace{0.1cm}
\noindent
\emph{\textbf{Proof}}.
These two identities are obviously equivalent by Proposition \ref{prop.theta}, (ii), and they follow immediately from the definition of $\G_x$ since $$\mathrm{Stab}_{\G(\K)}(gx) = g \mathrm{Stab}_{\G(\K)}(x)g^{-1}$$ for any non-Archimedean extension $\K/k$. \hfill $\Box$

\vspace{0.8cm}
\noindent \textbf{(2.1.2)} We have attached to each point $x$ of $\mathcal{B}^e(\G,k)$ a canonical (semi-)norm $\vartheta(x)$ on $k[\G]$. If $\G$ is split, we can give an explicit formula for these (semi-)norms.

Choose a maximal split torus $\T$ in $\G$ and let $\Phi = \Phi(\G,\T)$ denote the set of roots of $\G$ with respect to $\T$; choose also a Borel subgroup $\B$ of $\G$ containing $\T$ and let $\Phi^{+}$ denote the corresponding set of positive roots (those occurring in $\radu(\B)$).
Having fixed a total order on $\Phi^{+}$, the canonical map induced by multiplication
$$\prod_{\alpha \in \Phi^{+}} \U_{-\alpha} \times \T \times \prod_{\alpha \in \Phi^{+}} \U_{\alpha} \rightarrow \G$$
is an isomorphism onto an open subset $\Omega(\T,\B)$ of $\G$ which does not depend on the chosen ordering and is called the \emph{big cell} of $\G$ with respect to $(\T,\B)$.

Let $o$ denote a special point in $\mathcal{B}^e(\G,k)$. This point corresponds to a Demazure group scheme $\mathcal{G}$ over $k^{\circ}$ and we also fix a Chevalley basis of $\mathrm{Lie}(\mathcal{G},k^{\circ})$ (i.e., an \emph{integral} Chevalley basis of $\mathrm{Lie}(\G)(k)$~\cite{Steinberg}), which defines a collection of isomorphisms $z_{\alpha}: \mathbb{A}^1_k \tilde{\rightarrow} \U_{\alpha}$ for $\alpha \in \Phi$. We get therefore an isomorphism between the big cell $\Omega(\T,\B)$ and the spectrum of the $k$-algebra $k[\X^*(\T)][(\xi_{\alpha})_{\alpha \in \Phi}]$ of polynomials in the $\xi_{\alpha}$'s with coefficients in the group ring $k[\X^*(\T)]$ (the coordinate ring of $\T$). The open immersion $\Omega(\T, \B) \hookrightarrow \G$ corresponds to a $k$-homomorphism from $k[\G]$ to $k[\X^*(\T)][(\xi_{\alpha})_{\alpha \in \Phi}]$ and the (semi)norms $\vartheta(x)$ are conveniently described on the latter ring.

\begin{Prop} \label{prop.theta.explicit}
We assume that the group $\G$ is split and we use the notation introduced above.

\begin{itemize}
\item[(i)] The point $\vartheta(o)$ belongs to $\Omega(\T,\B)^{{\rm an}}$ and corresponds to the following multiplicative norm:
$$k[\X^*(\T)][(\xi_{\alpha})_{\alpha \in \Phi}] \rightarrow \mathbb{R}_{\geqslant 0}, \ \ \sum_{\chi \in \X^*(\T), \nu \in \mathbb{N}^{\Phi}} a_{\chi, \nu} \chi \ \xi^{\nu} \mapsto \max_{\chi, \nu} |a_{\chi, \nu}|.$$

\item[(ii)] If we use the point $o$ to identify the apartment $\A(\T,k)$ with $\V(\T) = \mathrm{Hom}_{\mathbf{Ab}}(\X^*(\T), \mathbb{R})$, the map $\V(\T) \rightarrow \G^{{\rm an}}$ induced by $\vartheta$ associates with $u \in \V(\T)$ the point of $\Omega(\T,\B)^{{\rm an}}$ defined by the multiplicative norm $$k[\X^*(\T)][(\xi_{\alpha})_{\alpha \in \Phi}] \rightarrow \mathbb{R}_{\geqslant 0}, \ \ \sum_{\chi \in \X^*(\T), \nu \in \mathbb{N}^{\Phi}} a_{\chi, \nu} \chi \ \xi^{\nu} \mapsto \max_{\chi, \nu} |a_{\chi, \nu}| \prod_{\alpha \in \Phi} \mathrm{e}^{\nu(\alpha) \langle u, \alpha \rangle}.$$
\end{itemize}
\end{Prop}

\vspace{0.1cm}
\noindent
\emph{\textbf{Proof}}.
(i) The chosen Chevalley basis $\mathbf{z} = (z_{\alpha})_{\alpha \in \Phi}$ of $\mathrm{Lie}(\mathcal{G}, k^{\circ})$ provides us with an integral model $\mathcal{B}_{\mathbf{z}}$ of $\B$ --- namely, the $k^{\circ}$-subgroup scheme of $\mathcal{G}$ generated by the $k^{\circ}$-torus $\mathcal{T} = \mathrm{Spec}(k^{\circ}[\X^*(\T)])$ and the unipotent $k^{\circ}$-groups $z_{\alpha}(\mathbb{A}^1_{k^{\circ}})$, $\alpha \in \Phi$ --- and the isomorphism $$\prod_{\alpha \in \Phi^{+}} \mathbb{A}^{1}_{k} \times \T \times \prod_{\alpha \in \Phi^{+}} \mathbb{A}^{1}_k \tilde{\rightarrow} \Omega(\T,\B) \subset \G$$ comes from a $k^{\circ}$-isomorphism $$\prod_{\alpha \in \Phi^{+}} \mathbb{A}^{1}_{k^{\circ}} \times \mathcal{T} \times \prod_{\alpha \in \Phi^{+}} \mathbb{A}^{1}_{k^{\circ}} \tilde{\rightarrow} \Omega(\mathcal{T}, \mathcal{B}_{z}) \subset \mathcal{G}$$ onto the big cell of $\mathcal{G}$ corresponding to $\mathcal{T}$ and $\mathcal{B}_{\mathbf{z}}$.

By definition, $\vartheta(o)$ is the unique point of $\G^{\mathrm{an}}$ contained in the affinoid domain $\G_o$ which the reduction map $\G_o = \mathcal{G}^{\rm an} \rightarrow \mathcal{G} \otimes_{k^{\circ}} \widetilde{k}$ (1.2.3) sends to the generic point of the target.
Since $\Omega(\mathcal{T}, \mathcal{B}_{\mathbf{z}})$ is an open subscheme of $\mathcal{G}$ meeting the special fibre of $\mathcal{G}$, the special fibre of $\Omega(\mathcal{T}, \mathcal{B}_{\mathbf{z}})$ contains the generic point of $\mathcal{G} \otimes_{k^{\circ}} \widetilde{k}$, and the affinoid space  $\Omega(\mathcal{T}, \mathcal{B}_{\mathbf{z}})^{\rm an}$ sits inside $\G_o$.
Therefore $\vartheta(o)$ is the unique point in $\Omega(\mathcal{T}, \mathcal{B}_{\mathbf{z}})^{\rm an}$ which reduces to the generic point of $\Omega(\mathcal{T}, \mathcal{B}_{\mathbf{z}}) \otimes_{k^{\circ}} \widetilde{k}$.
This means concretely that $\vartheta(o)$ is characterized by the following two conditions: for any $f \in k[\X^*(\T)][(\xi_{\alpha})_{\alpha \in \Phi}]$,
$$|f|(\vartheta(o)) \leqslant 1 \Longleftrightarrow f \in k^{\circ}[\X^*(\T)][(\xi_{\alpha})_{\alpha \in \Phi}] \ \ \textrm{ and } \ \ |f|(\vartheta(o)) < 1 \Longleftrightarrow f \textrm{ maps to $0$ in } \widetilde{k}[\X^*(\T)][(\xi_{\alpha})_{\alpha \in \Phi}].$$
From this, we immediately conclude that $$|f|(\vartheta(o)) = \max_{\chi, \nu} |a_{\chi, \nu}|$$ if $f = \sum_{\chi, \nu} a_{\chi, \nu} \chi \ \xi^{\nu}.$

\vspace{0.1cm}
(ii) For any $t \in \T(k)$ and any root $\alpha \in \Phi$, the element $t$ normalizes the root group $\U_{\alpha}$ and conjugation by $t$ induces an automorphism of $\U_{\alpha}$ which is just the homothety of ratio $\alpha(t) \in k^{\times}$ if we read it through the isomorphism $z_{\alpha}: \mathbb{A}^1_k \rightarrow \U_{\alpha}$.
We thus have a commutative diagram $$\xymatrix{\mathrm{Spec}(k[\X^*(\T)][(\xi_{\alpha})_{\alpha \in \Phi}]) \ar@{->}[r]^{\hspace{1cm} \sim} \ar@{->}[d]_{\tau} & \Omega(\T, \B) \ar@{->}[d]^{\mathrm{int}(t)} \\ \mathrm{Spec}(k[\X^*(\T)][(\xi_{\alpha})_{\alpha \in \Phi}]) \ar@{->}[r]_{\hspace{1cm} \sim}  & \Omega(\T, \B) ,}$$ where $\tau$ is induced by the $k[\X^*(\T)]$-automorphism $\tau^{*}$ of $k[\X^*(\T)][(\xi_{\alpha})_{\alpha \in \Phi}]$ mapping $\xi_{\alpha}$ to $\alpha(t)\xi_{\alpha}$ for any $\alpha \in \Phi$. It follows that $\vartheta(to) = t \vartheta(o)t^{-1}$ is the point of $\G^{\mathrm{an}}$ defined by the multiplicative norm on $k[\T][(\xi_{\alpha})_{\alpha \in \Phi}]$ mapping an element $f = \sum_{\chi, \nu} a_{\chi, \nu} \chi \ \xi^{\nu}$ to \begin{eqnarray*} |\tau^{*}(f)|(\vartheta(o)) & = & \left|\sum_{\chi, \nu} \left(a_{\chi, \nu} \prod_{\alpha \in \Phi} \alpha(t)^{\nu(\alpha)}\right) \chi \ \xi^{\nu} \right|(\vartheta(o)) \\ & = &  \max_{\chi, \nu} |a_{\chi, \nu}| \prod_{\alpha \in \Phi} |\alpha(t)|^{\nu(\alpha)} \\ & = & \max_{\chi, \nu} |a_{\chi,\nu}| \prod_{\alpha \in \Phi} \mathrm{e}^{\langle \alpha, \log|t| \rangle \nu(\alpha)}\end{eqnarray*} since, by definition, the group $\T(k)$ is mapped to $\V(\T)$ by sending $t$ to the linear form $\chi \mapsto \log|\chi(t)|$ on $\X^*(\T)$ (cf. 1.3.1).

To complete the proof, note that $\vartheta(x) = {\rm pr}_{\K/k}(\vartheta(x_{\K}))$ --- i.e., $\vartheta(x)$ is the restriction of $\vartheta(x_{\K})$ to $k[\T]$ --- for any non-Archimedean extension $\K/k$, and that any point $u$ of $\V(\T)$ belongs to the image of the map $\log |\cdot | : \T(\K) \rightarrow \V(\T \otimes_k \K) = \V(\T)$ for a suitable choice of $\K$.
\hfill $\Box$

\vspace{0.2cm} \begin{Rk} \label{rk.big-cell.charts} Let $\G$ be a connected Chevalley group scheme over a scheme $\Sr$ and let $\T$ denote a split maximal torus of $\G$. With each Borel subgroup $\B$ of $\G$ containing $\T$ is associated the open affine subscheme $\Omega(\T,\G)$ of $\G$ (the \emph{big cell}), whose definition is compatible with base change. When $\B$ runs over the set of Borel subgroups of $\G$ containing $\T$, the big cells $\Omega(\T,\G)$ cover the scheme ${\rm Bor}(\G)$. To prove this assertion, it is enough to check that those open subschemes cover each fiber $\G_{\overline{s}}$ over a geometric point $\overline{s}$ of $\Sr$. We are thus reduced to the case where $\Sr$ is the spectrum of an algebraic closed field, and then the conclusion follows easily from the refined Bruhat decomposition \cite[Expos\'e 13, Th\'eor\`eme 3 and Corollaire 1]{Bible}.

More generally, the same conclusion holds for the family of parabolic subgroups of a given type $t$ containing $\T$. Indeed, the morphism $\pi : {\rm Bor}(\G) \rightarrow {\rm Par}_t(\G)$ defined functor-theoretically by mapping each Borel subgroup to the unique parabolic subgroup of type $t$ containing it is surjective and maps each big cell $\Omega(\T,\B)$ onto the corresponding big cell of $\Omega(\T,\pi(\B))$.
\end{Rk}

\subsection{The canonical map $\vartheta: \mathcal{B}^e(\G, k) \rightarrow \G^{\mathrm{an}}$}


\begin{Prop} \label{prop.theta.properties} The map $\vartheta: \mathcal{B}^e(\G,k) \rightarrow \G^{\rm an}$ defined in Proposition \ref{prop.theta} enjoys the following properties.

\begin{itemize}
\item[(i)] This map is $\G(k)$-equivariant if we let the group $\G(k)$ act on $\G^{\rm an}$ by conjugation.

\item[(ii)] For any non-Archimedean extension $k'/k$, the natural diagram $$\xymatrix{\mathcal{B}^e(\G,k') \ar@{->}[r]^{\vartheta} & (\G \otimes_k k')^{{\rm an}} \ar@{->}[d]^{{\rm pr}_{k'/k}} \\ \mathcal{B}^e(\G,k) \ar@{->}[r]_{\vartheta} \ar@{->}[u] & \G^{{\rm an}}}$$ is commutative. Moreover, if $k'/k$ is a Galois extension, the upper arrow is ${\rm Gal}(k'/k)$-equivariant.

\item[(iii)] The map $\vartheta$ factors through the projection $\mathcal{B}^e(\G,k) \rightarrow \mathcal{B}(\G,k)$ and induces a continuous injection of $\mathcal{B}(\G,k)$ into $\G^{{\rm an}}$. Its restriction to any apartment of $\mathcal{B}(\G, k)$ is a homeomorphism onto a closed subspace of $\G^{{\rm an}}$.
If the field $k$ is locally compact, $\vartheta$ induces a homeomorphism between $\mathcal{B}(\G,k)$ and a closed subspace of $\G^{{\rm an}}$.
\end{itemize}
\end{Prop}

\vspace{0.1cm}
\noindent
\emph{\textbf{Proof}}.
(i) This assertion is Corollary \ref{cor.theta.equivariant}.

\vspace{0.1cm}
(ii) The first assertion follows immediately from the identity $\G_{x_{k'}} = \G_x \widehat{\otimes}_k k'$.

If $k'/k$ is a Galois extension, there is a natural action of the group $\Gamma = {\rm Gal}(k'/k)$ on $(\G \otimes_k k')^{\rm an}$ --- to an element $\gamma$ of $\Gamma$ corresponds the $k$-automorphism of $(\G \otimes_k k')^{\rm an}$ defined by ${\rm id} \otimes \gamma^{-1}$ at the level of the coordinate ring $k'[\G] = k[\G] \otimes_k k'$ --- and Galois-equivariance of $\vartheta$ amounts to the identity $(\G \otimes_k k')_{\gamma(x)} = \gamma(\G \otimes_k k')_x$ in $(\G \otimes_k k')^{\rm an}$. If $\iota : k' \rightarrow \K$ is any non-Archimedean extension, then $\gamma(\G \otimes_k k')_x(\K)$ consists by definition of elements $g$ in $\G(\K)$ which fix the image of $x$ in $\mathcal{B}(\G,\K)$ if we use the extension $u \circ \gamma$ to embed $\mathcal{B}(\G,k')$ into $\mathcal{B}(\G,\K)$, i.e., if we compose the embedding coming from $u$ with the automorphism of $\mathcal{B}(\G,k')$ induced by $\gamma$. Thus we have $\gamma(\G \otimes_k k')_x(\K) = (\G \otimes_k k')_{\gamma(x)}$, and therefore $\gamma(\G \otimes_k k')_x = (\G \otimes_k k')_{\gamma(x)}$.

\vspace{0.1cm}
(iii) For any two points $x, y \in \mathcal{B}^e(\G,k)$ such that $\vartheta(x) = \vartheta(y)$, $\G_x = \G_y$ and thus $\vartheta(x_{\K}) = \vartheta(y_{\K})$ for any non-Archimedean extension $\K/k$. By Proposition \ref{prop.transitivity}, we can choose $\K$ such that $y_{\K} = g x_{\K}$ for some $g \in \G(\K)$. Then $\vartheta(x_{\K}) = \vartheta(gx_{\K}) = g\vartheta(x_{\K})g^{-1}$ and thus $g \in (\Z\G_x)(\K)$ by Proposition \ref{prop.theta} (iii). Enlarging $\K$ if necessary, we may assume that $g$ belongs to $\Z(\K)\G_x(\K) = \G_x(\K)\Z(\K)$. Since $\G_x(\K) = \mathrm{Stab}_{\G(\K)}(x_{\K})$, $y_{\K} = g x_{\K} \in \Z(\K)x_{\K}$ and the points $x_{\K}$ and $y_{\K}$ belong therefore to the same fibre of $\mathcal{B}^e(\G,\K) \rightarrow \mathcal{B}(\G,\K)$. Relying on (ii), we have thus proved that the map $\vartheta$ factors through an injection $\mathcal{B}(\G,k) \rightarrow \G^{\mathrm{an}}$.

Given a maximal split torus $\Sr$ of $\G$, the continuity of $\vartheta$ is equivalent to the continuity of the map $$\G(k) \times \A(\Sr, k) \rightarrow \G^{\mathrm{an}}, \ \ (g,x) \mapsto g \vartheta(x)g^{-1}$$ with respect to the natural topology on the left hand side, for $\mathcal{B}(\G,k)$ is a topological quotient of $\G(k) \times \A(\Sr,k)$. Since the canonical map $\G(k) \times \G^{\mathrm{an}} \times \G(k) \rightarrow \G^{\mathrm{an}}, \ (g,x,h) \mapsto gxh$ is continuous, it remains to prove that the restriction of $\vartheta$ to the apartment $\A(\Sr, k)$ is continuous. In view of the previous assertion, there is no loss of generality in assuming that $\G$ is split, in which case the result is an obvious consequence of the explicit formula given in Proposition \ref{prop.theta.explicit}.

Relying again on this explicit formula, one now proves that $\vartheta(\A(\Sr,k))$ is closed in $\G^{\mathrm{an}}$. In view of (ii), it suffices to consider the case of a split group since the projection $(\G \otimes_k \K)^{\mathrm{an}} \rightarrow \G^{\mathrm{an}}$ is closed for any non-Archimedean extension $\K/k$; moreover, since $\vartheta$ factors through $\mathcal{B}(\G, k)$, we can also assume that $\G$ is semisimple. Consider now a sequence $(u_n)$ of points in $\V(\Sr) = {\rm Hom}_{\mathbf{Ab}}(\X^*(\T), \mathbb{R})$ such that the sequence $(\vartheta(u_n))$ converges in $\G^{\mathrm{an}}$. Then $e^{- \langle u_n, \alpha \rangle} = |\xi_{\alpha}|(\vartheta(u_n))$ converges in $\mathbb{R}_{\geqslant 0}$ for any root $\alpha \in \Phi$ and the limit belongs to $\mathbb{R}_{>0}$ since $|\xi_{\alpha}|(\vartheta(u_n))|\xi_{-\alpha}|(\vartheta(u_n))=  \mathrm{e}^{-\langle u_n, \alpha \rangle} \mathrm{e}^{\langle u_n, \alpha \rangle} = 1$ for any $n$. Thus we get  a map $u_{\infty}: \Phi \rightarrow \mathbb{R}$ which is obviously additive and extends therefore to a linear form on $\X^*(\T)$ since $\Phi$ spans the vector space $\X^*(\T) \otimes_{\mathbb{Z}} \mathbb{Q}$ (recall that we assumed $\G$ semisimple). The point $u_{\infty}$ is mapped to the limit of the sequence $(\vartheta(u_n))$ and thus $\vartheta(\A(\T,k))$ is closed.

We have proved that $\vartheta$ maps homeomorphically each apartment of $\mathcal{B}(\G,k)$ onto a closed subset of $\G^{\mathrm{an}}$. When the field $k$ is locally compact, this is true for the whole building $\mathcal{B}(\G,k)$. Indeed, if $\Sr$ is a maximal split torus of $\G$ and $x$ is a point in $\A(\Sr, k)$, the group $\G_x(k)$ is compact and $\mathcal{B}(\G,k)) = \G_x(k)\A(\Sr,k)$ (see (1.3.2)), hence $$\vartheta(\mathcal{B}(\G,k) = \G_x(k) \cdot \vartheta(\A(\Sr,k))$$ is a closed subspace of $\G^{\rm an}$ (the action of $\G(k)$ is by conjugation).
\hfill $\Box$

\begin{Rk} If $\G$ is split, the map $\vartheta$ introduced above coincides with the one defined by Berkovich in~\cite[5.4.4]{Ber1}. Indeed, with the notation of [loc. cit], $\mathbf{P} = \G_o$ where $o$ is the special point of $\mathcal{B}(\G,k)$ corresponding to the $k^{\circ}$-Demazure group scheme $\mathcal{G} \otimes_{\mathbb{Z}} k^{\circ}$ and, for any  $\lambda \in \mathrm{Hom}_{\mathbf{Ab}}(\X^*(\T), \mathbb{R}_{>0})$, $\mathbf{P}_{\lambda} = \G_{o - \log(\lambda)}$.

\end{Rk}

\subsection{The canonical map $\Theta: \mathcal{B}^e(\G,k) \times \mathcal{B}^e(\G,k) \rightarrow \G^{\mathrm{an}}$}

Given a point $x$ of $\mathcal{B}^e(\G,k)$ and a non-Archimedean extension $\K/k$, we always write $x$ instead of $x_{\K}$ in what follows.

\vspace{0.1cm}
The canonical map $\vartheta: \mathcal{B}^e(\G,k) \rightarrow \G^{\mathrm{an}}$ which we have defined above is equivariant with respect to the natural action of $\G(k)$ on $\G^{\mathrm{an}}$ by conjugation and therefore factors through the projection of $\mathcal{B}^e(\G,k)$ onto $\mathcal{B}(\G,k)$. It is in fact possible to embed equivariantly the whole extended building $\mathcal{B}^e(\G,k)$ into $\G^{\mathrm{an}}$ if we let the group $\G(k)$ act on $\G^{\mathrm{an}}$ by left translations. To be precise, we will use a canonical map $\Theta: \mathcal{B}^e(\G,k) \times \mathcal{B}^e(\G,k) \rightarrow \G^{\mathrm{an}}$ satisfying the following two conditions: for any point $o \in \mathcal{B}^e(\G,k)$, $\Theta(o,o) = \vartheta(o)$ and $\Theta(o,.): \mathcal{B}^e(\G,k) \rightarrow \G^{\mathrm{an}}$ is a $\G(k)$-equivariant embedding of $\mathcal{B}^e(\G,k)$ into $\G^{\mathrm{an}}$.

\vspace{0.1cm}
For any two points $x, y \in \mathcal{B}^e(\G,k)$, there exists by Proposition \ref{prop.transitivity} a
non-Archimedean extension $\K/k$ and an element $g \in \G(\K)$ such
that $y = gx$ in $\mathcal{B}^e(\G,\K)$.

One easily checks that the point $\mathrm{pr}_{\K/k}(g\vartheta(x))$
in $\G^{\an}$ depends neither on $\K$ nor on $g$. Indeed, if
$\K'$ is a non-Archimedean field extending $\K$ and if $g'$ is an
element of $\G(\K')$ such that $y=g'x$ in $\mathcal{B}^e(\G, \K')$, then
$g^{-1}g'x=x$ hence $g^{-1}g' \in \G_x(\K')$. Since $\G_x(\K')$ is
the subgroup of $\G(\K')$ fixing $\vartheta(x)$ in the natural
action of $\G(\K')$ on $(\G \otimes_k \K')^{\an}$ by left
translations (Proposition \ref{prop.theta} (iv)),
we have $g^{-1}g'\vartheta(x) = \vartheta(x)$ hence $g\vartheta(x) =
g'\vartheta(x)$ in $(\G \otimes_k \K')^{\an}$ and \begin{eqnarray*}
\mathrm{pr}_{\K'/k}(g'\vartheta(x)) & = &
\mathrm{pr}_{\K'/k}(g\vartheta(x)) \\ & = & \mathrm{pr}_{\K/k}
\mathrm{pr}_{\K'/\K}(g\vartheta(x)) \\ & = &
\mathrm{pr}_{\K/k}(g\mathrm{pr}_{\K'/\K}(\vartheta(x))) =
\mathrm{pr}_{\K/k}(g\vartheta(x)),\end{eqnarray*} for $\vartheta(x) =
\mathrm{pr}_{\K/k}(\vartheta(x))$.

\begin{Def} \label{def.Theta} For any two points $x, y$ in $\mathcal{B}^e(\G,k)$, we put $$\Theta(x,y) = {\rm pr}_{\K/k}(g\vartheta(x)),$$ where $g \in \G(\K)$ is such that $y = gx$ in $\mathcal{B}^e(\G,\K)$.
\end{Def}

\begin{Lemma} \label{lemma.Theta} Pick some points $x, x', y$ and $y'$ in $\mathcal{B}^e(\G,k)$. If
$\Theta(x,y) = \Theta(x',y')$ in $\G^{{\rm an}}$, then $\Theta(x,y) =
\Theta(x',y')$ in $(\G \otimes_k \K)^{{\rm an}}$ for any
non-Archimedean extension $\K/k$.
\end{Lemma}

\vspace{0.1cm}
\noindent
\emph{\textbf{Proof}}.
Let $x$ be a point in $\mathcal{B}(\G,k)$, $\K/k$ a non-Archimedean extension and $g$ a $\K$-point of $\G$. Given a point $z$ in $(\G \otimes_k \K)^{\rm an}$ whose image under the projection ${\rm pr}_{\K/k} : (\G \otimes_k \K)^{\rm an} \rightarrow \G^{\rm an}$ belongs to ${\rm pr}_{\K/k}(g\G_x)$, there exist a non-Archimedean extension $\K'/\K$ and a $\K'$-point $h$ of $\G$ localized in $z$ such that $h \in g\G_x(\K')$. The $\K'$-point $g^{-1}h$ is localized in $\G_x$, hence $g^{-1}z \in \G_x$ and $z \in g \G_x$. Therefore, $g\G_x = {\rm pr}_{\K/k}^{-1}({\rm pr}_{\K/k}(g\G_x))$ and, since $g\vartheta(x)$ is the only Shilov boundary point of $g\G_x$ by Proposition \ref{prop.theta} (i), $g\vartheta(x)$ is the only maximal point in ${\rm pr}_{\K/k}^{-1}({\rm pr}_{\K/k}(g\vartheta(x))) \subset g\G_x$, i.e., the only point at which each function $f \in \K[\G]$ reaches its supremum over ${\rm pr}_{\K/k}^{-1}({\rm pr}_{\K/k}(g\vartheta(x)))$.

We have thus proved that the point $g\vartheta(x)$ of $(\G \otimes_k \K)^{\rm an}$ is completely characterized by its image in $\G^{\rm an}$. The same argument applies more generally to $g\vartheta(x)h$ for any $g, h \in \G(\K)$.

\vspace{0.1cm} Consider now some points $x, y, x'$ and $y'$ in $\mathcal{B}(\G,k)$ such that $\Theta(x,y) = \Theta(x',y')$ and pick a non-Archimedean extension $\K/k$ such that $y = gx$, $x' = hx$ and $y' = jx$ with $g,h,j \in \G(\K)$. By definition, $\Theta(x,y)$ and $\Theta(x',y')$ are the images of $g\vartheta(x)$ and $$jh^{-1}\vartheta(x') = jh^{-1}(h\vartheta(x)h^{-1}) = j\vartheta(x)h^{-1}$$ respectively in $\G^{\rm an}$. Since those points are completely characterized by their images in $\G^{\rm an}$, we have $g\vartheta(x) = j\vartheta(x)h^{-1}$ and therefore the identity $\Theta(x,y) = \Theta(x',y')$ holds after any non-Archimedean extension of $k$. \hfill $\Box$

\vspace{0.2cm}
\begin{Prop} \label{prop.Theta} The map $\Theta: \mathcal{B}^e(\G,k) \times \mathcal{B}^e(\G,k)
\rightarrow \G^{{\rm an}}$ which we have just defined satisfies the
following properties.
\begin{itemize}
\item[(i)] For any points $x, y \in \mathcal{B}^e(\G,k)$ and any elements $g,h \in
\G(k)$, $$\Theta(gx,hy) = h\Theta(x,y)g^{-1}.$$
\item[(ii)] For any non-Archimedean extension $k'/k$, the natural diagram
$$\xymatrix{\mathcal{B}^e(\G,k') \times \mathcal{B}^e(\G,k') \ar@{->}[r]^{\hspace{1cm} \Theta} &
(\G \otimes_k k')^{{\rm an}} \ar@{->}[d]^{{\rm pr}_{k'/k}} \\
\mathcal{B}^e(\G,k) \times \mathcal{B}^e(\G,k) \ar@{->}[u] \ar@{->}[r]_{\hspace{1cm}
\Theta} & \G^{{\rm an}}}$$ is commutative.
\item[(iii)] Let the group $\G(k)$ act by left translations on $\G^{{\rm an}}$.
For any point $o$ in $\mathcal{B}^e(\G,k)$, the map $\Theta(o, \cdot)$
is a continuous and $\G(k)$-equivariant injection of $\mathcal{B}^e(\G,k)$
into $\G^{{\rm an}}$ which sends homeomorphically each apartment of
$\mathcal{B}^e(\G,k)$ onto a closed subset of $\G^{{\rm an}}$.

If the field $k$ is locally compact, the map $\Theta(o,\cdot)$ induces
a $\G(k)$-equivariant homeomorphism between $\mathcal{B}^e(\G,k)$ and a closed
subspace of $\G^{{\rm an}}$.
\end{itemize}
\end{Prop}

\vspace{0.1cm}
\noindent \emph{\textbf{Proof}}. (i) Consider a non-Archimedean extension $\K/k$ such that $y = jx$ for some $j \in \G(\K)$. We have $hy = hjg^{-1}gx$, hence $$\Theta(gx,hy) = {\rm pr}_{\K/k}(hjg^{-1}\vartheta(gx)) = {\rm pr}_{\K/k}(hj\vartheta(x)g^{-1})$$ and therefore $\Theta(gx,hy) = h\Theta(x,y)g^{-1}$ since the projection ${\rm pr}_{\K/k}$ is $\G(k)$-equivariant.

\vspace{0.1cm} (ii) This assertion follows immediately from the definition of $\Theta$.

\vspace{0.1cm} (iii) The map $\Theta(o, \cdot)$ is $\G(k)$-equivariant by (i). If $x$ and $y$ are two points of $\mathcal{B}^e(\G,k)$ such that $\Theta(o,x) = \Theta(o,y)$, the same equality holds after any non-Archimedean extension of $k$ by Lemma \ref{lemma.Theta}. Therefore, we may assume $x = go$ and $y = ho$ for some $g, h \in \G(k)$. It follows that $g\vartheta(o) = h\vartheta(o)$, hence $h^{-1}g$ belongs to $\G_o(k)$ by Proposition \ref{prop.theta.explicit} (iv) and $x=y$. Thus the map $\Theta(o,.)$ is injective.

In order to establish the continuity of $\Theta(o,\cdot)$, one may restrict to an apartment $\A(\Sr,k)$ containing $o$ since this map is equivariant and $\mathcal{B}^e(\G,k)$ is a topological quotient of $\G(k) \times \A(\Sr, k)$. We may also assume that $\G$ is split by (ii). Then, if $\N$ denotes the normalizer of $\Sr$ in $\G$ and if $\K/k$ is a (huge) non-Archimedean extension such that $|\K| = \mathbb{R}_{\geqslant 0}$, the group $\N(\K)$ acts transitively on $\A(\Sr, \K)$ and continuity of $\Theta(o,\cdot)$ is obvious since this map is induced by $$\N(\K) \rightarrow (\G \otimes_k \K)^{\rm an}, \ \  \ n \mapsto n\vartheta(o).$$ Existence of such an extension is established by transfinite induction for a well-ordering on $\mathbb{R}$; note that we could restrict to non-Archimedean extensions of $k$ containing a given enumerable family of extensions since sequential continuity of $\Theta(o, \cdot)$ on $\A(\Sr,k)$ is enough.

\vspace{0.1cm} If the field $k$ is locally compact, $\mathcal{B}^e(\G,k)$ is locally compact and the continuous bijection $\Theta(o,\cdot)$ is a homeomorphism onto a closed subset of $\G^{\rm an}$. \hfill $\Box$

\begin{Rk} If $\G$ is split and $o$ is a given special point of $\mathcal{B}^e(\G,k)$, the map $\Theta(o,\cdot)$ above coincides with the one defined by Berkovich in~\cite[5.4.2]{Ber1} starting with the $k^{\circ}$-Demazure group $\mathcal{G}_o$. Indeed, with the notations of [loc. cit], we have $\mathbf{P} = \G_o$, $\mathbf{p} = o$ and $\mathbf{t}_{\lambda} \ast \mathbf{P} = {\rm pr}_{\K/k}(t\G_o)$, where $\K/k$ is a non-Archimedean extension such that $\lambda$ takes values in $|\K^{\times}|$ and $t$ is an element of $\K$ satisfying $\lambda = \psi(t)$.
\end{Rk}

\subsection{Realizations of buildings in flag varieties}
\label{ss - buildings in flags}

\noindent \textbf{(2.4.1)} With each parabolic subgroup $\rP \in \mathrm{Par}(\G)(k)$ is associated a morphism $\lambda_{\rP}: \G \rightarrow \mathrm{Par}(\G)$ defined by the following condition: for any $k$-scheme $\Sr$, a point $g \in \G(\Sr)$ is mapped to the parabolic subgroup $\lambda_{\rP}(g) =  g (\rP \times_k \Sr) g^{-1}$ of $\G \times_k \Sr$. We recall that the image of $\lambda_{\rP}$ is the connected component $\mathrm{Par}_{t(\rP)}(\G)$ of $\mathrm{Par}(\G)$ which defines the type $t(\rP)$ of $\rP$, and that the morphism $\lambda_{\rP}$ identifies the scheme $\mathrm{Par}_{t(\rP)}(\G)$ with the quotient $\G/\rP$ (see 1.1.3).

\begin{Lemma} \label{lemma.theta.parabolic} For any parabolic subgroup $\rP \in {\rm Par}(\G)(k)$, the map $$\lambda_{\rP} \circ \vartheta: \mathcal{B}(\G,k) \rightarrow {\rm Par}(\G)^{{\rm an}}$$ depends only on the type $t$ of $\rP$ and is $\G(k)$-equivariant.
\end{Lemma}

\vspace{0.1cm}
\noindent
\emph{\textbf{Proof}}.
For any non-Archimedean extension $\K/k$, the following diagram $$\xymatrix{\mathcal{B}(\G,\K) \ar@{->}[r]^{\vartheta} & (\G \otimes_k \K)^{\mathrm{an}} \ar@{->}[d] \ar@{->}[r]^{\lambda_{\rP \otimes_k \K}} &  \mathrm{Par}(\G \otimes_k \K)^{\mathrm{an}} \ar@{->}[d] \\ \mathcal{B}(\G,k) \ar@{->}[u] \ar@{->}[r]_{\vartheta} & \G^{\mathrm{an}} \ar@{->}[r]_{\lambda_{\rP}} & \mathrm{Par}(\G)^{\mathrm{an}}}$$ is commutative by Proposition \ref{prop.theta.properties}. The vertical arrows are $\G(k)$-equivariant; we can therefore assume that the group $\G$ is split and it suffices to check that the restriction of the map $\lambda_{\rP} \circ \vartheta$ to the set of special points does not depend on the choice of the parabolic subgroup $\rP \in \mathrm{Par}_t(\G)(k)$ and is $\G(k)$-equivariant.

Let $o$ be a special point of $\mathcal{B}(\G,k)$ and let $\mathcal{G}$ denote the corresponding Demazure group scheme over $k^{\circ}$ with generic fibre $\G$. Since $\mathrm{Par}(\G)(k) = \mathrm{Par}(\mathcal{G})(k^{\circ})$, the group $\rP$ is the generic fibre of a parabolic subgroup $\mathcal{P}$ of $\mathcal{G}$ of type $t$ and the map $\lambda_{\rP}$ is induced by the map $\lambda_{\mathcal{P}}: \mathcal{G} \rightarrow \mathrm{Par}_t(\mathcal{G})$. If we let $r$ denote the reduction maps, it follows that the diagram $$\xymatrix{\G^{\mathrm{an}} \ar@{->}[r]^r \ar@{->}[d]_{\lambda_{\rP}} & \mathcal{G} \otimes_{k^{\circ}} \widetilde{k} \ar@{->}[d]^{\lambda_{\mathcal{P}} \otimes_{k^{\circ}} \widetilde{k}} \\ \mathrm{Par}(\G)^{\mathrm{an}} \ar@{->}[r]_{r} & \mathrm{Par}(\mathcal{G}) \otimes_{k^{\circ}} \widetilde{k}}$$ is commutative. Since the morphism $\lambda_{\mathcal{P}} \otimes_{k^{\circ}} 1\widetilde{k} = \lambda_{\mathcal{P} \otimes_{k^{\circ}} \widetilde{k}}$ is dominant, the generic point of $\mathcal{G} \otimes_{k^{\circ}} \widetilde{k}$ is mapped to the generic point of the connected component $\mathrm{Par}_t(\mathcal{\G}) \otimes_{k^{\circ}} \widetilde{k}$ of ${\rm Par}(\mathcal{G}) \otimes_{k^{\circ}} \widetilde{k}$ and therefore $\lambda_{\mathrm{\rP}} \circ \vartheta(o)$ is the unique point in $\mathrm{Par}(\G)^{\mathrm{an}}$ lying over the generic point of $\mathrm{Par}_t(\mathcal{G}) \otimes_{k^{\circ}} \widetilde{k}$. In particular, this point does not depend on the choice of $\rP \in \mathrm{Par}_t(\G)(k)$.

For any $g \in \G(k)$, we have $\lambda_{\rP}(\vartheta(o)g^{-1}) = \lambda_{g^{-1}\rP g}(\vartheta(o))$, hence $\lambda_{\rP}(\vartheta(o)g^{-1}) = \lambda_{\rP}(\vartheta(o))$. On the other hand, since the map $\lambda_{\rP}: \G^{\mathrm{an}} \rightarrow \mathrm{Par}(\G)^{\mathrm{an}}$ is $\G(k)$-equivariant when we let $\G(k)$ act by left translations on $\G^{\mathrm{an}}$, we get $\lambda_{\rP}(gx) = g\lambda_{\rP}(x)g^{-1}$ for any point $x$ in $\G^{\mathrm{an}}$. We obtain therefore $$(\lambda_{\rP} \circ \vartheta)(go) = \lambda_{\rP}(g\vartheta(o)g^{-1}) = g(\lambda_{\rP} \circ \vartheta)(o)g^{-1}$$ and this shows that the map $\lambda_{\rP} \circ \vartheta$ is $\G(k)$-equivariant.
\hfill $\Box$

\vspace{0.2cm}
\begin{Rk} \label{rk.Berkovich.comparison} While proving the lemma above, we have shown that $(\lambda_{\rP} \circ \vartheta)(gx) = \lambda_{\rP}(g\vartheta(x))$ for any element $g \in \G(k)$ and any point $x \in \mathcal{B}(\G,k)$. Since $g\vartheta(x) = g\Theta(x,x) = \Theta(x,gx)$, it follows that $$\lambda_{\rP} \circ \vartheta (gx) = \lambda_{\rP} \circ \Theta (x,gx).$$ Note that the right hand side makes it obvious that the map $\lambda_{\rP} \circ \vartheta$ is $\G(k)$-equivariant; moreover, this is also the definition adopted by Berkovich in~\cite[Sect. 5.5]{Ber1}, when $\G$ is split.
\end{Rk}

\vspace{0.2cm}
\begin{Def} \label{def.theta-t} For a $k$-rational type $t$, we denote by $\vartheta_t: \mathcal{B}(\G,k) \rightarrow {\rm Par}(\G)^{{\rm an}}$ the $\G(k)$-equivariant map defined by $\vartheta_t = \lambda_{\rP} \circ \vartheta$ for any $\rP \in {\rm Par}_t(\G)(k)$.
\end{Def}

\vspace{0.2cm}

\begin{Prop} \label{prop.theta-t.extension} For any $k$-rational type $t$ of $\G$ and any non-Archimedean extension $k'/k$, the diagram $$\xymatrix{\mathcal{B}(\G,k') \ar@{->}[r]^{\vartheta_t} & {\rm Par}(\G \otimes_k k')^{{\rm an}} \ar@{->}[d]^p\\ \mathcal{B}(\G,k) \ar@{->}[r]_{\vartheta_t} \ar@{^{(}->}[u]^{i} & {\rm Par}(\G)^{{\rm an}} }$$ in which $i$ denotes the canonical injection and $p$ the canonical projection, is commutative. \\
\indent Moreover, if $k'$ is a Galois extension of $k$, the upper horizontal arrow is ${\rm Gal}(k'/k)$ equivariant and the restriction of $p$ to the Galois-fixed point set in ${\rm Par}(\G \otimes_k k')^{\rm an}$ is injective.
\end{Prop}

\vspace{0.1cm}
\noindent
\emph{\textbf{Proof}}.
The first assertion follows immediately from Proposition \ref{prop.theta.properties} (ii), and from the commutativity of the diagram $$\xymatrix{(\G \otimes_k \K)^{\mathrm{an}} \ar@{->}[r]^{\hspace{-0.5cm} \lambda_{\rP \otimes_k \K}} \ar@{->}[d]_{p} & \mathrm{Par}(\G \otimes_k \K)^{\mathrm{an}} \ar@{->}[d]^{p} \\ \G^{\mathrm{an}} \ar@{->}[r]_{\lambda_{\rP}} & \mathrm{Par}(\G)^{\mathrm{an}}}$$ for any $\rP \in \mathrm{Par}_t(\G)(k)$.

\vspace{0.1cm} The second assertion follows from Galois-equivariance of $\vartheta$ and $\lambda_{\rP}$. The third assertion follows from the fact that each fibre of $p$ is a Galois orbit.
\hfill $\Box$

\vspace{0.8cm}
\noindent \textbf{(2.4.2)} We still consider a $k$-rational type $t$ of $\G$. Assuming that $\G$ is split, we give an explicit description of the map $\vartheta_t$, completely similar to the one in (2.1.2).

Let $\rP$ be a parabolic subgroup of $\G$ of type $t$ and pick a maximal split torus $\T$ of $\G$ contained in $\rP$. If we denote by $\rP^{\mathrm{op}}$ the parabolic subgroup of $\G$ opposite to $\rP$ with respect to $\T$, the morphism $$\radu(\rP^{\mathrm{op}}) \rightarrow \mathrm{Par}(\G), \ \ g \mapsto g\rP g^{-1}$$ (defined in terms of the functor of points) is an isomorphism onto an open subscheme of $\mathrm{Par}(\G)$ which we still denote by $\Omega(\T,\rP)$ by abuse of notation.

Let $\Phi(\G,\T)$ be the set of roots of $\G$ with respect to $\T$, pick a
special point $o$ in $\mathcal{B}(\G,k)$ and consider the corresponding
$k^{\circ}$-Chevalley group $\mathcal{G}$. The choice of an \emph{integral}
Chevalley basis in $\mathrm{Lie}(\G)(k)$ leads to an isomorphism of
$\radu(\rP^{\rm op})$ with the affine space
$$\prod_{\alpha \in \Psi} \U_{\alpha} \simeq \prod_{\alpha \in \Psi} \mathbb{A}^{1}_k,$$
where $\Psi = \Phi(\radu(\rP^{\rm op}),\T) = -\Phi(\radu(\rP),\T)$.

\begin{Prop} \label{prop.theta-t.explicit}
We assume that the group $\G$ is split and we use the notation introduced above.

\begin{itemize}
\item[(i)] The map $\vartheta_t$ sends the point $o$ to the point of
$\Omega(\T,\rP)^{{\rm an}}$ corresponding to the multiplicative (semi)norm
$$k\left[(\X_{\alpha})_{\alpha \in \Psi}\right] \rightarrow  \mathbb{R}_{\geqslant 0}, \ \
\sum_{\nu \in \mathbb{N}^{\Psi}} a_{\nu} \X^{\nu} \mapsto \max_{\nu} |a_{\nu}|.$$
\item[(ii)] Using the point $o$ to identify the apartment $\A(\T,k)$ with the vector space $\V(\T) = {\rm Hom}_{\mathbf{Ab}}(\X^*(\T),\mathbb{R})$,
the map $\V(\T) \rightarrow {\rm Par}(\G)^{{\rm an}}$ induced by $\vartheta_t$
associates with an element $u$ of $\V(\T)$
the point of $\Omega(\rP,\T)^{{\rm an}}$ corresponding to the mutiplicative
seminorm
$$k\left[(\X_{\alpha})_{\alpha \in \Psi}\right] \rightarrow  \mathbb{R}_{\geqslant 0}, \ \
\sum_{\nu \in \mathbb{N}^{\Psi}} a_{\nu} \X^{\nu} \mapsto \max_{\nu}
|a_{\nu}|\prod_{\alpha \in \Psi} e^{\nu(\alpha)\langle u, \alpha \rangle}.$$
\end{itemize}
\end{Prop}

\vspace{0.1cm}

\noindent \emph{\textbf{Proof}}.
We can argue exactly as for Proposition \ref{prop.theta.explicit}.
\hfill $\Box$

\vspace{0.2cm}

\begin{Cor}\label{cor.theta-t.norm}
For each point $x$ of $\mathcal{B}(\G,k)$, the seminorm $\vartheta_t(x)$ induces an extension of the absolute value of $k$ to the function field of ${\rm Par}_t(\G)$.
\end{Cor}

\vspace{0.2cm} This means that $\vartheta_t(x)$ is mapped to the generic point of ${\rm Par}_t(\G)$ by the canonical map $\rho : {\rm Par}_t(\G)^{\rm an} \rightarrow {\rm Par}(\G)$ (see 1.2.2).

\vspace{0.1cm}
\noindent
\emph{\textbf{Proof}}.
It suffices to prove this assertion when the group $\G$ is split. By the preceding proposition, $\vartheta_t(x)$ induces a multiplicative seminorm on the $k$-algebra $\A$ of any big cell $\Omega(\rP,\T)$ of $\Par_t(\G)$ which extends the absolute value of $k$ and satisfies the following condition: given any $f \in \A$, we have $|f|(\vartheta_t(x)) = 0$ if and only if $f=0$. Therefore, this seminorm is a norm and extends to an absolute value on the fraction field $\mathrm{Quot}(\A)$ of $\A$ extending the absolute value of $k$. Finally, since $\Omega(\rP,\T)$ is an affine open subset of the integral scheme $\Par_t(\G)$, the field $\mathrm{Quot}(\A)$ is nothing but the function field of $\Par_t(\G)$. \hfill $\Box$

\vspace{0.8cm}

\noindent \textbf{(2.4.3)} The map $\vartheta_t$ can be defined more generally for a type $t$ which is not $k$-rational, i.e., corresponds to a connected component ${\rm Par}_t(\G)$ of ${\rm Par}(\G)$ such that ${\rm Par}_t(\G)(k) = \varnothing$. The most important case is the type $t = \varnothing$ of Borel subgroups for a group $\G$ which is not quasi-split.

\vspace{0.1cm} Consider a finite Galois extension $k'/k$ splitting $\G$, set $\Gamma = {\rm Gal}(k'|k)$ and pick a type $t'$ of $\G \otimes_k k'$ over $t$, i.e., a connected component ${\rm Par}_{t'}(\G \otimes_k k')$ of ${\rm Par}(\G \otimes_k k') = {\rm Par}(\G) \otimes_k k'$ lying over ${\rm Par}_t(\G)$. Letting ${\rm pr}_{k'/k}$ denote the canonical projection ${\rm Par}_{t'}(\G \otimes_k k')^{\rm an} \rightarrow {\rm Par}(\G)^{\rm an}$, the map ${\rm pr}_{k'/k} \circ \vartheta_{t'}$ does not depend on the choice of $t'$ since $\Gamma$ acts transitively on the types of $\G \otimes_k k'$ lying over $t$ and we set $\vartheta_t = {\rm pr}_{k'/k} \circ \vartheta_{t'}$.

\vspace{0.1cm} One easily checks that proposition \ref{prop.theta-t.extension} holds in this more general situation.

\section{Compactifications of buildings}
\label{s - compactifications}

In this section we define, for a given reductive group $\G$ over a complete non-Archimedean field $k$, the Berkovich compactifications of the Bruhat-Tits building $\mathcal{B}(\G,k)$. If $k$ is a local field, these compactifications are defined by considering the maps $\vartheta_t : \mathcal{B}(\G,k) \to {\rm Par}_t(\G)^{\rm an}$ defined in the previous section (\ref{ss - buildings in flags}) and taking closures of their images; in general, we have first to restrict to apartments.
In any case, restricting the map $\vartheta_t$ to an apartment is the key point in order to obtain an injectivity criterion for $\vartheta_t$ and to analyse the corresponding compactification of $\mathcal{B}(\G,k)$. The latter space is described in terms of multiplicative seminorms on the coordinate ring of big cells of ${\rm Par}_t(\G)$ (Proposition \ref{prop.theta-t.explicit} and proof of Proposition \ref{prop.compact.extension}).

\vspace{0.1cm}
All types $t$ of parabolic subgroups considered in this section are $k$-\emph{rational}, i.e., correspond to a connected component ${\rm Par}_t(\G)$ of ${\rm Par}(\G)$ having a $k$-point; equivalently, $t$ corresponds to a conjugacy class of parabolic subgroups of $\G$. A similar construction can be made for any type $t$, maybe non-rational, since we have a well-defined map $\vartheta_t : \mathcal{B}(\G,k) \rightarrow {\rm Par}_t^{\rm an}$, but it will be shown in Appendix C that there is nothing to be gained since the corresponding compactification of $\mathcal{B}(\G,k)$ already occurs among compactifications associated with $k$-rational types.

\vspace{0.1cm} Given a $k$-rational type $t$ of parabolic subgroups of $\G$, we begin by introducing a class of parabolic subgroups of $\G$, which we call $t$-\emph{relevant} and which will later be useful to describe (the boundary components of) the compactification of type $t$ of $\mathcal{B}(\G,k)$ (see \ref{ss - relevant parabolics}).

\vspace{0.3cm}
\subsection{Reminder on quasisimple factors, and a warning}

Let $k$ be a field and $\G$ a connected semisimple $k$-group. There exist a unique (finite) family $(\G_i)_{i \in \I}$ of pairwise commuting smooth, normal and connected closed subgroups of $\G$, each of them quasi-simple, such that the product morphism $$\prod_{i \in \I} \G_i \rightarrow \G$$ is a central isogeny. The $\G_i$'s are the \emph{quasi-simple components} of $\G$. More generally, the quasi-simple components of a reductive $k$-group are the quasi-simple components of its derived subgroup.

\vspace{0.1cm} The isogeny $\prod_{i \in \I} \G_i \rightarrow \G$ induces an isomorphism of buildings
$$\prod_{i \in \I} \mathcal{B}(\G_i,k) = \mathcal{B}\big(\prod_{i \in \I} \G_i,k \big) \xymatrix{{} \ar@{->}[r]^{\sim} & {}} \mathcal{B}(\G,k)$$
and a $k$-isomorphism
$$\prod_{i \in \I}\mathrm{Par}(\G_i) = \mathrm{Par}\big(\prod_{i \in \I}\G_i\big) \xymatrix{{} \ar@{->}[r] & {}} \mathrm{Par}(\G).$$

\vspace{0.1cm} For each $k$-rational type $t$ of $\G$, the \emph{restriction} of $t$ to the quasi-simple component $\G_i$ is by definition the type $t_i$ of its parabolic subgroup $\rP \cap \G_i$, where $\rP$ is any element of $\mathrm{Par}_t(\G)(k)$. When no confusion seems likely to arise, we write $t$ instead of $t_i$.

We say that the restriction of $t$ to $\G_i$ is \emph{trivial} if $t_i$ is the type of the maximal parabolic subgroup $\G_i$, i.e., if any $\rP \in \mathrm{Par}_t(\G)(k)$ contains the full component $\G_i$. A $k$-rational type $t$ of $\G$ is \emph{non-degenerate} if its restriction to each isotropic quasi-simple component of $\G$ is non-trivial, i.e., if any $\rP \in \mathrm{Par}_t(\G)(k)$ induces a proper parabolic subgroup on each isotropic quasi-simple component of $\G$.

\begin{Def} \label{def.t-factor} For any $k$-rational type $t$, we let $\mathcal{B}_t(\G,k)$ denote the factor of $\mathcal{B}(\G,k)$ obtained by removing from the building each quasi-simple component on which the restriction of $t$ is trivial: $$\mathcal{B}_t(\G,k) = \prod_{\tiny \begin{array}{l} \hspace{0.5cm} i \in \I \\ t_i \textrm{ is non-trivial}\end{array}} \mathcal{B}(\G_i,k).$$ \end{Def}

\vspace{0.2cm}\emph{One word of caution about the notation $\mathcal{B}_t(G,k)$ and $\overline{\mathcal{B}}_t(\G,k)$ to be introduced in this section: the first one denotes the factor building of $\mathcal{B}(G,k)$ associated with the $k$-rational type $t$, the second a compactification of $\mathcal{B}_t(G,k)$ which still depends on $t$; for example, if $t$ and $t'$ are distinct non-degenerate $k$-rational types of $G$, then $\mathcal{B}_t(G,k) = \mathcal{B}_{t'}(G,k) = \mathcal{B}(G,k)$ but $\overline{\mathcal{B}}_t(\G,k) \neq \overline{\mathcal{B}}_{t'}(G,k)$.}

\vspace{0.3cm}
\subsection{Relevant parabolic subgroups}
\label{ss - relevant parabolics}

\vspace{0.3cm} \noindent
\emph{Reminder ([\textbf{SGA3}], Expos\'e XXVI, D\'efinition 4.4.2)} ---
Let $\Sr$ be a scheme and let $\G$ be a reductive $\Sr$-group scheme.
Two parabolic subgroups of $\G$, say $\rP$ and $\Qr$, are called \emph{osculatory} if $\rP \cap \Qr$ is a parabolic subgroup of $\G$.
This is equivalent to the following requirement: locally for the \'etale topology on $\Sr$, there exists a Borel subgroup of $\G$ simultaneously contained in $\rP$ and $\Qr$.

\vspace{0.5cm}

\noindent \textbf{(3.2.1)} Let $k$ be a field and $\G$ a reductive $k$-group. We consider a $k$-rational type $t$ of $\G$ and attach with each parabolic subgroup of $\G$ a closed subscheme of ${\rm Par}_t(\G)$.

\begin{Prop}
\label{prop.osc}
For any parabolic subgroup $\Qr$ of $\G$, the functor $$\left(\mathbf{Sch}/k\right)^{\mathrm{op}} \rightarrow \mathbf{Sets}, \ \Sr \mapsto \left\{\rP \in {\rm Par}_t(\G)(\Sr); \ \rP \ \textrm{and } \Qr \times_k \Sr \ \textrm{are osculatory} \right\}$$ is representable by a closed subscheme ${\rm Osc}_t(\Qr)$ of $\mathrm{Par}(\G)$, the \emph{osculatory subvariety} of $\Qr$ in ${\rm Par}_t(\G)$. This scheme is homogeneous under $\Qr$ and the morphism $\varepsilon_{\rP}: {\rm Osc}_t(\Qr) \rightarrow \mathrm{Par}(\Qr)$, defined functor-theoretically by $${\rm Osc}_t(\Qr)(\Sr) \rightarrow {\rm Par}(\Qr)(\Sr), \ \rP \mapsto \rP \cap \Qr,$$ is an isomorphism onto a connected component of ${\rm Par}(\Qr)$.
\end{Prop}

\vspace{0.1cm}
\noindent
\emph{\textbf{Proof}}.
Pick a parabolic subgroup $\Qr$ of $\G$ and first note that there exists a parabolic subgroup of $\G$ of type $t$ osculatory with $\Qr$. Indeed, given any parabolic subgroup $\rP$ of $\G$ of type $t$, $\Qr$ and $\rP$ contain minimal parabolic subgroups $\Qr_1$ and $\rP_1$ respectively; since any two minimal parabolic subgroups of $\G$ are conjugate in $\G(k)$~\cite[Expos\'e XXVI, Corollaire 5.7]{SGA3}, there exists $g \in \G(k)$ such that $\Qr_1 = g\rP_1 g^{-1}$ and therefore $g \rP g^{-1}$ is a parabolic subgroup of type $t$ osculatory with $\Qr$.

\vspace{0.1cm} Now we consider a parabolic subgroup $\rP$ of $\G$ of type $t$ osculatory with $\Qr$. For any $k$-scheme $\Sr$ and any parabolic subgroup $\rP' \in \mathrm{Par}_t(\G)(\Sr)$ osculatory with $\Qr \times_k \Sr$, the pairs $(\rP', \Qr \times_k \Sr)$ and $(\rP \times_k \Sr, \Qr \times_k \Sr)$ are conjugate \'etale locally over $\Sr$~\cite[Expos\'e XXVI, Corollaire 4.4.3]{SGA3}: there exist a covering \'etale morphism $\Sr' \rightarrow \Sr$ and an element $g \in \G(\Sr')$ such that $\rP' \times_{\Sr} \Sr' = g(\rP \times_k \Sr')g^{-1}$ and $\Qr \times_k \Sr' = g(\Qr \times_k \Sr')g^{-1}$. The last condition amounts to $g \in \Qr(\Sr')$ since $\Qr = \mathrm{Norm}_{\G}(\Qr)$. Set $\Sr'' = \Sr' \times_{\Sr} \Sr'$ and let $p_1, p_2: \Sr'' \rightarrow \Sr'$ denote the two canonical projections. The elements $g_1 = g \circ p_1$ and $g_2 = g \circ p_2$ of $\Qr(\Sr'')$ satisfy \begin{eqnarray*} g_1(\rP \times_k \Sr'') g_1^{-1} = p_1^{*}\left(g(\rP \times_k \Sr' )g^{-1}\right) & = & p_1^{*}(\rP' \times_\Sr \Sr') \\ & = & \rP' \times_{\Sr} \Sr'' \\ & = & p_2^* (\rP' \times_{\Sr} \Sr'') = p_2^{*}\left(g(\rP \times_k \Sr')g^{-1}\right) = g_2 (\rP \times_k \Sr'') g_2^{-1}, \end{eqnarray*} hence $g_2^{-1}g_1 \in \rP(\Sr'')$ since $\rP = \mathrm{Norm}_{\G}(\rP)$. In other words, the element $g$ of $\Qr(\Sr')$ defines a section of the quotient sheaf $\Qr/\rP \cap \Qr$ over $\Sr$ and we have thus proved that the natural morphism $\Qr \rightarrow \Par _t (\G), g \mapsto g \rP g^{-1}$ induces an isomorphism between the quotient sheaf $\Qr/\Qr \cap \rP$ and the subfunctor of $\Par _t (\G)$ consisting of parabolic subgroups of type $t$ of $\G$ osculatory with $\Qr$.

Finally, since $(\rP \cap \Qr)/{\rm rad}(\Qr)$ is a parabolic subgroup of $\Qr/{\rm rad}(\Qr)$, the quotient sheaf $$\Qr/\rP \cap \Qr \simeq \left(\Qr/\rad (\Qr)\right)/\left(\rP \cap \Qr /\rad (\Qr)\right)$$ is representable by a smooth and projective $k$-scheme, canonically isomorphic to the connected component of $\Par (\Qr) = \Par \left(\Qr/\rad (\Qr)\right)$ containing $\rP \cap \Qr$. We conclude that the same assertion holds for the functor of parabolic subgroups of $\G$ of type $t$ osculatory with $\Qr$. \hfill $\Box$

\vspace{0.3cm} \noindent \begin{Rk} \label{rk.osc.conjugation} For any parabolic subgroup $\Qr$ of $\G$ and any element $g$ of $\G(k)$, the $k$-automorphism $\operatorname{int}(g)$ of $\Par(\G)$ maps $\mathrm{Osc}_t(\Qr)$ onto $\mathrm{Osc}_t(g\Qr g^{-1})$. Indeed, given a $k$-scheme $\Sr$ and a parabolic subgroup $\rP \in \Par(\G)(\Sr)$ of $\G \times_k \Sr$ osculatory with $\Qr \times_k \Sr$, $g \rP g^{-1} \cap (g \Qr g^{-1} \times_k \Sr) = g (\rP \cap (\Qr \times_k \Sr)) g^{-1}$ is a parabolic subgroup of $\G \times_k \Sr$, hence $g\rP g^{-1}$ is osculatory with $g \Qr g^{-1}$.\end{Rk}

\vspace{0.3cm} \noindent \emph{\textbf{Notation}} --- Given a $k$-rational type $t$ of $\G$ and a parabolic subgroup $\Qr$ of $\G$, we still let the letter $t$ denote the type of the $k$-reductive group $\Qr_{\rm ss} = \Qr /\rad (\Qr)$ defined by the parabolic subgroup $(\rP~\cap~\Qr) /\rad (\Qr)$, where $\rP$ is any element of $\Par_{t} (\G)(k)$ osculatory with $\Qr$. Equivalently, the canonical morphism $\varepsilon_{\rP}: \mathrm{Osc}_{t}(\Qr) \rightarrow \Par (\Qr) = \Par(\Qr_{\rm ss})$ is an isomorphism onto the connected component $\Par_{t} (\Qr_{\rm ss})$ of $\Par (\Qr_{\rm ss})$.

\vspace{0.3cm} Recall that, if $\G$ a reductive group, $\Sr$ is a split maximal torus and $\rP$ is a parabolic subgroup containing $\Sr$, then $\rP^{\mathrm{op}}$ denotes the parabolic subgroup of $\G$ opposite to $\rP$ with respect to ${\rm Cent}_{\G}(\Sr)$ and the morphism $\radu (\rP^{\mathrm{op}})~\rightarrow~\Par(\G)$, functor-theoretically defined by $g \mapsto g\rP g^{-1}$, is an open immersion whose image $\Omega(\Sr,\rP)$ is the \emph{big cell} of $(\Sr, \rP)$ in $\Par(\G)$. Next proposition gives explicit equations defining an osculatory subvarieties in a big cell.

\vspace{0.3cm} \begin{Prop} \label{prop.osc.bigcell} Let $\rP$ and $\Qr$ be two osculatory parabolic subgroups of $\G$ containing a maximal split torus $\Sr$ and let $t$ denote the type of $\rP$. We let $\overline{\Qr}$ denote the reductive $k$-group $\Qr/{\rm rad}^{\rm u} (\Qr)$, $\overline{\Sr}$ the maximal split torus in $\overline{\Qr}$ induced by $\Sr$ and we set $\overline{\rP} = \rP \cap \Qr /{\rm rad}^{\rm u} (\Qr).$
\begin{itemize}
\item[(i)] The canonical isomorphism ${\rm Osc}_t(\Qr) \tilde{\rightarrow} {\rm Par}_t(\overline{\Qr})$ identifies the open subscheme ${\rm Osc}_t(\Qr) \cap \Omega(\Sr, \rP)$ of ${\rm Osc}_t(\Qr)$ and the big cell $\Omega(\overline{\Sr},\overline{\rP})$ of $(\overline{\Sr},\overline{\rP})$ in ${\rm Par}_t(\overline{\Qr})$.

\item[(ii)] Let $\Psi = \Phi(\radu (\rP^{\rm op}), \Sr)$ denote the set of roots of ${\rm rad}^{\rm u} (\rP^{\rm op})$ with respect to $\Sr$ and fix a total order on $\Psi$. The preimage of the closed subscheme ${\rm Osc}_t(\Qr)$ of $\mathrm{Par}_t(\G)$ under the immersion $$j: \prod_{\alpha \in \Psi} \U_{\alpha} \simeq \ {\rm rad}^{\rm u}(\rP^{\rm op}) \hookrightarrow \mathrm{Par}_t(\G)$$ is the closed subscheme defined by the equations $u_{\alpha} = 1$, where $\alpha$ runs over the complement of $\Phi(\Qr, \Sr)$ in $\Psi$.
\end{itemize}
\end{Prop}

\vspace{0.1cm}
\noindent
\emph{\textbf{Proof}}.
We first prove the second assertion.

\vspace{0.1cm}

(ii) Let $\Z$ be the closed subscheme of $\prod_{\alpha \in \Psi} \U_{\alpha}$ defined by the equations $u_{\alpha} =1$, $\alpha$ runing over the complement of $\Phi(\Qr, \Sr)$ in $\Psi$.

Both $\Z$ and $j^{-1}\mathrm{Osc}_t(\Qr)$ are integral (i.e., reduced and irreducible) closed subschemes of $\prod_{\alpha \in \Psi} \U_\alpha$: this is obvious for $\Z$ since $\U_\alpha$ is a smooth and geometrically irreducible $k$-scheme for any root $\alpha$; for $j^{-1}\mathrm{Osc}_t(\Qr)$, this follows from the fact that this scheme is isomorphic to a non-empty open subscheme of the integral $k$-scheme $\mathrm{Osc}_t(\Qr) \simeq \Par _t(\overline{\Qr})$. The canonical morphism $\prod_{\alpha \in \Psi} \U_{\alpha} \rightarrow \G$ maps $\Z$ into $\Qr$, hence $\Z \subset j^{-1}\mathrm{Osc}_t(\Qr)$; since these two closed subschemes of $\prod_{\alpha \in \Psi} \U_{\alpha}$ are integral, we are reduced to checking that they have the same dimension.

Let $\Lr$ denote the Levi subgroup of $\Qr$ containing $\T$ and write $$\Phi(\Qr,\Sr) = \Phi(\Lr, \Sr) \cup \Phi(\radu (\Qr), \Sr).$$

Since the parabolic subgroups $\rP$ and $\Qr$ are osculatory, $\radu (\rP^{\mathrm{op}}) \cap \radu (\Qr) = \{1\}$ and thus \begin{eqnarray*} \Phi(\Qr,\Sr) \cap \Psi & = & \Phi(\Qr, \Sr) \cap \Phi(\radu (\rP^{\rm op}),\Sr) \\ & = & \Phi(\Lr, \Sr) \cap \Phi(\radu (\rP^{\mathrm{op}}), \Sr) = \Phi(\Lr, \Sr) \cap \Psi. \end{eqnarray*}
It follows that the canonical projection $$\prod_{\alpha \in \Psi} \U_{\alpha} \rightarrow \prod_{\alpha \in \Phi(\Lr,\Sr) \cap \Psi} \U_{\alpha}$$ restricts to an isomorphism between $\Z$ and $\prod_{\alpha \in \Phi(\Lr,\Sr) \cap \Psi} \U_{\alpha}$. The subgroup $\Lr \cap \rP$ of $\Lr$ is parabolic and the set $-(\Phi(\Lr,\Sr) \cap \Psi)$ consists of roots of its unipotent radical with respect to $\Sr$; since the morphism $f: \Lr \rightarrow \overline{\Qr}$ induced by the canonical projection $\Qr \rightarrow \overline{\Qr} = \Qr/\radu (\Qr)$ is an isomorphism of reductive $k$-groups, we deduce that $f$ leads to an isomorphism between $\prod_{\alpha \in \Phi(\Lr,\Sr) \cap \Psi} \U_{\alpha}$ and the unipotent radical of $\overline{\rP}^{\mathrm{op}}$.
The conclusion is now obvious: since $\radu (\overline{\rP}^{\mathrm{op}})$ is isomorphic to an open dense subset of the irreducible $k$-scheme $\Par _t (\G)$, we have
\begin{eqnarray*} \dim \Z & = & \dim \radu (\overline{\rP}^{\mathrm{op}}) \\ & = & \dim \Par _t(\overline{\Qr}) = \dim \mathrm{Osc}_t (\Qr) \end{eqnarray*}
and therefore $\Z = j^{-1}\mathrm{Osc}_t(\Qr)$.
This proves (ii).

\vspace{0.1cm}

(i) We have just proved that the canonical isomorphism $\radu (\rP^{\mathrm{op}}) \tilde{\rightarrow} \Omega(\Sr,\rP)$ identifies the closed subschemes $\mathrm{Osc}_t(\Qr) \cap \Omega(\Sr, \rP)$ and $\radu (\rP^{\mathrm{op}}) \cap \Lr = \radu \left((\rP \cap \Lr)^{\mathrm{op}}\right)$. The canonical isomorphism $\xymatrix{\Lr \ar@{->}[r]^{\sim} & \overline{\Qr}}$ thus leads to a commutative diagram $$\xymatrix{\radu (\overline{\rP}^{\mathrm{op}}) \ar@{<-}[r]^{\sim} \ar@{^{(}->}[d] & \radu (\rP^{\mathrm{op}}) \cap \Lr  \ar@{^{(}->}[d] \ar@{->}[rd]^{\sim} & \\ \Par_t(\overline{\Qr}) \ar@{<-}[r]^{\sim} & \mathrm{Osc}_t(\Qr) & \mathrm{Osc}_t(\Qr) \cap \Omega(\Sr,\rP) \ar@{_{(}->}[l] }$$ and we deduce that the isomorphism $\mathrm{Osc}_t(\Qr) \tilde{\rightarrow} \Par_t(\overline{\Qr})$ identifies the open subscheme $\mathrm{Osc}_t(\Qr)~\cap~\Omega(\Sr,\rP)$ with the big cell $\Omega(\overline{\Sr},\overline{\rP})$. \hfill $\Box$

\vspace{0.3cm}
\begin{Ex}
\label{ex.osc.sl}
Let $\V$ be a $k$-vector space of dimension $d+1$ ($d \in \mathbb{N}$) and consider the semisimple $k$-group $\G = \SL(\V)$.
The types of $\G$ are in one-to-one correspondence with the types of flags of linear subspaces of $\V$ and we let $\delta$ denote the type corresponding to the flags $\left( \{ 0 \} \subset \Hr \subset \V \right)$ with $\dim (\Hr) = d$.

Recall that the $k$-scheme $\mathbb{P}(\V)$ represents the functor \begin{multline*} \left(\mathbf{Sch}/k\right)^{\mathrm{op}} \rightarrow \mathbf{Sets}, \ \Sr \mapsto \left\{ \textrm{isomorphism  classes of invertible quotients of } \V \otimes_k \Sr\right\} \\  \hspace{0.5cm} \simeq \left\{ \mathcal{O}_{\Sr}-\textrm{submodules of } \V \otimes_k \Sr, \textrm{ locally direct summands of rank } d\right\}.\end{multline*} There exists a unique $k$-isomorphism $\lambda: \mathbb{P}(\V) \rightarrow \Par_{\delta}(\G)$ such that, for any $k$-scheme $\Sr$, the map $\lambda(\Sr): \mathbb{P}(\V,\Sr) \rightarrow \Par_{\delta}(\G)(\Sr)$ sends an $\mathcal{O}_{\Sr}$-submodule of $\V \otimes_k \Sr$, locally a direct summand of rank $d$, to the parabolic subgroup of $\G \times_k \Sr$ stabilizing it.

For two flags $\F, \ \F'$ of linear subspaces in $\V$, the condition that their stabilizers are osculatory amounts to requiring that there exists a complete flag containing both $\F$ and $\F'$.

Let us now consider a parabolic subgroup $\Qr$ of $\G$, which is the stabilizer of a flag
$$\{ 0 \} = \V_0 \varsubsetneq \V_1 \varsubsetneq \ldots \varsubsetneq \V_r \varsubsetneq \V_{r+1} = \V.$$
A parabolic subgroup $\rP \in \mathrm{Par}_{\delta}(\G)(k)$, corresponding to a flag $(\{ 0 \} \subset \Hr \subset \V)$ with $\dim (\Hr) = d$, is osculatory with $\Qr$ if and only if the hyperplane $\Hr$ contains the linear subspace $\V_r$ and, since this holds more generally for any $k$-scheme $\Sr$ and any $\rP \in \Par_{\delta}(\G)(\Sr)$, the isomorphism $\lambda: \mathbb{P}(\V) \tilde{\rightarrow} \Par_{\delta}(\G)$ identifies the closed subscheme $\mathrm{Osc}_{\delta}(\Qr)$ of $\Par_{\delta}(\G)$ with the projective subspace $\mathbb{P}(\V/\V_r)$ of $\mathbb{P}(\V)$.
\end{Ex}

\vspace{0.8cm}
\noindent \textbf{(3.2.2)} The example above clearly shows that two different parabolic subgroups $\Qr, \ \Qr'$ of $\G$ may define the same closed subscheme $\mathrm{Osc}_t(\Qr) = \mathrm{Osc}_t(\Qr')$ in $\Par_t(\G)$. It turns out that there is a distinguished parabolic subgroup attached with each osculatory subvariety in ${\rm Par}_t(\G)$.

\vspace{0.3cm}
\begin{Prop} \label{prop.stab} Let $t$ denote a $k$-rational type of $\G$. For any parabolic subgroup $\Qr$ of $\G$, the set of parabolic subgroups $\Qr'$ of $\G$ satisfying $$\mathrm{Osc}_t(\Qr) = \mathrm{Osc}_t(\Qr')$$ has a maximal element.
\end{Prop}

\vspace{0.1cm}
\noindent
\emph{\textbf{Proof}}.
First note that the group functor $$\left(\mathbf{Sch}/k\right)^{\mathrm{op}} \rightarrow \mathbf{Sets}, \ \Sr \mapsto \mathrm{Stab}_{\G(\Sr)}(\mathrm{Osc}_t(\Qr \times_k \Sr))$$ is representable by a closed subscheme of $\G$. Indeed, the group $\G$ acts naturally on the Hilbert scheme $\mathcal{H}$ of the projective $k$-scheme $\Par_t(\G)$ and, for any $k$-scheme $\Sr$, the stabilizer of $\mathrm{Osc}_t(\Qr) \times_k \Sr$ in $\G(\Sr)$ is exactly the subgroup of $\G(\Sr)$ fixing the point $x \in \mathcal{H}(k) \subset \mathcal{H}(\Sr)$ defined by $\mathrm{Osc}_t(\Qr) \times_k \Sr$. It follows that our functor is represented by the fibre of the morphism $\G \rightarrow \mathcal{H}, \ g \mapsto g x$, over the point $x$. We let $\Pi$ denote the subgroup of $\G$ thus defined.

Since the subgroup $\Qr$ stabilizes $\mathrm{Osc}_t(\Qr)$, the inclusion $\Qr \subset \Pi$ is obvious.

Pick a finite Galois extension $k'/k$ splitting $\G$ and let us consider a purely inseparable extension $k''/k'$ such that the reduced $k''$-scheme $\rP'' = \left(\Pi \otimes_k k''\right)_{\mathrm{red}}$ underlying $\Pi \otimes_k k''$ is a smooth $k''$-group. Since $\Qr$ is smooth, $\Qr \otimes_k k''$ is reduced and therefore the closed immersion $\Qr \otimes_k k'' \hookrightarrow \Pi \otimes_k k''$ factors through $\rP'' \hookrightarrow \Pi \otimes_k k''$. This proves that the smooth $k''$-group $\rP''$ is a parabolic subgroup of $\G \otimes_k k''$ containing $\Qr \otimes_k k''$.

Since the $k'$-group $\G \otimes_k k'$ is split, there exists a parabolic subgroup $\rP'$ of $\G \otimes_k k'$ containing $\Qr \otimes_k k'$ such that $\rP'' = \rP' \otimes_{k'} k''$. Thanks to faithfully flat descent, the closed immersion $\rP' \otimes_k k'' = \rP'' \hookrightarrow \Pi \otimes_k k''$ comes from a closed immersion $\rP' \hookrightarrow \Pi \otimes_k k'$ and $\rP'$ is therefore the greatest parabolic subgroup of $\G \otimes_k k'$ containing $\Qr \otimes_k k'$ and contained in $\Pi \otimes_k k'$. It follows immediately from this description of $\rP'$ that this $k'$-group descends to a parabolic $k$-group $\rP$ of $\G$ containing $\Qr$ and contained in $\Pi$.

The identity $$\mathrm{Osc}_t(\rP) = \mathrm{Osc}_t(\Qr)$$ is a direct consequence of the inclusions $\Qr \subset \rP \subset \Pi$, for the first one implies ${\rm Osc}_t(\Qr) \subset {\rm Osc}_t(\rP)$ whereas the second gives ${\rm Osc}_t(\rP) \subset {\rm Osc}_t(\Qr)$ by the very definition of $\Pi$. Therefore, $\rP$ is the maximal element of the set $$\left\{\Qr' \in \Par(\G)(k) \ ; \ \mathrm{Osc}_t(\Qr') = \mathrm{Osc}_t(\Qr)\right\}.$$ \hfill $\Box$

\vspace{0.3cm} \begin{Def} \label{def.relevant} Let $t$ denote a $k$-rational type of $\G$. A parabolic subgroup $\Qr$ of $\G$ is said to be \emph{t-relevant} if it coincides with the maximal element of the set $$\left\{\Qr' \in {\rm Par}(\G)(k) \ ; \ {\rm Osc}_t(\Qr') = {\rm Osc}_t(\Qr)\right\}.$$
\end{Def}

\vspace{0.1cm}

It follows from the proof of Proposition \ref{prop.stab} that this condition amounts to requiring that $\Qr$ is the maximal parabolic subgroup of $\G$ stabilizing ${\rm Osc}_t(\Qr)$.

\vspace{0.3cm} \begin{Rk} \label{rk.minimal.t-rel} Each parabolic subgroup $\Qr$ of $\G$ is contained in a unique minimal $t$-relevant parabolic subgroup, namely the maximal parabolic subgroup stabilizing ${\rm Osc}_t(\Qr)$.
\end{Rk}

\vspace{0.3cm}
\begin{Ex}
\label{ex.relevant}
(i) Let us focus again on the example above: $\V$ is a finite dimensional $k$-vector space, $\G = \SL(\V)$ and $t = \delta$ is the type of flags $\left(\{ 0 \} \subset \Hr \subset \V\right)$ with $\operatorname{codim} (\Hr) = 1$. In this situation, the $\delta$-relevant parabolic subgroups of $\G$ are the stabilizers of flags $\left(\{ 0 \} \subset \W \subset \V\right)$ (we allow $\W = \{ 0 \}$ or $\W = \V$).

\vspace{0.1cm}
(ii) If the group $\G$ is quasi-split and $t = \varnothing$ is the type of Borel subgroups, then each parabolic subgroup of $\G$ is $\varnothing$-relevant. Indeed, for all parabolic subgroups $\rP, \ \Qr$ of $\G$ with $\Qr \varsubsetneq \rP$, there exists a Borel subgroup of $\G$ contained in $\rP$ but not in $\Qr$, hence $\mathrm{Osc}_\varnothing (\Qr) \neq \mathrm{Osc}_\varnothing (\rP)$ and therefore $\Qr$ is the maximal parabolic subgroup of $\G$ stabilizing $\mathrm{Osc}_\varnothing (\Qr)$.
\end{Ex}

\vspace{0.3cm}
\begin{Rk}
\label{rk.relevant.extension}
Let $t$ denote a $k$-rational type of $\G$ and consider a parabolic subgroup $\Qr$ of $\G$. If $\Qr$ is $t$-relevant, then for any extension $k'/k$ the parabolic subgroup $\Qr \otimes_k k'$ of $\G \otimes_k k'$ is $t$-relevant. Indeed, if $\rP'$ denotes the $t$-relevant parabolic subgroup of $\G \otimes_k k'$ stabilizing $\mathrm{Osc}_t(\Qr \otimes_k k')$ and $\Pi$ the subgroup of $\G$ stabilizing $\mathrm{Osc}_t(\Qr)$, then $\Pi \otimes_k k'$ is the stabilizer of $\mathrm{Osc}_t(\Qr \otimes_k k')$ in $\G \otimes_k k'$ for $\mathrm{Osc}_t(\Qr \otimes_k k') = \mathrm{Osc}_t(\Qr) \otimes_k k'$. As shown in the proof of Proposition \ref{prop.stab}, $\Qr$ is the maximal parabolic subgroup of $\G$ contained in $\Pi$ and $\Qr \otimes_k k'' = \left(\Pi \otimes_k k''\right)_{\mathrm{red}}$ for a convenient extension $k''/k'$, hence $\rP' \otimes_{k'} k'' = \Qr \otimes_k k''$ and therefore $\rP' = \Qr \otimes_k k'$.
\end{Rk}

\vspace{0.6cm}
\noindent \textbf{(3.2.3)} We give in section 3.3 a description of $t$-relevant parabolic subgroups of the semisimple $k$-group $\G$ in terms of its Dynkin diagram. As an immediate consequence, we will see that, if $\Qr$ is a parabolic subgroup of $\G$ and $\Qr'$ is the smallest $t$-relevant parabolic subgroup of $\G$ containing $\Qr$, then the semisimple group $\Qr/\rad (\Qr)$ is isogeneous to a quotient of the semisimple group $\Qr'/\rad(\Qr')$.

\subsection{Fans and roots}

We consider again in this paragraph an arbitrary field $k$ and a semisimple $k$-group $\G$. The basic notions on fans and their associated compactifications are collected in appendix B.

\vspace{0.4cm} \noindent \textbf{(3.3.1)}
Let $\Sr$ be a maximal split torus of $\G$ with character group $\X^*(\Sr) = \mathrm{Hom}_{k-\mathbf{Gr}}(\Sr, \mathbb{G}_{{\rm m},k})$ and let $\Phi = \Phi(\G, \Sr)$ denote the set of roots of $\G$ with respect to $\Sr$. Since it is more convenient to adopt multiplicative notation in order to compactify affine spaces, we let $$\Lambda(\Sr) = {\rm Hom}_{\mathbf{Ab}}(\X^*(\Sr), \mathbb{R}_{>0})$$ denote the multiplicative dual of $\X^*(\Sr)$. Each character $\chi \in \X^*(\Sr)$ defines a positive real function on $\Lambda(\Sr)$.

For any parabolic subgroup $\rP$ of $\G$ containing $\Sr$, the set $\Phi(\rP,\Sr)$ of roots of $\rP$ with respect to $\Sr$ is the subset of $\Phi(\G,\Sr)$ consisting of all roots $\alpha$ such that $\rP$ contains the root group $\U_{\alpha}$.

\vspace{0.1cm} We first recall that the set of parabolic subgroups of $\G$ containing $\Sr$ has a nice description in terms of cones in $\Lambda(\Sr)$ (Coxeter complex).

\vspace{0.3cm} \begin{Prop} \label{prop.weylfan} Let $\rP$ be a parabolic subgroup of $\G$ containing $\Sr$.
\begin{itemize}
\item[(i)] The subset $\mathfrak{C}(\rP)$ of $\Lambda(\Sr)$, defined by the condition $\alpha \leqslant 1$ for all $\alpha \in \Phi(\rP^{\rm op},\Sr) = - \Phi(\rP,\Sr)$, is a strictly convex polyhedral cone.
\item[(ii)] The cone $\mathfrak{C}(\rP)$ spans $\Lambda(\Sr)$ if and only if $\rP$ is minimal.
\item[(iii)] The faces of the cone $\mathfrak{C}(\rP)$ are the cones $\mathfrak{C}(\Qr)$, where $\Qr$ runs over the set of parabolic subgroups of $\G$ containing $\rP$.
\item[(iv)] For any parabolic subgroup $\Qr$ of $\G$ containing $\Sr$, $\mathfrak{C}(\rP) \cap \mathfrak{C}(\Qr)$ is the cone associated with the smallest parabolic subgroup of $\G$ containing both $\rP$ and $\Qr$.
Moreover, when $\Qr$ runs over the set of parabolic subgroups of $\G$ containing $\Sr$, the cones $\mathfrak{C}(\Qr)$ are pairwise distinct and they cover $\Lambda(\Sr)$.
\end{itemize}
\end{Prop}

\vspace{0.1cm}
\noindent
\emph{\textbf{Proof}}.
All the assertions above are well-known and follow immediately from the fact that the map $\rP \mapsto \Phi(\rP,\Sr)$ sets up an increasing one-to-one correspondence between parabolic subgroups of $\G$ containing $\Sr$ and closed and generating subsets of $\Phi$, i.e., subsets $\Psi$ of $\Phi$ satisfying the following two conditions: \begin{itemize} \item for all $\alpha, \ \beta \in \Psi$, $\alpha + \beta \in \Phi \Longrightarrow \alpha + \beta \in \Psi$; \item for any $\alpha \in \Phi(\G,\Sr)$, either $\alpha$ or $-\alpha$ belongs to $\Psi$.
\end{itemize}
(See~\cite[Expos\'e XXVI, Proposition 7.7]{SGA3}). The first condition amounts to $\Psi = \langle \Psi \rangle^+ \cap \Phi$, where $\langle \Psi \rangle^+$ denotes the semigroup spanned by $\Psi$ in $\X^*(\Sr)$, whereas the second one implies that $\Phi$ and $\Psi$ span the same subgroup of $\X^*(\Sr)$.  \hfill $\Box$

\vspace{0.2cm} \begin{Rk} \label{rk.par.roots} Given a parabolic subgroup $\rP$ of $\G$ containing $\Sr$, we have $$\Phi(\rP,\Sr) = \Phi(\Lr_{\rP},\Sr) \sqcup \Phi({\rm rad}^{\rm u}(\rP),\Sr),$$ where $\Lr_{\rP}$ is the Levi subgroup of $\rP$ containing ${\rm Cent}_{\G}(\Sr)$. The set $\Phi(\Lr_{\rP},\Sr)$ consists precisely of roots $\alpha \in \Phi(\G,\Sr)$ such that both $\alpha$ and $-\alpha$ belong to $\Phi(\rP,\Sr)$; geometrically, this characterization is equivalent to $$\Phi(\Lr_{\rP},\Sr) = \{\alpha \in \Phi(\rP,\Sr) \ | \ \alpha_{|\mathfrak{C}(\rP)} = 1\}.$$
\end{Rk}

\vspace{0.3cm} If $\rP$ is a minimal parabolic subgroup of $\G$ containing $\Sr$, the interior of the cone $\mathfrak{C}(\rP)$ is usually referred to as the \emph{Weyl chamber} of $\rP$ in $\Lambda(\Sr)$. This motivates the following definition.

\vspace{0.3cm} \begin{Def} \label{def.weylfan} The \emph{Weyl fan} on the vector space $\Lambda(\Sr)$ is the fan consisting of the cones $\mathfrak{C}(\rP)$, where $\rP \in {\rm Par}(\G)(k)$ and $\Sr \subset \rP$.
\end{Def}

\vspace{0.3cm} Now we consider a $k$-rational type $t$ of $\G$ and associate with it a new family of polyhedral cones in $\Lambda(\Sr)$. The cones of higher dimension are roughly speaking the "combinatorial neighborhoods" \ of all Weyl cones $\mathfrak{C}(\rP)$ with $\rP \in {\rm Par}_t(\G)(k)$ and $\Sr \subset \rP$.

\vspace{0.3cm}
\begin{Def} \label{def.t-cone} For any parabolic subgroup $\rP$ of $\G$ of type $t$ which contains $\Sr$, we let $\C_t(\rP)$ denote the union of all cones $\mathfrak{C}(\Qr)$ associated with the parabolic subgroups $\Qr$ of $\G$ satisfying $\Sr \subset \Qr \subset \rP$: $$\C_t(\rP) = \bigcup_{\tiny \begin{array}{l} \Qr \in {\rm Par}(\G)(k) \\ \hspace{2pt} \Sr \subset \Qr \subset \rP \end{array}} \mathfrak{C}(\Qr).$$
\end{Def}

Note that is suffices to consider \emph{minimal} parabolic subgroups $\rP_0$ satisfying $\Sr \subset \rP_0 \subset \rP$ in the above definition.

\vspace{0.3cm} In order to analyze this definition, we recall that with any $k$-rational type $t$ of $\G$ are associated two normal and semisimple subgroups $\G'$ and $\G''$ of $\G$, uniquely characterized by the following conditions: \begin{itemize} \item the canonical morphism $\G' \times \G'' \rightarrow \G$ is a central isogeny;
\item the restriction of $t$ to $\G'$ (to $\G''$, respectively) is non-trivial on any quasi-simple component of $\G'$ (is trivial, respectively).
\end{itemize}
The subgroup $\G'$ ($\G''$, respectively) is simply the product of quasi-simple components of $\G$ to which the restriction of $t$ is non-trivial (is trivial, respectively).


The groups $\Sr' = (\Sr \cap \G')^{\circ}$ and $\Sr'' = (\Sr \cap \G'')^{\circ}$ are maximal split tori in $\G'$ and $\G''$ respectively and $\Sr = \Sr'\Sr''$. The isogeny $\Sr' \times \Sr'' \rightarrow \Sr$ induces an injective homomorphism $\X^*(\Sr) \rightarrow \X^*(\Sr') \oplus \X^*(\Sr'')$ whose image has finite index, hence a canonical isomorphism $\xymatrix{\Lambda(\Sr') \oplus \Lambda(\Sr'') \ar@{->}[r]^{\sim} & \Lambda(\Sr).}$ Finally, the set $\Phi = \Phi(\G, \Sr)$ is the union of the two disjoint subsets $$\Sigma' = \{\alpha \in \Phi \ | \ \alpha_{|\Sr''}=1\} \ \textrm{ and } \ \Sigma'' = \{\alpha \in \Phi \ | \ \alpha_{|\Sr''} = 1\}$$ and the canonical projection $\X^*(\Sr) \rightarrow \X^*(\Sr'), \ \alpha \mapsto \alpha_{|\Sr'}$ ($\X^*(\Sr) \rightarrow \X^*(\Sr''), \ \alpha \mapsto \alpha_{|\Sr''}$, respectively) induces a bijection between $\Sigma'$ and $\Phi' = \Phi'(\G,\Sr')$ (between $\Sigma''$ and $\Phi''=\Phi(\G'',\Sr'')$, respectively).


\vspace{0.3cm} \begin{Lemma} \label{lemma.t-cone} Let $\rP$ be a parabolic subgroup of $\G$ of type $t$ containing $\Sr$.
\begin{itemize}
\item[(i)] The subset $\C_t(\rP)$ of $\Lambda(\Sr)$ is the convex polyhedral cone $\{\alpha \leqslant 1 \ ; \ \alpha \in \Phi({\rm rad}^{\rm u}(\rP^{\rm op}),\Sr)\}.$
\item[(ii)] The maximal linear subspace contained in $\C_t(\rP)$ is $\Lambda(\Sr'')$.
\item[(iii)] For any parabolic subgroup $\rP'$ of $\G$ of type $t$ containing $\Sr$, the cones $\C_t(\rP)$ and $\C_t(\rP')$ intersect along a common face.
\end{itemize}
\end{Lemma}

\emph{\textbf{Proof}}. Each parabolic subgroup considered in what follows contains the maximal split torus $\Sr$.

\vspace{0.1cm}
(i) Note that $\Phi(\radu (\rP^{\rm op}), \Sr) = -\Phi(\radu(\rP),\Sr)$ is precisely the subset of $\Phi(\rP^{\rm op}, \Sr)$ consisting of all roots $\alpha$ such that $-\alpha \notin \Phi(\rP^{\rm op}, \Sr)$. We set $\C = \{\alpha \leqslant 1 \ ; \  \alpha \in \Phi(\radu (\rP^{\rm op}), \Sr)\}$ and consider a minimal parabolic subgroup $\rP_0$. If $\rP_0 \subset \rP$, then $\Phi(\rP_0^{\rm op}, \Sr) \subset \Phi(\rP^{\rm op}, \Sr)$ and $\Phi(\radu (\rP^{\rm op}),\Sr) \subset \Phi(\radu (\rP_0^{\rm op}, \Sr) = \Phi(\rP_0^{\rm op},\Sr)$, hence $\mathfrak{C}(\rP_0) \subset \C$ and therefore $\C_t(\rP) \subset \C$.

\vspace{0.1cm}
If $\rP_0 \nsubseteq \rP$, then $\Phi(\rP_0^{\rm op},\Sr) \nsubseteq \Phi(\rP^{\rm op},\Sr)$ and thus there exists a root $\alpha \in \Phi(\rP_0^{\rm op},\Sr)$ such that $\alpha \notin \Phi(\rP^{\rm op},\Sr)$ and $-\alpha \in \Phi(\rP^{\rm op},\Sr)$. Since $\alpha < 1$ on the interior $\mathfrak{C}(\rP_0)^{\circ}$ of $\mathfrak{C}(\rP_0)$, it follows that $\mathfrak{C}(\rP_0)^{\circ}$ is disjoint from $\C$ and thus $$\bigcup_{\tiny \begin{array}{l} \rP_0 \textrm{ minimal} \\ \hspace{2pt} \rP_0 \nsubseteq \rP \end{array}} \mathfrak{C}(\rP_0) \subset \Lambda(\Sr) - \C^{\circ}.$$

We remark that the left hand side is exactly the complement of the interior of $\C_t(\rP)$, so that $\C^{\circ} \subset \C_t(\rP)^{\circ}$ and $\C \subset \C_t(\rP)$. We have thus proved $$\C_t(\rP) = \{\alpha \leqslant 1 \ ; \ \alpha \in \Phi(\radu (\rP^{\rm op}),\Sr)\}.$$

\vspace{0.1cm} (ii) We use the notation introduced before stating the proposition. We can write $\Phi(\rP,\Sr) = \Psi' \cup \Sigma''$, where $\Psi'$ is the closed and generating subset of $\Sigma'$ whose image under the bijection $\Sigma' \xymatrix{\ar@{->}[r]^{\sim}&} \Phi'$ is the set $\Phi(\rP', \Sr')$ of roots of the parabolic subgroup $\rP' = \rP \cap \G'$ of $\G'$ with respect to $\Sr'$. Since $\alpha_{|\Lambda(\Sr'')} = 1$ for each root $\alpha \in \Sigma'$, the cone $$\C_t(\rP) = \{\alpha \leqslant 1 \ ; \ \alpha \in \Phi(\radu (\rP^{\rm op}),\Sr)\} = \{\alpha \leqslant 1 \ ; \ \alpha \in (-\Psi') \ \textrm{ and } \alpha \notin \Psi'\}$$ contains the linear subspace $\Lambda(\Sr'')$ and it is enough to check that the cone  of $\Lambda(\Sr')$ defined by the conditions: $\alpha \leqslant 1$ for all $\alpha \in \Phi(\radu ({\rP'}^{\rm op}), \Sr')$, is strictly convex. Thus we are reduced to proving that, if the $k$-rational type $t$ is non-degenerate, then the cone $\C_t(\rP)$ is strictly convex.

\vspace{0.1cm} Let us assume that the cone $\C_t(\rP)$ is \emph{not} strictly convex and let $\Lr$ denote the maximal linear subspace of $\Lambda(\Sr)$ it contains. We let $\W$ denote the Weyl group of the root system $\Phi$. The subgroup $\W_{\rP}$ of $\W$ stabilizing the cone $\mathfrak{C}(\rP)$ acts simply transitively on the set of cones $\mathfrak{C}(\rP_0)$, where $\rP_0$ is a minimal parabolic subgroup contained in $\rP$. This subgroup stabilizes $\C_t(\rP)$, hence the linear subspace $\Lr$ by maximality.

Pick a minimal parabolic subgroup $\rP_0$ contained in $\rP$ and denote by $\Delta \subset \Phi$ the corresponding set of simple roots. We also equip $\X^*(\Sr)$ and $\Lambda(\Sr)$ with a $\W$-invariant scalar product.

By (i), we have $\alpha_{|\Lr} = 1$ for each root $\alpha \in \Delta \cap \Phi({\rm rad}^{\rm u}(\rP^{\rm op}),\Sr)$, i.e., each root $\alpha \in \Delta$ whose restriction to $\mathfrak{C}(\rP)$ is not identically equal to $1$. Since $\Delta$ spans a subgroup of finite index in $\X^*(\Sr)$, the set $\Gamma = \{\alpha \in \Delta \ | \ \alpha_{|\Lr} \neq 1\}$ is non-empty as $\Lr \neq \{1\}$. Pick a root $\beta$ in $\Gamma$. Since $\beta_{|\Lr} \neq 1$, $\beta_{|\mathfrak{C}(\rP)} = 1$ and thus the orthogonal reflection $w_{\beta}$ with respect to the hyperplane $\{\beta =1\}$ belongs to $\W_{\rP}$. For any root $\alpha$ in $\Delta - \Gamma$, $\alpha_{|\Lr} =1$ hence $$w_{\beta}(\alpha)_{|\Lr} = w_{\beta}(\alpha)_{|w_{\beta}(\Lr)} = \alpha_{|\Lr} = 1$$ and therefore $(\alpha|\beta) = 0$ since $w_{\beta}(\alpha) = \alpha - 2\frac{(\alpha|\beta)}{(\beta|\beta)}\beta$.

We have just proved that $\Gamma$ and $\Delta - \Gamma$ are orthogonal, which implies that $\Gamma$ contains a connected component of the Dynkin diagram of $\Phi$. Since moreover $\Gamma$ is contained in the subset $ \Phi(\Lr_{\rP}, \Sr) \cap \Delta = \{\alpha \in \Delta \ | \ \alpha_{|\mathfrak{C}(\rP)} = 1\}$ of $\Delta$ associated with the parabolic subgroup $\rP$, the latter contains therefore an quasi-simple component of $\G$ and thus the type $t=t(\rP)$ is trivial on this component.

\vspace{0.1cm} (iii) Let us consider two distinct parabolic subgroups $\rP$ and $\rP'$ of type $t$. The cones $\C_t(\rP)$ and $\C_t(\rP')$ have disjoint interiors, hence their intersection is contained in a proper face of each by convexity. Let $\F$ and $\F'$ denote the minimal faces of $\C_t(\rP)$ and $\C_t(\rP')$ containing $\C_t(\rP) \cap \C_t(\rP')$. We have $\C_t(\rP) \cap \C_t(\rP') = \F \cap \F'$ and this cone meets the interior of both $\F$ and $\F'$ by minimality.

Assume $\F \nsubseteq \F'$, hence $\F^{\circ} \nsubseteq \F'$. Since $\F^{\circ} \cap {\F'}^{\circ} \neq \varnothing$, it follows that $\F^{\circ}$ meets $\partial \F'$ and thus there exists a Weyl cone $\mathfrak{C}$ whose interior meets both $\F^{\circ}$ and $\partial \F'$. We have $\mathfrak{C} \subset \F'$ and $\mathfrak{C} \subset \partial \F'$ for both $\F$ and $\F'$ are a union of Weyl cones. Let $\rP_0$ and $\rP'_0$ be two minimal parabolic subgroups satisfying $\Sr \subset \rP_0 \subset \rP$ and $\Sr \subset \rP_0' \subset \rP'$, and such that $\mathfrak{C}$ is a common face of $\mathfrak{C}(\rP_0)$ and $\mathfrak{C}(\rP_0')$. There exists a unique element $w$ in the Weyl group such that $\mathfrak{C}(\rP_0') = w\mathfrak{C}(\rP_0)$. By construction, the cone $w\F$ is a face of $\C_t(\rP')$ whose interior meets $\mathfrak{C}$, hence the smallest face of $\C_t(\rP')$ containing $\mathfrak{C}$; since $\mathfrak{C} \subset \partial \F'$, we deduce $w\F \subset \partial \F'$ and therefore ${\rm dim}(\F') > {\rm dim}(\F)$.    

It is now easy to conclude. If $\F \neq \F'$, then one of the following three situations occurs: \begin{itemize}
\item[a)] $\F \nsubseteq \F'$ and $\F' \nsubseteq \F$; \item[b)] $\F \subset \F'$ and $\F' \nsubseteq \F$; \item[c)] $\F' \subset \F$ and $\F \nsubseteq \F'$.\end{itemize}
In each case, the discussion above leads to a contradiction: \begin{itemize}
\item[a)] we obtain ${\rm dim}(\F') > {\rm dim}(\F)$ and ${\rm dim}(\F) > {\rm dim}(\F')$; \item[b)] we obtain $\F \subset \F'$ and ${\rm dim}(\F) > {\rm dim}(\F')$; \item[c)] we obtain $\F' \subset \F$ and ${\rm dim}(\F') > {\rm dim}(\F)$.\end{itemize}
Therefore, $\F = \F'$ and the cones $\C_t(\rP)$ and $\C_t(\rP')$ do intersect along a common face. \hfill $\Box$

\vspace{0.2cm} \begin{Rk} \label{rk.strict.convexity} It follows from assertions (i) and (ii) above that $\rP$ induces a non-trivial parabolic subgroup of each isotropic quasi-simple factor of $\G$ if and only if the cone $\C_t(\rP) = \{\alpha \leqslant 1 \ ; \ \alpha \in \Phi({\rm rad}^{\rm u}(\rP^{\rm op}),\Sr)\}$ is strictly convex, hence if and only if $\Phi({\rm rad}^{\rm u}(\rP),\Sr) = -\Phi({\rm rad}^{\rm u}(\rP^{\rm op}), \Sr)$ spans $\X^*(\Sr) \otimes_{\mathbb{Z}} \mathbb{Q}$.
\end{Rk}

\vspace{0.3cm} When $\rP$ runs over the set of parabolic subgroups of $\G$ of type $t$ containing $\Sr$, the family of faces of the cones $\C_t(\rP)$ fulfills all requirements defining a fan except possibly strict convexity. In fact, assertion (ii) of Lemma \ref{lemma.t-cone} shows that this family is the preimage of a fan on $\Lambda(\Sr')$ under the canonical projection $\Lambda(\Sr) \rightarrow \Lambda(\Sr)/\Lambda(\Sr'') \simeq \Lambda(\Sr')$.

\vspace{0.2cm}
\begin{Def} \label{def.t-fan} For any $k$-rational type $t$, the prefan $\mathcal{F}_t$ of type $t$ on $\Lambda(\Sr)$ is the collection of all faces of the cones $\C_t(\rP)$, where $\rP$ runs over the set ${\rm Par}_t(\G)(k)$.
\end{Def}

\vspace{0.2cm} All the cones which occur in the prefan $\mathcal{F}_t$ can be described in terms of parabolic subgroups. Note that, for every parabolic subgroup $\Qr$ of $\G$ containing $\Sr$, the set $$\{\C \in \mathcal{F}_t \ | \ \mathfrak{C}(\Qr) \subset \C\}$$ is non-empty --- indeed, $\mathfrak{C}(\Qr) \subset \mathfrak{C}(\Qr') \subset \C_t(\rP)$ if $\Qr'$ is a minimal parabolic subgroup containing $\Qr$ and $\rP$ is the unique element of $\Par _t(\G)(k)$ containing $\Qr'$ --- and is stable under intersection by Lemma \ref{lemma.t-cone}, (iii). Hence the following definition makes sense.

\vspace{0.3cm} \begin{Def} \label{def.t-fan.bis} Given any parabolic subgroup $\Qr$ of $\G$ containing $\Sr$, we let $\C_t(\Qr)$ denote the smallest cone in $\mathcal{F}_t$ containing $\mathfrak{C}(\Qr)$.
\end{Def}

\vspace{0.3cm} \begin{Rk} 1. This definition coincides with Definition \ref{def.t-cone} if $\Qr$ is of type $t$.

2. For any parabolic subgroup $\Qr$ of $\G$ containing $\Sr$ and any cone $\C$ in $\mathcal{F}_t$, with $\mathfrak{C}(\Qr) \subset \C$, the following conditions are equivalent: \begin{itemize} \item $\C_t(\Qr)=\C$;  \item $\mathfrak{C}(\Qr)$ meets the interior of $\C$. \end{itemize}

In particular, since each cone $\C \in\mathcal{F}_t$ is the union of Weyl cones of parabolic subgroups of $\G$ containing $\Sr$, we see immediately that $\C = \C_t(\Qr)$ for a convenient $\Qr$: indeed, we just have to choose $\Qr$ such that the cone $\mathfrak{C}(\Qr)$ meets the interior of $\C$.

\vspace{0.1cm}
2. If $t$ is the type of a minimal parabolic subgroup, $\C_t(\Qr) = \mathfrak{C}(\Qr)$ for any parabolic subgroup $\Qr$ of $\G$ containing $\Sr$ and $\mathcal{F}_t$ is therefore nothing but the Weyl fan on $\Lambda(\Sr)$.
\end{Rk}

\vspace{0.5cm} \noindent \textbf{(3.3.2)} For any $k$-rational type $t$, we now relate the cones $\C_t(\rP)$ to $t$-relevant parabolic subgroups.

\vspace{0.1cm} Throughout this paragraph, $\Sr$ is a maximal split torus of $\G$ and we let $\Phi = \Phi(\G,\Sr)$ denote the set of roots of $\G$ with respect to $\Sr$. For any parabolic subgroup $\rP$ of $\G$ containing $\Sr$, we denote by $\Lr_{\rP}$ the Levi subgroup of $\rP$ containing $\mathrm{Cent}_{\G}(\rP)$. The set $\Phi(\rP,\Sr)$ of roots of $\rP$ with respect to $\Sr$ is the disjoint union of the subsets $\Phi(\Lr_{\rP},\Sr)$ and $\Phi(\radu (\rP),\Sr)$.

\vspace{0.2cm} The following proposition relates the combinatorial construction of 3.2.1 and the geometric viewpoint of 3.2. Moreover, it gives an explicit description of the cones in $\mathcal{F}_t$.

\begin{Prop} \label{prop.roots.fan} Let $\rP$ and $\Qr$ be two parabolic subgroups of $\G$ containing $\Sr$ and assume that $\rP$ is of type $t$.
\begin{itemize}
\item[(i)] $\rP$ and $\Qr$ are osculatory if and only if $\C_t(\Qr) \subset \C_t(\rP)$.
\item[(ii)] The cone $\C_t(\rP)$ is defined by the inequalities $\alpha \leqslant 1$, $\alpha \in \Phi({\rm rad}^{\rm u} (\rP^{\rm op}),\Sr)$.
\item[(iii)] If $\rP$ and $\Qr$ are osculatory, $\C_t(\Qr)$ is the polyhedral cone defined by the conditions $$\left\{ \begin{array}{ll} \alpha \leqslant 1, &  \alpha \in \Phi({\rm rad}^{\rm u} (\rP^{\rm op}), \Sr) \\ \alpha =1, & \alpha \in \Phi({\rm rad}^{\rm u} (\rP^{\rm op}), \Sr) \cap \Phi(\Lr_{\Qr}, \Sr).\end{array} \right.$$
\end{itemize}
\end{Prop}

\vspace{0.1cm}
\noindent
\emph{\textbf{Proof}}.
(i) The inclusions $\C_t(\Qr) \subset \C_t(\rP)$ and $\mathfrak{C}(\Qr) \subset \C_t(\rP)$ are equivalent and the latter amounts to the existence of a minimal parabolic subgroup $\rP_0$ of $\G$ containing $\Sr$ such that $\mathfrak{C}(\rP_0)$ contains both $\mathfrak{C}(\rP)$ and $\mathfrak{C}(\Qr)$, i.e., such that $\rP_0$ is simultaneously contained in $\rP$ and $\Qr$. Thus $\C_t(\Qr) \subset \C_t(\rP)$ if and only if the parabolic subgroups $\rP$ and $\Qr$ are osculatory.

\vspace{0.1cm} (ii) This assertion was proved in Lemma \ref{lemma.t-cone}, (i).

\vspace{0.1cm} (iii) We assume that the parabolic subgroups $\rP$ and $\Qr$ are osculatory and let $\F$ denote the face of the cone $\C_t(\rP)$ defined by the equations $\alpha = 1$ for all $\alpha \in \Phi(\radu (\rP^{\rm op}),\Sr) \cap \Phi(\Lr_{\Qr}, \Sr)$. Since the conditions $\alpha \in \Phi(\Lr_{\Qr}, \Sr)$ and $\alpha_{|\C(\Qr)} = 1$ are equivalent for any root $\alpha \in \Phi$, $\F$ is clearly the smallest face of $\C_t(\rP)$ containing $\mathfrak{C}(\Qr)$, and thus $\F = \C_t(\Qr)$. \hfill $\Box$

\vspace{0.3cm} \begin{Ex}
\label{ex-fan}
Let $d \geqslant 1$ be an integer and $\G$ the semisimple $k$-group $\SL (d+1)$. We consider the type $\delta$ corresponding to flags $\left(\{0 \} \subset \Hr \subset k^{d+1} \right)$, where $\Hr$ is a hyperplane in $k^{d+1}$.

Let $\T$ denote the torus of diagonal matrices and $\B$ the Borel subgroup of $\G$ consisting of upper triangular matrices.
 If $\chi_1, \ldots, \chi_{d+1}$ are the characters of $\T$ defined by $\chi_i \left(\mathrm{diag}(t_1, \ldots, t_{d+1}) \right) = t_i$, $1 \leqslant i \leqslant d+1$, the set of roots is $\Phi( \SL(d+1),  \T)= \{ \chi_i - \chi_j: i \neq  j\}$.

The simple roots associated with $\B$ are $\alpha_i = \chi_i - \chi_{i+1}$, where $1 \leqslant i \leqslant d$.


Let $\N$ be the normalizer of $\T$ in $\SL(d+1)$. Then the Weyl group $\N(k) / \T(k)$ can be identified with the symmetric group $\mathfrak{S}_{d+1}$.

The parabolic subgroup $\rP$ of $\G$ of type $\delta$ containing $\B$ consists of upper triangular block matrices with a $(d) \times (d)$ block in the top left hand corner and a $(1) \times (1)$ block in the bottom right hand corner.
By definition, $\C_\delta(\rP)$ is the union of all Weyl cones $\mathfrak{C}(\B')$, where $\B'$ is a Borel subgroup with $\T \subset \B' \subset \rP$. Any Borel group $\B'$ containing $\T$ is of the form $n \B n\inv$ for some $n \in \N(k)$. It is contained in $\rP$ if and only if $n$ is contained in $\rP(k)$, which is equivalent to the fact that the permutation $\sigma \in \mathfrak{S}_{d+1}$ induced by $n$ fixes $d+1$. Since $$\mathfrak{C}(\B) = \{\chi_{i+1} - \chi_i \leqslant 1 \  : \ i = 1, \ldots, d\}$$ we deduce
$$\C_\delta(P) =  \{ \chi_{d+1} - \chi_i \leqslant 1: i= 1, \ldots, d\}.$$
Note that $\Phi({\rm rad}^{\rm u}(\rP^{\rm op}),\T) = \{\chi_{d+1} - \chi_i : i = 1,\ldots, d\}$, so that we recover the description from Proposition \ref{prop.roots.fan}, (ii).

\vspace{0.1cm}
If $\Qr$ is a $\delta$-relevant parabolic containing $\B$, it consists of upper triangular block matrices with a $(r) \times (r)$ block in the top left hand corner and a $(d+1-r) \times (d+1-r)$ block in the bottom right corner for some $r\geq 1$, cf. Example \ref{ex.relevant}. Hence we find $$\Phi(\Lr_{\Qr},\T) = \{ \chi_i - \chi_j ; i\neq j\mbox{ and } i,j \leq r\} \cup \{\chi_i - \chi_j ; i\neq j \mbox{ and } i,j > r\}$$ and  $$\mathfrak{C}(\Qr) = \{ \chi_{r+1} - \chi_r \leqslant 1 \textrm{ and }  \chi_{i+1} - \chi_i = 1 \mbox{ for all } i \leq r-1 \mbox{ and all }i \geq r+1\}.$$ The face of $\C_\delta(\rP)$ containing $\mathfrak{C}(\Qr)$ is
$$\C_\delta(\Qr) = \{ \chi_{d+1} - \chi_i \leq 1 \mbox{ for all } i \leq r \mbox{ and } \chi_{d+1} - \chi_i= 1 \mbox{ for all } i \geq r+1\}.$$
\end{Ex}

\vspace{0.3cm} We go a little further and establish a characterization of $t$-relevant parabolic subgroups in terms of the Dynkin diagram of $\G$. This will allow us to compare in \cite{RTW2} the prefan $\mathcal{F}_t$ with the collection of cones defined in \cite{Wer2} from the viewpoint of \emph{admissibility}

\vspace{0.1cm}
We fix a minimal parabolic subgroup $\rP_0$ of $\G$ containing $\Sr$ and let $\Delta$ denote the corresponding set of simple roots in $\Phi$. The map $\rP \mapsto \Y_{\rP} = \Phi(\Lr_{\rP}, \Sr) \cap \Delta$ sets up an increasing bijection between the set of parabolic subgroups of $\G$ containing $\rP_0$ and the power set of $\Delta$. The inverse bijection associates with a subset $\Y$ of $\Delta$ the parabolic subgroup $\rP$ of $\G$ containing $\rP_0$ such that
$$\Phi(\rP,\Sr) = \left(\langle \Y \rangle \cap \Phi \right) \cup \Phi(\rP_0, \Sr).$$

Equivalently, $\Y_{\rP}$ is the subset of $\Delta$ defining the face $\mathfrak{C}(\rP)$ of the cone $\mathfrak{C}(\rP_0)$:
$$\mathfrak{C}(\rP) = \{x \in \mathfrak{C}(\rP_0) \ | \ \alpha(x)=1, \ \alpha \in \Y_{\rP}\}.$$

\vspace{0.1cm}
The group $\X^*(\Sr)$ is equipped with a $\W$-invariant scalar product $( \cdot | \cdot )$ and we agree to see each finite subset $\E$ of $\X^*(\Sr)$ as a \emph{graph} by introducing an edge between any two vertices $\alpha, \ \beta \in \E$ if $(\alpha | \beta) \neq 0$.

\vspace{0.3cm}
\begin{Prop} \label{prop.simpleroots.fan} Let $\rP$ denote the unique parabolic subgroup of $\G$ of type $t$ containing $\rP_0$. For any parabolic subgroup $\Qr$ of $\G$ containing $\rP_0$, the following conditions are equivalent:
\begin{itemize}
\item[(i)] the cones $\C_t(\Qr)$ and $\mathfrak{C}(\Qr)$ have the same dimension;
\item[(ii)] the linear subspace $\{\alpha  = 1 \ ; \ \alpha \in \Y_{\Qr}\}$ of $\Lambda(\Sr)$ is the support of a face of the cone $\C_t(\rP)$, namely of $\C_t(\Qr)$;
\item[(iii)] seeing $\Delta$ as a graph following the convention above, each connected component of $\Y_{\Qr}$ meets $\Y_{\Qr} - \Y_{\Qr} \cap \Y_{\rP}$.
\end{itemize}
\end{Prop}

\vspace{0.1cm}
\noindent
\emph{\textbf{Proof}}.
Equivalence of conditions (i) and (ii) follows immediately from the fact that the cone $\mathfrak{C}(\Qr)$ spans the linear subspace $\{\alpha =1 \ ; \ \alpha \in \Phi(\Lr_{\Qr},\Sr)\}$ of $\Lambda(\Sr)$.

\vspace{0.1cm} (iii) $\Longrightarrow$ (ii) We assume that each connected component of $\Y_{\Qr}$ meets $\Y_{\Qr} - \Y_{\Qr} \cap \Y_{\rP}$ and establish the inclusion $\Y_{\Qr} \subset \{\alpha \in \Phi \ | \ \alpha_{|\C_t(\Qr)}=1\}$. Since $\Y_{\Qr}$ generates $\Phi(\Lr_{\Qr}, \Sr) = \{\alpha \in \Phi \ | \ \alpha_{|\mathfrak{C}(\Qr)} = 1\}$, it will follow that $\C_t(\Qr)$ and $\mathfrak{C}(\Qr)$ generates the same linear subspace of $\Lambda(\Sr)$.

We pick $\alpha \in \Y_{\Qr}$ and, up to replacing $\alpha$ by $-\alpha$, we assume that $\alpha$ belongs to $\Phi(\rP^{\rm op}, \Sr) = \Phi(\Lr_{\rP}, \Sr) \cup \Phi({\rm rad}^{\rm u}(\rP^{\rm op}),\Sr)$. The case $\alpha \in \Phi({\rm rad}^{\rm u}(\rP^{\rm op}),\Sr)$ is trivial: indeed, $\alpha$ cuts out a face of the cone $\C_t(\rP)$ by Lemma \ref{lemma.t-cone} (i) and this face contains $\mathfrak{C}(\Qr)$ since $\alpha_{|\mathfrak{C}(\Qr)} = 1$.

\vspace{0.1cm}
We now address the case $\alpha \in \Phi(\Lr_{\rP}, \Sr)$, i.e., $\alpha \in \Y_{\rP}$. Our assumption implies the existence of a natural integer $d$ and of roots $\alpha_0, \ldots , \alpha_d$ satisfying \begin{itemize}
\item  $\alpha_0 \in \Y_{\Qr} - \Y_{\Qr} \cap \Y_{\rP}$ and $\alpha_d = \alpha$; \item $\alpha_i \in \Y_{\rP} \cap \Y_{\Qr}$ for any $i \in \{1, \ldots, d-1\}$; \item $(\alpha_i | \alpha_{i+1}) <  0$ for any $i \in \{1, \ldots, d-1\}$ and $(\alpha_i | \alpha_j) = 0$ if $|i-j| \geqslant 2$.
\end{itemize}

In this situation the root $\beta = r_{\alpha_{d-1}} \circ \ldots \circ r_{\alpha_1} (\alpha_0)$ (and $\beta = \alpha_0$ if $d= 0$) is given by \begin{eqnarray*} \beta & = & \alpha_0 - 2\frac{(\alpha_0|\alpha_1)}{(\alpha_1|\alpha_1)} \alpha_1 + \ldots + (-2)^{d-1} \frac{(\alpha_0 | \alpha_1) \ldots (\alpha_{d-2}|\alpha_{d-1})}{(\alpha_{1}|\alpha_{1}) \ldots (\alpha_{d-1}|\alpha_{d-1})} \alpha_{d-1} \\ & = & \alpha_0 + m_1 \alpha_1 + \ldots + m_{d-1} \alpha_{d-1} \end{eqnarray*} with $m_1, \ldots, m_d \in \mathbb{Z} - \{0\}$. Since $\alpha_0$ belongs to $\Delta - \Y_{\rP} \subset \Phi(\radu (\rP^{\rm op}), \Sr)$ and $\alpha_1, \ldots, \alpha_{d-1} \in \Phi(\Lr_{\rP}, \Sr)$, this root belongs to $\Phi(\radu (\rP^{\rm op}), \Sr)$ and therefore $\beta$ cuts out a face of $\C_t(\rP)$. Moreover, since all the roots $\alpha_0, \ldots, \alpha_{d-1}$ belong to $\Y_{\Qr}$, $\beta_{|\mathfrak{C}(\Qr)} = 1$ and therefore $\beta_{|\C_t(\Qr)} = 1$ since $\C_t(\Qr)$ is the smallest face of $\C_t(\rP)$ containing $\mathfrak{C}(\Qr)$.

Now we have \begin{eqnarray*} r_{\alpha_d}(\beta) & = & \beta - 2 \frac{(\alpha_d|\beta)}{(\alpha_d|\alpha_d)} \alpha_d \\ & = & \beta + (-2)^d \frac{(\alpha_0|\alpha_1) \ldots (\alpha_{d-1}|\alpha_d)}{(\alpha_1|\alpha_1) \ldots (\alpha_d | \alpha_d)} \alpha_d \\ & = & \beta + m\alpha\end{eqnarray*} with $m \in \mathbb{Z} - \{0\}$. As before, this root belongs to $\Phi({\rm rad}^{\rm u}(\rP^{\rm op}),\Sr)$ and is identically equal to $1$ on $\mathfrak{C}(\Qr)$, hence $r_{\alpha_d}(\beta)_{|\C_t(\Qr)} = 1$. Since $m \neq 0$, we finally reach our goal: $\alpha_{|\C_t(\Qr)} = 1$.

\vspace{0.1cm} (ii) $\Longrightarrow$ (iii) We prove the converse assertion. Let $\Y$ denote the union of the connected components of $\Y_{\Qr}$ which meet $\Y_{\Qr} - \Y_{\Qr} \cap \Y_{\rP}$ and set $\Y' = \Y_{\Qr} - \Y$. The proof of (iii) $\Longrightarrow$ (ii) shows that $\{\alpha = 1 \ ; \ \alpha \in \Y\}$ is the support of a face $\F$ of $\C_t(\rP)$, namely of the cone $\C_t(\Qr')$, where $\Qr'$ is the parabolic subgroup containing $\Sr$ associated with the subset $\Y'$ of $\Delta$. We have moreover $\mathfrak{C}(\Qr) = \F \cap \{\beta =1 \ ; \ \beta \in \Y'\}$.

Suppose that $\alpha$ is a root belonging to $\Y'$. On the one hand, the hyperplane $\Hr_{\alpha} = \{\alpha = 1 \}$ does not contain $\F$ since the subset $\Y \cup \{\alpha\}$ of $\Delta$ consist of linearly independent roots. On the other hand, orthogonality of $\alpha$ and $\Y$ implies that the cone $\F$ is invariant under the orthogonal reflection with respect to $\Hr_{\alpha}$. Thus, if $\Y'$ is non-empty, the cone $\mathfrak{C}(\Qr) = \F \cap \bigcap_{\alpha \in \Y'} \Hr_{\alpha}$ meets the interior of $\F$, hence $\F$ is the smallest face of $\C_t(\rP)$ containing $\mathfrak{C}(\Qr)$ and therefore $$\dim \mathfrak{C}(\Qr) \leqslant \dim \C_t(\Qr) - 1 < \dim \C_t(\Qr).$$ \hfill $\Box$

\vspace{0.2cm} \begin{Cor} \label{cor.simpleroots.fan} Let $\rP$ denote the unique parabolic subgroup of $\G$ of type $t$ containing $\rP_0$ and let $\Qr$ be a parabolic subgroup of $\G$ containing $\rP_0$. The linear subspace spanned by $\C_t(\Qr)$ is defined by the conditions $\alpha = 1$, where $\alpha$ runs over all connected components of $\Y_{\Qr}$ meeting $\Delta - \Y_{\rP}$.
\end{Cor}

\vspace{0.1cm} \noindent \textbf{\emph{Proof}}. This assertion was proved while establishing (iii) $\Longrightarrow$ (ii) above. \hfill $\Box$

\vspace{0.3cm} Here is finally our root-theoretic characterization of $t$-relevant parabolic subgroups. We still denote by $\rP_0$ denote a minimal parabolic subgroup of $\G$ containing $\Sr$.

\begin{Prop} \label{prop.roots.relevant} Let $\rP$ denote the parabolic subgroup of $\G$ of type $t$ containing $\rP_0$. For any parabolic subgroup $\Qr$ of $\G$ containing $\rP_0$, we let $\widetilde{\Y_{\Qr}}$ denote the union of connected components of $\Y_{\Qr}$ meeting $\Delta - \Y_{\rP}$. Then the following conditions are equivalent:
\begin{itemize}
\item[(i)] $\Qr$ is $t$-relevant;
\item[(ii)] for any root $\alpha \in \Delta$, $$(\alpha \in \Y_{\rP} \ \textrm{ and } \ \ \alpha \perp \widetilde{\Y_{\Qr}}) \Longrightarrow \alpha \in \Y_{\Qr}.$$
\end{itemize}
\end{Prop}

\vspace{0.1cm}
\noindent
\emph{\textbf{Proof}}.
By definition, the parabolic subgroup $\Qr$ is $t$-relevant if and only if it is maximal among all parabolic subgroups $\Qr'$ of $\G$ satisfying $\mathrm{Osc}_t(\Qr) = \mathrm{Osc}_t(\Qr')$. We can obviously restrict to parabolic subgroups $\Qr'$ containing $\Sr$, in which case we proved in Proposition \ref{prop.osc.bigcell}, (ii) that the latter condition amounts to $\Phi(\radu (\rP),\Sr) \cap \Phi(\Lr_{\Qr}, \Sr) = \Phi(\radu (\rP), \Sr) \cap \Phi(\Lr_{\Qr'}, \Sr)$, or equivalently to $\C_t(\Qr) = \C_t(\Qr')$ by application of Proposition \ref{prop.roots.fan}, (iii). It follows that the parabolic subgroup $\Qr$ is $t$-relevant if and only if, for any root $\alpha \in \Delta - \Y_{\Qr}$, the parabolic subgroup $\Qr_{\alpha}$ associated with the subset $\Y_{\Qr_{\alpha}} = \Y_{\Qr} \cup \{\alpha\}$ of $\Delta$ satisfies $\C_t(\Qr_{\alpha}) \varsubsetneq \C_t(\Qr)$.

We consider a root $\alpha$ in $\Delta - \Y_{\Qr}$ and let $\widetilde{\Y_{\Qr}}$ ($\widetilde{\Y_{\Qr_{\alpha}}}$, respectively) denote the union of the connected components of $\Y_{\Qr}$ (of $\Y_{\Qr_{\alpha}}$, respectively) meeting $\Phi - \Y_{\rP}$. The conditions $\C_t(\Qr) = \C_t(\Qr_{\alpha})$ and $\widetilde{\Y_{\Qr}} = \widetilde{\Y_{\Qr_{\alpha}}}$ are equivalent by Corollary \ref{cor.simpleroots.fan} and one immediately checks that the identity $\widetilde{\Y_{\Qr}} = \widetilde{\Y_{\Qr_{\alpha}}}$ amounts to orthogonality of $\alpha$ and $\widetilde{\Y_{\Qr}}$. Therefore, the parabolic subgroup $\Qr$ is $t$-relevant if and only if there is no root in $\Y_{\rP} - \Y_{\rP} \cap \Y_{\Qr}$ orthogonal to each connected component of $\Y_{\Qr}$ meeting $\Delta - \Y_{\rP}$. \hfill $\Box$

\vspace{0.3cm} \begin{Rk} \label{rk.relevant} 1. Letting $\widetilde{\Y_{\Qr}}$ denote the union of the connected components of $\Y_{\Qr}$ meeting $\Delta - \Y_{\rP}$, condition (ii) above is equivalent to the following requirement: the complement of $\Y_{\Qr}$ in $\Y_{\rP}$ consists only of roots whose distance to $\widetilde{\Y_{\Qr}}$ is at most one.

2. Given a parabolic subgroup $\Qr$ of $\G$ containing $\rP_0$, the smallest $t$-relevant parabolic subgroup of $\G$ containing $\Qr$ corresponds to the subset of $\Delta$ deduced from $\Y_{\Qr}$ by adjoining all the roots in $\Y_{\rP}$ which are orthogonal to each connected component of $\Y_{\Qr}$ meeting $\Delta - \Y_{\rP}$.

3. For any parabolic subgroup $\Qr$ of $\G$ containing $\Sr$, the smallest $t$-relevant parabolic subgroup of $\G$ containing $\Qr$ coincides with the largest parabolic subgroup $\Qr'$ of $\G$ containing $\Qr$ such that $\C_t(\Qr') = \C_t(\Qr)$.
\end{Rk}

\vspace{0.3cm}
\begin{Cor} \label{cor.roots.decomposition} For any parabolic subgroup $\Qr$ of $\G$ containing $\Sr$, both $$\{\alpha \in \Phi(\Lr_{\Qr},\Sr) \ | \ \alpha \textrm{ vanishes identically on } \C_t(\Qr)\}$$ and its complement are closed subsets of $\Phi(\Lr_{\Qr},\Sr)$.

Moreover, if we let $\Qr'$ denote the smallest $t$-relevant parabolic subgroup of $\G$ containing $\Qr$, then $\Phi(\Lr_{\Qr},\Sr) \subset \Phi(\Lr_{\Qr'},\Sr)$ and $$\{\alpha \in \Phi(\Lr_{\Qr'},\Sr)\ ; \ \alpha \textrm{ vanishes identically on } \C_t(\Qr')\} = \{\alpha \in \Phi(\Lr_{\Qr},\Sr)\ ; \ \alpha \textrm{ vanishes identically on } \C_t(\Qr)\}.$$
\end{Cor}

\vspace{0.1cm}
\noindent \emph{\textbf{Proof}}. Let $\Sigma$ denote the set of roots in $\Phi(\Lr_{\Qr},\Sr)$ which vanish identically on the cone $\C_t(\Qr)$; this is obviously a closed subset of $\Phi(\Lr_{\Qr},\Sr)$. We consider now a minimal parabolic subgroup $\rP_0$ of $\G$ containing $\Sr$ and contained in $\Qr$, and we let $\Delta$ denote the corresponding set of simple roots in $\Phi(\G,\Sr)$. By Corollary \ref{cor.simpleroots.fan}, $\Sigma \cap \Delta$ is an union of connected component of $\Phi(\Lr_{\Qr},\Sr) \cap \Delta$, thus $$ \Phi(\Lr_{\Qr},\Sr) \cap \Delta = (\Sigma \cap \Delta) \cup (\Sigma^{\mathrm{c}} \cap \Delta)$$ is a decomposition of $\Phi(\Lr_{\Qr},\Sr) \cap \Delta$ into mutually orthogonal subsets. It follows that $\Phi(\Lr_{\Qr},\Sr)$ is the disjoint union of the closed subsets $\R$ and $\R'$ respectively spanned by $\Sigma \cap \Delta$ and $\Sigma^{\mathrm{c}} \cap \Delta$. Since any root in $\Sigma$ is a linear combination of roots of $\Sigma \cap \Delta$, $\Sigma = \R$ and therefore $\Sigma^{\mathrm{c}} = \R'$ is closed.

\vspace{0.1cm} The second assertion follows immediately of Corollary \ref{cor.simpleroots.fan} and Remark \ref{rk.relevant}, 2.\hfill $\Box$

\vspace{0.3cm}

\begin{Ex}
\label{ex - relevant SL}
We use the notation of Example \ref{ex-fan}. The Dynkin diagram of ${\rm SL}(d+1)$ is the graph
$$\xymatrix{\circ  \ar@{-}[r]_{\hspace{-1.2cm} \alpha_{1}} & \circ \ar@{-}[r]_{\hspace{-1.4cm} \alpha_{2}} & \cdot \ \cdot \ \cdot &  \circ \ar@{-}[l]^{\hspace{1.4cm} \alpha_{d-1}} \ar@{-}[r]_{\hspace{1.4cm} \alpha_d} & \circ }.$$

For any proper parabolic subgroup $\Qr$ of $\G$ containing $\B$ and not contained in $\rP$, the only connected component of $\Y_{\Qr}$ meeting $\Delta - \Y_{\rP} = \{\alpha_{d}\}$ is $\widetilde{\Y_{\Qr}} = \{\alpha_{\ell+1}, \ldots, \alpha_{d}\}$, where $\ell$ is greatest index $i$ such that $\alpha_i \notin \Y_{\Qr}$. The roots in $\Y_{\rP} = \{\alpha_1, \ldots, \alpha_{d-1}\}$ which are orthogonal to $\widetilde{\Y_{\Qr}}$ are  $\alpha_1, \ldots, \alpha_{\ell-1}$. They are all contained in $\Y_{\Qr}$ if and only if
\[ \Y_{\Qr} = \Delta - \{\alpha_{\ell}\}, \]
or equivalently if $\Qr$ is the stabilizer of the linear subspace $\mathrm{Span}(e_1, \ldots, e_{\ell})$. Applying proposition \ref{prop.roots.relevant}, we thus recover the description of $\delta$-relevant parabolic subgroups of $\G = \SL (d+1)$ given in 3.2.1, Example \ref{ex.relevant}.
\end{Ex}

\subsection{Berkovich compactifications}

From now on, we work again under the assumptions of (1.3.4).

\vspace{0.2cm}
\noindent \textbf{(3.4.1)} Let $t$ denote a $k$-rational type of $\G$. We consider the central isogeny $\G' \times \G'' \rightarrow \G$ introduced after Definition \ref{def.t-cone}, which induces identifications $\mathcal{B}(\G,k) = \mathcal{B}(\G',k) \times \mathcal{B}(\G'',k)$, ${\rm Par}(\G) = {\rm Par}(\G') \times {\rm Par}(\G'')$ and ${\rm Par}_{t}(\G) = {\rm Par}_{t'}(\G')$, where $t'$ denotes the restriction of $t$ to $\G'$.

Moreover, we let $p'$ denote the canonical projection of $\Par(\G') \times \Par(\G'')$ on $\Par(\G')$ and $j$ the closed immersion $\Par(\G') \hookrightarrow \Par(\G)$ defined (functor-theoretically) by $\rP' \mapsto i(\rP' \times \G'')$.

\vspace{0.3cm} \begin{Lemma} \label{lemma.type.reduction} With the notation and convention introduced above, the diagram $$\xymatrix{\mathcal{B}(\G,k) \ar@{->}[r]^{\vartheta_t} \ar@{->}[d]_{p'} & \mathrm{Par}(\G)^{an} \ar@{<-}[d]^j \\ \mathcal{B}(\G',k) \ar@{->}[r]_{\vartheta_{t'}} & {\rm Par}(\G')^{an}}$$ is commutative.
\end{Lemma}

\vspace{0.1cm}
\noindent
\emph{\textbf{Proof}}.
Relying on Proposition \ref{prop.theta-t.extension}, it is enough to prove that this diagram is commutative after replacing $k$ by a non-Archimedean extension. Hence we can assume that $\G$ is split --- then $\G'$ and $\G''$ are also split --- and we may restrict to check that the maps $j \circ \vartheta_{t'} \circ p'$ and $\vartheta_t$ coincide on the set of special vertices of $\mathcal{B}(\G,k)$.

\vspace{0.1cm}
The diagram under consideration can be decomposed in four diagrams $$\xymatrix{\mathcal{B}(\G' \times_k \G'',k) \ar@{->}[r]^{\vartheta} \ar@{=}[d] & (\G' \times_k \G'')^{\mathrm{an}} \ar@{->}[d]^{i} \\ \mathcal{B}(\G,k) \ar@{->}[r]_{\vartheta} & \G^{\mathrm{an}}} \hspace{0.5cm} \xymatrix{\mathcal{B}(\G' \times_k \G'', k) \ar@{->}[r]^{\vartheta} \ar@{->}[d]^{p'} & (\G' \times_k \G'')^{\mathrm{an}} \ar@{->}[d]^{p'} \\ \mathcal{B}(\G',k) \ar@{->}[r]_{\vartheta} & {\G'}^{\mathrm{an}}}$$
$$\xymatrix{(\G' \times_k \G'')^{\mathrm{an}} \ar@{->}[r]^{\iota_{\rP' \times \G''}} \ar@{->}[d]_{i} & \Par (\G' \times_k \G'')^{\mathrm{an}} \ar@{<-}[d]^{i^{*}} \\ \G^{\mathrm{an}} \ar@{->}[r]_{\iota_{\rP}} & \Par(\G)^{\mathrm{an}}} \hspace{0.5cm} \xymatrix{(\G' \times_k \G'')^{\mathrm{an}} \ar@{->}[r]^{\iota_{\rP' \times \G''}} \ar@{->}[d]_{p'} & \Par(\G' \times \G'')^{\mathrm{an}} \ar@{->}[d]^{j} \\ {\G'}^{\mathrm{an}} \ar@{->}[r]_{\iota_{\rP'}} & \Par(\G')^{\mathrm{an}}}$$ where $\rP'$ and $\rP$ are elements of $\Par_{t'}(\G',k)$ and $\Par_t(\G)(k)$ respectively satisfying $i^{-1}(\rP) = \rP' \times_k \G''$. It suffices to check that each of these four diagrams is commutative.

This is obviously true for the last two.

Consider a special point $o$ in $\mathcal{B}(\G,k)$, whose associated $k^{\circ}$-Chevalley group we denote $\mathcal{G}$. We may find two $k^{\circ}$-Chevalley groups $\mathcal{G}'$ and $\mathcal{G}''$ with generic fibres $\G'$ and $\G''$ respectively, such that the isogeny $i: \G' \times_k \G'' \rightarrow \G$ extends to a $k^{\circ}$-isogeny $\mathcal{G}' \times_{k^{\circ}} \mathcal{G}'' \rightarrow \mathcal{G}$ (this follows from the equivalence between the category of split reductive groups over $k$ equipped with a splitting datum and the category of root data \cite[Expos\'e XXIII, Th\'eor\`eme 4.1]{SGA3}, together with the fact that any isogeny extends to an isogeny of splitting data \cite[Expos\'e XXII, Corollaire 4.2.3]{SGA3}). These Chevalley groups correspond to special points $o'$ and $o''$ in $\mathcal{B}(\G',k)$ and $\mathcal{B}(\G'',k)$, and the bijection between $\mathcal{B}(\G' \times_k \G'', k) = \mathcal{B}(\G',k) \times \mathcal{B}(\G'',k)$ and $\mathcal{B}(\G,k)$ induced by $i$ maps $(o',o'')$ to $o$. The commutativity of the first two diagrams now follows from the very definition of the map $\vartheta$  together with observation that the isogeny $\mathcal{G}' \times_{k^{\circ}} \mathcal{G}'' \rightarrow \mathcal{G}$ induces a finite morphism between special fibres and thus maps the generic point to the generic point. \hfill $\Box$

\vspace{0.3cm} Replacing the group $\G$ by the normal subgroup $\G'$ and the building $\mathcal{B}(\G,k)$ by its factor $\mathcal{B}(\G',k)$, we may use the lemma above to reduce the study of the map $\vartheta_t$ to the case of a non-degenerate $k$-rational type $t$, i.e., a $k$-rational type whose restriction to any quasi-simple component of $\G$ is non-trivial.

If $\G$ is split and $\Sr$ denotes a split maximal torus, a parabolic subgroup $\rP$ of $\G$ containing $\Sr$ is non-degenerate if and only if the set of roots of $\radu (\rP)$ with respect to $\Sr$ spans a subgroup of finite index in the character group $\X^*(\Sr)$ of $\Sr$ (cf. Remark \ref{rk.strict.convexity}).

\vspace{0.3cm}
\begin{Prop} \label{prop.theta-t.injective} If the $k$-rational type $t$ is non-degenerate, the map $$\vartheta_t: \mathcal{B}(\G,k) \rightarrow \mathrm{Par}(\G)^{\mathrm{an}}$$ is injective.
\end{Prop}

\vspace{0.1cm}
\noindent
\emph{\textbf{Proof}}.
By Proposition \ref{prop.theta-t.extension}, we may assume the group $\G$ to be split.

Given a split maximal torus $\Sr$ of $\G$, it follows from the explicit formula established in Proposition \ref{prop.theta-t.explicit} that the map $\vartheta_t$ is injective on the apartment $\A(\Sr,k)$. Indeed, having identified $\A(\Sr,k)$ with the vector space $\Lambda(\Sr)$ of real linear forms on $\X^*(\Sr)$, two linear forms $u, \ v \in \Lambda(\Sr)$ satisfying $\vartheta_t(u) = \vartheta_t(v)$ coincide on the subset $\Phi(\radu (\rP), \Sr)$ of $\X^*(\Sr)$, where $\rP$ denotes any parabolic subgroup of $\G$ of type $t$ containing $\Sr$. Since the type $t$ is non-degenerate, $\Phi(\radu (\rP), \Sr)$ spans $\X^*(\Sr) \otimes_{\mathbb{Z}} \mathbb{Q}$, hence $u=v$.

\vspace{0.1cm} Injectivity of $\vartheta_t$ on the whole building follows from the fact that any two points are contained in a common apartment.\hfill $\Box$

\vspace{0.8cm}
\noindent \textbf{(3.4.2)} We fix a $k$-rational type $t$ of (parabolic subgroups of) $\G$.

\vspace{0.1cm} For any maximal split torus $\Sr$ of $\G$, we let $\overline{\A}_t(\Sr,k)$ denote the closure of $\vartheta_t\left(\A(\Sr,k)\right)$ in $\Par(\G)^{\mathrm{an}}$ endowed with the induced topology. This is a compact topological space to which we refer as the \emph{compactified apartment} of type $t$ of $\Sr$.

We let $\overline{\mathcal{B}}_t(\G,k)$ denote the image of the map $$\G(k) \times \overline{\A}_t(\Sr,k) \rightarrow \Par(\G)^{\mathrm{an}}, \  (g,x) \mapsto gxg^{-1},$$ which we endow with the \emph{quotient} topology. Set-theoretically, $\overline{\mathcal{B}}_t(\G,k)$ is the union of all compactified apartments of type $t$ in $\Par(\G)^{\mathrm{an}}$.

\vspace{0.3cm} \begin{Def} \label{def.compact-t} The $\G(k)$-topological space $\overline{\mathcal{B}}_t(\G,k)$ is the \emph{Berkovich compactification of type $t$} of the building $\mathcal{B}(\G,k)$.
\end{Def}

\vspace{0.3cm} \begin{Rk} It is somehow incorrect to use the word "compactification" in this context for two reasons: \begin{itemize} \item if the type $t$ is degenerate, the map $\vartheta_t$ is not injective; \item if the field $k$ is not locally compact, the topological space $\overline{\mathcal{B}}_t(\G,k)$ is not compact.
\end{itemize}
However, the image of $\vartheta_t : \mathcal{B}(\G,k) \rightarrow \overline{\mathcal{B}}_t(\G,k)$ is obviously dense and we shall prove later (Proposition \ref{prop.compact.theta-t}) that this map is open.
\end{Rk}

\vspace{0.2cm} Functoriality with respect to the field extends to the compactifications.

\vspace{0.1cm} \begin{Prop} \label{prop.compact.extension} Let $k'/k$ be a non-Archimedean extension.

\begin{itemize}
\item[(i)] There exists a unique continuous map $\overline{\mathcal{B}}_t(\G,k) \rightarrow \overline{\mathcal{B}}_t(\G,k')$ extending the canonical injection of $\mathcal{B}(\G,k)$ into $\mathcal{B}(\G,k')$. This map is a $\G(k)$-equivariant homeomorphism onto its image.

\item[(ii)] If the field $k$ is locally compact, the image of $\overline{\mathcal{B}}_t(\G,k)$ in $\overline{\mathcal{B}}_t(\G,k')$ is closed.
\end{itemize}
\end{Prop}

\vspace{0.1cm}
\noindent
\emph{\textbf{Proof}}.
(i) There exists clearly at most one continuous extension $\overline{\mathcal{B}}_t(\G,k) \rightarrow \overline{\mathcal{B}}_t(\G,k')$ of the canonical injection $\mathcal{B}(\G,k) \hookrightarrow \mathcal{B}(\G,k')$ since the image of $\mathcal{B}(\G,k)$ in $\overline{\mathcal{B}}_t(\G,k)$ is dense.

For any maximal split torus $\Sr$ of $\G$, we set $\A_t(\Sr,k') = \vartheta_t(\A(\Sr,k))$ and let $\overline{\A}_t(\Sr,k')$ denote its closure in ${\rm Par}(\G \otimes_k k')^{\rm an}$. We recall that there exists a torus $\T$ of $\G$ satisfying the following conditions: \begin{itemize} \item $\T$ contains $\Sr$; \item $\T \otimes_k k'$ is a maximal split torus of $\G \otimes_k k'$ ; \item the injection of $\mathcal{B}(\G,k) \hookrightarrow \mathcal{B}(\G,k')$ maps $\A(\Sr,k)$ into $\A(\T,k') = \A(\T \otimes_k k',k')$. \end{itemize}
Equivalently, $\overline{\A}_t(\Sr,k')$ is the closure of $\vartheta_t(\A(\Sr,k)$ in $\overline{\A}(\T, k')$.

Relying on the commutativity of the diagram $$\xymatrix{\mathcal{B}(\G,k') \ar@{->}[r]^{\vartheta_t} \ar@{<-}[d] & \Par(\G)^{\mathrm{an}} \widehat{\otimes}_k k' \ar@{->}[d]^{{\rm pr}_{k'/k}} \\ \mathcal{B}(\G,k) \ar@{->}[r]_{\vartheta_t} & \Par(\G)^{\mathrm{an}}}.$$
(see Proposition \ref{prop.theta-t.extension}), it suffices to prove that the canonical projection ${\rm pr}_{k'/k}: \Par(\G)^{\mathrm{an}} \widehat{\otimes}_k k' \rightarrow \Par(\G)^{\mathrm{an}}$ induces a homeomorphism between $\overline{\A}_t(\Sr,k')$ and $\overline{\A}_t(\Sr, k)$, as well as a bijection between $\bigcup_{\Sr} \overline{\A}_t(\Sr,k')$ and $\bigcup_{\Sr}\overline{\A}_t(\Sr,k)$, and to consider the inverse bijection.

\vspace{0.1cm} We split the proof in two steps.

\vspace{0.1cm}
\emph{First step: $k'/k$ is a finite Galois extension}. Set $\Gamma = \mathrm{Gal}(k'/k)$. In this case, it suffices to note that the projection ${\rm pr}_{k'/k}$ induces a homeomorphism between the closed $\Gamma$-fixed point subspace $(\Par(\G)^{\mathrm{an}} \otimes_k k')^{\Gamma}$ of $\Par(\G)^{\mathrm{an}} \otimes_k k'$ and its image in $\Par(\G)^{\mathrm{an}}$, the latter being closed since the map ${\rm pr}_{k'/k}$ is closed. Since the maps $\mathcal{B}(\G,k) \hookrightarrow \mathcal{B}(\G,k')$ and $\vartheta_t: \mathcal{B}(\G,k') \rightarrow \Par(\G)^{\mathrm{an}} \otimes_k k'$ are $\Gamma$-equivariant, ${\rm pr}_{k'/k}$ induces therefore a homeomorphism between $\overline{\A}_t(\Sr,k')$ and $\overline{\A}_t(\Sr,k)$, as well as a bijection between $\bigcup_{\Sr}\overline{\A}_t(\Sr,k')$ and $\bigcup_{\Sr} \overline{\A}_t(\Sr,k)$.

\vspace{0.1cm}
\emph{Second step: the group $\G$ is split}. In this case, the result will follow from the construction of a continuous section $\sigma$ of ${\rm pr}_{k'/k}$ over $\overline{\A}_t(\Sr,k)$ mapping $\A_t(\Sr,k)$ onto $\A_t(\Sr,k')$. We rely on the explicit formula established in Proposition \ref{prop.theta-t.explicit} to define $\sigma$ and we use the notation introduced there. First note that each point $x$ of $\overline{\A}_t(\Sr,k)$ belongs to the open subset $\Omega(\rP,\Sr)^{\mathrm{an}}$ of $\Par(\G)^{\mathrm{an}}$ for a convenient choice of the parabolic subgroup $\rP$ of $\G$ containing $\Sr$ (cf. Remark \ref{rk.big-cell.charts}). Then this point corresponds to the multiplicative seminorm on the $k$-algebra $k\left[(\X_{\alpha})_{\alpha \in \Psi}\right]$ defined by $$f = \sum_{\nu \in \mathbb{N}^{\Psi}} a_{\nu} \X^{\nu} \mapsto |f|(x) = \max_{\nu} |a_{\nu}| \prod_{\alpha \in \Psi} |\X_{\alpha}|(x)^{\nu(\alpha)}$$ since the function $|f| - \max_{\nu} |a_{\nu}| \prod_{\alpha \in \Psi} |\X_{\alpha}|^{\nu(\alpha)}$ is continuous on $\Par(\G)^{\mathrm{an}}$ and vanishes identically on $\A_t(\Sr,k) = \vartheta_t(\A(\Sr,k))$. We define $\sigma(x)$ as the point in $\Omega(\rP,\Sr)^{\mathrm{an}} \widehat{\otimes}_{k} k'$ corresponding to the multiplicative seminorm on $k'[(\X_{\alpha})_{\alpha \in \Psi}]$ satisfying the same identity: $$\left| \sum_{\nu \in \mathbb{N}^{\Psi}} a_{\nu} \X^{\nu} \right|(\sigma(x)) = \max_{\nu} |a_{\nu}| \prod_{\alpha \in \Psi} |\X_{\alpha}|(x)^{\nu(\alpha)}.$$

The map $\sigma: \overline{\A}_t(\Sr,k) \rightarrow \Par(\G)^{\mathrm{an}} \widehat{\otimes}_k k'$ is a continuous section of the projection ${\rm pr}_{k'/k}$ mapping  $\A_t(\Sr,k)$ onto $\A_t(\Sr,k')$ by (2.4.2), Proposition \ref{prop.theta-t.explicit}.
Thus, the maps ${\rm pr}_{k'/k}$ and $\sigma$ induce continuous and mutually inverse bijections between the sets $\overline{\A}_t(\Sr,k')$ and $\overline{\A}_t(\Sr,k)$, and our assertion follows.

\vspace{0.1cm} (ii) For any maximal split torus $\Sr$ of $\G$ and any point $x$ in $\A(\Sr,k)$, $\mathcal{B}(\G,k) = \G_x(k)\A(\Sr,k)$ (see reminders of Bruhat-Tits theory in (1.3.3)), hence the image of $\mathcal{B}(\G,k)$ in $\overline{\mathcal{B}}_t(\G,k')$ is contained in the subspace $\F = \G_x(k) \overline{\A}_t(\Sr,k')$. If the field $k$ is locally compact, the topological group $\G_x(k)$ is compact and therefore $\F$, like $\overline{\A}_t(\Sr,k')$, is a closed subset of $\overline{\mathcal{B}}_t(\G,k')$. It follows that $\F$ contains the closure $\F'$ of $\mathcal{B}(\G,k)$ in $\overline{\mathcal{B}}_t(\G,k')$. Since $\F'$ contains $\overline{\mathcal{B}}_t(\G,k) = \G(k) \overline{\A}_t(\Sr,k')$, we see finally that $\overline{\mathcal{B}}_t(\G,k)$ is the closure of $\mathcal{B}(\G,k)$ in $\overline{\mathcal{B}}_t(\G,k')$. \hfill $\Box$

\vspace{0.3cm}
\vspace{0.1cm} \begin{Lemma} \label{lemma.compact.apartment} Let $\Sr$ be a maximal split torus of $\G$ and $x$ a point in the compactified apartment $\overline{\A}_t(\Sr,k)$. If there exists an element $g$ in $\G(k)$ such that $gx$ belongs to $\mathcal{B}_t(\G,k)$, then $x$ belongs the subspace $\A_t(\Sr,k) = \vartheta_t(\A(\Sr,k))$.
\end{Lemma}

\vspace{0.1cm}
\noindent
\emph{\textbf{Proof}}. We can restrict to a non-degenerate type by Lemma \ref{lemma.type.reduction}.

We first assume that $\G$ is split and rely in this case on the explicit formula in Proposition \ref{prop.theta-t.explicit}. There exists a parabolic subgroup $\rP$ of $\G$ containing $\Sr$ such that the point $x$ of $\Par(\G)^{\mathrm{an}}$ belongs to the big cell $\Omega(\rP,\Sr)^{\mathrm{an}}$ and corresponds to a multiplicative seminorm on the $k$-algebra $k[(\X_{\alpha})_{\alpha \in \Psi}]$, where $\Psi = \Phi({\rm rad}^{\rm u}(\rP^{\rm op}), \Sr)$. As noticed in the proof of Proposition \ref{prop.compact.extension}, the explicit formula for a point $x$ lying in the image of $\A(\Sr,k)$ holds more generally for any point in $\overline{\A}_t(\Sr,k)$: $$\left|\sum_{\nu} a_{\nu} \X^{\nu}\right|(x) = \max_{\nu}|a_{\nu}| \prod_{\alpha \in \Psi} |\X_{\alpha}|(x)^{\nu(\alpha)}.$$ One easily checks that the point $x$ belongs to (the image of) $\A(\Sr,k)$ if and only if $|\X_{\alpha}|(x)>0$ for any root $\alpha \in \Phi(\rP,\Sr)$, which amounts to requiring that this seminorm is in fact a norm, or equivalently that $|f|(x) > 0$ for any non-zero germ $f \in \mathcal{O}_{\Par(\G)^{\mathrm{an}}, x}$.

If there exists an element $g$ of $\G(k)$ such that the point $gx$ belongs to $\overline{\mathcal{B}}_t(\G,k) - \mathcal{B}(\G,k)$, then $gx$ belongs to $\overline{\A}_t(\Sr',k) - \A(\Sr',k)$ for some maximal torus $\Sr'$ of $\G$ and there exists therefore a non-zero germ $f' \in \mathcal{O}_{\Par(\G)^{\mathrm{an}}, gx}$ satisfying $|f'|(gx) = 0$. Then we have $|g^*f'|(x) = |f'|(gx) = 0$ and, since $g^*f'$ is a non-zero germ at $x$, the point $x$ belongs to $\overline{\A}_t(\Sr,k) - \A(\Sr,k)$.

\vspace{0.1cm} We now address the general case. Let $k'/k$ be finite Galois extension splitting $\G$ and consider a maximal torus $\T$ of $\G$ satisfying the following conditions: \begin{itemize} \item[(a)] $\T$ contains $\Sr$; \item[(b)] $\T \otimes_k k'$ is split; \item[(c)] the injection $\mathcal{B}(\G,k) \rightarrow \mathcal{B}(\G,k')$ identifies $\A(\Sr,k)$ with the Galois-fixed point set of $\A(\T,k') = \A(\T \otimes_k k', k')$. \end{itemize} It follows from Proposition \ref{prop.compact.extension} and continuity that the compactified apartment $\overline{\A}_t(\Sr,k)$ is identified with the Galois-fixed point set in $\overline{\A}_t(\T, k')$. If $x$ is a point of $\overline{\A}_t(\Sr,k)$ whose $\G(k)$ orbits meets $\mathcal{B}_t(\G,k)$, then we know that $x$ belongs to $\A(\T, k')$, and therefore to $\A(\Sr,k)$ since $x$ is Galois-fixed.
\hfill $\Box$

\vspace{0.3cm} \begin{Prop} \label{prop.compact.theta-t} The map $\vartheta_t: \mathcal{B}(\G,k) \rightarrow \overline{\mathcal{B}}_t(\G,k)$ is continuous, open and its image is dense.

This map is injective if and only if the type $t$ is non-degenerate. Finally, if the field $k$ is locally compact, the topological space $\overline{\mathcal{B}}_t(\G,k)$ is compact.
\end{Prop}

\vspace{0.1cm}
\noindent
\emph{\textbf{Proof}}.
Continuity of $\vartheta_t$ and density of its image follow immediately from the definition of $\overline{\mathcal{B}}_t(\G,k)$. By Proposition \ref{prop.theta-t.injective}, the map $\vartheta_t$ is injective if the type $t$ is non-degenerate; conversely, if $t$ is degenerate, then $\vartheta_t$ is not injective by Lemme \ref{lemma.type.reduction}.

It remains to check that this map is open. Let us consider the following commutative diagram $$\xymatrix{\G(k) \times \A(\Sr,k) \ar@{->}[d]_{\pi} \ar@{->}[r]^{\mathrm{id} \times \vartheta_t} & \G(k) \times \overline{\A}_t(\Sr,k) \ar@{->}[d]^{\pi} \\ \mathcal{B}(\G,k) \ar@{->}[r]_{\vartheta_t} & \overline{\mathcal{B}}_t(\Sr,k)}$$ associated with a maximal split torus $\Sr$ of $\G$, where the maps $\pi$ are defined by $\pi(g,x) = gx$. Given an open subset $\U$ in $\mathcal{B}(\G,k)$, $\V = (\mathrm{id} \times \vartheta_t)(\pi^{-1}(\U))$ is an open subset of $\G(k) \times \overline{\A}_t(\Sr,k)$. This is moreover a $\pi$-saturated subset, since any point $(g,x) \in \G(k) \times \overline{\A}_t(\Sr,k)$ such that $gx$ belongs to $\vartheta_t(\U)$ is contained in the image of $\mathrm{id} \times \vartheta_t$ by Lemma \ref{lemma.compact.apartment}, hence in $(\mathrm{id} \times \vartheta_t)(\pi^{-1}(\U))$. Since $\V$ meets each fibre of $\pi$ over $\vartheta_t(\U)$, $\V = \pi^{-1}(\vartheta_t(\U))$ and thus $\vartheta_t(\U)$ is open in $\overline{\mathcal{B}}_t(\G,k)$, for the map $\pi$ is open and surjective.

\vspace{0.1cm} If the field $k$ is locally compact, the space $\overline{\mathcal{B}}_t(\G,k)$ is compact by the same argument as for Proposition \ref{prop.compact.extension} (ii).\hfill $\Box$

\vspace{0.8cm}
\noindent \textbf{(3.4.3)} Let $t$ denote a $k$-rational type of $\G$ and $\Sr$ a maximal split torus. We prove in this paragraph that the compactified apartment $\overline{\A}_t(\Sr,k)$, defined as the closure of $\vartheta_t(\A(\Sr,k))$ in $\Par(\G)^{\mathrm{an}}$, coincides with the compactification of the apartment $\A(\Sr,k)$ associated with the prefan $\mathcal{F}_t$ on $\Lambda(\Sr)$.

\vspace{0.3cm}
\begin{Prop}  \label{prop.apartment.identification} The map $\vartheta_t: \A(\Sr,k) \rightarrow \mathrm{Par}_t(\G)^{\rm an}$ extends to a homeomorphism $$\overline{\A(\Sr,k)}^{\mathcal{F}_t} \xymatrix{{} \ar@{->}[r]^{\sim} & {}} \overline{\A}_t(\Sr,k).$$
\end{Prop}

\vspace{0.1cm}
\noindent
\emph{\textbf{Proof}}.
We split the proof in three steps, the third one consisting in Lemma \ref{lemma.reduction} below.\hfill $\Box$

\vspace{0.1cm}
\emph{First step: reduction to the split case}. We pick a finite Galois extension $k'/k$ splitting $\G$ and choose a maximal torus $\T$ of $\G$ satisfying the following conditions: \begin{itemize} \item[(a)] $\T$ contains $\Sr$; \item[(b)] $\T \otimes_k k'$ is split; \item[(c)] the injection $\mathcal{B}(\G,k) \rightarrow \mathcal{B}(\G,k')$ maps $\A(\Sr,k)$ into $\A(\T,k') = \A(\T \otimes_k k', k')$. \end{itemize} It follows from Proposition \ref{prop.compact.extension} that $\overline{\A}_t(\Sr,k)$ is identified with the closure of $\A(\Sr,k)$ in $\overline{\A}_t(\T,k')$.

Let $\rP$ denote a parabolic subgroup of $\G$ of type $t$ containing $\Sr$. Set $\T' = \T \otimes_k k'$, $\Sr' = \Sr \otimes_k k'$, $\rP' = \rP \otimes_k k'$ and let $\lambda$ denote the homomorphism $\X^*(\T') \rightarrow \X^*(\Sr') = \X^*(\Sr), \ \alpha \mapsto \alpha_{|\Sr'}$. The cone $\C_t(\rP')$ in $\Lambda(\T') = {\rm Hom}_{\mathbf{Ab}}(\X^*(\T'), \mathbb{R}_{>0})$ is defined by the inequalities $\alpha \leqslant 1$, $\alpha \in \Phi(\radu ({\rP'}^{\rm op}), \T')$ and $(\lambda^{\vee})^{-1} \C_t(\rP')$ is therefore the cone in $\Lambda(\Sr) = {\rm Hom}_{\mathbf{Ab}}(\X^*(\Sr),\mathbb{R}_{>0})$ defined by the inequalities $\alpha \leqslant 1$, $\alpha \in \Phi(\radu ({\rP'}^{\rm op}), \T')$. Since $\radu ({\rP'}^{\rm op}) = \radu (\rP^{\rm op}) \otimes_k k'$,
$$\Phi(\radu ({\rP}^{\rm op}), \Sr) \subset \lambda \left(\Phi(\radu ({\rP'}^{\rm op}), \T') \right) \subset \Phi(\radu ({\rP}^{\rm op}), \Sr) \cup \{0\},$$
hence $(\lambda^{\vee})^{-1}\C_t({\rP'}^{\rm op}) = \C_t(\rP)$. Thus, the prefan $\mathcal{F}_t$ on $\Lambda(\Sr)$ is the restriction of the prefan $\mathcal{F}_t$ on $\Lambda(\T')$, and consequently the canonical injection $\A(\Sr,k) \hookrightarrow \A(\T,k')$ extends to a homeomorphism between $\overline{\A(\Sr,k)}^{\mathcal{F}_t}$ and the closure of $\A(\Sr,k)$ in $\overline{\A(\T,k')}^{\mathcal{F}_t}$.

It follows from the discussion above that any homeomorphism $\varphi'$ between $\overline{\A(\T,k')}^{\mathcal{F}_t}$ and $\overline{\A}_t(\T,k')$ fitting into the commutative diagram $$\xymatrix{& \A(\T,k') \ar@{->}[ld]_{p}  \ar@{->}[rd]^{\vartheta_t} & \\ \overline{\A(\T,k')}^{\mathcal{F}_t} \ar@{->}[rr]_{\varphi'} & & \overline{\A}_t(\T,k')}$$ induces a homeomorphism $\varphi$ between $\overline{\A(\Sr,k)}^{\mathcal{F}_t}$ and $\overline{\A}_t(\Sr,k)$ fitting into the commutative diagram $$\xymatrix{& \A(\Sr,k') \ar@{->}[ld]_{p} \ar@{->}[rd]^{\vartheta_t} & \\ \overline{\A(\Sr,k')}^{\mathcal{F}_t} \ar@{->}[rr]_{\varphi'} & & \overline{\A}_t(\Sr,k')}$$ and it suffices therefore to prove the proposition when the group $\G$ is split.

\vspace{0.1cm} \emph{Second step: the split case}. We fix a special point $o$ in the apartment $\A(\Sr,k)$ with associated $k^{\circ}$-Chevalley group $\mathcal{G}$.

Let $\mathcal{S}$ denote the split $k^{\circ}$-torus with generic fibre $\Sr$. Any parabolic subgroup $\rP$ of $\G$ containing $\Sr$ extends uniquely to a parabolic subgroup $\mathcal{P}$ of $\mathcal{G}$ containing $\mathcal{S}$ and, if $\mathcal{P}^{\mathrm{op}}$ denotes the opposite parabolic subgroup with respect to $\mathcal{S}$, the morphism $\radu (\mathcal{P}^{\rm op}) \rightarrow \Par(\mathcal{G})$ defined functor-theoretically by $g \mapsto g \mathcal{P} g^{-1}$ is an isomorphism onto an affine open subscheme of $\Par(\mathcal{G})$, which we denote $\Omega_o(\rP,\Sr)$ and whose generic fibre is the big cell $\Omega(\rP,\Sr)$ of $\Par(\G)$. Equivalently, choosing a $k^{\circ}$-Chevalley basis of $\mathrm{Lie}(\G)$ and a total order on $\Psi = \Phi(\radu (\rP^{\rm op}), \Sr)$ allows us to identify $\Omega(\rP,\Sr)$ with the affine space $\mathrm{Spec}\left(k\left[(\X_{\alpha})_{\alpha \in \Psi}\right]\right)$, in which case $\Omega_o(\rP, \Sr)$ corresponds to the $k^{\circ}$-scheme $\mathrm{Spec}\left(k^{\circ}\left[(\X_{\alpha})_{\alpha \in \Psi}\right]\right)$. Finally, from the analytic point of view, the affine open subspace $\Omega_o(\rP,\Sr)$ of $\Par(\mathcal{G})$ determines an affinoid domain $\Omega_o(\rP,\Sr)^{\mathrm{an}}$ in $\Par(\G)^{\mathrm{an}}$ which, in the identification $\Omega(\rP,\Sr) \simeq \mathrm{Spec}\left(k\left[(\X_{\alpha})_{\alpha \in \Psi}\right]\right)$ above, is simply the affinoid domain of $\Omega(\rP,\Sr)^{\mathrm{an}}$ defined by the inequalities $|\X_{\alpha}| \leqslant 1$, $\alpha \in \Psi$.

When $\rP$ runs over the set of parabolic subgroups of $\G$ of type $t$ containing $\Sr$, the affine open subschemes $\Omega_o(\rP,\Sr)$ cover the connected component $\Par_t(\mathcal{G})$ of $\Par(\mathcal{G})$ and the affinoid domains $\Omega_o(\rP,\Sr)^{\mathrm{an}}$ cover therefore the connected component $\Par_t(\G)^{\mathrm{an}}$ of $\Par(\G)^{\mathrm{an}}$ (cf. Remark \ref{rk.big-cell.charts}).

\vspace{0.1cm}
Now we use the special point $o$ to identify the affine space $\A(\Sr,k)$ and the vector space $\V(\Sr) = {\rm Hom}_{\mathbf{Ab}}(\X^*(\Sr),\mathbb{R})$ and we identify the latter with $\Lambda(\Sr) = {\rm Hom}_{\mathbf{Ab}}(\X^*(\Sr), \mathbb{R}_{>0})$ via $$\V(\Sr) \times \X^*(\Sr) \rightarrow \mathbb{R}_{>0}, \ \ \ (u,\chi) \mapsto e^{\langle u, \chi \rangle}.$$ For any parabolic subgroup $\rP$ of $\G$ of type $t$ containing $\Sr$, the image of the map $\vartheta_t: \mathcal{B}(\G,k) \rightarrow \Par(\G)^{\mathrm{an}}$ is contained in the big cell $\Omega(\rP,\Sr)^{\mathrm{an}}$ and its restriction to the apartment $\A(\Sr,k)$ associates with an element $u$ of $\Lambda(\Sr)$ the multiplicative seminorm $$f = \sum_{\nu} a_{\nu} \X^{\nu} \mapsto \max_{\nu} |a_{\nu}| \prod_{\alpha \in \Psi} u(\alpha)^{\nu(\alpha)}$$
on the $k$-algebra $k\left[(\X_{\alpha})_{\alpha \in \Psi}\right]$ (Proposition \ref{prop.theta-t.explicit}). By Lemma \ref{lemma.t-cone}, (i), the polyhedral cone $\C_t(\rP)$ is the preimage of the affinoid domain $\Omega_o(\rP,\Sr)^{\mathrm{an}}$: $$\C_t(\rP) = \vartheta_t^{-1}\left(\Omega_o(\rP,\Sr)^{\mathrm{an}}\right) \cap \Lambda(\Sr).$$
Moreover, if we let $\langle \Psi \rangle^+$ denote the semigroup spanned by $\Psi$ in $\X^*(\Sr)$, the formula above allows us more generally to associate with any homomorphism of unitary monoids $u: \langle \Psi \rangle^+ \rightarrow [0,1]$ a multiplicative seminorm on $k\left[(\X_{\alpha})_{\alpha \in \Psi}\right]$ extending the absolute value of $k$. It follows that we get a continuous and injective map $$\omega_{t,\rP}: \overline{\C_t(\rP)} = \mathrm{Hom}_{\mathbf{Mon}}(\langle \Psi \rangle^+, [0,1]) \rightarrow \Omega_o(\rP,\Sr)^{\mathrm{an}}$$ which fits into the commutative diagram $$\xymatrix{\C_t(\rP) = \{u \in \Lambda(\Sr) \ | \ \alpha(u) \leqslant 1, \ \textrm{for all } \alpha \in \Psi\} \ar@{->}[r]^{\hspace{2.5cm} \vartheta_t} \ar@{->}[d]_{u \mapsto u_{|\Psi}} & \Omega_o(\rP,\Sr)^{\mathrm{an}} \\ \overline{\C_t(\rP)} = \mathrm{Hom}_{\mathbf{Mon}}(\langle \Psi \rangle^+, [0,1]) \ar@{->}[ru]_{\omega_{t,\rP}} & }.$$ If $\rP$ and $\rP'$ are two parabolic subgroups of $\G$ of type $t$ containing $\Sr$, the maps $\omega_{t,\rP}$ and $\omega_{t,\rP'}$ coincide on $\C_t(\rP) \cap \C_t(\rP')$, hence on $\overline{\C_t(\rP)} \cap \overline{\C_t(\rP')} = \overline{\C_t(\rP) \cap \C_t(\rP')}$. We thus get a continuous map $$\overline{\vartheta}_t: \overline{\A(\Sr,k)}^{\mathcal{F}_t} \rightarrow \Par(\G)^{\mathrm{an}}$$ extending $\vartheta_t$.

\vspace{0.1cm} Since the topological spaces $\overline{\A(\Sr,k)}^{\mathcal{F}_t}$ and $\Par(\G)^{\mathrm{an}}$ are compact, the continuous map $\overline{\vartheta}_t$ is proper and its image coincides with the closure $\overline{\A}_t(\Sr,k)$ of $\vartheta_t(\A(\Sr,k))$ in $\Par(\G)^{\mathrm{an}}$.

It only remains to prove that the map $\overline{\vartheta}_t$ is injective.  Since its restriction to any compactified cone $\overline{\C}$ is injective for $\C \in\mathcal{F}_t$, it suffices to check that any two points $x, \ y$ in $\overline{\A(\Sr,k)}^{\mathcal{F}_t}$ such that $\overline{\vartheta}_t(x) = \overline{\vartheta}_t(y)$ belong to the compactification of the same cone in $\mathcal{F}_t$; this is indeed the case by the lemma below. \hfill $\Box$

\vspace{0.3cm} Recall that $\G$ is assumed to be split. Using the notation introduced in the previous proof, let us consider the semisimple $\widetilde{k}$-group $\widetilde{\G}:= \mathcal{G} \otimes_{k^{\circ}} \widetilde{k}$ and the reduction map (1.2.3)
$$r_o : \Par(\G)^{\mathrm{an}} \rightarrow \Par(\widetilde{\G}).$$
Each parabolic subgroup $\Qr$ of $\G$ extends uniquely to a parabolic subgroup $\mathcal{Q}$ of $\mathcal{G}$ and $\widetilde{\Qr}:= \mathcal{Q} \otimes_{k^{\circ}} \widetilde{k}$ is a parabolic subgroup of $\widetilde{\G}$; moreover, if $\Qr$ contains $\Sr$, then $\widetilde{\Qr}$ contains $\widetilde{\Sr}$.

Note that with any parabolic subgroup $\Qr$ of $\G$ containing $\Sr$ we can associate: \begin{itemize} \item the polyhedral cone $\C_t(\Qr)$ in $\Lambda(\Sr)$ (Definition \ref{def.t-cone}), \item the integral closed subscheme $\mathrm{Osc}_t(\widetilde{\Qr})$ of $\Par(\widetilde{\G})$ (Proposition \ref{prop.osc}). \end{itemize}

Finally, for any polyhedral cone $\C$, we define the \emph{interior} $\mathrm{int}(\overline{\C})$ of the compactified cone $\overline{\C}$ as the complement of closures of all proper faces of $\C$: $$\mathrm{int}(\overline{\C}) = \overline{\C} - \bigcup_{\tiny \begin{array}{l} \F \subsetneq \C \\ \hspace{0.1cm} \mathrm{face} \end{array}} \overline{\F}.$$

\begin{Lemma} \label{lemma.reduction} For any parabolic subgroup $\Qr$ of $\G$ containing $\Sr$, the interior of the compactified cone $\overline{\C_t(\Qr)}$ is the preimage under the map $$r_o \circ \overline{\vartheta}_t: \overline{\A(\Sr,k)}^{\mathcal{F}_t} \rightarrow {\rm Par}(\widetilde{\G})$$ of the generic point of the irreducible closed subscheme ${\rm Osc}_t(\widetilde{\Qr})$.

In particular: any two points $x, \ y$ in $\overline{\A(\Sr,k)}^{\mathcal{F}_t}$ with $\overline{\vartheta}_t(x) = \overline{\vartheta}_t(y)$ belong to the compactification of the same cone in $\mathcal{F}_t$.
\end{Lemma}

\vspace{0.1cm}
\noindent
\emph{\textbf{Proof}}.
Taking into account the partition $$\overline{\A(\Sr,k)}^{\mathcal{F}_t} = \bigsqcup_{\C \in\mathcal{F}_t} \mathrm{int}(\overline{\C}),$$ it suffices to check that, for any parabolic subgroup $\Qr$ of $\G$ containing $\Sr$, the map $r_o \circ \overline{\vartheta}_t$ maps the interior of the compactified cone $\overline{\C_t(\Qr)}$ to the generic point of $\mathrm{Osc}_t(\widetilde{\Qr})$.

Let us fix a parabolic subgroup $\rP$ of $\G$ of type $t$ containing $\Sr$, set $\Psi = \Phi(\radu (\rP^{\rm op}), \Sr)$ and identify as above the $k^{\circ}$-scheme $\Omega_{o}(\rP,\Sr)$ with $\mathrm{Spec}\left(k^{\circ}\left[(\X_{\alpha})_{\alpha \in \Psi}\right]\right)$. The restriction of the map $r_o$ to the affinoid domain $$\Omega_o(\rP,\Sr)^{\mathrm{an}} = \left\{x \in \Omega(\rP,\Sr)^{\mathrm{an}} \ \left| |f|(x) \leqslant 1 \right. , \ \textrm{for all } f \in k^{\circ}\left[(\X_{\alpha})_{\alpha \in \Psi}\right] \right\} $$ of $\Omega(\rP,\Sr)^{\mathrm{an}}$ takes values in the affine open subset $\Omega(\widetilde{\rP}, \widetilde{\Sr}) = \Omega_o(\rP,\Sr) \otimes_{k^{\circ}} \widetilde{k}$ of $\Par(\widetilde{\G})$: given a point $x \in \Omega_o(\rP,\Sr)^{\mathrm{an}}$, the set of elements $f \in k^{\circ}\left[(\X_{\alpha})_{\alpha \in \Psi}\right]$ satisfying $|f|(x) < 1$ is a prime ideal containing the maximal ideal of $k^{\circ}$, hence its image in $\widetilde{k}\left[(\X_{\alpha})_{\alpha \in \Psi}\right]$ is a prime ideal and $r_o(x)$ is the point so defined in
$\Omega_o(\widetilde{\rP},\widetilde{\Sr}) \simeq \mathrm{Spec}\left(\widetilde{k}\left[(\X_{\alpha})_{\alpha \in \Psi}\right] \right)$.

Now we consider a parabolic subgroup $\Qr$ of $\G$ containing $\Sr$ and osculatory with $\rP$. By Proposition \ref{prop.roots.fan}, the interior of the compactified cone $\overline{\C_t(\Qr)}$ is the subspace of $\overline{\C_t(\rP)}$ defined by the following conditions: $$\left\{ \begin{array}{ll} \alpha = 1, & \alpha \in \Psi \cap \Phi(\Lr_{\Qr},\Sr) \\ \alpha < 1, & \alpha \in \Psi - \Psi \cap \Phi(\Lr_{\Qr},\Sr). \end{array} \right.$$
It follows that, for any point $x$ in $\mathrm{int}(\overline{\C_t(\Qr)})$, the set of elements $f \in k^{\circ}\left[(\X_{\alpha})_{\alpha}\right]$ satisfying $|f|(x) < 1$ is exactly the ideal generated by the maximal ideal of $k^{\circ}$ and the coordinates $\X_{\alpha}$ with $\alpha \in \Psi - \Psi \cap \Phi(\Lr_{\Qr}, \Sr)$. The point $r_o(x)$ is therefore the generic point of the closed subscheme of $\mathrm{Spec}\left(\widetilde{k}\left[(\X_{\alpha})_{\alpha \in \Psi}\right]\right)$ defined by the vanishing of the coordinates $\X_{\alpha}$ with $\alpha \in \Psi - \Psi \cap \Phi(\Lr_{\Qr},\Sr)$. Finally, since this closed subscheme is the intersection of the open subscheme $\Omega(\widetilde{\rP},\widetilde{\Sr})$ with the closed irreducible subscheme $\mathrm{Osc}_t(\widetilde{\Qr})$ (Proposition \ref{prop.osc.bigcell}, (ii)), $r_o(x)$ is nothing but the generic point of $\mathrm{Osc}_t(\widetilde{\Qr})$ and the proof is complete.
\hfill $\Box$

\section{Group action on the compactifications}
\label{s - group actions}

In this section, for a given reductive group $\G$ over a complete non-Archimedean field $k$ and a  given $k$-rational type $t$ of parabolic subgroups of $\G$,  we describe the Berkovich compactification $\overline{\mathcal{B}}_t(\G,k)$ of type $t$  of the Bruhat-Tits building $\mathcal{B}(\G,k)$.
This means that we describe the boundary components of $\overline{\mathcal{B}}_t(\G,k)$, which are in one-to-one correspondence with $t$-relevant parabolic subgroups as defined in the previous section from a geometric viewpoint in (3.2.1). A root-theoretic interpretation was given in (3.3.2).
It turns out that the boundary component of $\overline{\mathcal{B}}_t(\G,k)$ parameterized by a $t$-relevant parabolic subgroup $\Qr$ of $\G$ can naturally be identified with the Bruhat-Tits buildings of the semisimple quotient of the latter group (Theorem \ref{thm.stratification}). Natural fibrations between flag varieties induce $\G(k)$-equivariant maps between the corresponding compactifications, which we study.
Finally, we also describe the action of the group $\G(k)$ on $\overline{\mathcal{B}}_t(\G,k)$, which enables us to prove a mixed Bruhat decomposition (Proposition \ref{prop - mixed Bruhat dec}).
\vspace{0.3cm}

\subsection{Strata and stabilizers}

Throughout this section, we consider a semisimple $k$-group $\G$ and let $t$ denote a $k$-rational type of parabolic subgroups. We recall that, if $\Qr$ is a parabolic subgroup of $\G$, we still let $t$ denote the $k$-rational type of the parabolic subgroup $(\rP \cap \Qr)/\rad (\Qr)$ of the reductive group $\Qr/{\rm rad}(\Qr)$, where $\rP$ is any parabolic subgroup in $\Par_t(\G)(k)$ osculatory with $\Qr$ (see (3.2.1)).

\vspace{0.3cm} \noindent \textbf{(4.1.1)} For any parabolic subgroup $\rP$ of $\G$, we may use the canonical isomorphism $$\varepsilon_{\rP}: \mathrm{Osc}_t(\rP) \xymatrix{{}\ar@{->}[r]^{\sim} & {}} \Par_t(\rP) = \Par_t(\rP_{\mathrm{ss}})$$ described in Proposition \ref{prop.osc} to define the composite map $$\xymatrix{\mathcal{B}(\rP_{\mathrm{ss}},k) \ar@{->}[r]^{\vartheta_t} & \mathrm{Par}_t(\rP_{\mathrm{ss}})^{\mathrm{an}} \ar@{->}[r]^{\varepsilon_{\rP}^{-1}} & \mathrm{Osc}_t(\rP)^{\mathrm{an}} \ar@{^{(}->}[r] & \mathrm{Par}(\G)^{\mathrm{an}}}$$ and thus we get a continuous injection of the factor $\mathcal{B}_t(\rP_{\mathrm{ss}},k)$ of $\mathcal{B}(\rP_{\mathrm{ss}},k)$ into $\Par(\G)^{\mathrm{an}}$.

\vspace{0.1cm}
\begin{Thm} \label{thm.stratification} Let ${\rm Rel}_t(\G,k)$ denote the set of $t$-relevant parabolic subgroups of $\G$. When $\Qr$ runs over ${\rm Rel}_t(\G,k)$, the buildings $\mathcal{B}_t(\Qr_{\rm ss},k)$ define a stratification of $\overline{\mathcal{B}}_t(\G,k)$ into pairwise disjoint locally closed subspaces: $$\overline{\mathcal{B}}_t(\G,k) = \bigsqcup_{\Qr \in {\rm Rel}_t(\G,k)} \mathcal{B}_t(\Qr_{\rm ss},k).$$

For any $t$-relevant parabolic subgroup $\Qr$ of $\G$, the injection of $\mathcal{B}_t(\Qr_{\rm ss},k)$ into $\overline{\mathcal{B}}_t(\G,k)$ extends to a homeomorphism between the compactified building $\overline{\mathcal{B}}_t(\Qr_{\rm ss},k)$ and the closed subset $$\bigcup_{\tiny \begin{array}{l} \rP \in {\rm Rel}_t(\G,k) \\ \hspace{0.4cm} \rP \subset \Qr \end{array}} \mathcal{B}_t(\rP_{\rm ss},k)$$ of $\mathcal{B}_t(\G,k)$.
\end{Thm}

\vspace{0.2cm} We establish two lemmas before proving this theorem.

\vspace{0.2cm} \begin{Lemma} \label{lemma.stratification1} Let $\rP$ and $\Qr$ be two $t$-relevant parabolic subgroups of $\G$. For any $g \in \G(k)$, $g\mathcal{B}_t(\rP_{ss},k)g^{-1} \cap \mathcal{B}_t(\Qr_{ss},k) \neq \varnothing$ in $\mathrm{Par}(\G)^{{\rm an}}$ if and only if $g\rP g^{-1} = \Qr$. \\
\indent In particular: \begin{itemize} \item[(i)] if $\rP$ and $\Qr$ are distinct, $\mathcal{B}_t(\rP_{ss},k)$ and $\mathcal{B}_t(\Qr_{ss},k)$ are disjoint; \item[(ii)] for any points $x, y \in \mathcal{B}(\rP_{\rm ss}, k)$ and any $g \in \G(k)$, if $g \cdot x=y$ in ${\rm Par}_t(\G)^{\rm an}$, then $g \in \rP(k)$.\end{itemize}
\end{Lemma}

\vspace{0.1cm} \noindent \emph{\textbf{Proof}}. We recall that the $k$-analytic space $\X^{\rm an}$ associated with an algebraic $k$-scheme $\X$ is naturally equiped with a map $\rho : \X^{\rm an} \rightarrow \X$ (see preliminaries on Berkovich theory, 1.2.2). Let us consider the $\G(k)$-equivariant map $\rho: \Par(\G)^{\mathrm{an}} \rightarrow \Par(\G)$ defined at the end of (1.2.2) and pick $x$ in $\mathcal{B}(\rP_{\rm ss},k)$. By Corollary \ref{cor.theta-t.norm}, the map $\rho \circ \vartheta_t : \mathcal{B}(\rP_{ss},k) \rightarrow {\rm Par}_t(\rP_{ss})$ sends $x$ to the generic point of ${\rm Par}_t(\rP_{ss})$. It follows therefore from Proposition \ref{prop.osc} that our canonical embedding of $\mathcal{B}(\rP_{ss},k)$ into ${\rm Par}_t(\G)^{\rm an}$ maps $x$ to a point lying over the generic point of the integral scheme $\mathrm{Osc}_t(\rP)$. Since $g\mathrm{Osc}_t(\rP)g^{-1} = \mathrm{Osc}_t(g \rP g^{-1})$ for any $g \in \G(k)$ (see Remark \ref{rk.osc.conjugation}), the subsets $g\mathcal{B}_t(\rP_{ss},k)g^{-1}$ and $\mathcal{B}_t(\Qr_{ss},k)$ of $\Par(\G)^{\mathrm{an}}$ are non-disjoint if and only if the closed subschemes $\mathrm{Osc}_t(g\rP g^{-1}) = g \mathrm{Osc}_t(\rP)g^{-1}$ and $\mathrm{Osc}_t(\Qr)$ of $\Par(\G)$ coincide. Finally, since the parabolic subgroups $\rP$ and $\Qr$ are $t$-relevant, $g\rP g^{-1}$ and $\Qr$ are also $t$-relevant and the identity $\mathrm{Osc}_t(g\rP g^{-1}) = \mathrm{Osc}_t(\Qr)$ amounts to $g \rP g^{-1} = \Qr$. This completes the proof of our first assertion. Both (i) and (ii) are immediate consequences of what has been said.\hfill $\Box$

\vspace{0.2cm} Let $\Sr$ be a maximal split torus of $\G$ and $\Qr$ a parabolic subgroup of $\G$ containing $\Sr$. We let $\overline{\Sr}$ denote the image of $\Sr$ under the canonical projection $\Qr \rightarrow \Qr_{\mathrm{ss}}$ and define a map $j_{\Qr}: \A(\overline{\Sr},k) \rightarrow \overline{\A(\Sr,k)}^{\mathcal{F}_t}$ as follows:
\begin{itemize}
\item the apartment $\A(\overline{\Sr},k)$ is canonically isomorphic to the quotient of the apartment $\A(\Sr,k)$ by the linear subspace $\X^*(\overline{\Sr})^{\bot} = \langle \mathfrak{C}(\Qr) \rangle$ of $\Lambda(\Sr)$ (1.3.5);
\item by Proposition B.4, (iv), the quotient of $\A(\Sr,k)$ by the linear subspace $\langle \C_t(\Qr) \rangle$ is a stratum of $\overline{\A(\Sr,k)}^{\mathcal{F}_t}$;
\item since $\mathfrak{C}(\Qr) \subset \C_t(\Qr)$ by definition of the latter cone, $\langle \mathfrak{C}(\Qr) \rangle \subset \langle \C_t(\Qr) \rangle$ and thus the canonical projection of $\A(\Sr,k)/\langle \mathfrak{C}(\Qr)\rangle$ onto $\A(\Sr,k)/\langle \C_t(\Qr)\rangle$ leads to a map $$j_{\Qr} : \A(\overline{\Sr},k) = \A(\Sr,k)/\langle \mathfrak{C}(\Qr) \rangle \rightarrow \A(\Sr,k)/\langle \C_t(\Qr) \rangle \subset \overline{\A(\Sr,k)}^{\mathcal{F}_t}.$$
\end{itemize}

Note that this map may not be injective. Observe also that $\overline{\A(\Sr,k)}^{\mathcal{F}_t}$ is covered by images of the maps $j_{\Qr}$ when $\Qr$ runs over the set of parabolic subgroups containing $\Sr$.

\vspace{0.3cm}
\begin{Lemma} \label{lemma.stratification2} With the notation introduced above, the diagram $$\xymatrix{\overline{\A(\Sr,k)}^{\mathcal{F}_t} \ar@{<-}[d]_{j_{\Qr}} \ar@{->}[r]^{\overline{\vartheta}_t} & {\rm Par}(\G)^{{\rm an}} \\ \A(\overline{\Sr},k) \ar@{->}[r]_{\vartheta_t} & {\rm Par}(\Qr_{\rm ss})^{\rm an} \ar@{^{(}->}[u]_{\varepsilon_{\Qr}} }$$ is commutative.
\end{Lemma}

\vspace{0.1cm}
\noindent
\emph{\textbf{Proof}}.
We first reduce to the case of a split group by considering a finite Galois extension splitting $\G$ and a maximal torus $\T$ of $\G$ satisfying the following conditions: \begin{itemize} \item[(a)] $\T$ contains $\Sr$; \item[(b)] $\T \otimes_k k'$ is split; \item[(c)] the injection of $\mathcal{B}(\G,k)$ into $\mathcal{B}(\G,k')$ maps $\A(\Sr,k)$ in $\A(\T,k')$. \end{itemize}
We leave the details to the reader.

\vspace{0.1cm} Now we suppose that the group $\G$ is split. Fix a parabolic subgroup $\rP \in \Par_t(\G)(k)$ osculatory with $\Qr$ and let $\overline{\rP}$ denote the parabolic subgroup $(\rP \cap \Qr)/\rad(\Qr) \in \Par_t(\Qr_{\mathrm{ss}})(k)$. We choose a special point $o$ in $\A(\Sr,k)$ and let $\overline{o}$ denote its image under the canonical projection $\A(\Sr,k) \rightarrow \A(\overline{\Sr},k)$; this is a special point of $\A(\overline{\Sr},k)$, and we use $o$ and $\overline{o}$ as base points to identify $\A(\Sr,k)$ and $\A(\overline{\Sr},k)$ with $\Lambda(\Sr)$ and $\Lambda(\overline{\Sr})$ respectively.

Since the vector space $\Lambda(\overline{\Sr})$ is covered by the cones $\C_t(\overline{\rP})$ when $\rP$ runs over the set of parabolic subgroups $\rP \in \mathrm{Osc}_t(\Qr)(k)$ containing $\Sr$, it suffices to prove that the maps $\overline{\vartheta}_t \circ j$ and $\varepsilon_{\rP} \circ \vartheta_t$ coincide on $\C_t(\overline{\rP})$. Introducing as in the proof of Proposition \ref{prop.apartment.identification} the affinoid domains $\Omega_o(\rP,\Sr)$ and $\Omega_{\overline{o}}(\overline{\rP},\overline{\Sr})$ in $\Par(\G)^{\mathrm{an}}$ and $\Par_t(\Qr_{\mathrm{ss}})^{\mathrm{an}}$ respectively, $\varepsilon_{\Qr}$ identifies $\Omega_{\overline{o}}(\overline{\rP},\overline{\Sr})^{\mathrm{an}}$ with
$\Omega_o(\rP,\Sr)^{\mathrm{an}} \cap \mathrm{Osc}_t(\Qr)^{\mathrm{an}}$ by Proposition \ref{prop.osc.bigcell} and it remains to check that the diagram $$\xymatrix{\overline{\C_t(\rP)} \ar@{->}[r]^{\overline{\vartheta}_t} & \Omega_o(\rP,\Sr)^{\mathrm{an}} \\ \C_t(\overline{\rP}) \ar@{->}[u]^{j} \ar@{->}[r]_{\vartheta_t} & \Omega_{\overline{o}}(\overline{\rP},\overline{\Sr})^{\mathrm{an}} \ar@{->}[u]_{\varepsilon_{\Qr}}}$$ is commutative.

Set $\Psi = \Phi(\radu (\rP^{\rm op}), \Sr)$ and $\overline{\Psi} = \Phi({\rm rad}^{\rm u}(\overline{\rP}^{\rm op}),\overline{\Sr}) = \Psi \cap \Phi(\Lr_{\Qr},\Sr)$, and let $\langle \Psi \rangle^{+}$ and $\langle \overline{\Psi} \rangle^{+}$ denote the semigroups in $\X^*(\Sr)$ and $\X^*(\overline{\Sr})$ spanned by $\Psi$ and $\overline{\Psi}$ respectively. It follows easily from Proposition \ref{prop.roots.fan} that both semigroups $\langle \Psi \rangle^+ \cap \X^*(\overline{\Sr})$ and $\langle \overline{\Psi} \rangle^+$ span the same cone in $\X^*(\Sr) \otimes_{\mathbb{Z}} \mathbb{R}$. The proof of Proposition B.3, (i), shows that $j$ is the map\begin{eqnarray*} \C_t(\overline{\rP}) = \mathrm{Hom}_{\mathbf{Mon}}(\langle \overline{\Psi} \rangle^+ ,]0,1]) & \rightarrow & \mathrm{Hom}_{\mathbf{Mon}}(\langle \Psi \rangle^{+}, [0,1]) = \overline{\C_t(\rP)} \\ \ u & \mapsto & \widetilde{u} = \left\{\begin{array}{ll} u & \textrm{ on } \langle \overline{\Psi} \rangle^+ \\ 0 & \textrm{ on } \langle \Psi \rangle^+ - \langle \Psi \rangle^+ \cap \X^*(\overline{\Sr}). \end{array} \right. .\end{eqnarray*}
Once we have chosen a total order on $\Psi$, we may identify $\Omega_o(\rP,\Sr)$ and $\Omega_{\overline{o}}(\overline{\rP},\overline{\Sr})$ with the spectra of $k^{\circ}\left[(\X_{\alpha})_{\alpha \in \Psi}\right]$ and $k^{\circ}\left[(\X_{\alpha})_{\alpha \in \overline{\Psi}}\right]$ respectively. By Proposition \ref{prop.osc.bigcell}, $\varepsilon_{\Qr}$ is then the morphism deduced from the $k^{\circ}$-homomorphism $$k^{\circ}\left[(\X_{\alpha})_{\alpha \in \Psi}\right] \rightarrow  k^{\circ}\left[(\X_{\alpha})_{\alpha \in \overline{\Psi}}\right], \ \X_{\alpha} \mapsto \left\{\begin{array}{ll} \X_{\alpha} & \textrm{ if } \alpha \in \overline{\Psi} \\ 0 & \textrm{ if } \alpha \in \Psi - \overline{\Psi} \end{array} \right.$$ and, finally, the maps $\vartheta_t \circ j$ and $\varepsilon_{\Qr} \circ \vartheta_t$ both associate with a point $u \in \C_t(\overline{\rP}) = \mathrm{Hom}_{\mathbf{Mon}}(\langle \overline{\Psi} \rangle^+, ]0,1])$ the seminorm $$f = \sum_{\nu} a_{\nu} \X^{\nu}  \mapsto \max_{\nu} |a_{\nu}| \prod_{\alpha \in \Psi} \tilde{u}(\alpha)^{\nu(\alpha)}$$ on $k\left[(\X_{\alpha})_{\alpha \in \Psi}\right]$.

\vspace{0.3cm} \emph{\textbf{Proof of Theorem \ref{thm.stratification}}}. By the very definition of $\overline{\mathcal{B}}_t(\G,k)$ in (3.4.2), any point $x$ of this compactified building belongs to the compactified apartment $\overline{\A}_t(\Sr,k)$ of some maximal split torus $\Sr$ of $\G$. It follows from Lemma \ref{lemma.stratification2} that there exists a parabolic subgroup $\Qr$ such that $x \in \mathcal{B}(\Qr_{\rm ss},k)$. According to Remark \ref{rk.minimal.t-rel} and Lemma \ref{lemma.stratification1}, this $\Qr$ is unique if we assume it to be $t$-relevant. Conversely, if $\Qr$ is a parabolic subgroup of $\G$, any maximal split torus $\Sr'$ of $\Qr_{\mathrm{ss}}$ is the image of some maximal split torus $\Sr$ of $\Qr$ under the canonical projection $\Qr \rightarrow \Qr_{\mathrm{ss}}$ and $\overline{\A}_t(\Sr',k)$ is contained in $\overline{\A}_t(\Sr,k)$ by Lemma \ref{lemma.stratification2}. We have therefore $$\overline{\mathcal{B}}_t(\G,k) = \bigsqcup_{\Qr \in \mathrm{Rel}_t(\G,k)} \mathcal{B}_t(\Qr_{ss},k).$$

\vspace{0.1cm} Let $\Qr$ be a $t$-relevant parabolic subgroup of $\G$. Our injection of $\mathcal{B}_t(\Qr_{\mathrm{ss}},k)$ in $\Par(\G)^{\mathrm{an}}$ obviously extends to a continuous injection of $\overline{\mathcal{B}}_t(\Qr_{\mathrm{ss}},k)$ in $\Par(\G)^{\mathrm{an}}$ and, replacing $\G$ by $\Qr_{\mathrm{ss}}$ in what precedes, we get $$\overline{\mathcal{B}}_t(\Qr_{\mathrm{ss}},k) = \bigsqcup_{\tiny \begin{array}{l} \rP \in \mathrm{Rel}_t(\G,k) \\ \hspace{0.3cm} \Qr \subset \rP \end{array}} \mathcal{B}_t(\rP_{\mathrm{ss}},k).$$

\vspace{0.1cm} Now we check that $\mathcal{B}_t(\Qr_{\mathrm{ss}},k)$ is locally closed in $\overline{\mathcal{B}}_t(\G,k)$. Let us choose a maximal split torus $\Sr$ in $\Qr$ and consider the map $$\pi: \G(k) \times \overline{\A}_t(\Sr,k) \rightarrow \overline{\mathcal{B}}_t(\G,k), \ (x,g) \mapsto g.x:= gxg^{-1}$$ (conjugation takes place in ${\rm Par}(\G)^{\rm an}$). We pick a point $x$ in $\overline{\A}_t(\Sr,k)$ and let $\rP$ denote the $t$-relevant parabolic subgroup containing $\Sr$ such that $x$ is contained in the stratum $\overline{\A}_t(\Sr,k) \cap \mathcal{B}_t(\rP_{\mathrm{ss}},k)$ of $\overline{\A}_t(\Sr,k)$. For any element $g$ of $\G(k)$ such that $g.x$ belongs to $\mathcal{B}_t(\Qr_{\mathrm{ss}},k)$, we have $g\rP g^{-1} = \Qr$ by Lemma \ref{lemma.stratification1}. Since the parabolic subgroups $\Qr$ and $\rP$ both contain the maximal split torus $\Sr$, they are in fact conjugate under the Weyl group $\W$ of $(\G,\Sr)$. Hence there exists $n \in \G(k)$ normalizing $\Sr$ such that $n^{-1}\Qr n = \rP$, thus $(gn^{-1})\Qr (gn^{-1})^{-1} = \Qr$ and therefore $gn^{-1} \in \Qr(k)$. If we pick $n_1, \ldots, n_r$ in $\mathrm{Norm}_{\G}(\Sr)(k)$ lifting the elements of $\W$ and set $\Sigma = \overline{\A}_t(\Sr,k) \cap \mathcal{B}_t(\Qr_{\mathrm{ss}},k)$, then it follows that $$\pi^{-1}\left(\mathcal{B}_t(\Qr_{\mathrm{ss}},k)\right) = \bigcup_{i=1}^{r} \Qr(k) n_i \times (n_i^{-1}\Sigma).$$ Since this subset of $\G(k) \times \overline{\A}_t(\Sr,k)$ is locally closed, $\mathcal{B}_t(\Qr_{\mathrm{ss}},k)$ is a locally closed subspace of $\overline{\mathcal{B}}_t(\G,k)$.

One checks similarly that $\overline{\mathcal{B}}_t(\Qr_{\mathrm{ss}},k)$ is the closure of $\mathcal{B}_t(\Qr_{\mathrm{ss}},k)$ in $\overline{\mathcal{B}}_t(\G,k)$. \hfill $\Box$

\vspace{0.3cm}
\begin{Ex} \label{SL.stratification}
Let $\G$ be the group $\PGL(\V)$, where $\V$ is a vector space of dimension $d+1$ over a locally compact non-Archimedean field $k$.
Following Goldman and Iwahori~\cite{GoldmanIwahori}, the building $\mathcal{B}(\G, k)$ can be identified with the space of norms on $\V$ modulo scaling.
Let $\delta$ be the type of a stabilizer of a flag $(\{0\} \subset \Hr \subset \V)$ with ${\rm codim}(\Hr) = 1$.
It will be shown in a sequel to this article \cite{RTW2} that there exists a $\PGL(\V,k)$-equivariant homeomorphism $\iota$ from $\overline{\mathcal{B}}_{\delta}(\G,k)$ to the space of seminorms on $\V$ modulo scaling, thus extending the Goldman-Iwahori identification.
Let $\Qr$ be a $\delta$-relevant parabolic subgroup.
By Example \ref{ex.relevant}, $\Qr$ is the stabilizer of a flag $(\{0\} \subset \W \subset \V)$. Hence $\Qr_{ss}$ is isogenous to the product $\PGL(\W) \times \PGL(\V/\W)$. Since the type $\delta$ is trivial on  $\PGL(\W)$, the building $\mathcal{B}_\delta (\Qr_{ss}, k)$ coincides with $\mathcal{B}(\PGL(\V/\W), k)$.
In the above identification, $\iota$ identifies the stratum $\mathcal{B}_\delta (\Qr_{ss},k)$ with the set of seminorm classes on $\V$ with kernel $\W$.
\end{Ex}

\vspace{0.3cm}

\begin{Prop} Let $\K/k$ be a non-Archimedean extension.

For any $t$-relevant parabolic subgroup $\Qr$ of $\G$, $\Qr \otimes_k \K$ is a $t$-relevant parabolic subgroup of $\G \otimes_k \K$ and the canonical injection of $\mathcal{B}(\G,k)$ in $\mathcal{B}(\G,\K)$ extends continuously to an injection of $\overline{\mathcal{B}}_t(\G,k)$ in $\overline{\mathcal{B}}_t(\G,\K)$ which induces the canonical injection of $\mathcal{B}_t(\Qr_{\rm ss},k)$ in $\mathcal{B}_t(\Qr_{\rm ss},\K)$.
\end{Prop}

\vspace{0.1cm}
\noindent
\emph{\textbf{Proof}}.
 We have already proved that the canonical injection $\mathcal{B}_t(\G,k) \rightarrow \mathcal{B}_t(\G,\K)$ extends continuously to an injection of $\overline{\mathcal{B}}_t(\G,k)$ in $\overline{\mathcal{B}}_t(\G,\K)$ (Proposition \ref{prop.compact.extension}) and that, for any $t$-relevant parabolic subgroup $\Qr$ of $\G$, the parabolic subgroup $\Qr \otimes_k \K$ of $\G \otimes_k \K$ is still $t$-relevant (Remark \ref{rk.relevant.extension}).

\vspace{0.1cm} It remains to check that our map $\overline{\mathcal{B}}_t(\G,k) \rightarrow \overline{\mathcal{B}}_t(\G,\K)$ induces the canonical injection of $\mathcal{B}_t(\Qr_{\mathrm{ss}},k)$ in $\mathcal{B}_t(\Qr_{\mathrm{ss}},\K)$ for any $t$-relevant parabolic subgroup $\Qr$. The arguments are completely similar to those we used in order to prove Proposition \ref{prop.compact.extension}: we first reduce to the split case, then we rely on the explicit formula of proposition \ref{prop.theta-t.explicit}. \hfill $\Box$

\vspace{0.8cm} \noindent \textbf{(4.1.2)} We describe in this paragraph the subgroups of $\G$ naturally attached to strata of the compactified building $\overline{\mathcal{B}}_t(\G,k)$.

\vspace{0.3cm} \begin{Prop} \label{prop.strata.group} The natural action of $\G(k)$ on $\mathcal{B}_t(\G,k)$ extends uniquely to an action on $\overline{\mathcal{B}}_t(\G,k)$ and, for any $t$-relevant parabolic subgroup $\Qr$ of $\G$ and any element $g$ of $\G(k)$, $$g\mathcal{B}_t(\Qr_{\rm ss},k) = \mathcal{B}_t((g\Qr g^{-1})_{\rm ss},k).$$
\end{Prop}

\vspace{0.1cm}
\noindent
\emph{\textbf{Proof}}.
Given a maximal split torus $\Sr$ of $\G$, the map $\pi: \G(k) \times \overline{\A}_t(\Sr,k) \rightarrow \Par(\G)^{\mathrm{an}}, \ (g,x) \mapsto g.\vartheta_t(x)$ is equivariant with respect to the obvious actions of $\G(k)$. Since its image is precisely the subset $\overline{\mathcal{B}}_t(\G,k)$ of $\Par(\G)^{\mathrm{an}}$, this proves the first assertion. The second follows from Lemma \ref{lemma.stratification1}. \hfill $\Box$

\vspace{0.2cm}
\begin{Prop}
\label{prop.strata.groups}
Let $\Qr$ be a $t$-relevant parabolic subgroup of $\G$.
\begin{itemize}
\item[(i)] For any non-Archimedean extension $\K/k$,  the subgroup $\Qr(\K)$ of $\G(\K)$ is the stabilizer of the stratum $\mathcal{B}_t(\Qr_{\rm ss},\K)$ in $\overline{\mathcal{B}}_t(\G,\K)$.
\item[(ii)] There exists a largest smooth and connected closed subgroup $\R_{t}(\Qr)$ of $\G$ satisfying the following conditions: \begin{itemize} \item $\R_t(\Qr)$ is a normal subgroup of $\Qr$ containing the radical ${\rm rad} (\Qr)$;
\item for any non-Archimedean extension $\K/k$, the subgroup $\R_t (\Qr,\K)$ of $\G(\K)$ acts trivially on the stratum $\mathcal{B}_t(\Qr_{\rm ss},\K)$ of $\overline{\mathcal{B}}_t(\G,\K)$.\end{itemize}

The canonical projection $\Qr_{\rm ss} \rightarrow \Qr/\R_t(\Qr)$ identifies the buildings $\mathcal{B}_t(\Qr_{\rm ss},k)$ and $\mathcal{B}(\Qr/\R_t(\Qr),k)$.

\item[(iii)] For any two points $x, \ y$ in the stratum $\mathcal{B}_t(\Qr_{\rm ss},k)$, there exists a non-Archimedean extension $\K/k$ and an element $g$ of $\Qr(\K)$ such that $y = gx$.
\end{itemize}
\end{Prop}

\vspace{0.1cm}
\noindent
\emph{\textbf{Proof}}.
(i) Since any parabolic subgroup of $\G$ coincides with its normalizer in $\G$, this assertion follows from Lemma \ref{lemma.stratification1}.

\vspace{0.1cm} (ii) Let us consider the central isogeny $\Qr_{\rm ss}' \times \Qr_{\rm ss}'' \rightarrow \Qr_{\rm ss}$ associated with the type $t$ after Definition \ref{def.t-cone}; it identifies $\mathcal{B}(\Qr_{\rm ss}',k)$ and $\mathcal{B}_t(\Qr_{\rm ss},k)$. The preimage of $\Qr_{ss}''$ under the canonical projection of $\Qr$ onto $\Qr_{\mathrm{ss}}$ is a smooth closed subgroup of $\G$ normal in $\Qr$, and we let $\R_t(\Qr)$ denote its identity component. Since the formation of $\R_t(\Qr)$ commutes with arbitrary field extension, the subgroup $\R_t(\Qr,\K)$ of $\G(\K)$ acts trivially on the stratum $\mathcal{B}_t(\Qr_{\mathrm{ss}},\K)$ of $\overline{\mathcal{B}}_t(\G,\K)$ for any non-Archimedean extension $\K/k$.

Suppose now that $\R$ is a smooth and connected closed subgroup of $\G$ which is normal in $\Qr$ and contains the radical $\rad (\Qr)$. The group $\R' = \R/\rad (\Qr)$ is a smooth, connected and normal closed subgroup of $\Qr_{ss}$. By \cite[2.15]{BoTi} (see also \cite[Expos\'e 17]{Bible}), $\R'$ is the image of the product morphism $$\prod_{i \in \I} \Hr_i \rightarrow \G,$$ where $\{\Hr_i\}_{i \in \I}$ is the set of almost simple factors of $\Qr_{ss}$ contained in $\R'$. If the group $\R(k)$ acts trivially on $\mathcal{B}_t(\Qr_{ss},k)$, this is a fortiori the case for each $\Hr_i(k)$ and therefore $\Hr_i$ is contained in $\Qr_{ss}''$ by definition of $\Qr_{ss}'$ and $\Qr_{ss}''$. It follows that $\R$ is contained in $\R_t(\Qr)$.

\vspace{0.1cm} (iii) Consider two points $x$ and $y$ in the stratum $\mathcal{B}_t(\Qr_{\rm ss},k)$. Combining (ii) with Proposition \ref{prop.transitivity}, there exists a non-Archimedean extension $\K/k$ and a $\K$-point $g$ of $\Qr/\R_t(\Qr)$ mapping $x$ to $y$ in $\mathcal{B}_t(\Qr_{\mathrm{ss}},\K) = \mathcal{B}(\Qr/\R_t(\Qr),\K)$. Extending $\K$ if necessary, we may assume that $g$ is the image of a $\K$-point of $\Qr$ and the assertion follows.\hfill $\Box$

\vspace{0.1cm} \begin{Rk} \label{rk.trivial.action} 1. Note that, for any $t$-relevant parabolic subgroup $\Qr$ of $\G$, the group $\R_t(\Qr)(k)$ acts trivially on the whole analytic subspace $\mathrm{Osc}_t(\Qr)^{\mathrm{an}}$ of $\Par(\G)^{\mathrm{an}}$. Indeed, $\R_t(\Qr)$ acts trivially on the subscheme ${\rm Osc}_t(\Qr) \simeq {\rm Par}_t(\Qr_{ss}')$ of ${\rm Par}_t(\G)$ by construction.

2. The formation of $\R_t(\Qr)$ commutes with non-Archimedean field extension: $\R_t(\Qr \otimes_k \K) = \R_t(\Qr) \otimes_k \K$.
\end{Rk}

\vspace{0.3cm} Here is a root-theoretic description of the subgroup $\R_t(\Qr)$ of a $t$-relevant parabolic subgroup $\Qr$ of $\G$. We fix a maximal split torus $\Sr$ of $\Qr$ and let $\overline{\Sr}$ denote its image under the canonical projection $\Qr \rightarrow \Qr_{\mathrm{ss}}$. The canonical injection $\X^*(\overline{\Sr}) \rightarrow \X^*(\Sr)$ identifies the subset $\Phi(\Qr_{\mathrm{ss}},\overline{\Sr})$ of $\X^*(\overline{\Sr})$ with the subset $\Phi(\Qr,\Sr) - \Phi(\radu (\Qr), \Sr) = \Phi(\Lr_{\Qr},\Sr)$ of $\Phi(\Qr,\Sr)$ (where $\Lr_{\Qr}$ denotes the Levi subgroup of $\Qr$ containing $\mathrm{Cent}_{\G}(\Sr)$).

\begin{Prop} Let $\Lambda$ be the set of roots in $\Phi(\Lr_\Qr,\Sr)$ which do not vanish identically on the cone $\C_t(\Qr) \subset \Lambda(\Sr)$.
\begin{itemize}
\item[(i)] The quotient group $\R_t(\Qr)/{\rm rad} (\Qr)$ contains the anisotropic component of $\Qr_{\rm ss}$.
\item[(ii)] The isotropic component of $\R_t(\Qr)/{\rm rad} (\Qr)$ is the subgroup of $\Qr_{ss}$ generated by the images of the root groups $\U_{\alpha}$ for all $\alpha \in \Lambda$.
\item[(iii)] The subgroup $\R_t(\Qr)$ of $\Qr$ is the semi-direct product of ${\rm rad}^u (\Qr)$ by the subgroup of $\Lr_{\Qr}$ generated by the anisotropic component of $L_{\Qr}$, the subtorus of $\Sr$ cut out by the roots in $\Phi(\Lr_{\Qr},\Sr)$ and the root groups $\U_{\alpha}$, $\alpha \in \Lambda$.
\end{itemize}
\end{Prop}

\vspace{0.1cm}
\noindent
\emph{\textbf{Proof}}.
(i) This assertion is clear since $\R_t(\Qr)/\rad({\Qr})$ is the subgroup of $\Qr/\rad(\Qr)$ generated by the quasi-simple components on which the type $t$ restricts trivially (see proof of Proposition \ref{prop.strata.groups}, (ii)).

(ii) Since both $\Lambda$ and its complement are closed subsets of $\Phi(\Lr_{\Qr},\Sr)$ (Corollary \ref{cor.roots.decomposition}), $\Lambda$ is a union of simple components of the root system $\Phi(\Lr_{\Qr},\Sr)$. Let $\rP_0$ be a minimal parabolic subgroup of $\G$ containing $\Sr$ and contained in $\Qr$; we denote by $\Delta$ the corresponding set of simple roots in $\Phi(\G,\Sr)$ and by $\rP$ the parabolic subgroup of $\G$ of type $t$ which contains $\rP_0$. By Corollary \ref{cor.simpleroots.fan}, $\Lambda$ is the union of all connected components of $\Delta \cap \Phi(\Lr_{\Qr},\Sr)$ which do not meet $\Phi(\radu (\rP),\Sr) \cap \Delta$. This amounts to saying that $\Lambda$ is precisely the set of roots of the isotropic quasi-simple components of $\Qr_{\mathrm{ss}}$ on which the restriction of $t$ is trivial. Therefore it follows from the proof of Proposition \ref{prop.strata.groups}, (ii), that the isotropic component of $\R_t(\Qr)/\rad (\Qr)$ is precisely the normal subgroup of $\Qr_{\rm ss}$ corresponding to $\Lambda$. By~\cite[Corollaire 5.11]{BoTi}, the latter is generated by the root groups $\U_{\alpha}$ for all $\alpha \in \Lambda$.

\vspace{0.1cm} (iii) The group $\Qr$ (its radical $\rad (\Qr)$, respectively) is the semi-direct product of its unipotent radical $\radu (\Qr)$ by the Levi subgroup $\Lr_{\Qr}$ (by the radical of $\Lr_{\Qr}$, respectively). Let $\Hr$ denote the maximal anisotropic connected normal subgroup of $\Lr_{\Qr}$. The reductive group $\rad(\Lr_{\Qr})$ is the identity component of the center of $\Lr_{\Qr}$; it is a torus, generated by its anisotropic component $(\rad (\Lr_{\Qr}) \cap \Hr)^0$ and its maximal split subtorus~\cite[Proposition 1.8]{BoTi}. The latter is a subtorus of $\Sr$, namely the connected component of $$\bigcap_{\alpha \in \Phi(\Lr_{\Qr},\Sr)} {\rm ker}(\alpha).$$

The group $\R_t(\Qr)$ is the semi-direct product of $\radu (\Qr)$ by $\Lr_{\Qr} \cap \R_t(\Qr)$. It follows from (i) and (ii) that $\Lr_{\Qr} \cap \R_t(\Qr)$ is the subgroup of $\Lr_{\Qr}$ generated by $\Hr$, the subtorus of $\Sr$ cut out by the roots in $\Phi(\Lr_{\Qr}, \Sr)$ and the root groups $\U_{\alpha}$, $\alpha \in \Lambda$.
\hfill $\Box$

\vspace{0.3cm}

\begin{Ex}
\label{SL(V).trivial}
As in Example \ref{SL.stratification}, let $\G$ be the group $\PGL(\V)$, and let $\delta$ be the type of the stabilizer of a flag $(\{0\} \subset \Hr \subset \V)$ with ${\rm codim}(\Hr) = 1$.
Let $\T$ denote the torus  of diagonal matrices and $\B$ the Borel subgroup of $\G$ consisting of upper triangular matrices
(modulo center, of course), so that $\Hr$ is generated by $e_1, \ldots, e_d$ for a diagonal basis $e_1, \ldots, e_{d+1}$ of $\V$ with respect to $\T$.

Let $\Qr$ be the $\delta$-relevant parabolic subgroup induced by the stabilizer of the subspace $\W$ generated by $e_1, \ldots, e_r$ for some $1 \leq r \leq d+1$. Then, by Example \ref{SL.stratification}, the stabilizer $\R_\delta(\Qr)$ of the stratum $\mathcal{B}_\delta (\Qr_{ss},k)$ is the kernel of the natural map  $\Qr \rightarrow \PGL(\V/\W)$. It obviously contains the unipotent radical $\radu(\Qr)$. The natural morphism $\Lr_{\Qr} \rightarrow \PGL(\W) \times \PGL(\V/\W)$  maps $\R_\delta(\Qr) / \radu(\Qr)$ surjectively on the first factor $\PGL(\W)$. Its kernel is the subgroup of $\T$ given by all diagonal matrices with entries $(a, \ldots, a, b,\ldots, b)$, where $a$ appears $r$ times. This coincides with the subtorus of $\T$ cut out by $\Phi(\Lr_{\Qr},\T)$. Using Example \ref{ex-fan}, we find that a root $\alpha$ of $\Qr$ does not vanish identically on $\C_\delta(\Qr)$ if and only if $\alpha = \chi_i - \chi_j$ for $i\neq j$ and  $i,j \leq r$. The corresponding root groups are exactly the root groups in $\Lr_{\Qr}$ which are mapped to $\PGL(\W)$ under $\Lr_{\Qr} \rightarrow \PGL(\W) \times \PGL(\V/\W)$. Hence we recover the description of $\R_\delta(\Qr)$ in Proposition 4.10.
\end{Ex}

\vspace{0.5cm} \noindent \textbf{(4.1.3)} Now we extend our initial Theorem \ref{thm.affinoid.subgroups} to the compactified building $\overline{\mathcal{B}}_t(\G,k)$ by attaching with each point its stabilizer in $\G^{\mathrm{an}}$.

\begin{Thm} \label{thm.point.stabilizer}
For any point $x$ in $\overline{\mathcal{B}}_t(\G,k)$, there exists a unique geometrically reduced $k$-analytic subgroup ${\rm Stab}_{\G}^{\ t}(x)$ of $\G^{{\rm an}}$ such that, for any non-Archimedean extension $\K/k$, ${\rm Stab}^{\ t}_{\G}(x)(\K)$ is the subgroup of $\G(\K)$ fixing $x$ in $\overline{\mathcal{B}}_t(\G,\K)$.

\vspace{0.1cm}
Let $\Qr$ denote the $t$-relevant parabolic subgroup of $\G$ defining the stratum which contains $x$. The subgroup ${\rm Stab}^{\ t}_{\G}(x)$ is contained in $\Qr^{{\rm an}}$, it contains $\R_t(\Qr)^{{\rm an}}$ as a normal closed analytic subgroup and the canonical isomorphism $\Qr^{{\rm an}}/\R_t(\Qr)^{{\rm an}} \cong (\Qr/\R_t(\Qr))^{{\rm an}}$ identifies the quotient group ${\rm Stab}^{\ t}_{\G}(x)/\R_t(\Qr)^{{\rm an}}$ with the affinoid subgroup $(\Qr/\R_t(\Qr))_x$ of $(\Qr/\R_t(\Qr))^{{\rm an}}$
attached  by Theorem \ref{thm.affinoid.subgroups} to the point $x$ of $\mathcal{B}_t(\Qr_{\rm ss},k) = \mathcal{B}(\Qr/\R_t(\Qr),k)$.
\end{Thm}

\vspace{0.1cm}
\noindent
\emph{\textbf{Proof}}.
To the point $x$ of $\mathcal{B}_t(\Qr_{\mathrm{ss}},k) = \mathcal{B}(\Qr/\R_t(\Qr), k)$ corresponds by Theorem \ref{thm.affinoid.subgroups} a unique $k$-affinoid subgroup $(\Qr/\R_t(\Qr))_x$ of $(\Qr/\R_t(\Qr))^{\mathrm{an}}$ satisfying the following condition: for any non-Archimedean extension $\K/k$, $(\Qr/\R_t(\Qr))_x(\K)$ is the subgroup of $(\Qr/\R_t(\Qr))(\K)$ fixing the point $x$ in $\mathcal{B}(\Qr/\R_t(\Qr),\K)$. Using the canonical isomorphism $\Qr^{\mathrm{an}}/\R_t(\Qr)^{\mathrm{an}} \xymatrix{{} \ar@{->}[r]^{\sim} & {}} (\Qr/\R_t(\Qr))^{\mathrm{an}}$ to identify these analytic groups, we define $\mathrm{Stab}^{\ t}_{\G}(x)$ as the preimage of $(\Qr/\R_t(\Qr))_x$ under the canonical projection $\Qr^{\mathrm{an}} \rightarrow (\Qr/\R_t(\Qr))^{\mathrm{an}}$: $${\rm Stab}_{\G}^{\ t}(x) = \Qr^{\rm an} \times_{(\Qr/\R_t(\Qr))^{\mathrm{an}}} (\Qr/\R_t(\Qr))_x.$$ Since the morphism $\Qr \rightarrow \Qr/\R_t(\Qr)$ is smooth, $\mathrm{Stab}^{\ t}_{\G}(x)$ is a geometrically reduced $k$-analytic subgroup of $\Qr^{\mathrm{an}}$ which contains $\R_t(\Qr)^{\mathrm{an}}$ as a closed invariant analytic subgroup, and the quotient group $\Qr^{\mathrm{an}}/\R_t(\Qr)^{\mathrm{an}}$ is canonically isomorphic to the affinoid subgroup $(\Qr/\R_t(\Qr))_x$ of $(\Qr/\R_t(\Qr))^{\mathrm{an}}$. Moreover, for any non-Archimedean extension $\K/k$ and any element $g$ in ${\rm Stab}^{\ t}_{\G}(x)(\K)$, the action of $g$ on $\overline{\mathcal{B}}_t(\G,\K)$ stabilizes the stratum $\mathcal{B}_t(\Qr_{\mathrm{ss}},\K)$ and fixes the point $x$. The existence part of the proof is thus complete.

\vspace{0.1cm}
Uniqueness follows from the fact that two geometrically reduced analytic subgroups of $\G^{\mathrm{an}}$ having the same $\K$-points for any non-Archimedean extension $\K/k$ coincide.
 \hfill $\Box$

\vspace{0.3cm}
\begin{Prop}
\label{prop.generalized.parahoric}
Let $x$ be a point in $\overline{\mathcal{B}}_t(\G,k)$ and $\Qr$ the $t$-relevant parabolic subgroup of $\G$ such that $x$ belongs to the stratum $\mathcal{B}_t(\Qr_{\rm ss},k)$.
\begin{itemize}
\item[(i)] The group ${\rm Stab}^{\ t}_{\G}(x)(k)$ is Zariski dense in $\Qr$.
\item[(ii)] For any $g \in \G(k)$, ${\rm Stab}^{\ t}_{\G}(gx) = g{\rm Stab}^{\ t}_{\G}(x)g^{-1}$.
\end{itemize}
\end{Prop}

\vspace{0.1cm}
\noindent
\emph{\textbf{Proof}}.
(i) Fix a Levi subgroup $\Lr$ of $\Qr$. Since $\R_t(\Qr)$ contains $\radu (\Qr)$, the group $\mathrm{Stab}^{\ t}_{\G}(x)$ is the semi-direct product of the group $\radu (\Qr)^{\mathrm{an}}$ by the analytic subgroup $\mathrm{Stab}^{\ t}_{\G}(x) \cap \Lr^{\mathrm{an}}$. Therefore, $\mathrm{Stab}^{\ t}_{\G}(x)(k)$ is the semi-direct product of $\radu (\Qr)(k)$ by the subgroup $\mathrm{Stab}^{\ t}_{\G}(x) \cap \Lr(k)$ of $\Lr(k)$.

Let $\Sr_0$ denote the maximal split subtorus of $\rad (\Lr_{\Qr})$, $\Hr_0$ the anisotropic component of $\Lr$ and $(\Hr_i)_{i \in \I}$ the quasi-simple isotropic components of the derived subgroup $\D(\Lr)$ of $\Lr$. The product morphism $$\Sr_0 \times \Hr_0 \times \prod_{i \in \I} \Hr_i \rightarrow \Lr$$ is an isogeny. If we let $\J$ denote the subset of $\I$ consisting of indices $i \in \I$ such that the type $t$ is non-trivial on $\Hr_i$, then $\R_t(\Qr) \cap \Lr$ is the image of the subgroup $\Sr_0 \times \Hr_0 \times \prod_{i \in \J} \Hr_i$ and, for each index $i \in \I - \J$, $\mathrm{Stab}^{\ t}_{\G}(x) \cap \Hr_i^{\mathrm{an}}$ is the affinoid subgroup attached by Theorem \ref{thm.affinoid.subgroups} to the projection of $x$ on the factor $\mathcal{B}(\Hr_i,k)$ of $\mathcal{B}_t(\Qr_{\mathrm{ss}},k)$. It follows that the subgroup $\mathrm{Stab}^{\ t}_{\G}(x) \cap \Lr(k)$ of $\Lr(k)$ contains $\Sr_0(k), \ \Hr_0(k)$ and $\Hr_i(k)$ for each $i \in \J$, as well as a parahoric subgroup of $\Hr_i(k)$ for each $i \in \I - \J$.

The field $k$ is infinite as it carries a non-trivial absolute value. On the one hand, the groups $\Sr_0(k), \ \Hr_0(k)$ and $\Hr_i(k)$ are Zariski dense in the reductive groups $\Sr_0$, $\Hr_0$ and $\Hr_i$ respectively~\cite[Corollary 18.3]{Borel}; on the other hand, each parahoric subgroup of $\Hr_i(k)$ Zariski dense in $\Hr_i$ as well (Lemma \ref{Iwahori.Zariski.dense}) and therefore $\mathrm{Stab}^{\ t}_{\G}(x) \cap \Lr(k)$ is Zariski dense in $\Lr$.
It follows that $\mathrm{Stab}^{\ t}_{G}(x)(k)$ is Zariski dense in $\Qr$ since $\radu (\Qr)(k)$ is Zariski dense in $\radu (\Qr)$~\cite[Expos\'e XXVI, Cor. 2.7]{SGA3}.

(ii) This assertion is obvious.
\hfill $\Box$

\vspace{0.3cm}

\begin{Ex}
\label{SL_stabilizer}
In the setting of Example \ref{SL(V).trivial}, let $x$ be a point in the boundary component $\mathcal{B}_\delta(\Qr_{ss}, k)$.
Recall that $\mathcal{B}_\delta(\Qr_{ss},k)$ can be identified with $\mathcal{B}(\PGL(\V/\W),k)$.
We denote by $x$ also the corresponding point in $\mathcal{B}(\PGL(\V/\W),k)$.
Let $\phi: \Qr \rightarrow \PGL(\V/\W)$ be the natural map.
Then the preimage of the stabilizer of the point $x$ in $\PGL(\V/\W)$ under $\phi$ is equal to the stabilizer $\mbox{\rm Stab}^\delta_{\PGL(\V)}(x)(k)$ of $x$ in $\PGL(\V,k)$.
\end{Ex}

\vspace{0.5cm} \noindent \textbf{(4.1.4)} We will finally give an explicit description of the group $\mathrm{Stab}_{\G}^{\ t}(x)(k)$ for any point $x$ of $\overline{\mathcal{B}}_t(\G,k)$ by combining the theories of Borel-Tits and Bruhat-Tits. We consider a $t$-relevant parabolic subgroup $\Qr$ of $\G$ and pick a point $x$ in the stratum $\mathcal{B}_t(\Qr_{ss},k)$. We fix a maximal split torus $\Sr$ in $\G$ contained in $\Qr$ and such that $x$ belongs to the compactified apartment $\overline{\A}_t(\Sr,k)$. We set $\N = \mathrm{Norm}_{\G}(\Sr)$, $\Z = \mathrm{Cent}_{\G}(\Sr)$ and let $\Lr$ denote the Levi subgroup of $\Qr$ containing $\Z$. The Weyl group $\W$ of the root system $\Phi(\G,\Sr)$ acts on the set of parabolic subgroups containing $\Sr$ and the stabilizer $\W_{\Qr}$ of $\Qr$ in $\W$ is canonically isomorphic to the Weyl group of the root system $\Phi(\Lr,\Sr)$. Moreover, the normalizer of $\Sr$ in $\Qr$ is the subgroup $\N_{\Qr} = \N \cap \Qr$ and we have an exact sequence $$\xymatrix{1 \ar@{->}[r] & \Z(k) \ar@{->}[r] & \N_{\Qr}(k) \ar@{->}[r] & \W_{\Qr} \ar@{->}[r] & 1}.$$

We set $\Lr'' = \R_t(\Qr) \cap \Lr$ and let $\Lr'$ denote the semisimple subgroup of $\Lr$ generated by the isotropic quasi-simple components of $\Lr$ on which $t$ is non-trivial. Both the product morphism $\Lr' \times \Lr'' \rightarrow \Lr$ and the morphism $\pi : \Lr' \rightarrow \Qr/\R_t(\Qr)$ induced by the canonical projection of $\Qr$ onto $\Qr/\R_t(\Qr)$ are central isogenies. We set $\Sr' = (\Sr \cap \Lr')^{\circ}$ and $\Sr'' = (\Sr \cap \Lr'')^{\circ}$. The image $\overline{\Sr}$ of $\Sr'$ in $\Qr/\R_t(\Qr)$ is a maximal split torus of $\Qr/\R_t(\Qr)$ and the homomorphism $\pi^{*}: \X^*(\overline{\Sr}) \rightarrow \X^*(\Sr')$ identifies the root systems $\Phi(\Qr/\R_t(\Qr), \overline{\Sr})$ and $\Phi(\Lr', \Sr')$. Moreover, for any root $\alpha \in \Phi(\Lr', \Sr')$, $\pi$ induces an isomorphism between the root group $\U_{\alpha}$ in $\G$ and the corresponding root group $\overline{\U}_{\alpha}$ in $\Qr/\R_t(\Qr)$. We fix a special point in $\A(\Sr,k)$. Bruhat-Tits theory provides us with a decreasing filtration $\{\U_{\alpha}(k)_r\}_{r \in [-\infty,+\infty]}$ on the group $\U_{\alpha}(k)$ for each root $\alpha \in \Phi(\G,\Sr)$. We have $\U_{\alpha}(k)_{-\infty} = \U_{\alpha}(k)$, $\U_{\alpha}(k)_{+\infty} = \{1\}$ and, for any $r \in ]-\infty,+\infty[$, $\U_{\alpha}(k)_r$ is the subgroup of $\U_{\alpha}(k)$ which acts trivially on the half-space $\{\alpha \geqslant e^{-r} \}$ of $\overline{\A}(\Sr,k)$.

\vspace{0.1cm} Note that the decomposition $$\Phi(\Lr,\Sr) = \Phi(\Lr',\Sr') \cup \Phi(\Lr'', \Sr'')$$ is precisely the decomposition introduced after Definition \ref{def.t-cone}: $\Phi(\Lr'',\Sr'')$ is the union of all irreducible components of $\Phi(\Lr,\Sr)$ on which the type $t$ has trivial restriction whereas $\Phi(\Lr',\Sr')$ is the union of all irreducible components of $\Phi(\Lr, \Sr)$ on which the type $t$ has non-trivial restriction. The subgroups $\W'$ and $\W''$ of $\W_{\Qr}$ stabilizing $\Phi(\Lr',\Sr')$ and  $\Phi(\Lr'',\Sr'')$  respectively are canonically isomorphic to the Weyl groups of the latter root systems and $\W_{\Qr} = \W' \times \W''$.

\vspace{0.1cm}
The action of the group $\N(k)$ on the apartment $\A(\Sr,k)$ extends continuously to an action on the compactified apartment $\overline{\A}_t(\Sr,k)$: indeed, for any $n \in \N(k)$, the automorphism $\mathrm{int}(n)$ of $\mathrm{Par}(\G)^{\mathrm{an}}$ stabilizes the image of the equivariant map $\vartheta_t: \A(\Sr,k) \rightarrow \mathrm{Par}(\G)^{\mathrm{an}}$, hence induces an automorphism of its closure $\overline{\A}_t(\Sr,k)$ in $\mathrm{Par}(\G)^{\mathrm{an}}$. For any point $x$ of $\overline{\A}_t(\Sr,k)$, let $\N(k)_x:= \N(k) \cap \mathrm{Stab}_{\G}^{\ t}(x,k)$ be the subgroup of $\N(k)$ fixing $x$. We set analogously $\Z(k)_x:= \Z(k) \cap \mathrm{Stab}_{\G}^{\ t}(x,k)$ and define the \emph{local Weyl group} $\W_x$ as the image of $\N(k)_x$ in $\W$; we have therefore an exact sequence $$\xymatrix{1 \ar@{->}[r] & \Z(k)_x \ar@{->}[r] & \N(k)_x \ar@{->}[r] & \W_x \ar@{->}[r] & 1}.$$ Observe that, if $x$ belongs to the stratum $\mathcal{B}(\Qr_{\mathrm{ss}},k)$, then each element of $\N(k)_x$ stabilizes $\Qr$, thus $\N(k)_{x}$ is a subgroup of $\N_{\Qr}(k)$ by Lemma \ref{lemma.stratification1}. We also clearly have $\W'' \subset \W_x$, for $\N(k)_x$ contains the group $\mathrm{Norm}_{\Lr''}(\Sr'')(k)$, which is mapped onto $\W''$. It follows that $\W_x = \W'_x \times \W''$, where $\W'_x:= \W_x \cap \W'$.

\vspace{0.1cm}
Finally, $\overline{\A}_t(\Sr,k) \cap \mathcal{B}_t(\Qr_{\mathrm{ss}},k)$ is the apartment of $\Sr'$ in $\mathcal{B}(\Lr',k) = \mathcal{B}_t(\Qr_{\mathrm{ss}},k)$ and is canonically isomorphic to the quotient of $\A(\Sr,k)$ by the linear subspace $\langle \C_t(\Qr) \rangle$ of $\Lambda(\Sr)$ by Lemma \ref{lemma.stratification2}. The choice of an origin in $\A(\Sr,k)$ gives therefore an origin in this affine space, and each root $\alpha$ of $\Phi(\G,\Sr)$ belonging to the subset $\Phi(\Lr',\Sr')$ defines a function on $\overline{\A}_t(\Sr,k) \cap \mathcal{B}_t(\Qr_{\mathrm{ss}},k)$.

\vspace{0.1cm}
\begin{Thm} \label{thm.stab.k-points}
For any point $x$ in $\mathcal{B}_t(\Qr_{\rm ss},k) \cap \overline{\A}_t(\Sr,k)$, the group $\mathrm{Stab}_{\G}^{\ t}(x,k)$ is generated by the following subgroups of $\G(k)$:
\begin{itemize} \item $\N(k)_x$;
\item all $\U_{\alpha}(k)$ with $\alpha \in \Phi({\rm rad}^{\rm u}(\Qr),\Sr)$;
\item all $\U_{\alpha}(k)$ with $\alpha \in \Phi(\Lr'', \Sr'')$;
\item all $\U_{\alpha}(k)_{-\log \alpha(x)}$ with $\alpha \in \Phi(\Lr',\Sr')$.
\end{itemize}
\end{Thm}

\vspace{0.1cm}
\noindent
\emph{\textbf{Proof}}.
Since $\mathrm{Stab}_{\G}^{\ t}(x)$ is the semi-direct product of $\radu(\Qr)^{\mathrm{an}}$ by $\Lr^{\mathrm{an}} \cap \mathrm{Stab}_{\G}^{\ t}(x)$, $\mathrm{Stab}_{\G}^{\ t}(x)(k)$ is the semi-direct product of $\radu(\Qr,k)$ by $\Lr(k)_x:=\Lr(k) \cap \mathrm{Stab}_{\G}^{\ t}(x)(k)$ and it suffices to show that the latter group coincides with the subgroup $\F$ of $\Lr(k)$ generated by $\N(k)_x$, all $\U_{\alpha}(k)$ with $\alpha \in \Phi(\Lr'',\Sr'')$ and all $\U_{\alpha}(k)_{-\log \alpha(x)}$ with $\alpha \in \Phi(\Lr',\Sr')$. The inclusion $\F \subset \Lr(k)_x$ is obvious.

\vspace{0.1cm} Let us choose a minimal parabolic subgroup $\overline{\rP}_0$ in $\Qr/\R_t(\Qr)$ containing $\overline{\Sr}$. Its preimage $\rP_0$ under the isogeny $\pi: \Lr' \rightarrow \Qr/\R_t(\Qr)$ is a minimal parabolic subgroup of $\Lr'$ containing $\Sr'$.

\vspace{0.1cm} \emph{First step}. For any element $g$ in $\Lr(k)_x$, the element $\pi(g)$ of $\big(\Qr/\R_t(\Qr)\big)(k)$ belongs to $\big(\Qr/\R_t(\Qr)\big)_x(k)$, hence can be written as $\pi(g) = \overline{u}_- \overline{u}_+ \overline{n}$, with $$\overline{u}_- \in \radu(\overline{\rP}_0^{\mathrm{op}},k)_x = \big(\Qr/\R_t(\Qr)\big)_x(k) \cap \radu(\overline{\rP}_0^{\mathrm{op}},k),$$ $$\overline{u}_+ \in \radu(\overline{\rP}_0,k)_x = \big(\Qr/\R_t(\Qr)\big)_x(k) \cap \radu(\overline{\rP}_0,k)$$ and $$\overline{n} \in \overline{\N}_x(k) = \overline{\N}(k) \cap \big(\Qr/\R_t(\Qr)\big)_x(k)$$ (\cite{BT1a}, 7.1.4), where $\overline{\N}$ denotes the normalizer of $\overline{\Sr}$ in $\Qr/\R_t(\Qr)$. Since $\radu(\overline{\rP}_0^{\mathrm{op}})(k)_x$ and $\radu(\overline{\rP}_0)(k)$ are generated by the subgroups $\overline{\U}_{\alpha}(k)_{-\log \alpha(x)}$ with $\alpha \in - \Phi(\radu (\overline{\rP}_0),\overline{\Sr})$ and $\alpha \in \Phi(\radu (\overline{\rP}_0),\overline{\Sr})$ respectively, we may write $\overline{u}_- = \pi(u_-)$ and $\overline{u}_+ = \pi(u_+)$ with uniquely defined elements $$u_- \in \radu(\rP_0^{\mathrm{op}})(k)_x = {\big \langle} \U_{\alpha}(k)_{-\log \alpha(x)} \ ; \ \alpha \in -\Phi(\radu(\rP_0),\Sr') {\big \rangle}$$ and $$\ u_+ \in \radu(\rP_0)(k)_x = {\big \langle} \U_{\alpha}(k)_{-\log \alpha(x)}\ ; \ \alpha \in \Phi(\radu(\rP_0),\Sr') {\big \rangle}.$$  Thus, $h = (u_- u_+)^{-1}g$ is an element of $\mathrm{Stab}_{\G}^{\ t}(x)(k) \cap \Lr(k)$ whose image in $\Qr/\R_t(\Qr)$ normalizes $\overline{\Sr}$ and it follows that $h$ normalizes the torus $\Sr'$ in $\Lr$. Therefore we have: $$\mathrm{Stab}_{\G}^{\ t} (x)(k) \subset \F . (\mathrm{Norm}_{\Lr}(\Sr',k) \cap \mathrm{Stab}_{\G}^{\ t}(x)(k)).$$

\vspace{0.1cm} \textit{Second step}. The normalizer (centralizer, respectively) of $\Sr'$ in $\Lr$ is clearly the subgroup of $\Lr$ generated by $\Lr''$ and $\mathrm{Norm}_{\Lr'}(\Sr')$ (by $\Lr''$ and $\mathrm{Cent}_{\Lr'}(\Sr')$ respectively), hence $$\mathrm{Norm}_{\Lr}(\Sr')/\mathrm{Cent}_{\Lr}(\Sr') \simeq \mathrm{Norm}_{\Lr'}(\Sr')/\mathrm{Cent}_{\Lr'}(\Sr')$$ and $$\big(\mathrm{Norm}_{\Lr}(\Sr')/\mathrm{Cent}_{\Lr}(\Sr')
\big)(k) = \mathrm{Norm}_{\Lr}(\Sr')(k)/\mathrm{Cent}_{\Lr}(\Sr')(k)$$ is naturally isomorphic to the Weyl group $\W'$ of the root system $\Phi(\Lr',\Sr')$. Moreover, $\mathrm{Norm}_{\Lr}(\Sr) \subset \mathrm{Norm}_{\Lr}(\Sr')$, $\mathrm{Cent}_{\Lr}(\Sr) \subset \mathrm{Cent}_{\Lr}(\Sr')$ and the natural morphism $$\mathrm{Norm}_{\Lr}(\Sr)(k)/\mathrm{Cent}_{\Lr}(\Sr)(k) \rightarrow \mathrm{Norm}_{\Lr}(\Sr')(k)/\mathrm{Cent}_{\Lr}(\Sr')(k)$$ is the projection of the Weyl group $\W_{\Qr}$ onto its factor $\W'$. It follows that the group $$\mathrm{Norm}_{\Lr}(\Sr')(k)_x:=  \mathrm{Norm}_{\Lr}(\Sr')(k) \cap \mathrm{Stab}_{\G}^{\ t}(x)(k)$$ is an extension of the local Weyl group $\W'_{x}$ by $$\mathrm{Cent}_{\Lr}(\Sr')(k)_x:= \mathrm{Cent}_{\Lr}(\Sr)(k) \cap \mathrm{Stab}_{\G}^{\ t}(x)(k)$$ and, since the subgroup $\N(k)_x$ of $\mathrm{Norm}_{\Lr}(\Sr)(k)$ surjects onto $\W'_x$, the group $\mathrm{Norm}_{\Lr}(\Sr')(k)_x$ is generated by $\N(k)_x$ and $\mathrm{Cent}_{\Lr}(\Sr')(k)_x$. Therefore, $\mathrm{Stab}_{\G}^{\ t}(x)(k)$ is contained in the subgroup of $\G(k)$ generated by $\F$ and $\mathrm{Cent}_{\Lr}(\Sr')(k) \cap \mathrm{Stab}_{\G}^{\ t}(x)(k)$.

\vspace{0.1cm} \emph{Third step}.  The group $\Hr = \mathrm{Cent}_{\Lr}(\Sr') = \Lr''.\mathrm{Cent}_{\Lr'}(\Sr')$ is reductive, $\Sr$ is a maximal split torus and $\Phi(\Hr,\Sr) = \Phi(\Lr'',\Sr'')$. By Borel-Tits theory~\cite[Th\'eor\`eme 5.15]{BoTi}, the group $\Hr(k)$ is generated by the subgroups $\U_{\alpha}(k)$ with $\alpha \in \Phi(\Hr,\Sr)$ and by $\mathrm{Norm}_{\Hr}(\Sr)(k)$. Since the unipotent root group $\U_{\alpha}$ is contained in $\R_t(\Qr)$ for each root $\alpha \in \Phi(\Lr'',\Sr'')$, $\U_{\alpha}(k) \subset \mathrm{Stab}_{\G}^{\ t}(x)(k)$, and it follows that $\Hr(k) \cap \mathrm{Stab}_{\G}^{\ t}(x)(k)$ is generated by these unipotent subgroups and by $\mathrm{Norm}_{\Hr}(\Sr)(k) \cap \mathrm{Stab}_{\G}^{\ t}(x)(k)$. Therefore, $\mathrm{Stab}_{\G}^{\ t}(x)(k)$ is contained in the subgroup of $\G(k)$ generated by $\F$ and $\mathrm{Norm}_{\Hr}(\Sr)(k) \cap \mathrm{Stab}_{\G}^{\ t}(x)(k)$.

\vspace{0.1cm}  \emph{Fourth step}.  Finally, $$\mathrm{Norm}_{\Hr}(\Sr) = \mathrm{Norm}_{\Lr''}(\Sr'').\mathrm{Cent}_{\Lr'}(\Sr')$$ and $$\mathrm{Cent}_{\Hr}(\Sr) = \mathrm{Cent}_{\Lr''}(\Sr'').\mathrm{Cent}_{\Lr'}(\Sr') = \mathrm{Cent}_{\Lr}(\Sr),$$ hence $$\big(\mathrm{Norm}_{\Hr}(\Sr)/\mathrm{Cent}_{\Hr}(\Sr)\big)(k) = \mathrm{Norm}_{\Hr}(\Sr)(k)/\mathrm{Cent}_{\Hr}(\Sr)(k)$$ is naturally isomorphic to the Weyl group $\W''$ of the root system $\Phi(\Lr'',\Sr'')$ and the natural map $$\mathrm{Norm}_{\Hr}(\Sr)(k)/\mathrm{Cent}_{\Hr}(\Sr)(k) \rightarrow \mathrm{Norm}_{\Lr}(\Sr)(k)/\mathrm{Cent}_{\Lr}(\Sr)(k)$$ is the injection of $\W''$ into $\W_{\Qr}$. It follows that the group $$\mathrm{Norm}_{\Hr}(\Sr)(k)_x = \mathrm{Norm}_{\Hr}(\Sr)(k) \cap \mathrm{Stab}_{\G}^{\ t}(x)(k)$$ is an extension of the local Weyl group $\W''_x = \W''$ by \begin{eqnarray*} \mathrm{Cent}_{\Hr}(\Sr)(k)_x & = & \mathrm{Cent}_{\Hr}(\Sr)(k) \cap \mathrm{Stab}_{\G}^{\ t}(x)(k) \\  & = & \mathrm{Cent}_{\Lr}(\Sr)(k)_x \\ & = & \Z(k)_x. \end{eqnarray*} In particular, $\mathrm{Norm}_{\Hr}(\Sr)(k)_x$ is a subgroup of $\N(k)_x$, thus $\mathrm{Stab}_{\G}^{\ t}(x)(k) \subset \F$ and the proof is complete.
\hfill $\Box$

\vspace{0.3cm} The arguments given in the previous proof lead to an extension of Bruhat-Tits' definition of buildings to Berkovich compactifications. Together with the explicit description of the groups $\mathrm{Stab}_{\G}^{\ t }(x)(k)$ above, the next proposition will later allow us to compare Berkovich compactifications with the ones defined by the third author (see {\cite{RTW2}).

\vspace{0.3cm}
\begin{Cor} \label{cor.Berkovich.BT} Let $\Sr$ be a maximal split torus and let $x$ and $y$ be points in $\overline{\A}_t(\Sr,k)$. If there exists an element $g$ of $\G(k)$ such that $gx=y$ in $\overline{\mathcal{B}}_t(\G,k)$, then $y = nx$ for some element $n$ of $\N(k)$.

Consequently, the compactified building $\overline{\mathcal{B}}_t(\G,k)$ is the quotient of $\G(k) \times \overline{\A}_t(\Sr,k)$ by the following equivalence relation: $$(g,x) \sim (h,y) \ \Leftrightarrow \ \left(\exists n \in \N(k), \ y=nx \ \textrm{and} \ g^{-1}hn \in \mathrm{Stab}_{\G}^{\ t}(x)(k)\right).$$
\end{Cor}

\vspace{0.1cm}
\noindent
\emph{\textbf{Proof}}.
Let $\Qr$ and $\Qr'$ denote the $t$-relevant parabolic subgroups of $\G$ containing $\Sr$ such that $x \in \mathcal{B}_t(\Qr_{\mathrm{ss}},k)$ and $y \in \mathcal{B}_t(\Qr'_{\mathrm{ss}},k)$. The identity $gx=y$ implies $\Qr' = g\Qr g^{-1}$ (Lemma \ref{lemma.stratification1}) and thus there exists an element $n_1$ in $\N(k)$ such that $\Qr' = n_1 \Qr n_1^{-1}$. If we set $z = n_1^{-1}y$, then $y=n_1 z$ and $n_1^{-1}gx = z$, and therefore we may assume that the points $x$ and $y$ lie in the same stratum $\mathcal{B}_t(\Qr_{\mathrm{ss}},k)$ of $\overline{\mathcal{B}}_t(\G,k)$. This implies $g \in \Qr(k)$ by Lemma \ref{lemma.stratification1}.

\vspace{0.1cm}
Our final arguments are essentially the same as those given in the previous proof, the notation of which we use again here. Since $\Qr(k) = \radu(\Qr)(k).\Lr(k)$ and $\radu(\Qr)(k)$ acts trivially on $\mathcal{B}_t(\Qr_{\mathrm{ss}},k)$, we may assume that $g$ lies in $\Lr(k)$. Its image $\pi(g)$ in $\big(\Qr/\R_t(\Qr)\big)(k)$ satisfies $\pi(g)x=y$, hence there exists an element $\overline{n}$ of $\overline{\N}(k)$ such that $\overline{n}x=y$ and $\pi(g) \in \overline{n}\big(\Qr/\R_t(\Qr)\big)_x(k)$ (by the very definition of the building $\mathcal{B}(\Qr/\R_t(\Qr),k)$ in~\cite[7.4.1]{BT1a}). Relying on the decomposition $\big(\Qr/\R_t(\Qr)\big)_x(k) = \overline{\N}(k)_x \radu(\overline{\rP}_0^{\mathrm{op}})(k)_x \radu(\overline{\rP}_0)(k)_x$, we may find as in step 1 above unipotent elements $u_-$ and $u_+$ in $\mathrm{Stab}_{\G}^{\ t}(x)(k)$ such that $\pi(g(u_- u_+)^{-1})$ belongs to $\overline{\N}(k)$. If follows that $g (u_{-} u_+)^{-1}$ belongs to $$\mathrm{Norm}_{\Lr}(\Sr')(k) \subset \N(k) \R_t(\Qr)$$ by the last three steps above. We thus can write $g = ng'$ with $n \in \N(k)$ and $g' \in {\rm Stab}_{\G}^{\ t} (x)(k)$, hence $nx = ng'x = y$ and the first assertion of the lemma is established.

\vspace{0.1cm} The second assertion follows immediately from the first. \hfill $\Box$

\subsection{Natural fibrations between compactifications}

Natural morphisms between flag varieties induce fibrations between Berkovich compactifications of a building, which we now describe.

\vspace{0.2cm} \noindent \textbf{(4.2.1)} The set of types of parabolic subgroups of $\G$ is partially ordered as follows: given two types $t$ and $t'$, we set $t \leqslant t'$ if there exist $\rP \in {\rm Par}_t(\G)(k^{\rm a})$ and $\rP' \in {\rm Par}_t(\G)(k^{\rm a})$ with $\rP \subset \rP'$. The maximal type corresponds to the trivial parabolic subgroup $\G$ and the minimal one is given by Borel subgroups.

\vspace{0.1cm} Let $t$ and $t'$ be two types with $t \leqslant t'$. For any $k$-scheme $\Sr$ and any type $t$ parabolic subgroup $\rP$ of $\G \times_k \Sr$, there exists a unique type $t'$ parabolic subgroup $\rP'$ of $\G \times_k \Sr$ with $\rP \subset \rP'$. The map $$\pi_{t, \Sr}^{t'} : {\rm Par}_t(\G)(\Sr) \rightarrow {\rm Par}_{t'}(\G)(\Sr)$$ so defined is functorial with respect to $\Sr$, hence comes from a $k$-morphism $$\pi_t^{t'} : {\rm Par}_t(\G) \rightarrow {\rm Par}_{t'}(\G)$$ which obviously sits in a commutative diagram $$\xymatrix{& {\rm Par}_t(\G)^{\rm an} \ar@{->}[dd]^{\pi_t^{t'}} \\ \G \ar@{->}[ru]^{\lambda_{\rP'}} \ar@{->}[rd]_{\lambda_{\rP}} & \\ & {\rm Par}_{t'}(\G)^{\rm an}}$$ where $\rP \in {\rm Par}_t(\G)(k)$, $\rP ' \in {\rm Par}_{t'}(\G)(k)$, $\lambda_{\rP}(g) = g \rP g^{-1}$ and $\lambda_{\rP'}(g)=g\rP'g^{-1}$. This construction provides us with a continuous and $\G(k)$-equivariant map $$\pi_t^{t'} : \overline{\mathcal{B}}_t(\G,k) \rightarrow \overline{\mathcal{B}}_{t'}(\G,k)$$ such that $\pi_t^{t'} \circ \vartheta_t = \vartheta_{t'}$.

\vspace{0.1cm} \begin{Rk} \label{rk.maximal} Since each $k$-rational type dominates the type $t_{\rm min}$ of minimal parabolic subgroups of $\G$, we have a continuous, surjective and $\G(k)$-equivariant map $$\pi_{t_{\rm min}}^{t} : \overline{\mathcal{B}}_{t_{\rm min}}(\G,k) \rightarrow \overline{\mathcal{B}}_t(\G,k)$$ for each $k$-rational type $t$. Relying on on this observation, $\overline{\mathcal{B}}_{t_{\rm min}}(\G,k)$ is called the \emph{maximal compactification} of $\mathcal{B}(\G,k)$.
\end{Rk}

\vspace{0.2cm} \noindent \textbf{(4.2.2)} We restrict to $k$-\emph{rational} types in this paragraph. We fix two $k$-rational types $t$ and $t'$ with $t \leqslant t'$ and describe the map $\pi_{t}^{t'} : \overline{\mathcal{B}}_{t}(\G,k) \rightarrow \overline{\mathcal{B}}_{t'}(\G,k)$.

\begin{Lemma} For any parabolic subgroup $\Qr$ of $\G$, the stratum $\mathcal{B}_t(\Qr_{ss},k)$ is mapped onto the stratum $\mathcal{B}_{t'}(\Qr_{ss},k)$. Moreover, each $t'$-relevant subgroup is $t$-relevant.
\end{Lemma}

\vspace{0.1cm} \noindent \textbf{\emph{Proof}}. The first assertion follows from the fact that the morphism $\pi_t^{t'} : {\rm Par}_t(\G) \rightarrow {\rm Par}_{t'}(\G)$ maps the subscheme ${\rm Osc}_t(\Qr)$ onto the subscheme ${\rm Osc}_{t'}(\Qr)$. This is immediate in terms of functors: for any $k$-scheme $\Sr$ and any $\rP \in {\rm Osc}_t(\Qr)(\Sr)$, the subgroup $\pi_{t, \Sr}^{t'}(\rP) \cap (\Qr \times_k \Sr)$ of $\G \times_k \Sr$ contains the parabolic subgroup $\rP \cap (\Qr \times_k \Sr)$, hence is parabolic. We therefore have $\pi_{t,\Sr}^{t'}(\rP') \in {\rm Osc}_{t'}(\Qr)(\Sr)$, and the map ${\rm Osc}_t(\Qr) \rightarrow {\rm Osc}_{t'}(\Qr)$ is surjective since $\pi_{t}^{t'}$ is equivariant and both varieties are homogenous under $\Qr$ (Proposition \ref{prop.osc}}).

The stabilizer of ${\rm Osc}_{t'}(\Qr)$ contains the stabilizer of ${\rm Osc}_t(\Qr)$. If $\Qr$ is $t'$-relevant, then $\Qr = {\rm Stab}_{\G}\left({\rm Osc}_{t'}(\Qr)\right)$, hence $\Qr \subset {\rm Stab}_{\G}\left({\rm Osc}_{t}(\Qr)\right) \subset \Qr$ and therefore $\Qr$ is $t$-relevant.
\hfill $\Box$

\vspace{0.1cm}
\begin{Prop} \label{prop.fibrations} Let $\Qr$ be a $t'$-relevant parabolic subgroup of $\G$ and write $\Qr_{ss} = \Hr_1 \times \Hr_2$ (quasi-isogeny), where $\Hr_2$ is the largest semisimple factor of $\Qr_{ss}$ to which the restriction of $t'$ is trivial.
\begin{itemize}
\item[(i)] Let $t_2$ denote the restriction of $t$ to $\Hr_2$. We have $\mathcal{B}_t(\Qr_{ss},k) = \mathcal{B}(\Hr_1,k) \times \mathcal{B}_{t_2}(\Hr_2,k)$, $\mathcal{B}_{t'}(\Qr_{ss},k) = \mathcal{B}(\Hr_1,k)$, and the map $\pi_t^{t'}$ is the projection on the first factor.
\item[(ii)] The preimage of the stratum $\mathcal{B}_{t'}(\Qr_{ss},k)$ is the union of all strata $\mathcal{B}_t(\rP_{ss},k)$, where $\rP$ runs over the set of $t$-relevant parabolic subgroups of $\G$ contained in $\Qr$ and satisfying $\rP/{\rm rad}(\Qr) = \Hr_1 \times \rP_2$ with $\rP_2 \in {\rm Par}(\Hr_2)(k)$, hence it is homeomorphic to $\mathcal{B}(\Hr_1,k) \times \overline{\mathcal{B}}_{t_2}(\Hr_2,k)$.
\end{itemize}
\end{Prop}

\vspace{0.1cm} \noindent \textbf{\emph{Proof}}. (i) Let $t_1$ and $t'_1$ denote the restriction of $t$ and $t'$ respectively to $\Hr_1$. By construction, the restriction of $t'$ to each almost simple factor of $\Hr_1$ is non-trivial; since $t \leqslant t'$, this remark holds \emph{a fortiori} for the type $t$. The schemes ${\rm Osc}_t(\Qr)$, ${\rm Par}_t(\Qr_{ss})$ and ${\rm Par}_{t_1}(\Hr_1) \times {\rm Par}_{t_2}(\Hr_2)$ are canonically isomorphic; similarly, the schemes ${\rm Osc}_{t'}(\Qr)$, ${\rm Par}_{t'}(\Qr_{ss})$ and ${\rm Par}_{t'_1}(\Hr_1)$ are canonically isomorphic (Proposition \ref{prop.osc}). Moreover, the morphism $${\rm Par}_{t_1}(\Hr_1) \times {\rm Par}_{t_2}(\Hr_2) \rightarrow {\rm Par}_{t_1'}(\Hr_1)$$ induced by $\pi_t^{t'}$ is obviously the projection on the first factor composed by $\pi_{t_1}^{t_1'}$.

We have $\mathcal{B}_{t}(\Qr_{ss},k) = \mathcal{B}(\Hr_1,k) \times \mathcal{B}_{t_2}(\Hr_2,k)$, $\mathcal{B}_{t'}(\Qr_{ss},k) = \mathcal{B}(\Hr_1,k)$ and the restrictions of the maps $\vartheta_t$ and $\vartheta_{t'}$ to $\mathcal{B}(\Hr_1,k)$ coincide with the maps $\vartheta_{t_1}$ and $\vartheta_{t_1'}$ respectively by construction (cf. 4.1.1). Then the conclusion follows from commutativity of the diagram $$\xymatrix{& {\rm Par}_{t_1}(\Hr_1)^{\rm an} \ar@{->}[dd]^{\pi_{t_1}^{t_1'}} \\ \mathcal{B}(\Hr_1,k) \ar@{->}[ru]^{\vartheta_{t_1}} \ar@{->}[rd]_{\vartheta_{t_1'}} & \\ & {\rm Par}_{t'_1}(\Hr_1)^{\rm an}.}$$

\vspace{0.1cm} (ii) Given a $t$-relevant parabolic subgroup $\rP$ of $\G$, the stratum $\mathcal{B}_t(\rP_{ss},k)$ is mapped onto the statum $\mathcal{B}_{t'}(\rP_{ss},k)$. The latter coincides with $\mathcal{B}(\Qr_{ss},k)$ if and only if $\Qr$ is the smallest $t'$-relevant parabolic subgroup of $\G$ containing $\rP$, which amounts to saying that ${\rm Osc}_{t'}(\rP) = {\rm Osc}_{t'}(\Qr)$. In the isogeny between $\Hr_1 \times \Hr_2$ and $\Qr_{ss}$, $\rP/{\rm rad}(\Qr)$ corresponds to a parabolic subgroup $\rP_1 \times \rP_2$ of $\Hr_1 \times \Hr_2$, where $\rP_1 \in {\rm Par}(\Hr_1)(k)$ and $\rP_2 \in {\rm Par}(\Hr_2)(k)$. The condition above amounts to ${\rm Osc}_{t'}(\rP_1) = {\rm Par}_{t'}(\Hr_1)$, hence to $\rP_1 = \Hr_1$ by Lemma \ref{Lemma.osc} below.

\vspace{0.1cm} We have therefore \begin{eqnarray*}(\pi_t^{t'})^{-1}(\mathcal{B}_{t'}(\Qr_{ss},k)) & = & \bigcup_{\rP_2 \in {\rm Par}(\Hr_2)(k)} \mathcal{B}(\Hr_1,k) \times \mathcal{B}_{t_2}((\rP_2)_{ss},k) \\ & = & \mathcal{B}(\Hr_1,k) \times \overline{\mathcal{B}}_{t_2}(\Hr_2,k).\end{eqnarray*}
\hfill $\Box$

\vspace{0.2cm} \begin{Lemma} \label{Lemma.osc} Let $t$ denote a $k$-rational type of parabolic subgroups of $\G$ and assume that $t$ is non-degenerate (i.e., is non-trivial on each almost-simple component of $\G$). For any parabolic subgroup $\Qr$ of $\G$, the following conditions are equivalent: \begin{itemize}
\item[(i)] ${\rm Osc}_t(\Qr) = {\rm Par}_t(\G)$ ;
\item[(ii)] $\Qr = \G$.
\end{itemize}
\end{Lemma}

\vspace{0.1cm} \noindent \textbf{\emph{Proof}}. Consider a maximal split torus $\Sr$ of $\G$ contained in $\Qr$ and let $\rP$ denote a parabolic subgroup of $\G$ of type $t$, containing $\Sr$ and osculatory with $\Qr$. It follows from Proposition \ref{prop.osc.bigcell} that ${\rm Osc}_t(\Qr)$ and ${\rm Par}_t(\G)$ coincide if and only if $$\Phi({\rm rad}^{\rm u}(\rP^{\rm op}),\Sr) \subset \Phi(\Qr,\Sr).$$ Since $\rP$ and $\Qr$ are osculatory, ${\rm rad}^{\rm u}(\rP^{\rm op}) \cap {\rm rad}^{\rm u}(\Qr) = \{1\}$ and the latter condition is thus equivalent to $$\Phi({\rm rad}^{\rm u}(\rP^{\rm op}),\Sr) \subset \Phi(\Lr_{\Qr},\Sr),$$ where $\Lr_{\Qr}$ denotes the Levi subgroup of $\Qr$ containing ${\rm Cent}_{\G}(\Sr)$. Now, since $\rP$ induces a non-trivial parabolic subgroup on each almost-simple component of $\G$, $\Phi({\rm rad}^{\rm u}(\rP^{\rm op}),\Sr)$ spans a subgroup of finite index in $\X^*(\Sr)$ by Lemma \ref{lemma.t-cone}, hence $\Phi(\Lr_{\Qr},\Sr)$ spans $\X^*(\Sr) \otimes_{\mathbb{Z}} \mathbb{Q}$ and $\Qr = \G$. \hfill $\Box$

\subsection{The mixed Bruhat decomposition}

Let us choose as above a $k$-rational type of parabolic subgroups, say $t$, and let us consider the corresponding compactification $\mathcal{B}(\G,k) \rightarrow \overline{\mathcal{B}}_t(\G,k)$.

\vspace{0.1cm}
\noindent \emph{\textbf{Notation}}. We adopt the following conventions throughout this paragraph: for any stratum $\Sigma$ of $\overline{\mathcal{B}}_t(\G,k)$, we let $\rP_{\Sigma}$ denote the corresponding $t$-relevant parabolic subgroup of $\G$ and set $\R_{\Sigma} = \R_t(\rP_{\Sigma})$.
For any point $x$ of $\overline{\mathcal{B}}_t(\G,k)$, we let $\Sigma(x)$ denote the stratum --- possibly the building $\mathcal{B}_t(\G,k)$ --- containing $x$ and we set $\G_x = \mathrm{Stab}^{\ t}_{\G}(x)$.
\vspace{0.1cm}

\begin{Prop} \label{prop - mixed Bruhat dec}
Let $x$ and $y$ be any points in $\overline{\mathcal{B}}_t(\G,k)$.

\begin{itemize}
\item[(i)] There exists a maximal split torus $\Sr$ in $\G$ such that $x$ and $y$ lie in $\overline{\A}_t(\Sr,k)$.
\item[(ii)] The group $\G_x(k)$ acts transitively on the compactified apartments containing $x$.
\item[(iii)] Denoting by $\N$ the normalizer of $\Sr$ in $\G$, we have the following decomposition:
$$\G(k) = \G_x(k)\N(k)\G_y(k).$$
\end{itemize}
\end{Prop}

\vspace{0.3cm}
Let us start with the following statement.

\vspace{0.1cm}

\begin{Lemma}
\label{lemma - FoldingFixing}
Let $\overline{\A}$ be the compactified apartment associated with a maximal split torus $\Sr$ and let $\xi \in \overline{\A}$.

\begin{itemize}
\item[(i)] For any $x \in \A$, we have: $\G(k) = \G_\xi(k)\N(k)\G_x(k)$.
\item[(ii)] For any $\eta \in \overline{\mathcal{B}}_t(\G,k)$ such that $\Sigma(\eta) \cap \overline{\A}$ is an apartment in $\Sigma(\eta)$, there exists $g \in \rP_{\Sigma(\eta)} \cap \G_\xi (k)$ such that $g.\eta \in \overline{\A}$.
\end{itemize}
\end{Lemma}

\vspace{0.1cm}

\noindent \emph{\textbf{Proof of lemma}}. For any $\xi \in \overline{\mathcal{B}}_t(\G,k)$, there exists a point $\tilde\xi$ in the maximal compactification of $\mathcal{B}(\G,k)$ such that
${\rm Stab}_{\G(k)}(\tilde\xi) \subset \G_\xi(k)$ (see Remark \ref{rk.maximal}).
Therefore it is enough to work with the \emph{maximal} compactification.

\vspace{0.1cm}
(i) Let us denote by ${\rm H}$ the subset $\G_\xi(k)\N(k)\G_x(k)$.
We have to show that ${\rm H}=\G(k)$.
Let us denote by ${\rm M}$ the reductive Levi factor of $\rP_{\Sigma(\xi)}$ determined by $\Sr$.
For any vector chamber ${\rm D}$ in $\A$, we denote by ${\rm U}_{\rm D}^+(k)$ the unipotent group generated by all the corresponding positive root groups.
We choose ${\rm D}$ so that we have: ${\rm rad}^{\rm u}(\rP_{\Sigma(\xi)})(k) \subset {\rm U}_{\rm D}^+(k) \subset \rP_{\Sigma(\xi)}(k)$.

By definition, ${\rm H}$ contains $\N(k)\G_x(k)$ and by the Iwasawa decomposition~\cite[Prop. 7.3.1 (i)]{BT1a} we have: $\G(k) = {\rm U}_{\rm D}^+(k)\N(k)\G_x(k)$.
Therefore it remains to show that for any $u \in {\rm U}_{\rm D}^+(k)$, we have $u{\rm H} \subset {\rm H}$.
Let $u \in {\rm U}_{\rm D}^+(k)$ and $h \in {\rm H}$.
We write $u = u_\xi'v^+$ with $u_\xi' \in {\rm rad}^{\rm u}(\rP_{\Sigma(\xi)})(k)$ and
$v^+ \in {\rm M}(k) \cap{\rm U}_{\rm D}^+(k)$, and also $h = h_\xi n h_x$ with $h_\xi \in \G_\xi(k)$, $n \in \N(k)$ and $h_x \in \G_x(k)$.
For the factor $h_\xi$, we can write precisely: $h_\xi = u_\xi m_\xi$ with $u_\xi \in {\rm rad}^{\rm u}(\rP_{\Sigma(\xi)})(k)$ and $m_\xi \in {\rm M}_{\xi}(k)$.
Then we have:
$$uh = u'_\xi v^+ u_\xi m_\xi n h_x = (u'_\xi v^+u_\xi(v^+)^{-1})(v^+m_\xi)(nh_x).$$
Since ${\rm rad}^{\rm u}(\rP_{\Sigma(\xi)})(k)$ is normalized by ${\rm M}(k)$, the first factor $r_\xi = u'_\xi v^+u_\xi(v^+)^{-1}$ of the right hand-side belongs to $\G_\xi(k)$.
By Bruhat decomposition in ${\rm M}(k)$~\cite[Th. 7.3.4 (i)]{BT1a}, for any $\zeta \in \overline\A \cap \Sigma(\xi)$ we have
$v^+m_\xi = \ell_\xi n' \ell_\zeta$, with $\ell_\xi \in {\rm L}_{\xi}(k)$, $n' \in \N(k) \cap {\rm M}(k)$ and $\ell_\zeta \in {\rm L}_{\zeta}(k)$, where ${\rm L}$ denotes the semisimple Levi factor ${\rm L}=[{\rm M},{\rm M}]$.
Therefore, for any $\zeta \in \overline\A \cap \Sigma(\xi)$, we can write

\vspace{0.1cm}

\centerline{$uh = r_\xi \ell_\xi n' \ell_\zeta n h_x = (r_\xi\ell_\xi) (n'n) (n^{-1}\ell_\zeta n h_x)$.}

\vspace{0.1cm}

The fixed-point set in $\A$ of the bounded subgroup ${\rm L}_{\zeta}(k)$ is a non-empty intersection of root half-spaces, which is a fundamental domain for the action by translations of ${\rm S}(k) \cap {\rm L}(k)$ on $\A$~\cite[Prop. 7.6.4]{BT1a}.
We thus have the freedom to choose $\zeta \in \overline\A \cap \Sigma(\xi)$ so that $n^{-1}\ell_\zeta n$ fixes $x$.
For such a choice, we have: $r_\xi \ell_\xi \in \G_\xi(k)$, $n'n \in \N(k)$ and $n^{-1}\ell_\zeta n g_x \in \G_x(k)$, as required.

\vspace{0.1cm}

(ii) First, any point $\eta \in \mathcal{B}(\G,k)$ clearly satisfies the hypothesis in claim (ii).
Moreover if both $\xi$ and $\eta$ belong to the building $\mathcal{B}(\G,k)$, then the conclusion of the lemma follows from the facts that there is an apartment containing both of them~\cite[Th. 7.4.18 (i)]{BT1a} and that the stabilizer of $\xi$ in $\G(k)$ acts transitively on the apartments containing $\xi$~\cite[Cor. 7.4.9]{BT1a}.
We henceforth assume that $\xi$ and $\eta$ are not simultaneously contained in $\mathcal{B}(\G,k)$. 

If $\xi \in \A$, the conclusion follows from (i): there exist $g \in \G(k)$ and $\zeta \in \overline{\A} \cap \Sigma(\eta)$ such that $g.\zeta = \eta$; then we can write $g = g_\xi n g_\zeta$ with $g_\xi \in \G_\xi(k) \cap \rP_{\Sigma(\eta)}$, $n \in {\rm Stab}_{\G(k)}(\overline\A)$ and $g_\zeta \in \G_\zeta(k)$, which provides $\eta = g_\xi.(n.\zeta)$. To finish the proof, we assume that neither $\xi$ nor $\eta$ belong to the building $\mathcal{B}(\G,k)$ and argue by induction on the $k$-rank of $\G$.

First, we assume that this rank is equal to 1. Since $\eta \not\in \mathcal{B}(\G,k)$, then $\Sigma(\eta) = \{ \eta \}$ and the hypothesis that $\Sigma(\eta) \cap \overline{\A}$ is an apartment in $\Sigma(\eta)$ simply means that $\eta \in \overline{\A}$, so there is nothing to do.

We assume now that the $k$-rank of $\G$ is $\geqslant 2$ and we denote by $\Lr'$ the semisimple Levi factor of $\rP_{\Sigma(\eta)}$ determined by $\Sr$. Then there is a point $\zeta$ in the closure of $\Sigma(\eta) \cap \overline\A$ such that the stabilizer of $\zeta$ in ${\rm L}'(k)$ fixes $\xi$. To see this, recall that $\A' = \overline{\A} \cap \Sigma(\eta)$ is canonically isomorphic to the quotient of $\A$ by some linear subspace $\F$ and observe that, since we work with the maximal compactification, the projection $p : \A \rightarrow \A'$ extends continuously to a map $\overline{\A} \rightarrow \overline{\A'}$; indeed, the prefan on $\A$ deduced from the fan $\mathcal{F}'_{\varnothing}$ on $\A'$ consists of unions of cones occuring in $\mathcal{F}_{\varnothing}$ and thus $p$ extends to a map between $\overline{\A} = \overline{\A}^{\mathcal{F}_{\varnothing}}$ and $\overline{\A'} = \overline{\A'}^{\mathcal{F}'_{\varnothing}}$.
By the induction hypothesis, we can find $g \in \Lr'_{\zeta}(k)$ such that $g.\eta \in \Sigma(\eta) \cap \overline\A$, and the conclusion follows from this, since $\Lr'_{\zeta}(k) \subset \G_\xi(k)$.
\hfill $\Box$

\vspace{0.2cm}

Here is a proof in the case when the valuation is discrete.
We mention it because it is more geometric (using galleries).

\vspace{0.2cm}

\noindent \emph{\textbf{Second proof of lemma (discrete valuation)}}.
We argue by induction on the minimal length $\ell$ of a gallery in $\Sigma(\eta)$ connecting $\Sigma(\eta) \cap \overline{\A}$ to $\eta$ (by definition, such a gallery is a sequence $\mathfrak{a}_1, \mathfrak{a}_2, ... \, \mathfrak{a}_m$ of consecutively adjacent alcoves in $\Sigma(\eta)$, with $\overline{\mathfrak{a}_1}$ containing a codimension one face in $\overline{\A} \cap \Sigma(\eta)$ and $\eta \in \overline{\mathfrak{a}_m}$).
If $\ell=0$, we can simply take $g=1$.
We now assume that $\ell \geqslant 1$ and choose a corresponding gallery $\mathfrak{a}_1, \mathfrak{a}_2, ... \, \mathfrak{a}_\ell$ as above.
This codimension one face in $\overline{\mathfrak{a}_1} \cap \overline{\A}$ defines a wall in the Bruhat-Tits building $\Sigma(\eta)$, which itself defines a pair of opposite affine roots, say $\{ \pm \alpha \}$, in the root system of $(\rP_{\Sigma(\eta)}/\R_{\Sigma(\eta)})(k)$ with respect to $\Sr(k)$.
This defines a wall in the apartment $\A$ and at least one of the two closed root half-spaces of $\overline{\A}$ bounded by this wall, say the one defined by $\alpha$, contains $\xi$.
We therefore have $\U_\alpha (k) \subset \rP_{\Sigma(\eta)} (k) \cap \G_\xi(k)$.
Moreover, using Bruhat-Tits theory in the boundary stratum $\Sigma(\xi)$, there exists an element $u \in \U_\alpha (k) - \{1\}$ such that $u.\mathfrak{a}_1 \subset \overline{\A} \cap \Sigma(\xi)$
Applying $u$ to the minimal gallery $\mathfrak{a}_1, \mathfrak{a}_2, ... \, \mathfrak{a}_\ell$ and forgetting the first alcove, we see that the point $u.\eta$ can be connected to $\overline{\A} \cap \Sigma(\eta)$ by a gallery of length $\leqslant \ell-1$, so that we can apply our induction hypothesis to find an element $h \in \rP_{\Sigma(\eta)} \cap \G_\xi(k)$ such that $hu.\eta \in \overline{\A} \cap \Sigma(\eta)$.
Then we can finally take $g=hu$.
\hfill $\Box$

\vspace{0.2cm}

We can now proceed to the proof of the proposition.

\vspace{0.1cm}
\noindent
\emph{\textbf{Proof}}.
As a preliminary, we show that if $\Sigma$ and $\Sigma'$ are strata in the compactification $\overline{\mathcal{B}}_t(\G,k)$, then there exists an apartment $\A$ of the building $\mathcal{B}(\G,k)$ such that $\overline{\A}\cap \Sigma$ is an apartment of the building $\Sigma$ and $\overline{\A} \cap \Sigma'$ is an apartment of $\Sigma'$.
Indeed, let $\rP_\Sigma$ and $\rP_{\Sigma'}$ be the parabolic subgroup of $\G$ stabilizing $\Sigma$ and $\Sigma'$, respectively; it is enough to consider a maximal split torus $\Sr$ contained in $\rP_\Sigma \cap \rP_{\Sigma'}$.
The existence of such a maximal split torus (see \cite[20.7]{Borel}) corresponds to the fact that any two facets in the spherical building of $\rP$ are contained in an apartment.
The image, say $\Sr_\Sigma$ and $\Sr_{\Sigma'}$, of $\Sr$ by the canonical projection $\pi_\Sigma: \rP_\Sigma \twoheadrightarrow \rP_\Sigma/\R_\Sigma$ and $\pi_{\Sigma'}: \rP_{\Sigma'} \twoheadrightarrow \rP_{\Sigma'}/\R_{\Sigma'}$, respectively, is then a maximal split torus of the semisimple quotient  $\rP_\Sigma/\R_\Sigma$ (the semsimple quotient $\rP_{\Sigma'}/\R_{\Sigma'}$, respectively) and we have: $\overline{\A}(\Sr,k) \cap \Sigma = \A(\Sr_\Sigma, k)$ ($\overline{\A}(\Sr,k) \cap \Sigma' = \A(\Sr_{\Sigma'}, k)$, respectively).

\vspace{0.1cm}

(i) By the preliminary claim, there exists an apartment $\A$ such that $\overline{\A} \cap \Sigma(x)$ and $\overline{\A} \cap \Sigma(y)$ are apartments in $\Sigma(x)$ and $\Sigma(y)$, respectively.
Let us pick an auxiliary point $z \in \overline{\A} \cap \Sigma(x)$.
By Lemma \ref{lemma - FoldingFixing} with $\xi=z$ and $\eta=y$, we can find $g \in \rP_{\Sigma(y)} \cap \G_z(k)$ such that $g.y \in \overline{\A} \cap \Sigma(y)$ and by the same lemma with $\xi=g.y$ and $\eta=g.x$ we can find $h \in \rP_{\Sigma(\eta)} \cap \G_{g.y}(k)$ such that $hg.x \in \overline{\A} \cap \Sigma(\eta)$.
We finally have: $x, y \in g^{-1}h^{-1}\overline{\A}$.

\vspace{0.1cm}

(ii) The point $x$ lies in the closure $\overline{\A}_0$ of some apartment $\A_0$, which itself corresponds to a maximal split torus $\Sr_0$ in $\G$.
Let $\A$ be an arbitrary apartment such that $x \in \overline{\A}$; it corresponds to a maximal split torus $\Sr$ in $\G$.
Denoting by $\pi_{\Sigma(x)}$ the canonical map $\rP_{\Sigma(x)} \twoheadrightarrow \rP_{\Sigma(x)}/\R_{\Sigma(x)}$, we find that $\pi_{\Sigma(x)}(\Sr)$ is a maximal split torus of $\rP_{\Sigma(x)}/\R_{\Sigma(x)}$.
It follows from classical Bruhat-Tits theory~\cite[Cor. 7.4.9]{BT1a} that $\G_x(k)$ acts transitively on the apartments of the stratum $\Sigma(x)$ containing $x$, so there exists $g \in \G_x(k)$ such that $\pi_{\Sigma(x)}(g\Sr g^{-1})=\pi_{\Sigma(x)}(\Sr_0)$, meaning that both tori $\Sr_0$ and $g\Sr g^{-1}$ lie in the same algebraic $k$-group $\Sr_0 \ltimes \R_{\Sigma(x)}$.
The group $\R_{\Sigma(x)}(k)$ acts transitively (by conjugation) on the maximal $k$-split tori of $\Sr_0 \ltimes \R_{\Sigma(x)}$, so we can find $u \in \R_{\Sigma(x)}(k)$ such that
$(ug)\Sr(ug)^{-1}=\Sr_0$.
It remains to note that $ug \in \G_x(k)$ since $\R_{\Sigma(x)}(k)$ fixes $\Sigma(x)$ pointwise.

\vspace{0.1cm}

(iii) Let $g \in \G(k)$.
By (i), there exists a compactified apartment $\overline{\A}$ containing $ x$ and $y$ and a compactified apartment $\overline{\A'}$ containing $g \cdot y$ and $x$.
By (ii), we can find an element $h \in \G_x(k)$ such that $h\overline{\A} = \overline{\A'}$. Applying Corollary \ref{cor.Berkovich.BT} to the points $y$ and $h^{-1}g \cdot y$ in $\overline{\A}$, we get an element $n$  in $\N(k)$ satisfying $n^{-1}h^{-1}g \cdot y = y$, and therefore $$g = h n (n^{-1}h^{-1}g) \in \G_x(k) \N(k) \G_y(k).$$  
\hfill $\Box$

\vspace{0.2cm}

\begin{Rk}
Geometrically, the proof of (ii) can be described as follows.
Fix a reference apartment $\A_0$ whose closure in $\overline{\mathcal{B}}_t(\G,k)$ contains $x$, and pick an arbitrary apartment $\A$ with the same property.
First, we fold $\overline{\A} \cap \Sigma(x)$ onto $\overline{\A}_0 \cap \Sigma(x)$ by using actions of root group elements from the stabilizer of $x$ in the Levi factor of $\rP_\Sigma$ attached to $\overline{\A}_0 \cap \Sigma(x)$ (this transitivity property for actions of parahoric subgroups is, so to speak, "Bruhat-Tits theory in a stratum at infinity").
Then we use elements of the unipotent radical of $\rP_\Sigma$ to fold $\A$ onto $\A_0$.
\end{Rk}

\section*{Appendix A: on faithfully flat descent in Berkovich geometry}
\label{s - appendix}

In this first appendix, we develop the formalism of faithfully flat descent as introduced by Grothendieck~\cite[VIII]{SGA1}, in the context of Berkovich analytic geometry.
Some technicalities in connection with the Banach module or Banach algebra structures we consider have to be taken into account.
English references for the classical case from algebraic geometry are~\cite{Waterhouse} for affine schemes and~\cite{BLR} in general.

\vspace{0.3cm}

\noindent \textbf{(A.1)} Let $k$ denote a non-Archimedean field and let $\X = \mathcal{M}(\A)$ be a $k$-affinoid space.  For any non-Archimedean extension $\K/k$, the preimage of a $k$-affinoid domain $\D \subset \X$ under the canonical projection ${\rm pr}_{\K/k}: \X_{\K} = \X \widehat{\otimes}_k \K \rightarrow \X$ is a $\K$-affinoid domain in $\X_\K$ since the functor $\F_{{\rm pr}_{\K/k}^{-1}(\D)}$ is easily seen to be represented by the pair $(\A_{\D} \widehat{\otimes}_k \K, \varphi_{\D} \widehat{\otimes} \mathrm{id}_{\K})$. The converse assertion holds if the extension $\K/k$ is \emph{affinoid}, i.e., if $\K$ is a $k$-affinoid algebra, and this appendix is devoted to the proof of this fact.

\vspace{0.3cm}
\noindent \textbf{\emph{Proposition A.1}} --- \emph{Let $\X$ be a $k$-affinoid space and let $\K/k$ be an affinoid extension. A subset $\D$ of $\X$ is a $k$-affinoid domain if and only if the subset ${\rm pr}_{\K/k}^{-1}(\D)$ of $\X_{\K}$ is a $\K$-affinoid domain.}

\vspace{0.2cm}
\noindent \textbf{\emph{Lemma A.2}} --- \emph{Let $\K/k$ be a non-Archimedean extension.
The following conditions are equivalent:
\begin{itemize}
\item[(i)] the extension is affinoid;
\item[(ii)] there exist real positive numbers $r_1, \ldots, r_n$, linearly independent in $(\mathbb{R}_{>0}/|k^{\times}|)
\otimes_{\mathbb{Z}} \mathbb{Q}$ and such that the field $\K$ is a
finite extension of $k_\mathbf{r}$;
\item[(iii)] there exists a tower of non-Archimedean extensions
$$\K = \K(n) \supset \K(n-1) \supset \ldots \supset \K(1) \supset \K\{ 0 \}=k$$
such that $\K(i)/\K(i-1)$ is finite or $\K(i)=\K(i-1)_\mathbf{r}$ for some $\mathbf{r} \in \mathbb{R}_{>0} - |\K(i-1)^{\times}|^{\mathbb{Q}}$.
\end{itemize}}

\vspace{0.1cm}
\noindent
\emph{\textbf{Proof}}.
The implications (ii) $\Rightarrow$ (iii) $\Rightarrow$ (i) are obvious since each extension $\K(i)/\K(i-1)$ is affinoid.
The implication (i) $\Rightarrow$ (ii) is established in~\cite{Duc}.
\hfill $\Box$

\vspace{0.4cm}
\noindent \textbf{(A.2)} It seems adequate to begin by a brief review of faithfully flat descent in algebraic geometry (see also~\cite[VIII]{SGA1} and~\cite[\S 6]{BLR}).

\vspace{0.2cm} \emph{Faithfully flat descent in algebraic geometry.} If $\A$ is a ring, we let $\mathbf{Mod}(\A)$ denote the category of $\A$-modules. Any ring homomorphism $\varepsilon: \A \rightarrow \A'$ defines a functor $$\varepsilon^*: \mathbf{Mod}(\A) \rightarrow \mathbf{Mod}(\A'), \ \M \mapsto \varepsilon^*(\M) = \M \otimes_{\A} \A'.$$ The content of Grothendieck's faithfully flat descent theory is that the category $\mathbf{Mod(\A)}$ can be recovered from the category $\mathbf{Mod(\A')}$ if the homomorphism $\varepsilon$ is \emph{faithfully flat}, which is to say that the functor $\varepsilon^{*}$ is exact --- i.e., it commutes with taking kernels and images --- and faithful --- i.e., $\M\otimes_{\A}\A'=0$ if and only if $\M = 0$.

Consider the natural diagram $$\xymatrix{A \ar@{->}[r]^{\varepsilon} & \A' \ar@<1ex>[r]^{p_1 \hspace{0.4cm}} \ar@<-1ex>[r]_{p_2 \hspace{0.4cm}} & \A' \otimes_{\A} \A' \ar@<-0.3ex>[r]^{p_{13} \hspace{0.4cm}} \ar@<2ex>[r]^{p_{12}\hspace{0.4cm}} \ar@<-2ex>[r]_{p_{23}\hspace{0.4cm}} & \A' \otimes_{\A} \A' \otimes_{\A} \A',}$$ where the $\A$-linear maps are defined by $p_1(a) = a \otimes 1$, $p_2(a)= 1 \otimes a$ and $$p_{12}(a \otimes b) = a \otimes b \otimes 1, \ \ p_{23}(a \otimes b) = 1 \otimes a \otimes b, \ \ p_{13}(a\otimes b) = a \otimes 1 \otimes b,$$ so that $$p_{12} p_1 = p_{13} p_ 1 = q_1, \ \ p_{12} p_2 = p_{23}p_1 = q_2 \ \ \textrm{and} \ \ p_{23}p_2 = p_{13}p_2 = q_3,$$ where
$$q_1(a) = a \otimes 1 \otimes 1, \ \ q_2(a) = 1 \otimes a \otimes 1, \ \ \textrm{and} \ \ q_3(a)= 1 \otimes 1 \otimes a.$$
A \emph{descent datum} on an $\A'$-module $\M$ is an isomorphism of $\A' \otimes_{\A} \A'$-modules $\delta: p_1^{*}\M \quad \tilde{\rightarrow} \quad p_2^*M$ satisfying the following cocyle condition: $$p_{23}^{*}(\delta) \circ p_{12}^{*}(\delta) = p_{13}^{*}(\delta).$$ We denote by $\mathbf{Mod(\A')_{\mathbf{desc}}}$ the category whose objects are pairs $(\M,\delta)$ consisting of an $\A'$-module equipped with a descent datum and in which the morphisms between two objects $(\M, \delta)$ and $(\N, \delta')$ are the $\A'$-linear maps $\M \rightarrow \N$ compatible with descent data (in an obvious way). For any $\A$-module $\M$, the canonical isomorphism $p_1^*(\varepsilon^*\M) \xymatrix{{} \ar@{->}[r]^{\sim} & {}} p_2^*(\varepsilon^*\M)$ provides a descent datum $\delta_{\M}$ on the $\A'$-module $\varepsilon^{*}(\M)$ and the $\A'$-linear map $\varepsilon^{*}(\varphi): \varepsilon^{*}\M \rightarrow \varepsilon^{*}\N$ induced by an $\A$-linear map $\varphi: \M \rightarrow \N$ is automatically compatible with the descent data $\delta_{\M}, \ \delta_{N}$. Hence we get a functor $\underline{\varepsilon}^{*}: \mathbf{Mod(\A)} \rightarrow \mathbf{Mod(\A')_{\mathbf{desc}}}, \ \M \mapsto (\varepsilon^{*}\M, \delta_{\M})$.

\vspace{0.2cm}
\noindent \textbf{\emph{Theorem A.3}} --- \emph{The functor $\underline{\varepsilon}^{*}$ is an equivalence of categories.
Moreover there exists (up to a unique isomorphism) at most one descent datum on a given $A'$-module.}

This theorem follows readily from the next two statements:

\vspace{0.1cm}

(i) For any $\A$-module $\M$, the sequence $$\xymatrix{0 \ar@{->}[r] & \M \ar@{->}[r]^{\varepsilon_{\M}} & \varepsilon^{*}\M \ar@{->}[r]^{\varphi_{\M} \hspace{0.4cm}} & p_2^{*} (\varepsilon^{*}\M) ,}$$ where $\varepsilon_{\M} = \mathrm{id}_{\M} \otimes \varepsilon$ and $\varphi_{\M} = \mathrm{id}_{\M} \otimes p_2 - \delta_{\M} \circ (\mathrm{id}_{\M} \otimes p_1)$, is exact.

(ii) For any $\A'$-module $\M$ equipped with a descent datum $\delta$, let $\M_0$ be the kernel of the map $$\varphi_{\delta} = \mathrm{id}_{\M} \otimes p_2 - \delta \circ (\mathrm{id}_{\M} \otimes p_1): \M \longrightarrow p_2^{*}\M;$$ then the canonical map $\varepsilon^{*}\M_0 = \M_0 \otimes_{\A} \A' \rightarrow \M$, which is automatically compatible with the descent data $\delta$ and $\delta_{\M_0}$, is an isomorphism.

\vspace{0.1cm}

\textit{First step} ---
We begin by assuming that the homomorphism $\varepsilon$ admits a section $\sigma$. Defining the map $\tau: \A' \otimes_\A \A' \rightarrow \A'$ by $\tau(a \otimes b) = \sigma(a)b$ and setting $\sigma_\M = \mathrm{id}_\M \otimes \sigma$ and $\tau_\M = \mathrm{id}_{\M} \otimes \tau$, we have $\sigma_\M \circ \varepsilon_\M = \mathrm{id}_\M$ and $\tau_\M \circ \varphi_\M = \mathrm{id}_{\varepsilon^{*}\M} - \varepsilon_\M \circ \sigma_\M$, hence the sequence (i) is exact.

The descent datum $\delta$ induces an isomorphism $\tau^{*}(\delta)$ between the $\A'$-modules $\tau^{*}(p_1^{*}\M) = \varepsilon^{*}(\sigma^*\M)$ and $\tau^{*}(p_2^{*}\M) = \M$. Thanks to the cocycle condition satisfied by $\delta$, this isomorphism is compatible with the descent data $\delta$ and $\delta_{\sigma^*\M}$; in view of (i), it induces therefore an $\A$-module isomorphism between $\sigma^* \M$ and $\M_0$. Hence $(\M,\delta)$ is canonically isomorphic to $(\varepsilon^{*}\M_0, \delta_{M_0})$ and descent data on $\A'$-modules are therefore unique up to a unique isomorphism.

\textit{Second step} --- We now rely on faithful flatness of $\A'$ over $\A$ to deduce the general case from the first step. Indeed, the first assertion is true if and only if the sequence is exact after applying $\varepsilon^{*}$ (``assertion $\varepsilon^{*}$(i)'') whereas the second assertion is true if and only if the canonical map $\varepsilon^{*}(\varepsilon^{*}\M_0) \rightarrow \varepsilon^{*}\M$ is an isomorphism (``assertion $\varepsilon^{*}$(ii)''). Thanks to the associativity of tensor product and to the canonical identification $\M \otimes_{\A}' \A' = \M$ for any $\A'$-module $\M'$, assertion $\varepsilon^{*}$(i) is exactly assertion (i) if we consider the morphism $p_1: \A' \rightarrow \A' \otimes_{\A} \A'$ and the $\A'$-module $\varepsilon^{*}\M$ instead of the morphism $\varepsilon: \A \rightarrow \A'$ and the $\A$-module $\M$. By the same argument, $\varepsilon^{*}(\delta)$ is a descent datum on the $\A' \otimes_{A} \A'$-module $p_1^{*}\M' = \varepsilon^{*}(\varepsilon^{*}\M)$ with respect to the morphism $p_1$ and, since $\varepsilon^{*}\M_0$ is the kernel of $\varepsilon^{*}(\varphi_{\delta}) = \varphi_{\varepsilon^{*}(\delta)}$, assertion $\varepsilon^{*}$(ii) is precisely assertion (ii) if we consider the morphism $p_1$ and the $\A' \otimes_{\A} \A'$-module $p_1^{*}\M$ instead of the the morphism $\varepsilon$ and the $\A'$-module $\M$. But assertions $\varepsilon^{*}$(i) and $\varepsilon^{*}$(ii) are true since the morphism $p_1$ has an obvious section; assertions (i) and (ii) are therefore true and the theorem is proved.
\hfill $\Box$

\vspace{0.2cm}
\noindent \textbf{\emph{Remark A.4}} --- 1. It is worth recalling that faithfully flat descent includes Galois descent as a special case. Indeed, if $\Lr/\K$ is a finite Galois extension with group $\G$, the map
$$\Lr \otimes_\K \Lr \xymatrix{{} \ar@{->}[r]^{\sim} & {}} \prod_{g \in \G} \Lr, \  \ a \otimes b \mapsto (g(a)b)_g$$
is by definition an isomorphism of $\K$-algebras and, if $\M$ is an $\Lr$-module,
\begin{itemize}
\item an $\Lr \otimes_{\K} \Lr$-isomorphism $\delta: p_1^{*} \M \tilde{\rightarrow} p_2^{*}\M$ is nothing but a collection $(\delta_{g})_{g \in \G}$ of $\K$-automorphisms of $\M$ such that $\delta_g(ax) = g(a)\delta_g(x)$ for any $a \in \Lr$, $g \in \G$ and $x \in \M$;
\item $\delta$ is a descent datum, i.e., it satisfies the cocycle condition, if and only if $\delta_{gh} = \delta_g \circ \delta_h$ for any $g, h \in \G$.
\end{itemize}

In other words, a descent datum on an $\Lr$-module $\M$ is nothing but an action of $\G$ on $\M$ via semilinear automorphisms. Moreover, if $\delta = (\delta_g)_{g \in \G}$ is a descent datum on $\M$, then $\mathrm{Ker}(\varphi_{\delta})$ is the $\K$-module consisting of all elements $x$ in $\M$ such that $\delta_g(x)=x$ for any $g \in \G$.

2. Faithfully flat descent applies equally well to algebras: indeed, the functor $\underline{\varepsilon}: \mathbf{Mod(\A)} \rightarrow \mathbf{Mod(\A')_{desc}}$ obviously induces an equivalence between the subcategories $\mathbf{Alg(\A)}$ and $\mathbf{Alg(\A')_{desc}}$ if we restrict ourselves to descent data which are isomorphisms of $\A' \otimes_{\A} \A'$-algebras.

\vspace{0.4cm} \noindent \textbf{(A.3)} \emph{non-Archimedean field extensions.} We consider now a non-Archimedean extension $\K/k$ and we adapt the algebraic arguments above to the functor $$\mathbf{BMod}(k) \rightarrow \mathbf{BMod(\K)}, \ \ \M \mapsto \M \widehat{\otimes}_k \K.$$ Working with completed tensor products instead of standard tensor products requires only minor modifications as soon as one knows that this functor is exact on the subcategory $\mathbf{BMod^{\mathrm{st}}(k)}$; this nontrivial fact is due to L. Gruson~\cite{Gru}.

\noindent \textbf{\emph{Lemma A.5}} --- \emph{Let $\K/k$ be a non-Archimedean extension.
\begin{itemize}
\item[(i)] The functor $$\varepsilon^{*}: \mathbf{BMod}(k) \rightarrow \mathbf{BMod}(\K), \ M \mapsto M \widehat{\otimes}_k \K$$ transforms strict exact sequences of $k$-modules into strict exact sequences of $\K$-modules.
\item[(ii)] For any Banach $k$-module $\M$, the canonical homomorphism $\M \rightarrow \varepsilon^{*}\M$ is an isometric injection. In particular, the functor $\varepsilon^{*}$ is faithful.
\item[(iii)] A sequence of Banach $k$-modules is strict and exact if and only if it is strict and exact after applying $\varepsilon^{*}$.
\end{itemize}}

\vspace{0.1cm}
\noindent
\emph{\textbf{Proof}}.
(i) This is proved by Gruson in~\cite[Sect. 3]{Gru} and the argument goes as follows.

Let $\xymatrix{0 \ar@{->}[r] & \M' \ar@{->}[r]^{u} & \M \ar@{->}[r]^{v} & \M'' \ar@{->}[r] & 0}$ be a short exact and strict sequence of Banach $k$-modules; modifying norms in their equivalence classes if necessary, we can assume that both $u$ and $v$ are isometric. The sequence $$\xymatrix{0 \ar@{->}[r] & \M'\widehat{\otimes}_k \N \ar@{->}[r]^{u \widehat{\otimes} \mathrm{id}_{\N}} & \M \widehat{\otimes}_k \N \ar@{->}[r]^{v \widehat{\otimes} \mathrm{id}_{\N}} & \M'' \widehat{\otimes}_k \N \ar@{->}[r] & 0}$$ is obviously exact and isometric if $\N$ is a finite dimensional Banach $k$-module, since $\N$ is then the direct sum of a finite number of copies of $k$. Having proved that any Banach $k$-module $\N$ is the limit of a direct system $(\N_{\bullet})$ of finite dimensional Banach $k$-modules, one gets a short exact and isometric sequence $$\xymatrix{0 \ar@{->}[r] & \M'\widehat{\otimes}_k \N_{\bullet} \ar@{->}[r]^{u \widehat{\otimes} \mathrm{id}_{\N_{\bullet}}} & \M \widehat{\otimes}_k \N_{\bullet} \ar@{->}[r]^{v \widehat{\otimes} \mathrm{id}_{\N_{\bullet}}} & \M'' \widehat{\otimes}_k \N_{\bullet} \ar@{->}[r] & 0}$$ of direct systems of Banach $k$-modules. In this situation, taking limits preserves exactness as well as norms and we conclude from the commutativity of completed tensor products with limits that the sequence $$\xymatrix{0 \ar@{->}[r] & \M'\widehat{\otimes}_k \N \ar@{->}[r]^{u \widehat{\otimes} \mathrm{id}_{\N}} & \M \widehat{\otimes}_k \N \ar@{->}[r]^{v \widehat{\otimes} \mathrm{id}_{\N}} & \M'' \widehat{\otimes}_k \N \ar@{->}[r] & 0}$$ is exact and isometric.

(ii) Pick a direct system $\M_{\bullet}$ of finite dimensional Banach $k$-modules with limit $\M$. Since the assertion is obvious as long as $\M$ is decomposable, we get an isometric exact sequence of direct systems $\xymatrix{0 \ar@{->}[r] & \M_{\bullet} \ar@{->}[r] & \M_{\bullet} \widehat{\otimes}_k \K}$ and, taking limits, we conclude that the canonical homomorphism $\M \rightarrow \M \widehat{\otimes}_k \K$ is an isometric injection.

(iii) If a bounded $k$-linear map $u: \M \rightarrow \N$ between Banach $k$-modules is strict, then the bounded $\K$-linear map $u_{\K} = u \widehat{\otimes}_k \K$ is strict thanks to the exactness property of the functor $\cdot \ \widehat{\otimes}_k \K$ on $\mathbf{BMod}^{\mathbf{st}}(k)$. Conversely, consider the commutative diagram $$\xymatrix{0 \ar@{->}[r] & \M/\mathrm{ker}(u) \ar@{->}[r] \ar@{->}[d]_{u} & \left(\M/\mathrm{ker}(u)\right) \widehat{\otimes}_k \K  = \left(\M \widehat{\otimes}_k \K \right)/\mathrm{ker}(u_{\K}) \ar@{->}[d]^{u_{\K}} \\ 0 \ar@{->}[r] & \N \ar@{->}[r] & \N \widehat{\otimes}_k \K},$$ in which rows are exact and isometric; if the map $u_{\K}$ is strict, then so is $u$ and the conclusion follow from (i) and (ii).
\hfill $\Box$

\vspace{0.3cm} The definition of a descent datum is formally the same as in the algebraic situation.

\vspace{0.2cm}
\noindent \textbf{\emph{Proposition A.6}} --- \emph{Let $\K/k$ be an extension of non-Archimedean fields. The functor $$\mathbf{BMod}(k) \rightarrow \mathbf{BMod}(\K)_{\mathbf{desc}}, \ \ \M \mapsto (\M, \delta_{\M})$$ is an equivalence of categories.}

\vspace{0.1cm}
\noindent
\emph{\textbf{Proof}}.
By the same general arguments as in the proof of Theorem A.3, the proposition follows from the next two assertions.

(i) For any Banach $k$-module $\M$, the sequence $$(\Sr) \ \ \xymatrix{0 \ar@{->}[r] & \M \ar@{->}[r]^{\varepsilon_{\M}} & \varepsilon^{*}\M \ar@{->}[r]^{\varphi_{\M}} & p_2^{*}\M }$$ is strict and exact.

(ii) For any $\K'$-module $\M$ equipped with a descent datum $\delta$, let $\M_0$ be the kernel of the map $$\varphi_{\delta} = \mathrm{id}_{\M} \otimes p_2 - \delta \circ (\mathrm{id}_{\M} \otimes p_1): \M \longrightarrow p_2^{*}\M.$$ Then the canonical map $\varepsilon^{*}\M_0 = \M_0 \widehat{\otimes}_{k} \K \rightarrow \M$, which is automatically compatible with the descent data $\delta$ and $\delta_{\M_0}$, is a (strict) isomorphism.

As in the algebraic situation above, these assertions are true as soon as $\varepsilon$ is any morphism of Banach $k$-algebras admitting a section; they are therefore true if one substitutes the field extension $\varepsilon: k \rightarrow \K$ and the Banach $k$-module $\M$ (the Banach $\K$-module with descent datum $(\M, \delta)$, respectively) by the morphism $p_1: \K \rightarrow \K \widehat{\otimes}_k \K$ and the Banach $\K$-module $\varepsilon^{*}\M$ (the Banach $\K \widehat{\otimes}_k \K$-module with descent datum $(p_1^{*}\M, p_1^{*}(\delta))$, respectively). Thanks to the associativity of completed tensor product and to the canonical identification $\M \widehat{\otimes}_\K \K = \M$ for any Banach $\K$-module $\M$, the new sequences  relative to $p_1: \K \rightarrow \K \widehat{\otimes}_k \K$ are exactly the ones obtained by applying the functor $\varepsilon^{*}$ to the former sequences, relative to $\varepsilon: k \rightarrow \K$. Therefore (i) and (ii) follow from Lemma A.5. \hfill $\Box$

\vspace{0.2cm}
The following slightly more precise result will be useful in the study of maps between compactifications.

\vspace{0.2cm}
\noindent \textbf{\emph{Proposition A.7}} --- \emph{Let $\K/k$ be an extension of non-Archimedean fields and let $\M$ be a Banach $\K$-module equipped with a descent datum $\delta$. If $\delta$ is an isometry, then the canonical isomorphism $\mathrm{Ker}(\varphi_{\delta}) \widehat{\otimes}_k \K \tilde{\rightarrow} \M$ is an isometry.}

\vspace{0.1cm}
\noindent
\emph{\textbf{Proof}}.
If $\A$ is a Banach ring and if we let $\mathbf{BMod}_1(A)$ denote the subcategory of $\mathbf{BMod}(A)$ in which morphisms are required to be contractions (i.e., to have norm at most one), then a morphism between two Banach $\A$-modules is an isometric isomorphism if and only if it is an isomorphism in the category $\mathbf{BMod}_1(\A)$. According to this observation, our assertion will follow from descent theory for the categories $\mathbf{BMod}_1(k)$ and $\mathbf{BMod}_1(\K)$ instead of $\mathbf{BMod}(k)$ and $\mathbf{BMod}(\K)$. Since the canonical morphisms $\varepsilon$, $p_1$, $p_2$, $p_{12}$, $p_{23}$ and $p_{13}$ are contractions, we can apply the same arguments as in the proposition above to deduce that, indeed, the functor $\varepsilon^{*}$ defines an equivalence between the categories $\mathbf{BMod}_1(k)$ and $\mathbf{BMod}_1(\K)$. \hfill $\Box$

\vspace{0.2cm} Finally, if the non-Archimedean extension $\K/k$ is affinoid, then affinoid algebras behave well under descent.

\vspace{0.2cm}
\noindent \textbf{\emph{Proposition A.8}} --- \emph{Let $\K/k$ be an affinoid extension. A Banach $k$-algebra $\A$ is $k$-affinoid if and only if the Banach $\K$-algebra $\A \widehat{\otimes}_k \K$ is $\K$-affinoid.}

\vspace{0.1cm}
\noindent
\emph{\textbf{Proof}}.
When $\K = k_r$ with $r \in \mathbb{R}_{>0} - |k^{\times}|^{\mathbb{Q}}$, this statement is~\cite[Corollary 2.1.8]{Ber1}. The proof given there works more generally for any affinoid extension $\K/k$ once it has been noticed that $\K$ contains a dense and finitely generated $k$-subalgebra. \hfill $\Box$

\vspace{0.2cm}
\noindent \textbf{\emph{Corollary A.9}} --- \emph{Let $\K/k$ be a non-Archimedean extension. The functor from Banach $k$-algebras to Banach $\K$-algebras equipped with descent data is an equivalence of categories. Moreover, if the extension is affinoid, this functor maps $k$-affinoid algebras onto $\K$-affinoid algebras.}

\vspace{0.1cm}
\noindent
\emph{\textbf{Proof}}.
By the same argument as in Remark A.4, it follows from Proposition A.6 that a Banach $\K$-algebra $\A_{\K}$ with a descent datum comes from a Banach $k$-algebra $\A$. Moreover, in view of the previous proposition, $\A$ is a $k$-affinoid algebra if $\A_{\K}$ is a $\K$-affinoid algebra and if the extension $\K/k$ is affinoid. \hfill $\Box$

\vspace{0.1cm} We can now go back to our main technical descent result.

\vspace{0.3cm} \emph{Proof of Proposition A.1.} Let $\D$ be a subset of $\X$ such that $\D' = {\rm pr}_{\K/k}^{-1}(\D)$ is a $\K$-affinoid domain in $\X' = \X \widehat{\otimes}_k \K$ and denote by $(\A_{\D'}, \varphi')$ a pair representing the functor $\F_{\D'}: \mathbf{Aff(\K)} \rightarrow \mathbf{Sets}$. Denoting as above by $p_1$ and $p_2$ the two canonical maps from $\K$ to $\K \widehat{\otimes}_k \K$ as well as the corresponding projections $\X' \times_{\X} \X' = \X \widehat{\otimes}_k \K \widehat{\otimes}_k \K \rightarrow \X'$, $$\mathrm{Hom}_{\K \widehat{\otimes}_k \K}(p_i^{*}\A, \B) = \mathrm{Hom}_{\K}(\A, \B^{(i)})\ \ \ \ (i \in \{1,2\})$$ for any Banach $\K$-algebra $\A$ and any Banach $\K \widehat{\otimes}_k \K$-algebra $\B$, where $\B^{(i)}$ stands for $\B$ seen as a $\K$-algebra via the map $p_i$. Hence the pair $(p_i^{*}\A_{\D'}, p_i^{*}(\varphi_{V'})$ represents the functor $\F_{p_i^{-1}(\D')}$. Since $p_1 \circ {\rm pr}_{\K/k} = p_2 \circ {\rm pr}_{\K/k}$, we have $p^{-1}(\D') = p_2^{-1}(\D')$ and thus there exists an isomorphism of Banach $\K \widehat{\otimes}_k \K$-algebras $$\delta: \xymatrix{p_1^{*}\A_{V'} \ar@{->}[r]^{\sim} & p_2^{*}\A_{\D'}}$$ such that $\delta \circ p_1^{*}(\varphi_{\D'}) = p_2^{*}(\varphi_{\D'}) \circ \delta_{A}$. If we let as above $q_1$, $q_2$ and $q_3$ denote the three canonical projections from $\X' \times_{\X} \X' \times_{\X} \X'$ onto $\X'$, then $q_1^{-1}(\D') = q_2^{-1}(\D') = q_3^{-1}(\D')$ and it follows that $\delta$ satisfies the cocyle condition defining descent data. Hence $\delta$ is a descent datum on $\A_{\D'}$. One checks similarly that the map $\varphi_{\D'}: \A \rightarrow \A_{V'}$ is compatible with descent data.

Corollary 9 applies here, and thus we get a $k$-affinoid algebra
$\A_{\D}$ together with a bounded $k$-homomorphism $\varphi: \A
\rightarrow \A_{\D}$ which induce $\A_{\D'}$ and $\varphi_{\D'}$ after
base-change to $\K$. It also follows from this corollary
that, for any affinoid $k$-algebra $\B$, a bounded $k$-morphism
$\varphi: \A \rightarrow \B$ factors through $\varphi_{\D}$ if and
only if the morphism $\varepsilon^{*}(\varphi)$ factors through
$\varepsilon^{*}(\varphi_{\D}) = \varphi_{\D'}$; since this last
condition is equivalent to the inclusion of
${\rm pr}_{\K/k}^{-1}(\mathrm{im}(^a\varphi)) = \mathrm{im}(^{a}
(\varepsilon^{*}(\varphi)))$ into ${\rm pr}_{\K/k}^{-1}(\D)=\D'$, we deduce
from the surjectivity of the map ${\rm pr}_{\K/k}$ that $\varphi$ factors
through $\varphi_{\D}$ if and only if $\mathrm{im}(^{a}\varphi)$ is
contained in $\D$, i.e., if and only if $\varphi \in \F_{\D}(\B)$. Hence the
pair $(\A_{\D}, \varphi_{\D})$ represents the functor $\F_{\D}$, which completes the proof. \hfill $\Box$

\vspace{0.4cm} \noindent \textbf{(A.4)} We conclude this section with a technical result which follows easily from Proposition A.7. The norm of a Banach $k$-algebra $\A$ is said to be \emph{universally multiplicative} if, for any non-Archimedean extension $\K/k$, the norm of the Banach $\K$-algebra $\A \widehat{\otimes}_k \K$ is multiplicative.

\vspace{0.2cm}
\noindent \textbf{\emph{Lemma A.10}} --- \emph{Let $k'/k$ be a finite Galois extension and let $\A'$ be a Banach $k'$-algebra equipped with a descent datum $\delta$. We denote by $\A$ the Banach $k$-algebra such that $(\A',\delta) \simeq (\A \otimes_k k', \delta_{\A})$.
\begin{itemize}
\item[(i)] If the norm of $\A'$ is multiplicative, then the descent datum is an isometry.
\item[(ii)] If the norm of $\A'$ is universally multiplicative, then the norm of $\A$ is universally multiplicative.
\end{itemize}}

\vspace{0.1cm}
\noindent
\emph{\textbf{Proof}}.
(i) By definition, the descent datum $\delta$ is an isomorphism of Banach $k' \otimes_k k'$-algebras $p_1^{*}\A' \tilde{\rightarrow} p_2^{*}\A'$ satisfying the natural cocycle condition, where $p_1$ and $p_2$ are the canonical homomorphisms $k' \rightarrow k' \otimes_k k'$. Since $k'/k$ is a finite Galois extension, $k' \otimes_k k'$ is isometric to the product of a $[k':k]$ copies of $k'$ and thus the Banach $k$-algebra $p_1^{*}\A' = \A' \otimes_{k',p_1} (k' \otimes k')$ is isometric to the product of $[k':k]$ copies of $\A'$. The same argument applies also to $p_2^{*}\A'$.

Now, observe that the norm of $\A'$ coincides with the spectral norm since it is multiplicative.
This remains true for the product of a finite number of copies of $\A'$ since the induced norm is power-multiplicative and therefore the norms on the Banach $k'$-algebras $p_1^{*}\A'$ and $p_2^{*}\A'$ coincide with the spectral norms. Since any homomorphism of Banach algebras lowers the spectral (semi-)norms, isomorphisms are isometries with respect to the spectral (semi-)norms and we conclude that our descent datum $\delta$ is an isometry.

(ii) By construction, we have a canonical isometric monomorphism $\A \hookrightarrow \A'$ and saying that the descent datum is isometric amounts to saying that the induced isomorphism $\A \otimes_k k' \rightarrow \A'$ is an isometry.

Consider now a non-Archimedean extension $\K/k$ and pick a non-Archimedean field $\K'$ extending both $k'$ and $\K$. By assumption, the norm on $\A' \widehat{\otimes}_{k'} \K'$ is multiplicative. Thanks to the canonical isometric monomorphism $(\A \widehat{\otimes}_k \K) \widehat{\otimes}_{\K} \K' \simeq \A \widehat{\otimes}_k \K'$, it suffices to show that the norm of $\A \widehat{\otimes}_k \K'$ is multiplicative to deduce that the norm of $\A \widehat{\otimes}_k \K$ is multiplicative. Since $\A \widehat{\otimes}_k \K'$ is isometric to $(\A \otimes_k k') \widehat{\otimes}_k' \K'$, the conclusion follows from the isometry $\A \otimes_k k' \simeq \A'$ and our initial assumption. \hfill $\Box$

\section*{Appendix B: on fans}
\label{s - fans}

This second appendix deals with the technicalities useful to compactify vector spaces by means of the notion of a fan.
We use it in the case when the fan comes from Lie theory, that is when the ambient space is the Coxeter complex of a spherical root system, in which roots are seen as linear forms.

\vspace{0.3cm}

\noindent \textbf{(B.1)} Let $\M$ be a free abelian group of finite rank. We equip the abelian group $\Lambda = \mathrm{Hom}_{\mathbf{Ab}}(\M,\mathbb{R}_{>0})$ with the structure of a real vector space by setting $\lambda.\varphi = \varphi^{\lambda}$ for any $\lambda \in \mathbb{R}, \ \varphi \in \Lambda$.

\vspace{0.1cm} A (rational) \emph{polyhedral cone} is a subset of $\Lambda$ defined by a finite number of inequalities $\varphi \leqslant 1$ with $\varphi \in \M$. A \emph{face} of a polyhedral cone $\C$ is the intersection of $\C$ with a hyperplane $\{\varphi =1 \}$, where $\varphi$ is an element of $\M$ such that $\varphi_{|\C} \leqslant 1$. The cone $\C$ is \emph{strictly convex} if it contains no line.

\vspace{0.1cm}
For each strictly convex polyhedral cone $\C$ in $\Lambda$, $$\Sr_{\C} = \{\varphi \in \M \ | \varphi(u) \leqslant 1 \textrm{ for all } u \in \C\}$$ is a semigroup in $\M$ which spans $\M$ as a group and which is finitely generated (\emph{Gordan's Lemma}). Besides, $$ \C  = \{u \in \Lambda \ | \varphi(u) \leqslant 1 \ \textrm{ for all } \varphi \in \Sr_{\C}\}.$$ If $(\varphi_{i})_{i \in \I}$ is a set of generators of the semigroup $\Sr_{\C}$, each face $\F$ of $\C$ can be described by equalities $\varphi_i = 1$ with $i$ running over a subset of $\I$. Since $\Sr_{\C}$ is finitely generated, the set of faces of $\C$ is therefore \emph{finite}.

\vspace{0.2cm}
\noindent \textbf{\emph{Remark B.1}} --- Let $\C$ be a strictly convex polyhedral cone and consider a face $\F$ of $\C$. If $\F \neq \C$, there exists by definition an element $\varphi$ of $\Sr_{\C}$ such that $\varphi_{\F} = 1$ and $\varphi_{\C - \F} < 1$. Moreover, for any $\psi \in \M$ whose restriction to $\F$ is $1$, one can find a natural number $n$ such that $(n\varphi + \psi)_{\F} \leqslant 1$ on $\C$: indeed, on can find such a number so that $n\varphi + \psi$ is not greater than $1$ on any given ray (a one dimensional face) of $\C$ and, since the set of rays is finite, there exists an uniform $n$.

\vspace{0.3cm}
\noindent \textbf{(B.2)} A \emph{fan} on $\Lambda$ is a finite family $\mathcal{F}$ of polyhedral cones satisfying the following conditions: \begin{itemize} \item each cone is strictly convex; \item the union of all these cones covers $\Lambda$; \item for each cones $\C, \C' \in\mathcal{F}$, $\C \cap \C'$ is a face of $\C$ and $\C'$; \item each face of a cone $\C \in\mathcal{F}$ belongs to $\mathcal{F}$.
\end{itemize}

\vspace{0.1cm} To any fan $\mathcal{F}$ on the vector space $\Lambda$ corresponds a compactification $\overline{\Lambda}^{\mathcal{F}}$ of $\Lambda$ which we now describe.

Letting $\mathbf{Mon}$ denote the category of unitary monoids, the \emph{canonical compactification} of a polyhedral cone $\C$ is defined as the set $$\overline{\C} = \mathrm{Hom}_{\mathbf{Mon}}(\Sr_{\C},[0,1])$$ of all morphisms of unitary monoids $\Sr_{\C} \rightarrow [0,1]$, equipped with the coarsest topology for which each evaluation map $\overline{\C} \rightarrow [0,1], \ u \mapsto \varphi(u)$, is continuous, where $\varphi \in \Sr_{\C}$. This topological space is compact since it can be canonically identified with a closed subspace of the product space $[0,1]^{\Sr_{\C}}$. The canonical map $$\C \rightarrow \overline{\C}, \ \ u \mapsto (\varphi \mapsto \varphi(u))$$ identifies $\C$ homeomorphically with the open subset $\mathrm{Hom}_{\mathbf{Mon}}(\Sr_{\C}, ]0,1])$ of $\mathrm{Hom}_{\mathbf{Mon}}(\Sr_{\C}, [0,1])$ (that this subset is open follows from the finite generation of $\Sr_{\C}$).

\vspace{0.2cm} \noindent \textbf{\emph{Lemma B.2}} --- \emph{Let $\C$ be a strictly convex polyhedral cone and $\F$ a face of $\C$.
\begin{itemize}
\item[(i)] There exists a unique continuous map $\overline{\F} \rightarrow \overline{\C}$ extending the inclusion $\F \hookrightarrow \C$. This map is a homeomorphism between $\overline{\F}$ and the closure of $\F$ in $\overline{\C}$.
\item[(ii)] Let $\Sr_{\C}^{\F}$ denote the subset of $\Sr_{\C}$ consisting of those elements $\varphi$ such that $\varphi_{|\F} = 1$ and let $\langle \F \rangle$ denote the linear subspace of $\Lambda$ generated by $\F$. The set $$\C_{\F} = \{u \in \overline{\C} \ | \ \varphi(u) > 0 \ \textrm{ for all } \varphi \in \Sr_{\C}^{\F} \ \textrm{and } \varphi(u) = 0 \ \textrm{for all } \varphi \in \Sr_{\C} - \Sr_{\C}^{\F}\}$$ is canonically identified with a strictly convex polyhedral cone in the vector space $$\Lambda/\langle \F \rangle = {\rm Hom}_{\mathbf{Mon}}(\Sr_{\C}^{\F}, \mathbb{R}_{>0}).$$
\item[(iii)] If we let $\mathcal{F}^{\bullet}$ denote the set of faces of $\C$, $$\overline{\C} = \bigsqcup_{\F \in \mathcal{F}^{\bullet}} \C_{\F}.$$
\item[(iv)] For any cone $\C \in\mathcal{F}$ and any faces $\F, \F'$ of $\C$, $$\C_{\F} \cap \overline{\F'} = \left\{\begin{array}{ll} \F'_{\F} & \textrm{ if } \F \subset \F' \\ \varnothing & \textrm{ otherwise} \end{array} \right.$$ in $\overline{\C}$.
\end{itemize}}

\vspace{0.1cm}
\noindent
\emph{\textbf{Proof}}.
(i) To the inclusion $\F \subset \C$ corresponds an inclusion $\Sr_{\C} \subset \Sr_{\F}$, hence a natural continuous map $$i: \overline{\F} = \mathrm{Hom}_{\mathbf{Mon}}(\Sr_{\F}, [0,1]) \rightarrow \mathrm{Hom}_{\mathbf{Mon}}(\Sr_{\C},[0,1]) = \overline{\C}$$ extending the inclusion of $\F = \mathrm{Hom}_{\mathbf{Mon}}(\Sr_{\F}, ]0,1])$ into $\C = \mathrm{Hom}_{\mathbf{Mon}}(\Sr_{\C}, ]0,1])$. If the latter is strict, injectivity of $i$ follows from Remark B.1: with the notation introduced there, if $u, \ v \in \overline{\F}$ have the same restriction to $\Sr_{\C}$, then $\varphi(u)^n\psi(u) = \varphi(v)^n\psi(v)$ and thus $\psi(u)=\psi(v)$ since $\varphi(u) = \varphi(v) =1$.

The topological spaces $\overline{\F}$ and $\overline{\C}$ being compact, the continuous injection $i$ is a homeomorphism onto its image and $i(\overline{\F})$ is the closure of $i(\F)$ in $\overline{\C}$ since $\F$ is dense in $\overline{\F}$.

\vspace{0.1cm} (ii) We have $$\langle \F \rangle = \{u \in \Lambda \ | \varphi(u)=1 \ \textrm{ for all } \varphi \in \Sr_{\C}^{\F}\}$$ and the canonical map $\Lambda = \mathrm{Hom}_{\mathbf{Mon}}(\Sr_{\C}, \mathbb{R}_{>0}) \rightarrow  \mathrm{Hom}_{\mathbf{Mon}}(\Sr_{\C}^{\F}, \mathbb{R}_{>0})$ deduced from the inclusion of $\Sr_{\C}^{\F}$ into $\Sr_{\C}$ induces a linear isomorphism between the vector spaces $\Lambda/\langle \F \rangle$ and $\mathrm{Hom}_{\mathbf{Mon}}(\Sr_{\C}^{\F}, \mathbb{R}_{>0})$.

If $\N$ denotes the subgroup of $\M$ consisting of all elements $\varphi$ such that $\varphi_{|\F}=1$ and if $\W = \mathrm{Hom}_{\mathbf{Ab}}(\N,\mathbb{R}_{>0})$, then $\N$ is free of finite rank and $\Sr_{\C}^{\F}$ is canonically isomorphic to the semigroup in $\N$ associated with the strictly convex polyhedral cone $p(\C)$ of $\W$, where $p$ denotes the canonical projection of $\Lambda$ on $\W$. Thus $\Sr_{\C}^{\F}$ is finitely generated by Gordan'sLemma. Besides, it follows immediately from the definition of $\Sr_{\C}^{\F}$ that it contains the sum of two elements of $\Sr_{\C}$ if and only if it contains both summands. One deduces from this last property that, for any $u \in \mathrm{Hom}_{\mathbf{Mon}}(\Sr_{\C}^{\F}, ]0,1])$, the map $\tilde{u}$ from $\Sr_{\C}$ to $[0,1]$ defined by $$\varphi(\tilde{u}) = \left\{ \begin{array}{ll} \varphi(u) & \textrm{ if } \varphi \in \Sr_{\C}^{\F} \\ 0 & \textrm{ otherwise}\end{array}\right.$$ is a morphism of unitary monoids, hence defines a point in $\C_{\F}$. We thus get a homeomorphism between $\C_{\F}$ and the polyhedral cone $\mathrm{Hom}_{\mathbf{Mon}}(\Sr_{\C}^{\F}, ]0,1])$ in $\Lambda/\langle \F \rangle$.

\vspace{0.1cm} (iii) Let us consider a point $u$ in $\overline{\C}$. We let $\Sigma$ denote the set of all $\varphi \in \Sr_{\C}$ such that $\varphi(u) >0$ and $\F$ the subset of $\C$ defined by the conditions $\varphi = 1$, $\varphi \in \Sigma$. Then $\F$ is a face of $\C$ and $\Sigma \subset \Sr_{\C}^{\F}$. If we pick $\varphi_1, \ldots, \varphi_r$ in $\Sr_{\C}$ such that $\C \cap \{\varphi_i=1\}$ are the different faces of codimension one of $\C$ containing $\F$, then: \begin{itemize} \item for any $\varphi \in \Sr_{\C}^{\F}$, there exists an integer $n \geqslant 1$ such that $n\varphi$ belongs to $\mathbb{N}\varphi_1 + \ldots + \mathbb{N}\varphi_r$; \item for any $i \in \{1, \ldots, r\}$, there exists an element $\varphi$ in $\Sigma$ and an integer $n \geqslant 1$ such that $n\varphi = n_1\varphi_1 + \ldots n_r \varphi_r$ with $n_1, \ldots, n_r \in \mathbb{N}$ and $n_i \geqslant 1$.
\end{itemize}

Since the sum of two elements of $\Sr_{\C}$ belongs to $\Sigma$ if and only if both summands belong to $\Sigma$, the last property implies $\varphi_1, \ldots, \varphi_r \in \Sigma$ and then the identity $\Sigma = \Sr_{\C}^{\F}$ follows from the first one.

Finally, the point $u$ belongs to the cone $\C_{\F}$ and thus $\overline{\C} = \bigsqcup_{\F \in \mathcal{F}^{\bullet}} \C_{\F}$.

\vspace{0.1cm} (iv) We have $\F \subset \F'$ if and only if $\Sr_{\C}^{\F'} \subset \Sr_{\C}^{\F}$. If $\Sr_{\C}^{\F'} \nsubseteq \Sr_{\C}^{\F}$, there exists $\varphi \in \Sr_{\C}^{\F'}$ which does not belong to $\Sr_{\C}^{\F}$ and therefore $\C_{\F} \cap \overline{\F'} = \varnothing.$

If $\Sr_{\C}^{\F'} \subset \Sr_{\C}^{\F}$, then \begin{eqnarray*} \C_{\F} \cap \overline{\F'} & = & \left\{ u \in \mathrm{Hom}_{\mathbf{Mon}}(\Sr_{\C}, [0,1]) \ \left| \begin{array}{l} \varphi(u) = 0 \textrm{ for any } \varphi \notin \Sr_{\C}^{\F} \textrm{ and } \varphi(u) > 0 \textrm{ for any } \varphi \in \Sr_{\C}^{\F} \\ \hspace{2cm} \varphi(u) = 1 \textrm{ for any } \varphi \in \Sr_{\C}^{\F'} \end{array}\right.\right\}\\ & = & \left\{u \in {\rm Hom}_{\mathbf{Ab}}(\Sr_{\F'},[0,1]) \ \left| \begin{array}{l} \varphi(u) = 0 \textrm{ for any } \varphi \notin \Sr_{\C}^{\F} \\  \varphi(u) >0  \textrm{ for any } \varphi \in \Sr_{\C}^{\F'} \end{array} \right. \right\} \\ & = & \F'_{\F}, \end{eqnarray*} for $\Sr_{\F'}$ is the subgroup of $\M$ generated by $\Sr_{\C}$ and $-\Sr_{\C}^{\F'}$. \hfill $\Box$

\vspace{0.5cm}
Consider now a fan $\mathcal{F}$ on the vector space $\Lambda$. We deduce from the first assertion in the lemma above that the compactified cones $\{\overline{\C}\}_{\mathcal{F}^{\bullet}}$ glue together to define a compact topological space $\overline{\Lambda}^{\mathcal{F}}$ containing $\Lambda$ as a dense open subset. Indeed, it is enough to define $\overline{\Lambda}^{\mathcal{F}}$ as the quotient of the compact topological space $\bigsqcup_{\C \in\mathcal{F}} \overline{\C}$ by the following equivalence relation: two points $x \in \overline{\C}$ and $y \in \overline{\C'}$ are equivalent if and only if there exists a cone $\C'' \in\mathcal{F}$ contained in $\C$ and $\C'$, as well as a point $z \in \overline{\C''}$ mapped to $x$ and $y$, respectively, under the canonical injections of $\overline{\C''}$ into $\overline{\C}$ and $\overline{\C'}$, respectively. The quotient space $\overline{\Lambda}^{\mathcal{F}}$ is compact since we have glued together a finite number of compact spaces along closed subspaces.

Each compactified cone $\overline{\C}$ embeds canonically into $\overline{\Lambda}^{\mathcal{F}}$ and, since $\C' = \overline{\C'} \cap \C$ for any cones $\C, \ \C' \in\mathcal{F}$ satisfying $\C' \subset \C$, the natural map $\bigsqcup_{\C \in\mathcal{F}} \C \rightarrow \overline{\Lambda}^{\mathcal{F}}$ factors through the canonical projection $\bigsqcup_{\C \in\mathcal{F}} \C \rightarrow \Lambda$ and induces therefore a homeomorphism between $\Lambda$ and a dense open subset of $\overline{\Lambda}^{\mathcal{F}}$.

\vspace{0.2cm} \noindent \textbf{\emph{Proposition B.3}} --- \emph{Let us consider a fan $\mathcal{F}$ on the vector space $\Lambda$.
\begin{itemize}
\item[(i)] For any cone $\C \in\mathcal{F}$, there exists a canonical homeomorphism $i^{\C}$ between the vector space $\Lambda/\langle \C \rangle$ and a locally closed subset $\Sigma_{\C}$ of $\overline{\Lambda}^{\mathcal{F}}$. \\
The set $\mathcal{F}_{\C}$ of cones in $\mathcal{F}$ containing $\C$ induces a fan on the vector space $\Lambda/\langle \C \rangle$ and the map $i^{\C}$ extends to a homeomorphism between the associated compactification of $\Lambda/\langle \C \rangle$ and the closure of $\Sigma_{\C}$ in $\overline{\Lambda}^{\mathcal{F}}$.
\item[(ii)] The family $\{\Sigma_{\C}\}_{\C \in\mathcal{F}}$ is a stratification of $\overline{\Lambda}^{\mathcal{F}}$ into locally closed subspaces: $$\overline{\Lambda}^{\mathcal{F}} = \bigsqcup_{\C \in\mathcal{F}} \Sigma_{\C} \ \hspace{0.5cm} \ \textrm{and} \ \hspace{0.5cm} \ \overline{\Sigma_{\C}} = \bigsqcup_{\tiny \begin{array}{l} \C' \in\mathcal{F} \\ \C \subset \C'\end{array}} \Sigma_{\C'}.$$
\item[(iii)] The action of $\Lambda$ on itself by translations extends to an action of $\Lambda$ on $\overline{\Lambda}^{\mathcal{F}}$ by homeomorphisms stabilizing each stratum and, via the identification $i^{\C}: \Lambda/\langle \C \rangle \tilde{\rightarrow} \Sigma_{\C}$, the action induced on the stratum $\Sigma_{\C}$ is the action of $\Lambda$ on $\Lambda/\langle \C \rangle$ by translations.
\item[(iv)] A sequence $(p_n)$ of points in $\Lambda$ converges to a point of $\overline{\Lambda}^{\mathcal{F}}$ belonging to the stratum $\Sigma_{\C}$ if and only if the following two conditions hold: \begin{itemize} \item almost all points $p_n$ lie in the union of the cones $\C' \in\mathcal{F}$ containing $\C$; \item for any cone $C' \in\mathcal{F}_{\C}$ and any element $\varphi$ of $\Sr_{\C'}$, the sequence $(\varphi(p_n))$ converges in $[0,+\infty[$ and $$\lim \varphi(p_n) = 0 \Longleftrightarrow \varphi \notin \Sr_{\C'}^{\C}.$$ \end{itemize}
\end{itemize}}

\vspace{0.1cm}
\noindent
\emph{\textbf{Proof}}.
(i) For any cone $\C$ in $\mathcal{F}$, the quotient vector space $\Lambda/\langle \C \rangle$ is canonically isomorphic to $\mathrm{Hom}_{\mathbf{Ab}}(\M^{\C}, \mathbb{R}_{>0})$, where $\M^{\C}$ denotes the subgroup of $\M$ consisting in elements $\varphi$ satisfying $\varphi_{|\C} = 1$. For any cone $\C'$ in $\mathcal{F}$ containing $\C$, the semigroup $\Sr_{\C'} \cap \M^{\C} = \Sr_{\C'}^{\C}$ is finitely generated and spans $\M^{\C}$. If we let $\mathcal{F}_{\C}$ denote the set of all cones $\C' \in\mathcal{F}$ containing $\C$ and $p$ the canonical projection of $\Lambda$ onto $\Lambda/\langle \C \rangle$, it follows that the polyhedral cones $p(\C') = \mathrm{Hom}_{\mathbf{Mon}}(\Sr_{\C'}^{\C}, ]0,1])$, $\C' \in\mathcal{F}_{\C}$, define a fan on the vector space $\Lambda/\langle \C \rangle$.

For any cone $\C' \in\mathcal{F}$, extension by zero on $\Sr_{\C'} - \Sr_{\C'}^{\C}$ provides us with a map $$i^{\C}_{\C'}: p(\C') = \mathrm{Hom}_{\mathbf{Mon}}(\Sr_{\C'}^{\C}, ]0,1]) \rightarrow \mathrm{Hom}_{\mathbf{Mon}}(\Sr_{\C'}, [0,1]) = \overline{\C'} \subset \overline{\Lambda}^{\mathcal{F}}.$$ This map is a homeomorphism onto the locally closed subspace $$\C'_{\C} = \left\{u \in \overline{\C'} \ | \ \varphi(u) = 0 \textrm{ for any } \varphi \in \Sr_{\C'} - \Sr_{\C'}^{\C} \textrm{ and } \varphi(u) >0 \textrm{ for any } \varphi \in \Sr_{\C'}^{\C}\right\}.$$ Moreover, for any face $\C''$ of $\C'$ containing $\C$, $\Sr_{\C'} \subset \Sr_{\C''}$ and the natural diagram $$\xymatrix{p(\C') = \mathrm{Hom}_{\mathbf{Mon}}(\Sr_{\C'}^{\C}, ]0,1]) \ar@{->}[r]^{i_{\C'}^{\C}} \ar@{<-}[d] & \mathrm{Hom}_{\mathbf{Mon}}(\Sr_{\C'}, [0,1]) = \overline{\C'} \ar@{<-}[d] \\ p(\C'') = \mathrm{Hom}_{\mathbf{Mon}}(\Sr_{\C''}^{\C},]0,1]) \ar@{->}[r]_{i_{\C''}^{\C}} & \mathrm{Hom}_{\mathbf{Mon}}(\Sr_{\C''},[0,1]) = \overline{\C''}}$$ is commutative. Therefore there exists a unique map $i^{\C}$ from $\Lambda/\langle \C \rangle$ to $\overline{\Lambda}^{\mathcal{F}}$ whose restriction to each cone $p(\C')$, $\C' \in\mathcal{F}_{\C}$, coincides with $i_{\C'}^{\C}$. Let $\Sigma_{\C}$ denote the union of all cones $\C'_{\C}$ with $\C' \in\mathcal{F}_{\C}$.

Thanks to the gluing conditions (iv) of Lemma B.2, the map $i^{\C}$ is a homeomorphism between the vector space $\Lambda/\langle \C \rangle$ and the subspace $\Sigma_{\C}$ of $\overline{\Lambda}^{\mathcal{F}}$. Since $$\Sigma_{\C} \cap \overline{\C''} = \left\{\begin{array}{ll} \C''_{\C} & \textrm{ if } \C \subset \C'' \\ \varnothing & \textrm{ otherwise}\end{array} \right.$$
is a locally closed subspace of $\overline{\C''}$ for any cone $\C'' \in\mathcal{F}$, $\Sigma_{\C}$ is a locally closed subspace of $\overline{\Lambda}^{\mathcal{F}}$.

Finally, if $\C'$ is a cone in $\mathcal{F}_{\C}$, the closure of $\C'_{\C}$ in $\overline{\Lambda}^{\mathcal{F}}$ is canonically homeomorphic to the canonical compactification of this cone by lemma 2, (i) and it follows that the map $i^{\C}: \Lambda/\langle \C \rangle \rightarrow \Sigma_{\C} \subset \overline{\Lambda}^{\mathcal{F}}$ extends to a homeomorphism between the compactification of $\Lambda/\langle \C \rangle$ coming from the fan $\{p(\C')\}_{\C' \in\mathcal{F}_{\C}}$ and the closure of $\Sigma_{\C}$ in $\overline{\Lambda}^{\mathcal{F}}$.

\vspace{0.1cm} (ii) Given two cones $\C, \ \C'$ in $\mathcal{F}$ such that $\Sigma_{\C} \cap \Sigma_{\C'} \neq \varnothing$, we can pick $\Gamma$ and $\Gamma'$ in $\mathcal{F}$ with $\C \subset \Gamma$, $\C' \subset \Gamma'$ and $\Gamma_{\C} \cap \Gamma'_{\C'} \neq \varnothing$. Since $$\overline{\Gamma} \cap \Gamma'_{\C'} = \left\{\begin{array}{ll} \Gamma_{\C'} & \textrm{ if } \C' \subset \Gamma \\ \varnothing & \textrm{ otherwise}\end{array}\right.$$ (Lemma B.2, (iv)), we deduce $\C' \subset \Gamma$, hence $\Gamma_{\C} \cap \Gamma_{\C'} \neq \varnothing$ and, finally, $\C' = \C$. Thus the locally closed subspaces $\Sigma_{\C}$, $\C \in\mathcal{F}$, are pairwise disjoint.

Moreover, for any cones $\C, \ \C' \in\mathcal{F}$ with $\C \subset \C'$, $$\overline{\C_{\C'}} = \bigcup_{\tiny \begin{array}{l} \C'' \in\mathcal{F} \\ \C \subset \C'' \subset \C'\end{array}} \C'_{\C''}$$ and therefore $$\overline{\Sigma_{\C}} = \bigcup_{\tiny \begin{array}{l} \C' \in\mathcal{F} \\ \C \subset \C'\end{array}} \overline{\C'_{\C}} = \bigcup_{\tiny \begin{array}{l} \C', \ \C'' \in\mathcal{F} \\ \C \subset \C'' \subset \C'\end{array}} \C'_{\C''} = \bigcup_{\tiny \begin{array}{l} \C'' \in\mathcal{F} \\ \C \subset \C''\end{array}} \Sigma_{\C''}.$$

\vspace{0.1cm} (iii) Let us pick a vector $v \in \Lambda$ and consider the unique map $t_{v}: \overline{\Lambda}^{\mathcal{F}} \rightarrow \overline{\Lambda}^{\mathcal{F}}$ fulfilling the following requirement: for any cone $\C \in\mathcal{F}$, $t_v(\Sigma_{\C}) \subset \Sigma_{\C}$ and the map $(i^{C})^{-1} \circ t_v \circ i^{\C}: \Lambda/\langle \C \rangle \rightarrow \Lambda/\langle \C \rangle$ is the translation by the vector $v$. Given a cone $\C \in\mathcal{F}$, note that the union of all strata $\Sigma_{\C'}$, $\C' \subset \C$, is naturally homeomorphic to $\mathrm{Hom}_{\mathbf{Mon}}(\Sr_{\C}, \mathbb{R}_{\geqslant 0})$. This observation allows us to make the restriction of $t_v$ to $\overline{\C} = \mathrm{Hom}_{\mathbf{Mon}}(\Sr_{\C},[0,1])$ explicit: for any point $u \in \overline{\C}$, $t_v(u)$ is the point of $\overline{\Lambda}^{\mathcal{F}}$ corresponding to the morphism of unitary monoids $$\Sr_{\C} \rightarrow \mathbb{R}_{\geqslant 0}, \ \varphi \mapsto \varphi(v)\varphi(u).$$ Clearly, the restriction of the map $t_v$ to each compactified cone $\overline{\C}$ is continuous and therefore this map is continuous.
The map $\Lambda \times \overline{\Lambda}^{\mathcal{F}} \rightarrow \overline{\Lambda}^{\mathcal{F}}, \ (v,u) \mapsto t_v(u)$ provides $\overline{\Lambda}^{\mathcal{F}}$ with an action of $\Lambda$ by homeomorphisms; this action stabilizes each stratum, on which it induces the natural action of $\Lambda$ by translation via the identifications $i^{\C}: \Lambda/\langle \C \rangle \simeq \Sigma_{\C}$.

\vspace{0.1cm} (iv) We consider a sequence $(p_n)$ of points in $\Lambda$.

Let us first assume that this sequence converges to a point $p$ in $\overline{\Lambda}^{\mathcal{F}}$ belonging to the stratum $\Sigma_{\C}$. We consider a cone $\C' \in\mathcal{F}$ containing $\C$ and an element $\varphi$ in $\Sr_{\C'}$.

There is at least one cone $\C''$ in $\mathcal{F}_{\C}$ containing infinitely many $p_n$, since  $\bigcup_{\C'' \in\mathcal{F}_{\C}} \overline{\C''}$ contains the stratum $\Sigma_{\C}$. We pick one of them. For any $\psi \in \Sr_{\C''}$, $$\lim_{p_n \in \C''} \psi(p_n) = \psi(p)$$ and $\psi(p) = 0$ if and only if $\psi \notin \Sr_{\C''}^{\C}$. Since $\C$ is a common face of the cones $\C'$ and $\C''$, there exists an element $\psi \in \Sr_{\C''}^{\C}$ such that $\varphi + \psi$ belongs to $\Sr_{\C''}$. Moreover, $\varphi + \psi \in \Sr_{\C''}^{\C}$ if and only if $\varphi \in \Sr_{\C''}^{\C}$ as $(\varphi + \psi)_{|\C} = \varphi_{|\C}$. Since $\psi(p) > 0$, it follows that $$\lim_{p_n \in \C''} \varphi(p_n) = \varphi(p),$$ and $\varphi(p)=0$ if and only if $\varphi \notin \Sr_{\C'}^{\C}$.

\vspace{0.1cm} Let us now assume that the sequence $(p_n)$ is eventually contained in the union of all cones $\C' \in\mathcal{F}$ containing $\C$ and that, for any $\C' \in\mathcal{F}_{\C}$ and any $\varphi \in \Sr_{\C'}$, the sequence $(\varphi(p_n))$ converges in $[0,+\infty[$, with $\lim \varphi(p_n) = 0$ if and only if $\varphi \in \Sr_{\C'} - \Sr_{\C'}^{\C}$.

Given any cone $\C' \in\mathcal{F}_{\C}$ containing an infinite number of terms of the sequence $(p_n)$, we define a map $p_{\C'}: \Sr_{\C'} \rightarrow [0,1]$ by setting $\varphi(p_{\C'}) = \lim_{p_n \in \C'} \varphi(p_n)$ for all $\varphi \in \Sr_{\C'}$. This is obviously a morphism of unitary monoids, hence a point in $\overline{\C'}$, and it follows from our assumption that $p_{\C'}$ belongs to $\Sigma_{\C}$.

If $\C'$ and $\C''$ are two cones in $\mathcal{F}_{\C}$, both of them containing infinitely many $p_n$, then $\C' \cap \C''$ is a cone in $\mathcal{F}_{\C}$ and obviously $p_{\C'} = p_{\C' \cap \C''} = p_{\C''}$. Thus the sequence $(p_n)$ converges in $\overline{\Lambda}^{\mathcal{F}}$ to a point of $\Sigma_{\C}$. \hfill $\Box$

\vspace{0.3cm} \noindent \textbf{(B.3)} More generally, property (iii) in the proposition above allows us to compactify any affine space $\A$ under the vector space $\Lambda$. Let $\sim$ denote the usual equivalence relation on $\A \times \Lambda$: $(a,v) \sim (a',v')$ if $a+v = a'+v'$. The structural map $\A \times \Lambda \rightarrow \A, \ (a,v) \mapsto a+v$ induces a homeomorphism between the quotient space $\A \times \Lambda /\sim$ and $\A$. Embedding $\A \times \Lambda$ in $\A \times \overline{\Lambda}^{\mathcal{F}}$, one checks that the closure of the equivalence relation $\sim$ is an equivalence relation $\sim'$ which we can easily make explicit:

\vspace{0.1cm}
\begin{center}$(a,x) \sim' (a',y)$ if and only if $x$ and $y$ are contained in the same stratum $\Sigma_{\C}$ of $\overline{\Lambda}^\mathcal{F}$ \\
and there exists some $v \in \Lambda$ such that $y= t_v(x)$ and $a'+v \in a + \langle \C \rangle$.\end{center}
\vspace{0.1cm}

Then we define $\overline{\A}^{\mathcal{F}}$ to be the quotient topological space $\A \times \overline{\Lambda}^{\mathcal{F}}/\sim'$.

\vspace{0.3cm}
\noindent \textbf{\emph{Proposition B.4}} --- \emph{Let $\A$ be an affine space under the vector space $\Lambda$ and let $\mathcal{F}$ be a fan on $\Lambda$.
\begin{itemize}
\item[(i)] The topological space $\overline{\A}^{\mathcal{F}}$ is compact and the canonical map $\A \rightarrow \overline{\A}^{\mathcal{F}}$ is a homeomorphism onto a dense open subset of $\overline{\A}^{\mathcal{F}}$.
\item[(ii)] For any point $a \in \A$, the map $\Lambda \rightarrow \A, \ v \mapsto a +v$ extends uniquely to a homeomorphism between $\overline{\Lambda}^{\mathcal{F}}$ and $\overline{\A}^{\mathcal{F}}$.
\item[(iii)] For any vector $v \in \Lambda$, the translation $\A \rightarrow \A, \ a \mapsto a + v$ extends uniquely to an automorphism of the topological space $\overline{\A}^{\mathcal{F}}$.
\item[(iv)] The topological space $\overline{\A}^{\mathcal{F}}$ is stratified into affine spaces: $$\overline{\A}^{\mathcal{F}} = \bigsqcup_{\C \in\mathcal{F}} \A/\langle \C \rangle.$$
\end{itemize}}

\vspace{0.1cm}
\noindent
\emph{\textbf{Proof}}.
(i) and (ii) The topological space $\overline{\A}^{\mathcal{F}}$ is Hausdorff because the equivalence relation $\sim'$ is closed. Since the equivalence relation $\sim$ on $\A \times \Lambda$ is closed as well, $\A \times \Lambda$ is invariant under $\sim'$ and its image in $\overline{\A}^{\mathcal{F}}$ is thus a dense open subset.

Let us pick a point $a$ in $\A$ and check that the canonical projection $p: \A \times \overline{\Lambda}^{\mathcal{F}} \rightarrow \overline{\A}^{\mathcal{F}}$ induces a homeomorphism between $\{a\} \times \overline{\Lambda}^{\mathcal{F}}$ and $\overline{\A}^{\mathcal{F}}$. Since $\{a\} \times \overline{\Lambda}^{\mathcal{F}}$ is compact and $\overline{\A}^{\mathcal{F}}$ is Hausdorff, the continuous map $p: \{a\} \times \overline{\Lambda}^{\mathcal{F}} \rightarrow \overline{\A}^{\mathcal{F}}$ is closed; its image is a closed subset of $\overline{\A}^{\mathcal{F}}$ containing the dense open subset $p(\A \times \Lambda)$, thus $p(\{a\} \times \overline{\Lambda}^{\mathcal{F}}) = \overline{\A}^{\mathcal{F}}$ and therefore $\overline{\A}^{\mathcal{F}}$ is compact. Finally, given two points $x,\ y \in \overline{\Lambda}^{\mathcal{F}}$ with $(a,x) \sim' (a,y)$, we may choose sequences $(x_n)$ and $(y_n)$ of points in $\Lambda$ converging to $x$ and $y$ respectively and satisfying $(a,x_n) \sim (a, y_n)$ for all $n$; then we have $x_n = y_n$ for all $n$, hence $x=y$, for the topological space $\overline{\Lambda}^{\mathcal{F}}$ is Hausdorff. Thus the map $p: \{a\} \times \overline{\Lambda}^{\mathcal{F}} \rightarrow \overline{\A}^{\mathcal{F}}$ is a homeomorphism.

\vspace{0.1cm}

(iii) This assertion follows immediately from Proposition B.3, (iii).

\vspace{0.1cm}

(iv) Let $\C$ be a cone in $\mathcal{F}$ and $0_{\C}$ denote the origin of the stratum $\Sigma_{\C} \simeq \Lambda/\langle \C \rangle$ in $\overline{\Lambda}^{\mathcal{F}}$. The map $\A \rightarrow \A \times \overline{\Lambda}^{\mathcal{F}}, \ a \mapsto (a,0_{\C})$ is $\Lambda$-equivariant and induces a homeomorphism between the quotient affine space $\A/\langle \C \rangle$ and a locally closed subspace of $\overline{\A}^{\mathcal{F}}$ which we can also describe as the image of $\A \times \Sigma_{\C}$ under the canonical projection $p$. Relying on Proposition B.3, it follows from (ii) that the locally closed subspaces of this kind  define a stratification of $\overline{\A}^{\mathcal{F}}$:

\vspace{0.1cm}
\centerline{$\displaystyle \overline{\A}^{\mathcal{F}} = \bigsqcup_{\C \in\mathcal{F}} \A/\langle \C \rangle$.}
\hfill $\Box$

\vspace{0.3cm} \noindent \textbf{\emph{Remark B.5}} --- More generally, one can define a \emph{prefan} on the real vector space $\Lambda$ as the preimage $\mathcal{F}$ of a fan $\mathcal{F}'$ on a quotient space $\Lambda' = \Lambda/\Lambda_0$. It consists of rational polyhedral cones in $\Lambda$ satisfying all the defining conditions of a fan but strict convexity, since each cone contains the vector subspace $\Lambda_0$. If $\A$ is an affine space under $\Lambda$, one agrees on defining $\overline{\A}^{\mathcal{F}}$ as the compactification $\overline{\A'}^{\mathcal{F}'}$ of $\A' = \A/\Lambda_0$ with respect to the fan $\mathcal{F}'$.

\section*{Appendix C: on non-rational types} This last appendix deals with non-rational types of parabolic subgroups and with the corresponding compactifications of a building.
We consider a semisimple linear group $\G$ over a non-Archimedean field $k$ and recall that a \emph{type} $t$ of parabolic subgroups of $\G$ is by definition a connected component of the $k$-scheme ${\rm Par}(\G)$, which we denote by ${\rm Par}_t(\G)$. If $\G$ is split, then types are in one-to-one correspondence with $\G(k)$-conjugacy classes of parabolic subgroups of $\G$, and for any type $t$, ${\rm Par}_t(\G)$ is isomorphic to $\G/\rP$, where $\rP$ is any parabolic subgroup of $\G$ defining a $k$-point in ${\rm Par}_t(\G)$. In general, a type $t$ is said to be $k$-\emph{rational} if the component ${\rm Par}_t(\G)$ has a $k$-point. The most important example is the type $\varnothing$ of Borel subgroups of $\G$: the scheme ${\rm Bor}(\G) = {\rm Par}_{\varnothing}(\G)$ is a geometrically connected component of ${\rm Par}(\G)$, and the type $\varnothing$ is $k$-rational if and only if $\G$ has a Borel subgroup, i.e., if and only if $\G$ is quasi-split.

\vspace{0.1cm} Let $t$ be any type. The construction of (3.4.1) makes sense even if $t$ is non-rational: we just consider the map $\vartheta_t : \mathcal{B}(\G,k) \rightarrow {\rm Par}(\G)^{\rm an}$ defined in (2.4.3), take the closure $\overline{\A}_t(\Sr,k)$ of the image of some apartment $\A(\Sr,k)$ and define the compactified building $\overline{\mathcal{B}}_t(\G,k)$ as the topological quotient of $\G(k) \times \overline{\A}_t(\Sr,k)$ under the equivalence relation induced by the map $$\G(k) \times \overline{\A}_t(\Sr,k) \rightarrow {\rm Par}(\G)^{\rm an}, \ \ \ (g,x) \mapsto g \cdot x = gxg^{-1}.$$
Equivalently, $\overline{\mathcal{B}}_t(\G,k)$ is the closure of $\mathcal{B}(\G,k)$ in the compactified building $\overline{\mathcal{B}}_{t'}(\G,k')$, where $k'/k$ is a finite extension splitting $\G$ and $t'$ denotes a type of $\G \otimes_k k'$ dominating $t$.

\vspace{0.1cm}
Our aim is to show that there exists a $k$-\emph{rational} type $t'$ such that $\overline{\mathcal{B}}_t(\G,k) \cong \overline{\mathcal{B}}_{t'}(\G,k)$.

\vspace{0.1cm} Let $\rP_0$ be a minimal parabolic subgroup of $\G$. By Proposition \ref{prop.osc} and Galois descent, the functor $$(\mathbf{Sch}/k)^{\rm op} \rightarrow \mathbf{Sets}, \ \ \ \Sr \mapsto \{\rP \in {\rm Par}_t(\G)(\Sr) \ | \ \rP \textrm{ and } \rP_0 \times_k \Sr \textrm{ are osculatory}\}$$ is representable by a closed and smooth subscheme ${\rm Osc}_t(\rP_0)$ of ${\rm Par}_t(\G)$, homogeneous under $\rP_0$ and such that, for any finite Galois extension $k'/k$, $${\rm Osc}_t(\rP_0) \otimes_k k' = \bigcup_{t' \in \I} {\rm Osc}_{t'}(\rP_0 \otimes_k k'),$$ where $\I$ is the set of types of $\G \otimes_k k'$ dominating $t$. One proves as in Proposition \ref{prop.stab} the existence of a largest parabolic subgroup $\Qr_0$ of $\G$ stabilizing ${\rm Osc}_t(\rP_0)$. The conjugacy class of $\Qr_0$ does not depend on the initial choice of $\rP_0$ since minimal parabolic subgroups of $\G$ are conjugate under $\G(k)$, hence defines a $k$-rational type $\tau$.

\vspace{0.2cm} \noindent \textbf{\emph{Example}} ---  1. If $t$ is $k$-rational, then $\Qr_0$ is the unique parabolic subgroup of $\G$ of type $t$ containing $\rP_0$. Indeed, let $\rP$ be the parabolic subgroup of type $t$ containing $\rP_0$. Since ${\rm Osc}_t(\rP_0)$ is homogeneous under $\rP_0$, this scheme is reduced to the closed point $\rP$ of ${\rm Par}_t(\G)$ and thus $\Qr_0 = \rP$. We have therefore $\tau = t$ if $t$ is $k$-rational.


2. If $\G$ is quasi-split, then $\rP_0$ is a Borel subgroup of $\G$ and $\tau$ is the largest $k$-rational type dominated by $t$.

3. If $t = \varnothing$ is the type of Borel subgroups, then $\tau$ is the minimal $k$-rational type: $\tau = t_{\rm min}$. Indeed, if $k'/k$ is a finite Galois extension splitting $\G$, then ${\rm Osc}_{\varnothing}(\rP_0) \otimes_k k' = {\rm Osc}_{\varnothing}(\rP_0 \otimes_k k')$, $\Qr_0 \otimes_k k'$ is the largest parabolic subgroup of $\G \otimes_k k'$ stabilizing ${\rm Osc}_{\varnothing}(\rP_0 \otimes_k k')$ and $\Qr_0 \otimes_k k' = \rP_0 \otimes_k k'$ since $\rP_0 \otimes_k k'$ is $\varnothing$-relevant. It follows that $\Qr_0 = \rP_0$ by Galois descent.

\vspace{0.2cm} \noindent \textbf{\emph{Proposition}} --- \emph{With the notation above, we have $\overline{\mathcal{B}}_t(\G,k) \cong \overline{\mathcal{B}}_{\tau}(\G,k)$.}

\vspace{0.2cm} We first prove this result for the type $t = \varnothing$ of Borel subgroups, in which case $\tau = t_{\rm min}$ is the type of minimal parabolic subgroups of $\G$.

\vspace{0.2cm}
\noindent \textbf{\emph{Lemma 1}} --- \emph{The projection $$\pi_{\varnothing}^{t_{\rm min}} : {\rm Bor}(\G)^{\rm an} = {\rm Par}_{\varnothing}(\G)^{\rm an} \rightarrow {\rm Par}_{t_{\rm min}}(\G)^{\rm an}$$ induces an homeomorphism between $\overline{\mathcal{B}}_{\varnothing}(\G,k)$ and $\overline{\mathcal{B}}_{t_{\rm min}}(\G,k)$.}

\vspace{0.1cm} \noindent \textbf{\emph{Proof}}. Consider a finite Galois extension $k'/k$ splitting $\G$. It follows easily from results of 4.2 and Galois equivariance that the projection $\pi_{\varnothing}^{t_{\rm min}}$ induces a map $\overline{\mathcal{B}}_{\varnothing}(\G,k) \rightarrow \overline{\mathcal{B}}_{t_{\rm min}}(\G,k)$ satisfying the following condition: for any parabolic subgroup $\Qr$ of $\G$, the preimage of the stratum $\mathcal{B}_{t_{\rm min}}(\Qr_{ss},k)$ is $\mathcal{B}(\Hr_1,k) \times \overline{\mathcal{B}}_{\varnothing}(\Hr_2,k)$, where $\Hr_1$ and $\Hr_2$ are the semi-simple normal and connected subgroups of $\Qr_{ss}$ to which the restrictions of $t_{\rm min}$ are non-degenerate and trivial respectively. Since $t_{\rm min}$ is the minimal $k$-rational type, this implies that $\Hr_2$ has no non-trivial parabolic subgroup, hence is anisotropic over $k$. It follows that $\overline{\mathcal{B}}_{\varnothing}(\Hr_2,k) = \mathcal{B}_{\varnothing}(\Hr_2,k)$ is a point and that the map $\pi_{\varnothing}^{t_{\rm min}} : \overline{\mathcal{B}}_{\varnothing}(\G,k) \rightarrow \overline{\mathcal{B}}_{t_{\rm min}}(\G,k)$ is bijective. This is clearly a homeorphism. \hfill $\Box$

\vspace{0.2cm} We now prove the proposition at the level of apartments.

\vspace{0.1cm}
\noindent \textbf{\emph{Lemma 2}} --- \emph{For any maximal split torus, $\overline{\A}_t(\Sr,k) \cong \overline{\A}_{\tau}(\Sr,k)$.}

\vspace{0.1cm}
\noindent \textbf{\emph{Proof}}. We fix a finite Galois extension $k'/k$ splitting $\G$ and set $\Gamma = {\rm Gal}(k'|k)$. We still denote by $t$ a type of $\G \otimes_k k'$ dominating $t$. Let $\T$ be a maximal torus of $\G$ containing $\Sr$ and satisfying the following conditions: \begin{itemize} \item $\T' = \T \otimes_k k'$ is split; \item the injection $\mathcal{B}(\G,k) \hookrightarrow \mathcal{B}(\G,k')$ maps $\A(\Sr,k)$ into $\A(\T', k')$. \end{itemize}  It follows from the definition of the map $\vartheta_t$ in (2.4.3) that $\overline{\A}_t(\Sr,k)$ can be identified with the closure of $\A(\Sr,k)$ in $\overline{\A}_t(\T',k')$. By Proposition \ref{prop.apartment.identification}, we are reduced to checking that the prefans $\mathcal{F}_t$ and $\mathcal{F}_{\tau}$ on the vector space $\Lambda(\T')$ have the same restriction to $\Lambda(\Sr)$, i.e., that $$\C_t(\rP) \cap \Lambda(\Sr) = \C_{\tau}(\rP) \cap \Lambda(\Sr)$$ for any parabolic subgroup $\rP$ of $\G$ containing $\Sr$. It is enough to consider \emph{minimal} parabolic subgroups of $\G$ containing $\Sr$.

\vspace{0.1cm} So let $\rP_0$ be a minimal parabolic subgroup of $\G$ containing $\Sr$ and denote as above by $\Qr_0$ the largest parabolic subgroup of $\G$ stabilizing ${\rm Osc}_t(\rP_0)$. We write $\rP_0$ and $\Qr_0$ for $\rP_0 \otimes_k k'$ and $\Qr_0 \otimes_k k'$ respectively, and we recall that $\tau$ is by definition the type of $\Qr_0$. Let $\B$ be a Borel subgroup of $\G \otimes_k k'$ satisfying  $\T' \subset \B \subset \rP_0$ and let $\rP$ denote the unique parabolic subgroup of $\G \otimes_k k'$ of type $t$ containing $\B$. We have $\C_t(\rP_0) = \C_t(\Qr_0)$ since ${\rm Osc}_t(\rP_0) = {\rm Osc}_t(\Qr_0)$, and $\C_{\tau}(\rP_0) = \C_{\tau}(\Qr_0)$ since $\Qr_0$ is of type $\tau$ and contains $\rP_0$. Recall that $$\C_t(\rP) = \{\alpha \leqslant 1, \ \ \textrm{for all } \alpha \in \Phi({\rm rad}^{\rm u}(\rP^{\rm op}),\T')\}, \ \ \C_{\tau}(\Qr_0) = \{\alpha \leqslant 1, \ \ \textrm{ for all } \alpha \in \Phi({\rm rad}^{\rm u}(\Qr_0^{\rm op}),\T') \}$$ and $$\C_t(\Qr_0) = \{\alpha \leqslant 1, \ \textrm{ for all } \alpha \in \Phi({\rm rad}^{\rm u} (\rP^{\rm op}), \T')\} \cap \langle \C_t(\Qr_0)\rangle,$$ where $$\langle \C_t(\Qr_0)\rangle =  \{\alpha = 1, \ \textrm{ for all } \alpha \in \Phi(\Lr_{\Qr_0^{\rm op}},\T') \cap \Phi({\rm rad}^{\rm u}(\rP^{\rm op}),\T')\}$$ is the linear subspace spanned by $\C_t(\Qr_0)$ (see Proposition \ref{prop.roots.fan}, (iii)). Since $\B \subset \rP$ and $\B \subset \Qr_0$, the Weyl cone $\mathfrak{C}(\B)$ is contained in both $\C_t(\rP)$ and $\C_{\tau}(\Qr_0)$ and therefore these two cones have overlapping interiors. This observation has the following consequence: for any root $\alpha \in \Phi({\rm rad}^{\rm u}(\rP^{\rm op}),\T')$, the cones $\C_t(\rP)$ and $\C_{\tau}(\Qr_0)$ cannot lie on both sides of the hyperplane $\{\alpha =1\}$, hence $\C_{\tau}(\Qr_0)$ is not contained in the half-space $\{\alpha \geqslant 1\}$ since $\alpha \leqslant 1$ on $\C_t(\rP)$. This implies that $(-\alpha)$ does not belong to $\Phi({\rm rad}^{\rm u}(\Qr_0^{\rm op}),\T')$ or, equivalently, $\alpha \in \Phi(\Qr_0^{\rm op},\T')$. Thus we get $\Phi({\rm rad}^{\rm u}(\rP^{\rm op}),\T') \subset \Phi(\Qr_0^{\rm op},\T')$ and the inclusion $\langle \C_t(\Qr_0) \rangle \cap \C_{\tau}(\Qr_0) \subset \C_t(\Qr_0)$ follows immediately. Since $$\Lambda(\Sr) \subset \langle \mathfrak{C}(\rP_0)\rangle \subset \langle \C_t(\rP_0) \rangle = \langle \C_t(\Qr_0) \rangle,$$ the inclusion $$\C_{\tau}(\Qr_0) \cap \Lambda(\Sr) \subset \C_{t}(\Qr_0) \cap \Lambda(\Sr)$$ is established.

\vspace{0.1cm} Conversely, consider a root $\alpha \in \Phi({\rm rad}^{\rm u}(\Qr_0^{\rm op}),\T')$. The inclusion $\Phi({\rm rad}^{\rm u}(\Qr_0^{\rm op}),\T') \subset \Phi(\rP^{\rm op}, \T')$ being proved as above, $\alpha$ belongs either to $\Phi({\rm rad}^{\rm u}(\rP^{\rm op}), \T')$ or to $\Phi(\Lr_{\rP^{\rm op}},\T')$. In the first case, $\alpha \leqslant 1$ on $\C_t(\Qr_0)$ and thus $\alpha \leqslant 1$ on $\C_t(\Qr_0) \cap \Lambda(\Sr)$.

We address now the case $\alpha \in \Phi(\Lr_{\rP^{\rm op}}, \T')$. Note that $\C_t(\Qr_0) \cap \{\alpha = 1\}$ is a union of Weyl cones and assume that there exists a point $x \in \C_t(\Qr_0)^{\circ} \cap \Lambda(\Sr)$ such that $\alpha(x) = 1$. This point belongs to the interior of some Weyl cone $\mathfrak{C}$ contained in $\C_t(\Qr_0) \cap \{\alpha = 1\}$. Since $\mathfrak{C}^{\circ} \cap \Lambda(\Sr) \neq \varnothing$, this cone corresponds to a parabolic subgroup $\Qr_1$; moreover, we have $\C_t(\Qr_1) = \C_t(\Qr_0)$, for $\mathfrak{C} \cap \C_t(\rP_0)^{\circ} = \mathfrak{C} \cap \C_t(\Qr_0)^{\circ} \neq \varnothing$. It follows that $\Qr_1 \subset \Qr_0$, because $\Qr_0$ is by definition the largest parabolic subgroup of $\G$ such that $\mathfrak{C}(\Qr_0)$ meets the interior of $\C_t(\rP_0)$, hence $\mathfrak{C}(\Qr_0) \subset \mathfrak{C}(\Qr_1)$ and $\alpha = 1$ on $\mathfrak{C}(\Qr_0)$. This last condition amounts to $\alpha \in \Phi(\Lr_{\Qr_0},\T') = \Phi(\Lr_{\Qr_0^{\rm op}},\T')$ and thus leads to a contradiction since we assumed $\alpha \in \Phi({\rm rad}^{\rm u}(\Qr_0^{\rm op}),\T')$. We have therefore $\alpha <1$ or $\alpha > 1$ on $\C_t(\Qr_0)^{\circ} \cap \Lambda(\Sr)$ by convexity, hence $\alpha \leqslant 1$ or $\alpha \geqslant 1$ on $\C_t(\Qr_0) \cap \Lambda(\Sr)$. Since $\alpha$ belongs to $\Phi({\rm rad}^{\rm u}(\Qr_0^{\rm op}),\T')$, we have $\alpha < 1$ on the interior of $\mathfrak{C}(\Qr_0) \subset \C_t(\Qr_0)$ and therefore $\alpha \leqslant 1$ on $\C_t(\Qr_0) \cap \Lambda(\Sr)$.

We have thus proved that each root $\alpha \in \Phi({\rm rad}^{\rm u}(\Qr_0^{\rm op}),\T')$ satisfies $\alpha \leqslant 1$ on $\C_t(\Qr_0) \cap \Lambda(\Sr)$, hence $$\C_t(\Qr_0) \cap \Lambda(\Sr) \subset \C_{\tau}(\Qr_0)$$ and, finally, $$\C_t(\Qr_0) \cap \Lambda(\Sr) = \C_{\tau}(\Qr_0) \cap \Lambda(\Sr).$$ \hfill $\Box$

\vspace{0.2cm} \noindent \textbf{\emph{Proof of Proposition}}. Identifying $\overline{\mathcal{B}}_{\varnothing}(\G,k)$ and $\overline{\mathcal{B}}_{t_{\rm min}}(\G,k)$ by Lemma 1, we have two $\G(k)$-equivariant and continuous maps $$\xymatrix{& \overline{\mathcal{B}}_{\varnothing}(\G,k) \ar@{->}[rd]^{\pi_{\varnothing}^{t}} \ar@{->}[ld]_{\pi_{t_{\rm min}}^{\tau}} & \\ \overline{\mathcal{B}}_{\tau}(\G,k) & & \overline{\mathcal{B}}_t(\G,k).}$$ Consider two points $x, y$ in $\overline{\mathcal{B}}_{\varnothing}(\G,k)$ and let $\overline{\A}_{\varnothing}(\Sr,k)$ be a compactified apartment containing both of them (Proposition \ref{prop - mixed Bruhat dec}, (i)). The conditions $\pi_{\varnothing}^{t}(x) = \pi_{\varnothing}^{t}(y)$ and $\pi_{t_{\rm min}}^{\tau}(x) = \pi_{t_{\rm min}}^{\tau}(y)$ amount to saying that $x$ and $y$ have the same image in the compactified apartments $\overline{\A}_{t}(\Sr,k)$ and $\overline{\A}_{\tau}(\Sr,k)$ respectively, hence are equivalent by Lemma 2. It follows that the diagram above can be completed by a $\G(k)$-equivariant homeomorphism  $$\xymatrix{\overline{\mathcal{B}}_{\tau}(\G,k)   \ar@{->}[r]^{\sim} &  \overline{\mathcal{B}}_t(\G,k)}.$$  \hfill $\Box$

\bibliographystyle{amsalpha}
\bibliography{rtw12}

\providecommand{\bysame}{\leavevmode\hbox to3em{\hrulefill}\thinspace}
\providecommand{\MR}{\relax\ifhmode\unskip\space\fi MR }
\providecommand{\MRhref}[2]{%
  \href{http://www.ams.org/mathscinet-getitem?mr=#1}{#2}
}
\providecommand{\href}[2]{#2}
\begin{thebibliography}{KMRT98}

\bibitem[EGA]{EGA}
Alexander Grothendieck, \emph{{\'E}l\'ements de g\'eom\'etrie alg\'ebrique,
  chapitres {I-IV}}, Publ. Math. IH\'ES (1960-67), vol. 4, 8, etc., r\'edig\'es
  avec la collaboration de Jean Dieudonn\'e.

\bibitem[SGA1]{SGA1}
\bysame, \emph{S\'eminaire de g\'eom\'etrie alg\'ebrique du {B}ois {M}arie ---
  {SGA1}. {R}ev\^etements \'etales et groupe fondamental}, Documents
  math\'ematiques, vol.~3, Soci\'et\'e math\'ematique de France, 2003.

\bibitem[SGA3]{SGA3}
Michel Demazure and Alexander Grothendieck (eds.), \emph{Sch\'emas en groupes.
  {S}\'eminaire de g\'eom\'etrie alg\'ebrique du {B}ois {M}arie 1962/64 ({SGA
  3})}, Lecture Notes in Mathematics, vol. 151-153, Springer, 1970.

\bibitem[AC]{BourbakiAlgCo56}
Nicolas Bourbaki, \emph{{A}lg\`ebre commutative 5-6}, \'El\'ements de
  math\'ematique, Springer, 2007.

\bibitem[INT]{Integration78}
\bysame, \emph{Int\'egration 7-8}, \'El\'ements de math\'ematique, Springer,
  2007.

\bibitem[TS]{BBKspectral}
\bysame, \emph{{T}h\'eories spectrales 1-2}, \'El\'ements de math\'ematique,
  Springer, 2007.

\bibitem[AB08]{AbramenkoBrown}
Peter Abramenko and Kenneth~S. Brown, \emph{Buildings: theory and
  applications}, Graduate Texts in Mathematics, vol. 248, Springer, 2008.

\bibitem[BB66]{BailyBorel}
Walter~L. Baily and Armand Borel, \emph{Compactification of arithmetic
  quotients of bounded symmetric domains}, Ann. of Math. \textbf{84} (1966),
  442--528.

\bibitem[BC91]{BoutotCarayol}
Jean-Fran{\rm \c c}ois Boutot and Henri Carayol, \emph{Uniformisation
  {$p$}-adique des courbes de {S}himura : les th\'eor\`emes de \v {C}erednik et
  de {D}rinfel'd}, Ast\'erisque (1991), no.~196-197, 45--158.

\bibitem[Ber90]{Ber1}
Vladimir~G. Berkovich, \emph{Spectral theory and analytic geometry over
  non-archimedean fields}, Mathematical Surveys and Monographs, vol.~33,
  American Mathematical Society, 1990.

\bibitem[Ber93]{Ber2}
\bysame, \emph{\'{E}tale cohomology for non-archimedean analytic spaces}, Publ.
  Math. IH\'ES \textbf{78} (1993), 5--161.

\bibitem[Ber94]{Ber4}
\bysame, \emph{Vanishing cycles for formal schemes}, Invent. Math.
  \textbf{115} (1994), 539--571.

\bibitem[Ber98]{BerkoICM}
\bysame, \emph{{$p$}-adic analytic spaces}, Proceedings of the International
  Congress of Mathematicians, Vol. II (Berlin, 1998), 1998, pp.~141--151.

\bibitem[BGR84]{BGR}
Siegfried Bosch, Ulrich G{\"u}ntzer, and Reinhold Remmert,
  \emph{Non-archimedean analysis}, Grundlehren der Mathematischen
  Wissenschaften, vol. 261, Springer, 1984.

\bibitem[BJ06]{BorelJi}
Armand Borel and Lizhen Ji, \emph{Compactifications of symmetric and locally
  symmetric spaces}, Mathematics: Theory \& Applications, Birkh\"auser, 2006.

\bibitem[BLR90]{BLR}
Siegfried Bosch, Werner L{\"u}tkebohmert, and Michel Raynaud, \emph{N\'eron
  models}, Ergebnisse der Mathematik und ihrer Grenzgebiete (3), vol.~21,
  Springer, 1990.

\bibitem[Bor91]{Borel}
Armand Borel, \emph{Linear algebraic groups}, Graduate Texts in Mathematics,
  vol. 126, Springer, 1991.

\bibitem[BS73]{BorelSerre}
Armand Borel and Jean-Pierre Serre, \emph{Corners and arithmetic groups},
  Comment. Math. Helv. \textbf{48} (1973), 436--491, Appendice: Arrondissement
  des vari\'et\'es \`a coins, par A. Douady et L. H\'erault.

\bibitem[BT65]{BoTi}
Armand Borel and Jacques Tits, \emph{Groupes r\'eductifs}, Publ. Math. IH\'ES
  \textbf{27} (1965), 55--150.

\bibitem[BT72]{BT1a}
Fran{\rm \c c}ois Bruhat and Jacques Tits, \emph{Groupes r\'eductifs sur un
  corps local, {I}. {D}onn\'ees radicielles valu\'ees}, Publ. Math. IH\'ES
  \textbf{41} (1972), 5--251.

\bibitem[BT84]{BT1b}
\bysame, \emph{Groupes r\'eductifs sur un corps local, {II}. {S}ch\'emas en
  groupes. {E}xistence d'une donn\'ee radicielle valu\'ee}, Publ. Math. IH\'ES
  \textbf{60} (1984), 197--376.

\bibitem[Che95]{ChevalleyBBK}
Claude Chevalley, \emph{Certains sch\'emas de groupes semi-simples},
  S\'eminaire Bourbaki, Vol.\ 6, Exp.\ No.\ 219, Soc. Math. France, Paris,
  1995, pp.~219--234.

\bibitem[Che05]{Bible}
\bysame, \emph{Classification des groupes alg\'ebriques semi-simples},
  Collected works, vol.~3, Springer, 2005.

\bibitem[Del71]{DeligneShimura}
Pierre Deligne, \emph{Travaux de {S}himura}, S\'eminaire {B}ourbaki, 23i\`eme
  ann\'ee (1970/71), {E}xp. 389, Springer, 1971, pp.~123--165. Lecture Notes in
  Math., Vol. 244.

\bibitem[Dem65]{DemazureThese}
Michel Demazure, \emph{Sch\'emas en groupes r\'eductifs}, Bull. Soc. Math.
  France \textbf{93} (1965), 369--413.

\bibitem[DG70a]{DemazureGabriel}
Michel Demazure and Pierre Gabriel, \emph{Groupes alg\'ebriques. {T}ome {I}:
  {G}\'eom\'etrie alg\'ebrique, g\'en\'eralit\'es, groupes commutatifs}, Masson
  \& Cie, \'Editeur, Paris, 1970, avec un appendice {\it Corps de classes
  local} par Michiel Hazewinkel.

\bibitem[Duc08]{Duc}
Antoine Ducros, \emph{Triangulations et cohomologie \'etale sur une courbe
  analytique $p$-adique}, J. Algebraic Geom. \textbf{17} (2008), 503-575.

\bibitem[Fur63]{Furst}
Hillel Furstenberg, \emph{A {P}oisson formula for semi-simple {L}ie groups},
  Ann. of Math. \textbf{77} (1963), 335--386.

\bibitem[GI63]{GoldmanIwahori}
Oscar Goldman and Nagayoshi Iwahori, \emph{The space of p-adic norms}, Acta
  Math. \textbf{109} (1963), 137--177.

\bibitem[GJT98]{GJT}
Yves Guivarc'h, Lizhen Ji, and John~Ch. Taylor, \emph{Compactifications of
  symmetric spaces}, Progress in Mathematics, vol. 156, Birkh\"auser, 1998.

\bibitem[GR06]{GuiRem}
Yves Guivarc'h and Bertrand R{\'e}my, \emph{Group-theoretic compactification of
  {B}ruhat-{T}its buildings}, Ann. Sci. \'Ecole Norm. Sup. \textbf{39} (2006),
  871--920.

\bibitem[Gru66]{Gru}
Laurent Gruson, \emph{Th\'eorie de {F}redholm $p$-adique}, Bulletin de la
  S.M.F. \textbf{94} (1966), 67--95.

\bibitem[KMRT98]{Involutions}
Max-Albert Knus, Alexander Merkurjev, Markus Rost, and Jean-Pierre Tignol,
  \emph{The book of involutions}, American Mathematical Society Colloquium
  Publications, vol.~44, American Mathematical Society, 1998, With a preface in
  French by J. Tits.

\bibitem[Lan96]{La}
Erasmus Landvogt, \emph{A compactification of the {B}ruhat-{T}its building},
  Lecture Notes in Mathematics, vol. 1619, Springer, 1996.

\bibitem[Mar91]{Margulis}
Grigory~A. Margulis, \emph{Discrete subgroups of semisimple {L}ie groups},
  Ergebnisse der Mathematik und ihrer Grenzgebiete (3), vol.~17, Springer,
  1991.

\bibitem[Pra01]{Prasad}
Gopal Prasad, \emph{Galois-fixed points in the {B}ruhat-{T}its building of a
  reductive group}, Bull. Soc. Math. France \textbf{129} (2001), 169--174.

\bibitem[PY02]{PrasadYu}
Gopal Prasad and Jiu-Kang Yu, \emph{On finite group actions on reductive groups
  and buildings}, Invent. Math. \textbf{147} (2002), 545--560.

\bibitem[Rou77]{RousseauHDR}
Guy Rousseau, \emph{Immeubles des groupes r\'eductifs sur les corps locaux},
  U.E.R. Math\'ematique, Universit\'e Paris XI, Orsay, 1977, Th{\`e}se de
  doctorat, Publications Math{\'e}matiques d'Orsay.

\bibitem[Rou08]{RousseauGrenoble}
\bysame, \emph{Euclidean buildings}, G\'eom\'etries \`a courbure n\'egative ou
  nulle, groupes discrets et rigidit\'es (A.~Parreau L.~Bessi\`eres and
  B.~R\'emy, eds.), S\'eminaires et Congr\`es, no.~18, Soci\'et\'e
  math\'ematique de France, 2008.

\bibitem[RTW2]{RTW2}
Bertrand R\'emy, Amaury Thuillier and Annette Werner, \emph{Bruhat-Tits theory from Berkovich's point of view. II. Satake compactifications of buildings}, preprint February 2009.

\bibitem[Sat60]{Sa}
Ichiro Satake, \emph{On representations and compactifications of symmetric
  {R}iemannian spaces}, Ann. of Math. \textbf{71} (1960), 77--110.

\bibitem[Spr98]{Springer}
Tonny~A. Springer, \emph{Linear algebraic groups}, second ed., Progress in
  Mathematics, vol.~9, Birkh\"auser, 1998.

\bibitem[Ste68]{Steinberg}
Robert Steinberg, \emph{Lectures on {C}hevalley groups}, Yale University, New
  Haven, Conn., 1968, Notes prepared by John Faulkner and Robert Wilson.

\bibitem[Tit75]{TitsICM}
Jacques Tits, \emph{On buildings and their applications}, Proceedings of the
  International Congress of Mathematicians (Vancouver, B. C., 1974), Vol. 1,
  1975, pp.~209--220.

\bibitem[Tit79]{TitsCorvallis}
\bysame, \emph{Reductive groups over local fields}, Automorphic forms,
  representations and $L$-functions (Armand Borel and William~A. Casselman,
  eds.), Proc. Symp. Pure Math., vol. XXXIII, part 1, AMS, 1979, pp.~29--69.

\bibitem[Tit86]{TitsCome}
\bysame, \emph{Immeubles de type affine}, Buildings and the geometry of
  diagrams (Como, 1984), Lecture Notes in Math., vol. 1181, Springer, 1986,
  pp.~159--190.

\bibitem[Wat79]{Waterhouse}
William~C. Waterhouse, \emph{Introduction to affine group schemes}, Graduate
  Texts in Mathematics, vol.~66, Springer, 1979.

\bibitem[Wer04]{Wer1}
Annette Werner, \emph{Compactification of the {B}ruhat-{T}its building of {{\rm
  PGL}} by seminorms}, Math. Z. \textbf{248} (2004), 511--526.

\bibitem[Wer07]{Wer2}
\bysame, \emph{Compactifications of {B}ruhat-{T}its buildings associated to
  linear representations}, Proc. Lond. Math. Soc. \textbf{95} (2007), 497--518.

\end{thebibliography}

\vspace{1cm}

\vspace{0.5cm}
\begin{flushleft} \textit{Bertrand R\'emy and Amaury Thuillier} \\
Universit\'e de Lyon \\
Universit\'e Lyon 1 -- CNRS \\
Institut Camille Jordan -- UMR 5208 \\
43 boulevard du 11 novembre 1918 \\
F-69622 Villeurbanne cedex \\
\vspace{1pt}
$\{$remy; thuillier$\}$@math.univ-lyon1.fr
\end{flushleft}

\vspace{0.1cm}
\begin{flushleft}
\textit{Annette Werner} \\
Institut f\"ur Mathematik \\
Goethe-Universit\"at Frankfurt \\
Robert-Maier-Str., 6-8 \\
D-60325 Frankfurt-am-Main \\

\vspace{1pt}
werner@mathematik.uni-frankfurt.de
\end{flushleft}
\end{document}